\documentclass[12pt]{article}

\usepackage[utf8]{inputenc}
\usepackage{amsmath,amssymb,amsthm,geometry,nicefrac,enumerate,hyperref,comment,bbm,mathrsfs,xcolor,color,verbatim,amsfonts,graphicx,bm,dsfont,listings,xparse,etoolbox}
\usepackage{textcomp}
\usepackage{float}

\usepackage{tikz}
\usetikzlibrary{matrix,chains,positioning,decorations.pathreplacing,arrows}
\usetikzlibrary{shapes,fadings,arrows.meta}
\tikzset{font={\fontsize{9pt}{12}\selectfont}}
\usepackage{adjustbox}

\DeclareFontEncoding{LS1}{}{}
\DeclareFontSubstitution{LS1}{stix}{m}{n}
\DeclareMathAlphabet{\mathscr}{LS1}{stixscr}{m}{n}

\usepackage{datetime}
\usepackage[inline]{enumitem}
\usepackage[capitalise,sort]{cleveref}

\crefname{figure}{Figure}{Figures}
\crefname{subsection}{Subsection}{Subsections}
\crefname{enumi}{item}{items}
\crefname{equation}{}{}
\usepackage{cite}

\usepackage[english]{babel}

\usepackage{subfig}

\newcommand{\with}{\curvearrowleft}%"applied with ..." arrow
\newcommand{\descritic}{\cfadd{def:strict_saddle}\text{descending }}

\counterwithin{figure}{section}

\theoremstyle{plain}
\newtheorem{lemma}{Lemma}[section]
\newtheorem{remark}[lemma]{Remark}

\newtheorem{theorem}[lemma]{Theorem}
\newtheorem{definition}[lemma]{Definition}
\newtheorem{corollary}[lemma]{Corollary}
\newtheorem{proposition}[lemma]{Proposition}
\newtheorem{setting}[lemma]{Setting}

\newtheorem{prop}[lemma]{Proposition}
\newtheorem{cor}[lemma]{Corollary}

\theoremstyle{definition}

\numberwithin{equation}{section}

\DeclareFontEncoding{LS1}{}{}
\DeclareFontSubstitution{LS1}{stix}{m}{n}
\DeclareMathAlphabet{\mathscr}{LS1}{stixscr}{m}{n}

\definecolor{cyan(process)}{rgb}{0.0, 0.72, 0.92}
\definecolor{cinnamon}{rgb}{0.82, 0.41, 0.12}
\definecolor{darkred}{rgb}{0.55, 0.0, 0.0} 
\definecolor{persianblue}{rgb}{0.11, 0.22, 0.73}
\hypersetup{colorlinks=true, linkcolor={darkred}, urlcolor=persianblue}

\usepackage{mleftright}
\usepackage{mathtools}

\DeclarePairedDelimiter{\br}{[}{]}
\DeclarePairedDelimiter{\cu}{\{}{\}}
\DeclarePairedDelimiter{\rbr}{(}{)}
\DeclarePairedDelimiter{\abs}{\lvert}{\rvert}

\newcommand{\scalar}[1]{\left<{#1}\right>}
\newcommand{\R}{\mathbb{R}}
\newcommand{\N}{\mathbb{N}}
\newcommand{\E}{\mathbb{E}}

\renewcommand{\P}{\mathbb{P}}

\newcommand{\Z}{\mathbb{Z}}

\newcommand{\dens}{\mathfrak{p}}
\newcommand{\bD}{\mathbb{D}}

\newcommand{\width}{\mathfrak{h}}

\newcommand{\ssum}{\textstyle\sum}

\newcommand{\tint}{\textstyle\int}

\newcommand{\norm}[1]{\lVert #1 \rVert}

\newcommand{\normmm}[1]{{\left\vert\kern-0.25ex\left\vert\kern-0.25ex\left\vert #1
    \right\vert\kern-0.25ex\right\vert\kern-0.25ex\right\vert}}
\newcommand{\Abs}[1]{\big| #1 \big|}

\newcommand{\eps}{\varepsilon}
\newcommand{\cA}{\mathcal{A}}
\newcommand{\cB}{\mathcal{B}}

\newcommand{\cD}{\mathcal{D}}

\newcommand{\cF}{\mathcal{F}}
\newcommand{\cG}{\mathcal{G}}

\newcommand{\cJ}{\mathcal{J}}

\newcommand{\cM}{\mathcal{M}}
\newcommand{\cN}{\mathcal{N}}
\newcommand{\cO}{\mathcal{O}}

\newcommand{\cS}{\mathcal{S}}

\newcommand{\fA}{\mathfrak{A}}

\newcommand{\fC}{\mathfrak{C}}
\newcommand{\fD}{\mathfrak{D}}

\newcommand{\fG}{\mathfrak{G}}
\newcommand{\fH}{\mathfrak{H}}

\newcommand{\fJ}{\mathfrak{J}}

\newcommand{\fM}{\mathfrak{M}}
\newcommand{\fN}{\mathfrak{N}}

\newcommand{\fa}{\mathfrak{a}}
\newcommand{\fb}{\mathfrak{b}}
\newcommand{\fc}{\mathfrak{c}}
\newcommand{\fd}{\mathfrak{d}}

\newcommand{\fg}{\mathfrak{g}}
\newcommand{\fh}{\mathfrak{h}}
\newcommand{\fj}{\mathfrak{j}}
\newcommand{\fk}{\mathfrak{k}}

\newcommand{\fn}{\mathfrak{n}}

\newcommand{\fp}{\mathfrak{p}}

\newcommand{\ft}{\mathfrak{t}}

\newcommand{\fw}{\mathfrak{w}}

\newcommand{\bfc}{\mathbf{c}}

\newcommand{\bfm}{\mathbf{m}}
\newcommand{\bfn}{\mathbf{n}}

\newcommand{\bft}{\mathbf{t}}

\newcommand{\bfT}{\mathbf{T}}

\newcommand{\scrC}{\mathscr{C}}

\newcommand{\scrR}{\mathscr{R}}

\newcommand{\scra}{\mathscr{a}}
\newcommand{\scrb}{\mathscr{b}}
\newcommand{\scrc}{\mathscr{c}}

\newcommand{\scri}{\mathscr{i}}
\newcommand{\scrj}{\mathscr{j}}

\newcommand{\scrt}{\mathscr{t}}
\newcommand{\scru}{\mathscr{u}}
\newcommand{\scrv}{\mathscr{v}}
\newcommand{\scrw}{\mathscr{w}}
\newcommand{\scrx}{\mathscr{x}}
\newcommand{\scry}{\mathscr{y}}
\newcommand{\scrz}{\mathscr{z}}
\NewDocumentEnvironment{cproof}{m}
{\begin{proof}[Proof of \cref{#1}]}%
{\noindent The proof of \cref{#1} is thus complete.
\end{proof}}

\NewDocumentEnvironment{cproof2}{m}
{\begin{proof}[Proof of \cref{#1}]}%
{\noindent This completes the proof of \cref{#1}.
\end{proof}}

\ExplSyntaxOn

\NewDocumentCommand{\enum}{ O{;} m o }
 {
  \my_enum:nnn { #1 } { #2 } { #3 }
 }

\seq_new:N \l__my_enum_seq
\tl_new:N \l__my_enum_item_tl
\prop_const_from_keyval:Nn \l__verbs
{
show = shows ,
imply = implies ,
demonstrate = demonstrates ,
prove = proves ,
establish = establishes ,
ensure = ensures ,
assure = assures
}
\int_new:N \l__number_of_args

\cs_new_protected:Nn \my_enum:nnn
 {
  \seq_set_split:Nnn \l__my_enum_seq { #1 } { #2 }
  \seq_remove_all:Nn \l__my_enum_seq {}
  \int_set_eq:NN \l__number_of_args { \seq_count:N \l__my_enum_seq }
  \seq_use:Nnnn \l__my_enum_seq { ~and~ } { ,~ } { ,~and~ }
  \IfNoValueTF{#3}{}{
    \space
    \int_compare:nNnTF{ \l__number_of_args } < {2}{ \prop_item:Nn \l__verbs {#3} }{ #3 }
  }
 }

\seq_new:N \g_cflist_loaded
\seq_new:N \g_cflist_pending

\NewDocumentCommand{\cfadd}{ m }
{
  \seq_if_in:NnF \g_cflist_loaded { #1 } {
    \seq_if_in:NnF \g_cflist_pending { #1 } {
      \seq_gput_right:Nn \g_cflist_pending { #1 }
    }
  }
}

\NewDocumentCommand{\cfconsiderloaded}{ m }{
  \seq_gput_right:Nn \g_cflist_loaded {#1}
}

\NewDocumentCommand{\cfremove}{ m }
{
  \seq_gremove_all:Nn \g_cflist_pending { #1 }
}

\NewDocumentCommand{\cfload}{ o }
{
  \seq_if_empty:NTF \g_cflist_pending {\unskip} {
    (cf.\ \cref{\seq_use:Nn \g_cflist_pending {,}})\IfValueTF{#1}{#1~}{\unskip}
    \seq_gconcat:NNN \g_cflist_loaded \g_cflist_loaded \g_cflist_pending
    \seq_gclear:N \g_cflist_pending
  }
}

\NewDocumentCommand{\cfclear} {} {
  \seq_gclear:N \g_cflist_loaded
  \seq_gclear:N \g_cflist_pending
}

\NewDocumentCommand{\cfout}{ o }
{
  \seq_if_empty:NTF \g_cflist_pending {\unskip} {
    (cf.\ \cref{\seq_use:Nn \g_cflist_pending {,}})\IfValueTF{#1}{#1~}{\unskip}
    \seq_gclear:N \g_cflist_pending
  }
}

\NewDocumentCommand{\ifnocf} { m } {
  \seq_if_empty:NT \g_cflist_pending { #1 }
}

\ExplSyntaxOff
\ExplSyntaxOn

\bool_new:N \g_noteobserve

\NewDocumentCommand{\nobs}{}{
  \bool_if:nTF { \g_noteobserve } {
    \bool_gset_false:N \g_noteobserve
    note~
  } {
    \bool_gset_true:N \g_noteobserve
    observe~
  }
}

\NewDocumentCommand{\Nobs}{}{
  \bool_if:nTF { \g_noteobserve } {
    \bool_gset_false:N \g_noteobserve
    Note~
  } {
    \bool_gset_true:N \g_noteobserve
    Observe~
  }
}

\ExplSyntaxOff

\ExplSyntaxOn

\bool_new:N \g_hencetherefore

\NewDocumentCommand{\hence}{}{
  \bool_if:nTF { \g_hencetherefore } {
    \bool_gset_false:N \g_hencetherefore
    hence~
  } {
    \bool_gset_true:N \g_hencetherefore
    therefore~
  }
}

\NewDocumentCommand{\Hence}{}{
  \bool_if:nTF { \g_hencetherefore } {
    \bool_gset_false:N \g_hencetherefore
    Hence,~we~obtain~
  } {
    \bool_gset_true:N \g_hencetherefore
    Therefore,~we~obtain~
  }
}

\ExplSyntaxOff

\NewDocumentEnvironment {athm} {m m} {%
\begin{#1}\label{#2}\global\def\loc{#2}%
}{%
\end{#1}%
}

\NewDocumentEnvironment {adef} {m} {%
\begin{definition}\label{#1}\global\def\loc{#1}%
}{%
\end{definition}%
}

\NewDocumentEnvironment{aproof} {} {%
\begin{proof}[Proof~of~\cref{\loc}]%
}{%
\end{proof}%
}

\ExplSyntaxOff

\newcommand{\Exists}{\exists\,}
\newcommand{\Forall}{\forall\,}
\newcommand{\andq}{\text{and} \quad}
\newcommand{\qandq}{\quad \text{and} \quad}
\newcommand{\qqandqq}{\qquad \text{and} \qquad}

\newcommand{\loss}{\mathcal{L}}
\newcommand{\lossapp}{\mathfrak{L}}

\renewcommand{\d}{\mathrm{d}} %for integration "dx"
\renewcommand{\emptyset}{\varnothing}
\newcommand{\Hs}{\operatorname{Hess}}

\newcommand{\dimension}{\mathfrak{d}}
\newcommand{\degerror}{\kappa}
\newcommand{\0}{\scra}
\newcommand{\1}{\scrb}

\newcommand{\const}{\bfc}
\newcommand{\batchsize}{\mathbf{m}}
\newcommand{\indexset}{\fJ}
\newcommand{\Lip}{\operatorname{Lip}}
\newcommand{\ceil}[1]{ \left\lceil #1 \right\rceil \cfadd{def:ceiling}}
\newcommand{\floor}[1]{ \lfloor #1 \rfloor \cfadd{def:ceiling}}
\newcommand{\indicator}[1]{\mathbbm{1}_{\smash{#1}}}

\newcommand{\realization}[1] {\mathcal{N} ^{ #1  }}

\renewcommand{\emptyset}{\varnothing}
\newcommand{\act}{\mathbb{A}}

\DeclarePairedDelimiter{\spro}{\langle}{\rangle}

\renewcommand{\d}{ \mathrm{d}}

\definecolor{codegreen}{rgb}{0,0.6,0}
\definecolor{codegray}{rgb}{0.5,0.5,0.5}
\definecolor{codepurple}{rgb}{0.58,0,0.82}
\definecolor{backcolour}{rgb}{0.95,0.95,0.92}

\lstdefinestyle{mystyle}{
    backgroundcolor=\color{backcolour},   
    commentstyle=\color{codegreen},
    keywordstyle=\color{magenta},
    numberstyle=\tiny\color{codegray},
    stringstyle=\color{codepurple},
    basicstyle=\ttfamily\footnotesize,
    breakatwhitespace=false,         
    breaklines=true,                 
    captionpos=b,                    
    keepspaces=true,                 
    numbers=left,                    
    numbersep=5pt,                  
    showspaces=false,                
    showstringspaces=false,
    showtabs=false,                  
    tabsize=2
}

\crefname{listing}{Source code}{Source codes}

\makeatletter
\ExplSyntaxOn
\seq_new:N \g_abbrs
\prop_new:N \g_abbr_counts
\tl_new:N \l_abbr_count_tl

\ExplSyntaxOn

\bool_new:N \g_forexample

\NewDocumentCommand{\eg}{ o }{
	\IfValueT{#1}{
		\str_if_eq:noTF {fe} {#1} {
			\bool_gset_true:N \g_forexample
		} {\bool_gset_false:N \g_forexample}
	}
	\bool_if:nTF { \g_forexample } {
		\bool_gset_false:N \g_forexample
		for~example
	}{
		\bool_gset_true:N \g_forexample
		for~instance
	}
}

\NewDocumentCommand{\abbr}{m m O{#1} m m O{#4} m}{
	\expandafter\newcommand\csname#3\endcsname[1][]{
		\seq_if_in:NnTF \g_abbrs {#1} {
			\prop_get:NnN \g_abbr_counts {#1} \l_abbr_count_tl
			\prop_gput:Nnx \g_abbr_counts {#1} {\int_eval:n {\l_abbr_count_tl + 1}}
			\hyperref[#1]{#7}
		} {
			\seq_gput_left:Nn \g_abbrs {#1}
			\prop_gput:Nnn \g_abbr_counts {#1} {1}
			\expandafter\gdef\csname#1@def\endcsname{#2}
			\phantomsection\label{#1}
			\str_if_eq:nnTF{##1}{}{\emph{#2}}{##1}~(\hyperref[#1]{#7})
		}
	}
	\expandafter\newcommand\csname#6\endcsname[1][]{
		\seq_if_in:NnTF \g_abbrs {#1} {
			\prop_get:NnN \g_abbr_counts {#1} \l_abbr_count_tl
			\prop_gput:Nnx \g_abbr_counts {#1} {\int_eval:n {\l_abbr_count_tl + 1}}
			\hyperref[#1]{#4}
		} {
			\expandafter\gdef\csname#1@def\endcsname{#5}
			\seq_gput_left:Nn \g_abbrs {#1}
			\prop_gput:Nnn \g_abbr_counts {#1} {1}
			\phantomsection\label{#1}
			\str_if_eq:nnTF{##1}{}{\emph{#5}}{##1}~(\hyperref[#1]{#4})
		}
	}
}

\ExplSyntaxOff

\makeatother
\abbr{Adamax}{adaptive moment estimation maximum SGD}{Adamaxs}{Adaptive moment estimation maximum SGD}{Adamax}
\abbr{Adam}{adaptive moment estimation SGD}{Adams}{adaptive moment estimation SGD}{Adam}
\abbr{Nadam}{Nesterov accelerated adaptive moment estimation SGD}{Nadams}{Nesterov accelerated adaptive moment estimation SGD}{Nadam}
\abbr{RMSprop}{root mean square propagation SGD}{RMSprops}{the root mean square propagation}{RMSprop}
\abbr{ANN}{artificial neural network}{ANNs}{artificial neural networks}{ANN}
\abbr{DNN}{deep artificial neural network}{DNNs}{deep artificial neural networks}{DNN}
\abbr{GD}{Gradient descent}{GDs}{gradient descent}{GD}
\abbr{GF}{gradient flow}{GFs}{gradient flow}{GF}
\abbr{SGD}{stochastic gradient descent}{SGDs}{stochastic gradient descent}{SGD}
\abbr{iid}{independent and identically distributed}{i.i.d}{independent and identically distributed}{i.i.d.}
\abbr{DL}{Deep learning}{DLs}{Deep learning}{DL}
\abbr{AI}{artificial intelligence}{AIs}{artificial intelligence}{AI}
\abbr{LLM}{large language model}{LLMs}{large language models}{LLM}
\abbr{PDE}{partial differential equation}{PDEs}{partial differential equations}{PDE}
\abbr{as}{almost surely}{ass}{almost surely}{a.s.}
\abbr{pas}{$\P$-almost surely}{pass}{almost surely}{$\P$-a.s.}
\abbr{ReLUs}{rectified linear unit}{ReLU}{rectified linear unit}
\makeatother

\begin{document}

\title{Convergence to good non-optimal critical points\\ in the training of neural networks: Gradient descent optimization 
with one random initialization overcomes\\ all bad non-global local minima with high probability}

\author{Shokhrukh Ibragimov$^{1}$, Arnulf Jentzen$^{2,3}$, and Adrian Riekert$^{4}$\bigskip\\
\small{$^1$ Applied Mathematics: Institute for Analysis and Numerics,}\vspace{-0.1cm}\\
\small{University of M\"unster, Germany; e-mail: \texttt{sibragim}\textcircled{\texttt{a}}\texttt{uni-muenster.de}}\smallskip\\
\small{$^2$ School of Data Science and School of Artificial Intelligence,} \vspace{-0.1cm}\\
\small{The Chinese University of Hong Kong, Shenzhen, China; e-mail: \texttt{ajentzen}\textcircled{\texttt{a}}\texttt{cuhk.edu.cn}}\smallskip\\
\small{$^3$ Applied Mathematics: Institute for Analysis and Numerics,}\vspace{-0.1cm}\\
\small{University of M\"unster, Germany; e-mail: \texttt{ajentzen}\textcircled{\texttt{a}}\texttt{uni-muenster.de}}\smallskip\\
\small{$^4$ Applied Mathematics: Institute for Analysis and Numerics,}\vspace{-0.1cm}\\
\small{University of M\"unster, Germany; e-mail: \texttt{ariekert}\textcircled{\texttt{a}}\texttt{uni-muenster.de}}}

\date{\today}
%\date{\currenttime \, \today}

\maketitle
%\pagebreak
\begin{abstract}
Gradient descent (GD) optimization methods for the training of artificial neural networks (ANNs) belong nowadays to the most heavily employed computational schemes in the digital world. Despite the compelling success of such methods, it remains an open problem of research to provide a rigorous theoretical justification for the success of GD optimization schemes in the training of ANNs. The main difficulty is the issue that the optimization risk landscapes associated to ANNs usually admit many non-optimal critical points (saddle points as well as non-global local minima) whose risk values are strictly larger than the optimal risk value (the smallest possible risk value). It is a key contribution of this article to overcome this obstacle in certain simplified shallow ANN training situations where only the bias parameters but not the weight parameters of the ANNs are trained, where certain special normal random initializations are employed, and where the architecture of the considered shallow ANNs just consists of a one-dimensional input layer, a multi-dimensional hidden layer, and a one-dimensional output layer. Specifically, in such simplified ANN training scenarios we prove that the gradient flow (GF) dynamics with only one random initialization overcomes with high probability all bad non-global local minima (all non-global local minima whose risk values are much larger than the risk value of the global minima) and converges with high probability to a good critical point (a critical point whose risk value is very close to the optimal risk value of the global minima). This analysis allows us to establish convergence in probability \emph{to zero} of the risk value of the GF trajectories with convergence rates as the ANN training time and the width of the ANN increase to infinity. We expect that the reason for the convergence of the risk to zero established in this work in the above outlined simplified shallow ANN training situations (with high probability overcoming all bad non-global local minima and convergence to good critical points) is not only the reason of convergence in the simplified shallow ANN training setups studied analytically in this work but is also representative for other practically relevant situations in the training of deep ANNs with many hidden layers. In particular, we complement the analytical findings of this work with extensive numerical simulations for the standard stochastic GD optimization method as well as the popular Adam optimizer in the training of shallow and deep ANNs: the results of all these numerical simulations strongly suggest that with high probability the considered GD optimization method overcomes all bad non-global local minima, does not converge to a global minimum, but does converge to a good non-optimal critical point whose risk value is very close to the optimal risk value.
\end{abstract}

%\pagebreak 
\tableofcontents
%==============================================================================%
%------------------------------------------------------------------------------%
%=================================-----Section-----============================%
%------------------------------------------------------------------------------%
%==============================================================================%
%\newpage
\section{Introduction}

\GD\ optimization methods for the training of \ANNs\ belong these days to the most heavily employed numerical schemes in the digital world. In particular, such computational techniques are significantly employed in a multiplicity of applications ranging from text processing and speech recognition to computer vision, just to name a few (cf., e.g., Le Cun et al.~\cite{LeCunBengioHinton2015} and Bottou et al.~\cite[Section 2]{Bottou2018optimization}).

Despite the ubiquitous presence of \GD\ optimization methods in the training of \ANNs, it remains a fundamental open problem of research to provide a rigorous theoretical justification for the practical success of such methods. \GD\ optimization methods (by which we mean the negative gradient flow dynamics, deterministic \GD\ schemes, as well as stochastic \GD\ schemes) walk along the negative gradient with the aim to minimize the risk function of the considered supervised learning problem. The fundamental challenge, which arises in the mathematical analysis of such methods, is the fact that there are many non-optimal critical points (many saddle points as well as many non-global local minimum points) in the optimization landscape (cf., e.g.,~\cite{ShokhrukhLocalMinima,Safran2018spurious,Swirszcz2017local}) whose risk values are strictly larger than the optimal risk value (the smallest possible risk value) and that the probability of the domain of attraction of the non-global local minimum points does not vanish (see, e.g., Cheridito et al.~\cite{FlorianEscapeSaddles2022}).

\def\layersep{4cm}
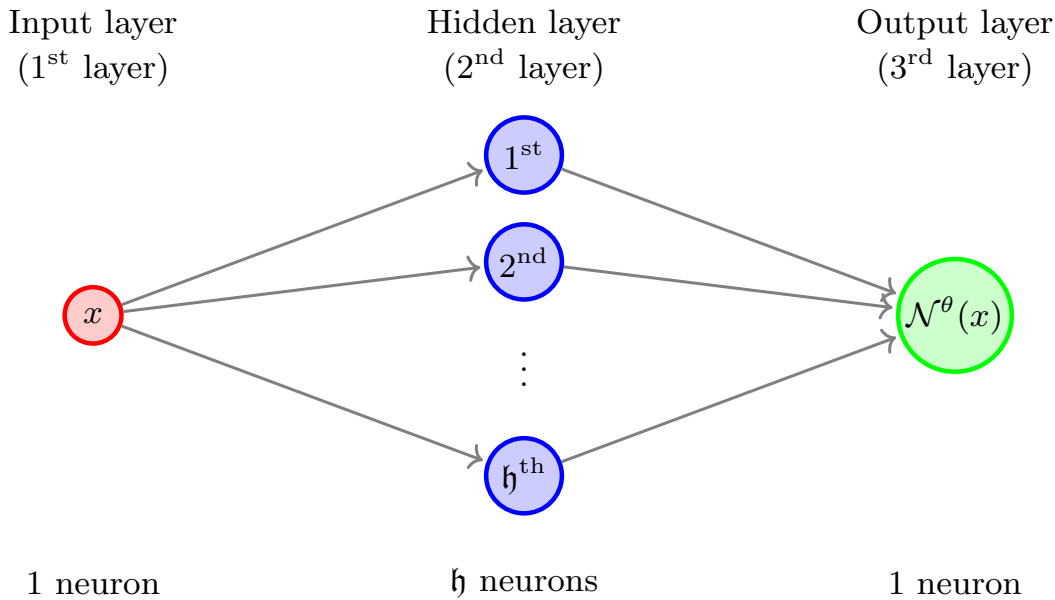
\begin{figure}[H]
\centering
\begin{adjustbox}{width=\textwidth}
\begin{tikzpicture}[shorten >=1pt,->,draw=black!50, node distance=\layersep]
\tikzstyle{every pin edge}=[<-,shorten <=1pt]
\tikzstyle{input neuron}=[very thick, circle,draw=red, fill=red!20, minimum size=15pt,inner sep=0pt]
\tikzstyle{output neuron}=[very thick, circle, draw=green,fill=green!20,minimum size=30pt,inner sep=0pt]
\tikzstyle{hidden neuron}=[very thick, circle,draw=blue,fill=blue!20,minimum size=20pt,inner sep=0pt]
\tikzstyle{annot} = [text width=9em, text centered]
\tikzstyle{annot2} = [text width=4em, text centered]

%----------Neuron(s) of input layer----------
\node[input neuron] (I) at (0,-1) {$x$};

%----------Neuron(s) of hidden layer----------
\path[yshift = 1.5cm]
%node[hidden neuron] (H0-\name) at (\layersep, -\y cm) {};
node[hidden neuron](H0-1) at (\layersep, -1 cm) {$1^{\text{st}}$};
\path[yshift = 1.5cm]
node[hidden neuron](H0-2) at (\layersep, -2 cm) {$2^{\text{nd}}$};
\path[yshift = 1.5cm]
node(H0-dots) at (\layersep, -2.9 cm) {\vdots};
\path[yshift = 1.5cm]
node[hidden neuron](H0-3) at (\layersep, -4 cm) {$\width^{\text{th}}$};

%----------Neuron(s) of output layer----------
\path[yshift = 1.5cm]
node[output neuron](O) at (2*\layersep,-2.5 cm) {$\cN^{\theta}(x)$};   

%----------Arrows from Input to Hidden layer----------
\foreach \dest in {1,2,3}
\path[line width = 0.8 ] (I) edge (H0-\dest);

%----------Arrows from Hidden to Output layer----------
\foreach \source in {1,2,3}
\path[line width = 0.8 ] (H0-\source) edge (O);

% Annotate the layers
\node[annot,above of=H0-1, node distance=1cm, align=center] (hl) {Hidden layer\\($2^{\text{nd}}$ layer)};
\node[annot,left of=hl, align=center] {Input layer\\ ($1^{\text{st}}$ layer)};
\node[annot,right of=hl, align=center] {Output layer\\($3^{\text{rd}}$ layer)};

\node[annot2,below of=H0-3, node distance=1cm, align=center] (sl) {$\width$ neurons};
\node[annot2,left of=sl, align=center] {$1$ neuron};
\node[annot2,right of=sl, align=center] {$1$ neuron};
\end{tikzpicture}
\end{adjustbox}
\caption{Graphical illustration of the \ANN\ architecture considered 
in \cref{theorem_intro_2} and \cref{theorem_intro}: 
\ANNs\ with a $ 1 $-dimensional input layer, $ \width $-dimensional hidden layer, and $ 1 $-dimensional output layer.}
\label{figure_shallow_nns_general_introduction}
\end{figure}

It is a key contribution of this work to overcome this difficulty in certain simplified training situations where the considered \ANNs\ are shallow with just one hidden layer and a one-dimensional input (see also \cref{figure_shallow_nns_general_introduction} for a graphical illustration of the analyzed \ANN\ architecture), where certain special random initialization schemes are employed, and where only the \emph{bias parameter} of the \ANNs\ but not the \emph{weight parameters} of the \ANNs\ are trained. Specifically, in such simplified \ANN\ training situations we prove under one random initialization
\begin{enumerate}[label = (\Roman*)]
\item
\label{intro_item_I}
that the \GF\ dynamics overcomes with high probability all \emph{bad} non-global local minimum points (overcomes all non-global local minimum points whose risk values are far away from the optimal risk value of the global minimum points) and
\item
\label{intro_item_II}
that the \GF\ dynamics does converge with high probability to a \emph{good} critical point (does converge with high probability to a critical point whose risk value is very close to the optimal risk value of the global minimum points).
\end{enumerate}
This analysis allows to conclude that with high probability we have that the risk of the \GF\ trajectory $(\Theta_t)_{t \in [0,\infty)}$ converges \emph{to zero} as, both, the training time $t \geq 0$ \emph{and} the width $\width \in \N = \{1, 2, 3, \dots \}$ of the \ANN\ converge to infinity.

In \cref{theorem_intro_2} below we disclose such a result for \ANNs\ with the clipping activation function $\R \ni x \mapsto \max\{ \min\{x, 1\}, 0\} \in \R$ (see \cref{theorem_intro_2:eq:risk}) and in \cref{theorem_intro} below we present such a result for \ANNs\ with the \ReLU\ activation function $\R \ni x \mapsto \max\{ x, 0 \} \in \R$ (see \cref{eq:theorem_intro_risk_function}). We now first present the precise statement of \cref{theorem_intro_2}, thereafter, we comment on the statement and the mathematical objects in \cref{theorem_intro_2}, and, thereafter, we finally present \cref{theorem_intro}.

\begin{samepage}
\cfclear
\begin{theorem}\label{theorem_intro_2}
Let $\0 \in \R$, $ \1 \in (\0, \infty)$, let $f, \dens \in C(\R, \R)$ be piecewise analytic, assume for all $ x \in \R $ that $\dens(x) \ge 0$ and $\dens^{ -1 }( \R \backslash \{ 0 \} ) = ( \0, \1 ) $, assume that $ f $ is strictly increasing, let $ ( \Omega, \cF, \P ) $ be a probability space, for every $ \width \in \N $ let $ \fd_\width = 3 \width + 1 $, let $ \loss_{ \width } \colon \R^{ \fd_{ \width } } \to \R $ satisfy for all $ \theta = ( \theta_1, \dots, \theta_{\dimension_{ \fh } } ) \in \R^{ \dimension_{ \width } } $ that
\begin{equation}\label{theorem_intro_2:eq:risk}
\textstyle \loss_{ \width }( \theta ) = \int_{\R} (f(x) - \theta_{\fd_\width} - \sum_{j = 1}^{\width} \theta_{2 \width + j} \max\cu{\min\cu{\theta_{\fh + j} x + \theta_{j}, 1}, 0})^2 \dens(x) \, \d x,
\end{equation}
let $\cG_\width = ( \cG_\width^1, \dots, \cG_\width^{ \fd_\width } )\colon \R^{ \fd_\width } \to \R^{ \fd_\width }$ satisfy for all $ k \in \{ 1, 2, \dots, \fd_\width \} $, $\theta = ( \theta_1, \dots, \theta_{ \dimension_{ \fh } } ) \in \R^{ \fd_\width }$ that $\cG_\width^k(\theta) = ( \frac{ \partial^{ - } }{ \partial \theta_k }\loss_\width )( \theta ) \indicator{ \N \backslash ( \width, 3 \width]}( k )$, and let $\Theta^\width = (\Theta^{\width, 1}, \ldots, \Theta^{\width, \fd_\width}) \colon [0, \infty) \times \Omega \to \R^{\fd_\width}$ be a stochastic process with continuous sample paths, assume for all $\width \in \N$, $t \in [0, \infty)$ that
\begin{equation}\label{theorem_intro_2:eq:assgf}
\textstyle \P\bigl(\Theta_t^\width = \Theta_0^\width - \int_0^t \cG_\width (\Theta_s^\width) \, \d s \bigr) = 1,
\end{equation}
assume for all $ \width \in \N $ that $ \Theta_0^{ \width, 1 }, \dots, \Theta_0^{ \width, \fd_\width } $ are independent, let $ \alpha \in (\nicefrac{ 3 }{ 4 }, 1 ) $, $ \beta \in ( \alpha + 2, \infty ) $, assume for all $\width \in \N$, $k \in \{1, 2, \ldots, \width\}$ that $\width^{- \beta} \Theta_0^{\width, k}$ is standard normal, and assume for all $\width \in \N$, $k \in \{1, 2, \ldots, \width\}$, $x \in \R$ that $\P(\width^{- \beta} \Theta_0^{\width, \width + k} \leq x) = \P(\width^{\alpha} \Theta_0^{\width, 2 \width + k} \le x) = [\frac{2}{\pi}]^{1/2} \int_0^x \exp(- \frac{y^2}{2}) \, \d y$ and $\E[\abs{\Theta_0^{\width, \fd_\width}}^2] < \infty $. Then there exist $\fc, \degerror \in (0, \infty)$ and random variables $\fC_{\fh} \colon \Omega \to \R$, $\fh \in \N$, such that for all $\fh \in \N$ it holds that $\limsup_{t \to \infty} \E[\loss_{\fh}(\Theta_t^{\fh})] \leq \fc \fh^{- \degerror}$ and
\begin{equation}\label{eq:theorem_2_intro_statement}
\textstyle \P\bigl(\Forall t \in (0, \infty) \colon \loss_{ \fh }(\Theta_t^{ \fh } ) \le \fc \fh^{- \degerror} + \fC_{\fh} t^{- 1} \bigr) \ge 1 - \fc \fh^{- \degerror}.
\end{equation}
\end{theorem}
\end{samepage}
\cref{theorem_intro_2} is a direct consequence of \cref{theo:main:conv} in \cref{sec:convergence_of_GF_with_random_init_clipping} below and in \cref{sec:anns_with_clipping} below we present the detailed proof of \cref{theo:main:conv}. Next let us add some comments on the statement and the mathematical objects appearing in \cref{theorem_intro_2}.

The real numbers $ \0, \1 \in \R $ with $ \0 < \1 $ in \cref{theorem_intro_2} specify the borders of the support $[\0, \1]$ of the probability distribution of the input data of the considered supervised learning problem. The function $ f \colon \R \to \R $ in \cref{theorem_intro_2} is the target function (the function which we intend to approximate via \ANNs) of the considered supervised learning problem and the function $ \dens \colon \R \to \R $ in \cref{theorem_intro_2} is an unnormalized density function of the probability distribution of the input data of the considered supervised learning problem. Both $ f $ and $ \dens $ are assumed to be continuous and piecewise analytic, by which we mean that there exist $ K \in \N $, $  \tau_1, \dots, \tau_K \in \R$ with
\begin{equation}
\tau_1 < \tau_2 < \dots < \tau_K = \infty
\end{equation}
such that for all $ k \in \{ 2, 3, \dots, K \} $ we have that $ f|_{ [\tau_{ k - 1 }, \tau_k ] } $, $f | _ { (- \infty , \tau_1 ]}$, $f| _ { [ \tau_K, \infty ) }$ and $ \dens|_{ [ \tau_{ k - 1 }, \tau_k ] }$, $\dens | _ { (- \infty , \tau_1 ]}$, $\dens | _ { [ \tau_K, \infty ) }$ are analytic functions (cf.\ \cref{def:piecewise_analytic} in \cref{sec:anns_with_clipping} below). 
This notion generalizes the setting of a piecewise polynomial function which we considered in our previous article \cite{JentzenRiekert2021ExistenceGlobMin}.
In addition, we also assume in \cref{theorem_intro_2} that the target function $ f $ is strictly increasing.

In \cref{theorem_intro_2} we study the \GF\ dynamics with one random initialization and the triple $ ( \Omega, \mathcal{F}, \P ) $ in \cref{theorem_intro_2} serves as the underlying probability space to specify this random initialization. In \cref{theorem_intro_2} we analyze \emph{shallow fully-connected feedforward \ANNs} with $ 1 $ neuron on the input layer (with an $ 1 $-dimensional input layer), with $ \width \in \N $ neurons on the hidden layer (with an $\width$-dimensional hidden layer), and with $ 1 $ neuron on the output layer (with an $ 1 $-dimensional output layer). There are thus $ \width $ weight parameters and $ \width $ bias parameters to describe the affine linear transformation between the $ 1 $-dimensional input layer and the $ \width $-dimensional hidden layer and there are thus $\width$ weight parameters and $ 1 $ bias parameter to describe the affine linear transformation between the $ \width $-dimensional input layer and the $ 1 $-dimensional output layer. Altogether there are thus
\begin{equation}
\width + \width + \width + 1 = 3 \width + 1
\end{equation}
real parameters to describe the considered \ANNs\ and it is thus precisely the role of the natural numbers $ \fd_{ \width } = 3 \width + 1 $ for $ \width \in \N $ to specify the number of real parameters in the considered \ANNs. We also refer to \cref{figure_shallow_nns_general_introduction} for a graphical illustration of the \ANN\ architecture used in \cref{theorem_intro_2}.

In \cref{theorem_intro_2:eq:risk} we consider the realization of \ANNs\ with the clipping activation function, which is for every $ \theta = (\theta_1, \dots, \theta_{ \fd_{ \width } } ) \in \R^{ \fd_\width } $ given by
\begin{equation}
\textstyle \R \ni x \mapsto \theta_{\fd_\width} + \sum_{j = 1}^{\width} \theta_{2 \width + j} \max\cu{\min\cu{\theta_{\fh + j} x + \theta_{j}, 1}, 0} \in \R.
\end{equation}
\Nobs that the clipping function is a continuous sigmoidal function in the sense of Cybenko~\cite{Cybenko1989} and, in particular, admits the universal approximation property (cf., e.g.,~\cite{Cybenko1989,Hornik1989,LeshnoLinPinkusSchocken1993}). Similarly to the more commonly used \ReLU\ activation function, the clipping function is piecewise linear and almost everywhere differentiable. Thus it is possible to define a suitable generalized gradient in an analogous way. In the setting considered in this article, the risk function is actually continuously differentiable due to the integral appearing in \cref{theorem_intro_2:eq:risk}. In addition, the clipping function has the convenient property that it is bounded and its derivative has bounded support, which makes it easier to establish boundedness of the considered \GF\ trajectories in some cases. For further uses of the clipping activation for \ANNs\ we refer, e.g., to \cite{ChoiWang2018PACT} and \cite{SakrChoi2018}.

The functions $ \loss_{ \width } \colon \R^{ \fd_{ \width } } \to \R $ for $ \width \in \N $ in \cref{theorem_intro_2} are the risk functions which specify the squared $ L^2 $-distance between the target function $ f \colon \R \to \R $ and its \ANN\ approximations. The functions $\cG_\width = ( \cG_\width^1, \dots, \cG_\width^{ \fd_\width } ) \colon \R^{\fd_\width} \to \R^{\fd_\width} $ for $ \width \in \N $ in \cref{theorem_intro_2} describes the generalized gradient functions that are used in the training process. Note that the indicator functions $\N \ni k \mapsto \indicator{ \N \backslash ( \width, 3 \width ] }( k ) \in \{ 0, 1 \}$ for $ \width \in \N $ in the description of the generalized gradient functions $\cG_{ \width }\colon \R^{ \fd } \to \R^{ \fd }$, $ \width \in \N $, ensures that only the bias parameters but not the weight parameters are moved in the training process. We would like to point out that, even though only the bias parameters are trained in \cref{theorem_intro_2}, there are still a large number of \emph{non-global local minimum points in the optimization landscape} in the setup of \cref{theorem_intro_2}.

The stochastic processes $\Theta^{ \width } = ( \Theta^{ \width, 1 }, \dots, \Theta^{ \width, \fd_\width } ) \colon [0, \infty) \times \Omega \to \R^{ \fd_{ \width } }$ for $ \width \in \N $ in \cref{theorem_intro_2} are the randomly initialized solutions of the \GF\ differential equations in \cref{theorem_intro_2:eq:assgf} and the real numbers $ \alpha \in ( \nicefrac{ 3 }{ 4 }, 1 ) $, $ \beta \in (\alpha + 2, \infty ) $ in \cref{theorem_intro_2} are two parameters for the distribution of the \GF\ trajectories $ ( \Theta^{ \width } )_{t \in [0,\infty) } $ at initial time $ t = 0 $. In \cref{eq:theorem_2_intro_statement} in \cref{theorem_intro_2} we establish with high probability that the risk of the \GF\ trajectory $\loss_{ \width }( \Theta^{ \width }_t ) $ converges with the convergence rates $ \kappa $ and $ 1 $ to zero as the training time $t$ and the width of the \ANN\ $ \width $ converge to infinity. In analogy to traditional numerical analysis, say, numerical discretization analyses for partial differential equations (PDEs), the convergence rate $ \kappa $ in \cref{eq:theorem_2_intro_statement} can be regarded as a \emph{spatial approximation order} and the convergence rate $ 1 $ in \cref{eq:theorem_2_intro_statement} can be regarded as a \emph{temporal approximation order}.

In the next result, \cref{theorem_intro} below, we establish a similar but in a certain sense more restrictive result than \cref{theorem_intro_2} above for the \ReLU\ activation function $ \R \ni x \mapsto \max\{ x, 0 \} \in \R $ (see \cref{eq:theorem_intro_risk_function} below) instead of for the clipping activation function $\R \ni x \mapsto \max\{ \min\{ x, 1 \}, 0\} \in \R$ (see \cref{theorem_intro_2:eq:risk} above).

\begin{samepage}
\cfclear
\begin{theorem}\label{theorem_intro}
Let $\0 \in \R$, $\1 \in (\max\{\0, 0\}, \infty)$, let $f \in C^2(\R, \R)$, $\fp \in C(\R, \R)$ be piecewise analytic, assume for all $x \in (\0, \1)$ that $\min\{f'(x), \allowbreak f''(x)\} \allowbreak > f(\0) = 0 \le \fp(x)$ and $\fp^{-1}(\R \backslash \{0\}) = (\0, \1)$, let $(\Omega, \cF, \P)$ be a probability space, for every $\fh \in \N$ let $\dimension_{\fh} = 3 \fh + 1$, let $\loss_{\fh} \colon \R^{\dimension_{\fh}} \to \R$ satisfy for all $\theta = (\theta_1, \ldots, \theta_{\dimension_{\fh}}) \in \R^{\dimension_{\fh}}$ that
\begin{equation}\label{eq:theorem_intro_risk_function}
\textstyle \loss_{\fh}(\theta) = \int_{\R} (f(x) - \theta_{\dimension_{\fh}} - \sum_{j = 1}^{\fh} \theta_{2 \fh + j} \max\{\theta_{\fh + j} x + \theta_{j}, 0\})^2 \fp(x) \, \d x,
\end{equation}
let $\cG_{\fh} = (\cG_{\fh}^1, \ldots, \cG_{\fh}^{\dimension_{\fh}}) \colon \R^{\dimension_{\fh}} \to \R^{\dimension_{\fh}}$ satisfy for all $k \in \{1, 2, \ldots, \dimension_{\fh}\}$, $\theta = (\theta_1, \ldots, \theta_{\dimension_{\fh}}) \in \R^{\dimension_{\fh}}$ that $\cG_{\fh}^{k}(\theta) = (\frac{\partial^{-}}{\partial \theta_k} \loss_{\fh})(\theta) \mathbbm{1}_{\N \cap (0, \fh]}(k)$, and let $\Theta^{\fh} = (\Theta^{\fh, 1}, \ldots, \Theta^{\fh, \dimension_{\fh}}) \colon [0, \infty) \times \Omega \to \R^{\dimension_{\fh}}$ be a stochastic process with continuous sample paths, assume for all $\fh \in \N$, $t \in [0, \infty)$ that 
\begin{equation}
\label{eq:def_GF_dynamics2}
\textstyle \P\bigl(\Theta_t^{\fh} = \Theta_0^{\fh} - \int_0^{t} \cG_{\fh}(\Theta_s^{\fh}) \, \d s\bigr) = 1,
\end{equation}
assume for all $\fh \in \N$ that $\Theta_0^{\fh, 1}, \ldots, \Theta_0^{\fh, \dimension_{\fh}}$ are independent, let $\alpha \in (0, 1)$, $\beta \in (0, \nicefrac{\alpha}{3})$, assume for all $\fh \in \N$, $k \in \N \cap (0, \fh] \cup \{\dimension_{\fh}\}$ that $\fh^{\alpha} \Theta_0^{\fh, k}$ is standard normal, and assume for all $\fh \in \N$, $k \in \N \cap (\fh, 3 \fh]$, $x \in \R$ that $\P(\fh^{\beta} \Theta_0^{\fh, k} \le x) = \int_{0}^{x} [\frac{2}{\pi}]^{1/2} \exp(-\frac{y^2}{2}) \, \d y$. Then there exist $\fc, \degerror \in (0, \infty)$ and random variables $\fC_{\fh} \colon \Omega \to \R$, $\fh \in \N$, such that for all $\fh \in \N$ it holds that $\limsup_{t \to \infty} \E[\loss_{\fh}(\Theta_t^{\fh})] \leq \fc \fh^{- \degerror}$ and
\begin{equation}
\textstyle \P\bigl(\Forall t \in (0, \infty) \colon \loss_{\fh}(\Theta_t^{\fh}) \le \fc \fh^{- \degerror} + \fC_{\fh} t^{-1}\bigr) \ge 1 - \fc \fh^{- \degerror}.
\end{equation}
\end{theorem}
\end{samepage}

\cref{theorem_intro} is a direct consequence of \cref{cor:main} in \cref{sec:ANNs_with_ReLU_final_result} below and in \cref{sec:anns_with_relu} we present a detailed proof of \cref{cor:main}. The statement of \cref{theorem_intro} is similar to the statement of \cref{theorem_intro_2} but \cref{theorem_intro} considers the \ReLU\ activation function instead of the clipping activation function and somehow imposes more restrictive assumptions. In particular, we note that in \cref{theorem_intro} \emph{only the inner bias parameters} (the $ \width $ bias parameters between the $1$-dimensional input layer and the $ \width $-dimensional hidden layer) but not the outer bias parameter (the bias parameter between the $ \width $-dimensional hidden layer and the $ 1 $-dimensional output layer) are trained while in \cref{theorem_intro_2} \emph{all bias parameters} are trained.

A result similar to \cref{theorem_intro_2} and \cref{theorem_intro}, respectively, 
has recently been established independently in Gentile \& Welper~\cite{GentileWelper2022} 
where only the inner bias parameters of a shallow neural network with one-dimensional input and \ReLU\ activation are trained by \GF, where the inner weight parameters are assumed to be $1$, where the outer weight parameters are initialized randomly and uniformly in $\{ - 1 , 1 \}$,
 and where the true risk with respect to a uniform input distribution is considered. 
While the results in \cite{GentileWelper2022} are 
regarding the probability distribution of the input data 
and the values of the inner weight parameters
more restrictive than ours, the target function in \cite{GentileWelper2022} 
is not assumed to be strictly increasing as in 
\cref{theorem_intro_2} and \cref{theorem_intro} 
and the article \cite{GentileWelper2022} even establishes exponential convergence of the \GF\ with high probability. 
The arguments of proof of the main convergence result in \cite{GentileWelper2022} 
follow a completely different approach than the arguments in the proofs of \cref{theorem_intro_2} and \cref{theorem_intro}. 
Specifically, the main convergence result in \cite{GentileWelper2022} is proved by analyzing the neural tangent kernel operator 
on the infinite-dimensional Hilbert space $ L^2( [ -1, 1 ] ) $ and various Sobolev subspaces. 

The results from \cite{GentileWelper2022} have recently been generalized in \cite{Welper2024,Welper2024_JML}
to the setting of deep \ANNs\ in a situation where only the second but last layer is trained while the remaining parameters are not changed after initialization.
These works also allow a multidimensional input domain,
but require a suitable coercivity property of the neural tangent kernel which has so far only been verified numerically.

Let us also comment on some of the shortcomings of \cref{theorem_intro_2} and \cref{theorem_intro}, respectively. 
First, we would like to point out that \cref{theorem_intro_2} and \cref{theorem_intro} only establish convergence to zero for the time-continuous \GF\ dynamics (see \cref{theorem_intro_2:eq:assgf} and \cref{eq:def_GF_dynamics2}) -- which is, of course, not implementable on a computer -- but not for \emph{time-discrete \GD\ optimization methods}, in particular, not for the practically relevant \emph{stochastic \GD\ method}. Moreover, instead of merely using a large number of input-output data pairs, \cref{theorem_intro_2} and \cref{theorem_intro} also exploit the explicit knowledge of the probability distribution of the input data in the gradient dynamics (see \cref{theorem_intro_2:eq:risk}--\cref{theorem_intro_2:eq:assgf} and \cref{eq:theorem_intro_risk_function}--\cref{eq:def_GF_dynamics2}). 
However, in our subjective point of view, the key difficulty in the mathematical analysis of the training of \ANNs\ is the issue to overcome with high probability all bad non-global local minimum points which is already fully present in the case of time-continuous \GF\ processes and we expect that it will be possible to extend \cref{theorem_intro_2} and \cref{theorem_intro}, respectively, to cover also these more practically relevant situations  by using, among other things, the findings in Dereich \& Kassing~\cite{DereichKassing2024}. 
Another key shortcoming of \cref{theorem_intro_2} and \cref{theorem_intro}, respectively, is the fact that only the (inner) bias parameters of the \ANNs\ are trained but not the weight parameters and this shortcoming seems indeed to be extremely hard to relax. Nonetheless, we would like to point out that the optimization landscapes associated to \cref{theorem_intro_2} and \cref{theorem_intro} still admit a large number of bad non-global local minimum points whose risk values are far away from the optimal risk value of the global minimum points. 

A further restriction in \cref{theorem_intro_2} and \cref{theorem_intro} is the assumption of a piecewise analytic target function and input density.
It should be pointed out that this assumption already provides a generalization compared to our previous works \cite{EberleJentzenRiekertWeiss2023,JentzenRiekert2021ExistenceGlobMin}, where we assumed that both these functions are piecewise polynomial.
One might wonder whether one can combine our results with approximation theory for continuous functions by polynomials or analytic functions to obtain convergence also for more general target functions.
However, this seems to be difficult to achieve since changing the target function or input density also changes the entire risk function and, consequently, the \GF\ dynamics.
In particular, the Kurdyka-\L ojasiewicz inequality, on which our proofs rely, would not necessarily hold in that case.
It is known that in the absence of this inequality the \GF\ might not converge, even for smooth objective functions (cf., e.g., \cite{AbsilMahonyAndrews2005}).
Although such an example is, to the best of our knowledge, currently not know for the training of \ANNs, it is not inconceivable that this can happen for certain target functions without structural restrictions similar to the ones in our theorem.

Furthermore, we assume in \cref{theorem_intro_2} and \cref{theorem_intro} that the target function $f$ is strictly increasing and, in the situation of \cref{theorem_intro}, also strictly convex.
We conjecture that it may be possible to relax these assumptions, for instance, using the fact that for a piecewise analytic target function $f$ the input interval can always be divided into a finite number of subintervals on each of which $f$ is either strictly increasing, strictly decreasing or constant. One could then try to separately analyze in each of these regions how the activity intervals at the limiting critical point behave. This leads to some technical difficulties at the boundaries of each subinterval, but might otherwise be doable by employing similar arguments as in this work. Numerical simulation results in \cref{subsec_num_sim_ReLU_all_params_1_hidden_layer_random_normal} and \cref{subsec_num_sim_ReLU_all_params_2_hidden_layers_random_normal} also strongly suggest that one can obtain similar results as in \cref{theorem_intro}, without assuming the monotonicity and the convexity of the target function $f$, by initializing $\Theta$ standard random normally instead of as in \cref{theorem_intro}, and by training all \ANN\ parameters instead of training only inner biases.

An additional shortcoming of \cref{theorem_intro_2} and \cref{theorem_intro}, respectively,
is the assumption of a one-dimensional input domain, and this assumption also seems to be difficult to relax with the current method of proof.
Indeed, our proof relies on the analysis of \emph{activity intervals}, i.e., the respective intervals where each neuron has a non-constant output.
For a multi-dimensional input these intervals would have to be replaced by corresponding polyhedral subsets of the input domain, which would make the corresponding arguments much more technically involved.

Nevertheless, we emphasize that the main contribution of our manuscript is to establish convergence of the risk value not via convergence to global minima, but instead via convergence to \emph{good critical points}.
In fact, in our accompanying article \cite{JentzenRiekert2024} we demonstrate that gradient-based optimization methods do (with high probability) indeed not converge to global minima in the \ANN\ optimization landscape; except for the overparameterized regime, where the target function can be represented exactly by an \ANN.
We emphasize that the results of \cite{JentzenRiekert2024} hold in a much more general setting allowing, in particular, more realistic situations where the input domain is high-dimensional and all parameters of the considered \ANNs\ are trained.
In the additional follow-up work \cite{HannibalJentzenThang2024} this non-convergence result has been further generalized to the case of deep \ANNs\ with an arbitrary number of hidden layers and an arbitrary input dimension.
In this sense, overcoming all bad local minima and proving convergence to good critical points is therefore the best result one can generally hope for.
It is precisely the subject of this work to rigorously establish such convergence results in the considered one-dimensional \ANN\ optimization problem, but we are confident that analogous results can be shown also in more general multi-dimensional setups.
Roughly speaking, our approach in this work is to simplify the training scenario for \ANNs\ in a way
\begin{itemize}
\item so that, on the one hand, the optimization problem becomes significantly \emph{simplified} which enabled us to establish \cref{theorem_intro_2} and \cref{theorem_intro} and
\item so that, on the other hand, the optimization problem is still \emph{difficult} enough to provide a complicated optimization landscape with a large number of bad non-global local minimum points.
\end{itemize}
We conjecture that the central reasons/arguments for the convergence results in \cref{theorem_intro_2} and \cref{theorem_intro} developed in \cref{sec:anns_with_clipping} (in the situation of \cref{theorem_intro_2}) and \cref{sec:anns_with_relu} (in the situation of \cref{theorem_intro}) -- with high probability \emph{overcoming all bad non-global local minimum points} and \emph{convergence to good critical points}; see \cref{intro_item_I,intro_item_II} above -- is not only the reason for the convergence in the simplified shallow \ANN\ training situations in \cref{theorem_intro_2} and \cref{theorem_intro} but is also the reason for the success of \GD\ optimization methods in the usual non-simplified training situations of \ANNs. 
In particular, in \cref{sec:numerical_simulations} we present a series of numerical simulation results for the \emph{standard stochastic \GD\ optimization method} as well as for the popular \emph{\Adam\ optimizer} (see \cite{DBLP:journals/corr/KingmaB14}) for several \ANN\ architectures (including shallow \ANNs\ with just one hidden layer and deep \ANNs\ with two or three hidden layers) and each of these numerical simulations suggests that with high probability the considered \GD\ optimization method \emph{overcomes all bad non-global local minimum points}, \emph{does not converge to a global minimum point}, but \emph{does converge to a good non-optimal critical point} whose risk value is very close to the optimal risk value of the global minimum points.

\textbf{Literature overview} -- Let us also add some comments on the literature on \SGD\ optimization methods as well as on the training of \ANNs. In the situation of a convex objective function, convergence of \GD\ and \SGD\ type optimization methods to the global minimum point with explicit convergence rates can be proved; see, e.g.,~\cite{BachMoulines2013,JentzenKuckuckNeufeldVonWurstemberger2021,BachMoulines2011,Nesterov2015,Nesterov2004}. For non-convex objective functions, generally one cannot even ensure convergence of the iterates of \GD\ methods; see Absil et al.~\cite{AbsilMahonyAndrews2005} for a counterexample. Nevertheless, in the case where the objective function satisfies a {\L}ojasiewicz type inequality the findings, e.g., in \cite{AbsilMahonyAndrews2005,AttouchBolte2009,BolteDaniilidis2006,DereichKassing2024} show convergence to a critical point for \GF\ and \GD\ type processes under the assumption that the process does not diverge to infinity. In the context of \ANNs, such convergence results have also been established in \cite{Davis2018stochastic,EberleJentzenRiekertWeiss2023,JentzenRiekert2021ExistenceGlobMin} and in this work we use similar arguments based on a {\L}ojasiewicz inequality for the risk function as in Eberle et al.~\cite{EberleJentzenRiekertWeiss2023}. However, the critical limit points in the above mentioned references might well be bad saddle points or bad local minimum points whose risk values are far away from the optimal risk values of the global minimum points. Thus, in order to establish convergence to a \emph{good} (non-optimal) critical point one also needs to employ more specific arguments related to the structure of the risk landscape in the training of \ANNs\ as performed in this work in the simplified \ANN\ training situations in \cref{theorem_intro_2} and \cref{theorem_intro}.

It should also be noted that establishing non-divergence or boundedness of \GF\ and \GD\ processes in the training of \ANNs, which is necessary to apply {\L}ojasiewicz type arguments, seems to be a difficult problem (cf.\ Eberle et al.~\cite{EberleNormalizedGF2022}). In particular, as the counterexamples in Gallon et al.~\cite{GallonJentzen2022} show, such a proof can only be possible under certain regularity assumptions on the target function. In the setting considered in this article, we can show boundedness of the considered \GF\ trajectories since only the bias parameters are trained.

In general, the problem of \emph{ruling out convergence to saddle points} seems to be more accessible than the problem of overcoming bad local minimum points. Indeed, sufficient conditions which ensure that the convergence of \GD\ type optimization methods to saddle points can be excluded have been revealed, e.g., in~\cite{LeePanageasRecht2019,LeeJordanRecht2016,PanageasPiliouras2017,Panageas2019firstorder}. Concretely, these results show that \GD\ processes can (up to a Lebesgue zero set of initial values) not converge to \emph{strict} saddles, i.e., saddle points at which the Hessian matrix of the objective function admits at least one strictly negative eigenvalue. They thus require the objective function to be at least twice continuously differentiable. For a similar result regarding time-continuous \GF\ processes we refer to Bah et al.~\cite{BahRauhutTerstiege2021}. In fact, a key step in our proofs consists in verifying the necessary regularity conditions and applying \cite[Theorem~28]{BahRauhutTerstiege2021} in our setting to exclude convergence of \GF\ trajectories to certain bad saddle points. In the context of \ANNs, the fact that \GD\ processes avoid certain saddle points also has recently been established in a special situation in Cheridito et al.~\cite{FlorianEscapeSaddles2022}. More generally, the results in Zhang et al.~\cite{Zhang2022embedding} suggest that strict saddle points appear regularly in the risk landscape related to \ANNs.

On the other hand, with high probability \emph{ruling out convergence to bad local minimum points} remains a challenging problem in essentially all practically relevant cases. A particular setup where positive results have been established is the so-called overparametrized regime. In this situation one considers the empirical risk, which is computed with respect to a finite set of input-output data pairs, and assumes that the number of neurons is much larger than the number of data pairs. Under such a condition, convergence of \GD\ type optimization methods to a global minimum point with high probability has been demonstrated, e.g., in~\cite{AroraDuHuLiWang2019,DuZhaiPoczosSingh2018,EMaWu2020,JentzenKroeger2021,LiLiang2019,RotskoffEijnden2018,ZhangMartensGrosse2019} for the case of shallow \ANNs\ and, e.g., in~\cite{AllenzhuLiLiang2019,AllenzhuLiSong2019,DuLeeLiWangZhai2019,SankararamanDeXuHuangGoldstein2020,ZouCaoZhouGu2019} for the case of deep \ANNs.
The limiting behavior of \GD\ processes associated to neural networks as the width increases to infinity has also been investigated, e.g., in Chizat \& Bach~\cite{ChizatBach2018,ChizatBach2020}, Jacot et al.~\cite{JacotGabrielHongler2020}, and Wojtowytsch~\cite{Wojtowytsch2020}.

Some further progress has been made in the case of training \ANNs\ for particularly simple target functions. Specifically, in the case of a constant target function the articles \cite{Riekert2021ConvergenceConstTF,RiekertJentzenConvergenceConstTargetfn2021,DNNReLUarXiv} establish convergence of \GF\ and \GD\ processes to a global minimum point. In fact, although the risk function is not convex one can show that all critical points are global minimum points in this simplified situation. In the case of an affine linear one-dimensional target function a complete characterization of all local minimum points and saddle points has been established in Cheridito et al.~\cite{Cheridito2022Landscape}.

A further approach is to show that certain global minima of the risk function satisfy suitable regularity properties and conclude convergence of \GD\ optimization methods provided that the initial value is contained in a suitable region of attraction. In order to obtain convergence with high probability, a large number of independent random initializations is then required. Such local convergence results can, e.g., be found in Chatterjee~\cite{Chatterjee2022}, Fehrman et al.~\cite{FehrmanGessJentzen2020}, and \cite{JentzenRiekert2022Piecewise,JentzenRiekert2021ExistenceGlobMin}. By contrast, it is a key strength of the results in this article that \emph{a single random initialization is sufficient to guarantee convergence}.

Finally, more detailed overviews and further references regarding gradient-based optimization methods can be found, e.g., in Bercu \& Fort~\cite{BercuFort2013}, Bottou et al.~\cite{Bottou2018optimization}, E et al.~\cite{EMaWojtowytschWu2020}, Fehrman et al.~\cite[Section~1.1]{FehrmanGessJentzen2020}, \cite[Section~1]{JentzenKuckuckNeufeldVonWurstemberger2021}, Ruder~\cite{Ruder2017overview}, and the references mentioned therein.

The remainder of this work is organized as follows. In \cref{sec:anns_with_clipping} we study the training of \ANNs\ with the clipping activation function and present the detailed proof of \cref{theorem_intro_2} above and \cref{theo:main:conv} below, respectively. 
In \cref{sec:anns_with_relu} we study the training of \ANNs\ with the \ReLU\ activation function and provide the detailed proof of \cref{theorem_intro} above and \cref{cor:main} below, respectively. Finally, in \cref{sec:numerical_simulations} we present several numerical simulation results which complement the analytical findings of this work.

%==============================================================================%
%------------------------------------------------------------------------------%
%=================================-----Section-----============================%
%------------------------------------------------------------------------------%
%==============================================================================%
\section{Artificial neural networks (ANNs) with clipping activation}
\label{sec:anns_with_clipping}

In this section we study \ANNs\ with the clipping activation function.
The main result of this section is \cref{theo:main:conv} below, which in particular implies \cref{theorem_intro_2} from the introduction. 
Broadly speaking, the proof of \cref{theo:main:conv} consists of three parts:

\begin{enumerate} [label = (\Roman*)]
	\item We show in \cref{cor:limit_point} that the considered \GF\ process converges to a critical point and that the risk values converge with rate $\cO ( t^{-1} )$. We also analyze properties of such critical points of the risk function.
	\item Secondly, we establish in \cref{theo:conv:bias:deterministic} an upper bound for the risk value at the limit point, provided that the weight parameters satisfy suitable estimates.
	\item Thirdly, we prove in \cref{lem:random:init} that the required conditions on the weights in \cref{theo:conv:bias:deterministic} are satisfied with high probability for our random initialization.
\end{enumerate}

\subsection{The risk function and its derivatives}

We first introduce the \ANNs\ with clipping activation we consider throughout this section. Afterwards we establish some basic regularity properties of the risk function.

\subsubsection{Mathematical setup for bias training of clipping ANNs}

In \cref{setting:ann:clip} we describe the mathematical framework 
which we frequently employ throughout this section.

\begin{setting} 
\label{setting:ann:clip}
Let $ \0 \in \R $, $ \1 \in (\0 , \infty) $, $ \width \in \N $, $v_1, v_2, \dots, v_{ \width }, w_1, w_2, \dots, w_\width \in (0, \infty)$,
$f, \dens, \fc \in C(\R, \allowbreak \R)$ satisfy for all $x \in \R$ that $\dens(x) \ge 0$, $ \dens^{-1}( (0, \infty) ) = (\0, \1) $, and 
$\fc( x ) = \min\allowbreak\cu{\max \allowbreak \cu{x, \allowbreak 0}, \allowbreak 1}$, let $\Lip \colon C(\R, \R) \to [0, \infty]$ satisfy for all $F \in C(\R, \R)$ that
\begin{equation}
\textstyle 
\Lip(F) = \sup_{x, y \in [\0, \1], x \neq y}\frac{\abs{F(x) - F(y)}}{\abs{x - y}},
\end{equation}
for every $i \in \{ 1, 2, \dots, \width \}$ let $\psi_i \colon \R^{\width + 1} \to \R$ satisfy for all $\theta = (\theta_1, \ldots, \theta_{\width + 1}) \in \R^{\width + 1}$ that $\psi_i(\theta) = \allowbreak - [ w_i ]^{ - 1 } \allowbreak \theta_i$, for every $\theta = (\theta_1, \ldots, \theta_{\width + 1}) \in \R^{ \width + 1 }$, $i \in \{1, 2, \dots, \width\}$ let $I_i^{\theta} \subseteq \R$ satisfy $I_i^{\theta} \allowbreak = ( \psi_i(\theta), \allowbreak \psi_i(\theta) \allowbreak + [ w_i ]^{ - 1 } ) \cap ( \0, \1 )$, for every $\theta = (\theta_1, \dots, \theta_{ \width + 1 } ) \in \R^{ \width + 1}$ let $\realization{ \theta } \colon \R \to \R$ satisfy for all $x \in \R$ that
\begin{equation}
\label{eq:realization_function_ANNs_clipping}
\textstyle \realization{ \theta }( x ) = \theta_{ \width + 1 } + \sum_{ i = 1 }^{ \width } v_i \fc( w_i x + \theta_i ),
\end{equation}
and let 
$\loss \colon \R^{\width + 1} \to \R$ satisfy for all $\theta \in \R^{\width + 1}$ that
\begin{equation}
\label{eq:risk_function_ANNs_clipping}
\textstyle \loss( \theta ) = \int_{ \0 }^{ \1 } (\realization{ \theta }(x) - f(x))^2 \dens(x) \, \d x.
\end{equation}
\end{setting}

In the following we add some brief explanatory comments 
regarding some of the mathematical objects 
appearing in \cref{setting:ann:clip}. 
As in the introduction, 
we also consider 
in \cref{setting:ann:clip} 
\emph{shallow \ANNs} with $ 1 $ neuron on the input layer, 
with $ \width \in \N $ neurons on the hidden layer, 
and with $ 1 $ neuron on the output layer 
(see \cref{figure_shallow_nns_general_introduction}). 
The function 
$ \fc \colon \R \to \R $
in \cref{setting:ann:clip} 
is used as a short notation for the clipping function 
\begin{equation}
  \R \ni x \mapsto \min\{ \max\{ x, 0 \}, 1 \} \in \R 
\end{equation}
and serves as the activation function 
of the considered \ANNs. 
In \cref{setting:ann:clip} 
the \emph{inner weight parameters} 
of the considered \ANNs\ are fixed positive real numbers 
denoted by 
$ w_1, w_2, \dots, w_{ \width } \in (0,\infty) $, 
the \emph{outer weight parameters} of the considered \ANNs\ are fixed positive real numbers 
denoted by 
$ v_1, v_2, \dots, v_{ \width } \in (0,\infty) $, 
and the $ \width $ \emph{inner bias parameters} 
and the \emph{outer bias parameter} are trainable parameters 
and thus appear as arguments 
$ 
  \theta_1, \theta_2, \dots, \theta_{ \width + 1 } 
  \in \R 
$ 
of the risk function 
$  
  \loss = 
  ( \loss( \theta_1, \dots, \theta_{ \width + 1 } ) )_{ 
    ( \theta_1, \dots, \theta_{ \width + 1 } ) 
    \in \R^{ \width + 1 } 
  } 
  \colon \R^{ \width + 1 } \to \R 
$  
specified 
in \cref{eq:risk_function_ANNs_clipping}.

For every \ANN\ parameter vector 
$ 
  \theta = ( \theta_1, \dots, \theta_{ \width + 1 } ) \in \R^{ \width + 1 } 
$ 
we have that the function 
$ \realization{ \theta } \colon \R \to \R $ 
in \cref{eq:realization_function_ANNs_clipping}
specifies the \emph{realization function} 
of the considered \ANN. 
For every \ANN\ parameter vector 
$ \theta = ( \theta_1, \dots, \theta_{ \width + 1 } ) \in \R^{ \width + 1 } $
and every $ i \in \{ 1, 2, \dots, \width \} $ 
we have that 
the $ i $-th summand in the sum on the right hand side 
of \cref{eq:realization_function_ANNs_clipping} 
coincides with the function 
\begin{equation}
  \R \ni x \mapsto 
  v_i \fc( w_i x + \theta_i )
  \in \R
\end{equation}
and specifies the contribution 
of the $ i $-th neuron on the hidden layer. 
Observe in \cref{setting:ann:clip} 
that for every \ANN\ parameter vector $\theta = ( \theta_1, \dots, \theta_{ \width + 1 } ) \in \R^{ \width + 1 } $ and every $ i \in \{ 1, 2, \dots, \width \} $, $ x \in \R $ we have that $  0 < \fc( w_i x + \theta_i ) < 1$ if and only if 
\begin{equation}
\textstyle \psi_i( \theta ) = - \theta_i [w_i]^{-1} < x < [w_i]^{-1}  - \theta_i [w_i]^{-1} = \psi_i( \theta ) + [w_i]^{- 1}
  .
\end{equation}
For every \ANN\ parameter vector 
$ \theta = ( \theta_1, \dots, \theta_{ \width + 1 } ) \in \R^{ \width + 1 } $
and every $ i \in \{ 1, 2, \dots, \width \} $
we thus have that the set 
$
  I^{ \theta }_i \subseteq \R
$ 
in \cref{setting:ann:clip} 
describes precisely 
the largest open subset of 
$ ( \0, \1 ) $ 
on which the function 
\begin{equation}
  ( \0, \1 ) \ni x \mapsto 
  v_i \fc( w_i x + \theta_i )
  \in \R
\end{equation}
is not constant 
and, in this sense, 
the intervals 
$
  I^{ \theta }_i \subseteq \R
$ 
for 
$
  i \in \{ 1, 2, \dots, \width \}
$
and 
$ \theta \in \R^{ \width + 1 } $ 
may be referred as 
\emph{activity intervals}.

In the next subsection, 
\cref{ssub:activity_intervals}, 
we collect a few elementary properties 
of such activity intervals. 
Thereafter, in \cref{ssub:differentiability_risk_clipping} 
we provide explicit representations 
for the first and second order partial derivatives of the risk function 
$
  \loss \colon \R^{ \width + 1 } \to \R
$ 
from \cref{setting:ann:clip}.

\subsubsection{Properties of activity intervals}
\label{ssub:activity_intervals}
\begin{prop}
\label{prop:active:intervals}
Assume \cref{setting:ann:clip} and let 
$ 
  \theta = ( \theta_1, \dots, \theta_{ \width + 1 } ) \in \R^{ \width + 1 } 
$.  
Then
\begin{enumerate}[label = (\roman*)]
\item
\label{prop:active:intervals:item1} 
it holds 
for all $ i \in \{ 1, 2, \dots, \width \} $
that 
$
  I_i^{ \theta } 
  = \cu{ x \in ( \0, \1 ) \colon 0 < w_i x + \theta_i < 1 } 
$,

\item
\label{prop:active:intervals:item2} 
it holds 
for all $ i \in \{ 1, 2, \dots, \width \} $
that 
$
  ( 
    ( 
      I_i^{ \theta } \neq \emptyset 
    )
    \Longleftrightarrow 
    ( 
      \0 - [ w_i ] ^{ - 1 } < \psi_i( \theta ) < \1 
    ) 
  )
$, 
and 

\item
\label{prop:active:intervals:item3} 
it holds for all 
$ \cJ \subseteq \{ 1, 2, \ldots, \width \} $, 
$ x, y \in ( \cap_{ i \in \cJ } I_i^{ \theta } ) $ 
with $ y - x \ge 0 $ that
\begin{equation}
\textstyle 
  \realization{ \theta }( y ) 
  -
  \realization{ \theta }( x ) 
  \ge 
  ( \sum_{ i \in \cJ } v_i w_i ) 
  ( y - x ) 
  .
\end{equation}
\end{enumerate}
\end{prop}

\begin{cproof}{prop:active:intervals}
\Nobs that the fact that 
for all 
$ i \in \{ 1, 2, \dots, \width \} $
it holds that 
\begin{equation}
  I_i^{ \theta } 
  = 
  ( 
    - \theta_i [ w_i ]^{ - 1 }, 
    ( 1 - \theta_i ) [ w_i ]^{ - 1 } 
  ) 
  \cap ( \0, \1 )
\end{equation}
establishes \cref{prop:active:intervals:item1,prop:active:intervals:item2}. 
Next \nobs that the fact that 
$ \fc \colon \R \to \R $ 
is non-decreasing assures 
that for all 
$ \cJ \subseteq \{ 1, 2, \dots, \width \} $, 
$ x, y \in ( \cap_{ i \in \cJ } I_i^{ \theta } ) $ 
with $ x \le y $ 
it holds that
\begin{equation}
\begin{split}
\textstyle \realization{\theta} ( y ) - \realization{\theta}(x) 
& \textstyle = \br[\big]{\theta_{ \width + 1 } + \ssum_{j=1}^\width v_j \fc ( w_j y + \theta_j)} - \br[\big]{\theta_{ \width + 1 } +  \ssum_{j=1}^\width v_j \fc ( w_j x + \theta_j)} \\
& \textstyle = \sum_{j = 1}^\width v_j \br{\fc ( w_j y + \theta_j ) - \fc ( w_j x + \theta_j ) } \\
& \textstyle \ge \sum_{j \in \cJ } v_j \br{\fc ( w_j y + \theta_j ) - \fc ( w_j x + \theta_j)}
\\
& \textstyle = \sum_{j \in \cJ} v_j \br{(w_j y + \theta_j ) - ( w_ j x + \theta_j)} = \rbr[\big]{\sum_{j \in \cJ }  v_j w_j } ( y - x ).
\end{split}
\end{equation}
This establishes \cref{prop:active:intervals:item3}.
\end{cproof}

\subsubsection{First and second order differentiability properties of the risk function}
\label{ssub:differentiability_risk_clipping}

We next show in \cref{prop:risk:gradient} that the considered risk function is once continuously differentiable. The proof is similar to our previous works \cite{ShokhrukhLocalMinima,JentzenRiekert2022Piecewise} and uses the facts that the clipping activation $\fc$ is almost everywhere differentiable and that the weights are nonzero.

\begin{prop}[First order partial derivatives of the risk function] 
\label{prop:risk:gradient}
Assume \cref{setting:ann:clip}. Then
\begin{enumerate} [label = (\roman*)]
\item
\label{prop:risk:gradient:item1} 
it holds that 
$ \loss \in C^1( \R^{ \width + 1 }, \R ) $,

\item 
\label{prop:risk:gradient:item2} 
it holds for all 
$ 
  \theta = ( \theta_1, \dots, \theta_{ \width + 1 } ) \in \R^{ \width + 1 } 
$, 
$ 
  i \in \{ 1, 2, \dots, \width \}
$ 
that
\begin{equation} \label{prop:risk:gradient:item2:eq}
\begin{split}
\textstyle 
  ( \frac{ \partial }{ \partial \theta_i } \loss )( \theta ) 
& 
  \textstyle 
  = 2 v_i \int_{ \0 }^{ \1 } 
  ( \realization{ \theta }( x ) - f( x ) ) 
  \dens(x) 
  \indicator{ (0, 1) }( w_i x + \theta_i ) 
  \, \d x 
\\
& \textstyle 
  = 2 v_i \int_{I_i^\theta} (\realization{\theta} (x) - f(x)) \dens(x) \, \d x
  ,
\end{split}
\end{equation}
and
\item
\label{prop:risk:gradient:item3} 
it holds for all 
$ 
  \theta = ( \theta_1, \dots, \theta_{ \width + 1 } ) \in \R^{ \width + 1 } 
$
that
\begin{equation} \label{prop:risk:gradient:item3:eq}
\textstyle 
  ( 
    \frac{\partial}{\partial \theta_{\width + 1}} 
    \loss 
  )( \theta ) 
  = 2 \int_{\0}^{\1} (\realization{\theta}(x) - f(x)) \dens(x) \, \d x.
\end{equation}
\end{enumerate}
\end{prop}
\begin{cproof}{prop:risk:gradient}
\Nobs that for all 
$ \theta = ( \theta_1, \dots, \theta_{ \width + 1 } ) \in \R^{ \width + 1 } $, 
$ i \in \{ 1, 2, \dots, \width \} $, 
$ 
  x \in [\0, \1] \allowbreak \backslash 
  \cu{ \psi_i (\theta), \allowbreak \psi_i (\theta) \allowbreak + [w_i]^{-1}} 
$ 
it holds that $\R \ni \vartheta \mapsto (\realization{(\theta_1, \ldots, \theta_{i-1}, \vartheta, \theta_{i+1}, \ldots, \theta_{\width + 1})}(x) - f(x))^2 \dens(x) \in \R$ is differentiable at $\theta_i$ with
\begin{equation}
\textstyle \frac{\partial}{\partial \theta_{i}} \rbr[\big]{(\realization{\theta}(x) - f(x))^2 \dens (x)}
= 2 v_i (\realization{\theta}(x) - f(x)) \dens(x) \indicator{(0, 1)} (w_i x + \theta_i).
\end{equation}
Combining this with \cite[Corollary~2.3]{ShokhrukhLocalMinima} establishes \cref{prop:risk:gradient:item2:eq}. We obtain \cref{prop:risk:gradient:item3:eq} in the same way. Furthermore, \nobs that, e.g., \cite[Lemma 2.6]{JentzenRiekert2022Piecewise} implies for all $i \in \{1, 2, \ldots, \width + 1\}$ that $\R^{\width + 1} \allowbreak \ni \theta \mapsto \rbr{\frac{\partial}{\partial \theta_i} \loss}(\theta) \in \R$ is continuous.
\end{cproof}

Next we show that the risk function is even twice continuously differentiable. This follows from the Leibniz integral rule and the assumption that $\dens^{-1} ( \R \backslash \cu{0} ) = ( \0 , \1 )$, which means that no regularity problems can occur at the boundary.

\begin{prop}[Second order partial derivatives of the risk function]  
\label{prop:risk:hess}
Assume \cref{setting:ann:clip}. Then
\begin{enumerate} [label = (\roman*)]

\item
\label{item1:prop:risk:hess} 
it holds that 
$ \loss \in C^2( \R^{ \width + 1 }, \R ) $,

\item
\label{item2:prop:risk:hess} 
it holds for all $\theta = ( \theta_1, \dots, \theta_{ \width + 1 } ) \in \R^{ \width + 1 } $, $i, j \in \{ 1, 2, \ldots, \width \}$ that
\begin{equation}
\begin{split}
\textstyle (\frac{ \partial^2 }{ \partial \theta_i \partial \theta_j } \loss )( \theta ) & \textstyle = 2 v_i v_j \int_{ \R } \dens(x) \indicator{ (0, 1) }( w_i x + \theta_i ) 
  \indicator{ (0, 1) }( w_j x + \theta_j ) 
  \, \d x 
\\
& \textstyle \quad - 2 v_i (w_i)^{- 1} [\mathbbm{1}_{\{i\}}(j)] \br[\big]{( \realization{ \theta }(x) - f(x) ) \dens(x)}_{x = \min\cu{\max\cu{\0, \psi_i(\theta)}, \1}}^{x = \min \cu{\max \cu{\0, \psi_i(\theta) + [w_i]^{-1}}, \1}} \\
& \textstyle = 2 v_i v_j \int_{ I_i^{ \theta } \cap I_j^{ \theta } } \dens(x) \, \d x - 2 v_i ( w_i )^{- 1} [\mathbbm{1}_{\{i\}}(j)] \br[\big]{(\realization{ \theta }( x ) - f( x ) ) \dens(x)}_{x = \psi_i( \theta ) }^{ x = \psi_i( \theta ) + [ w_i ]^{ - 1 } },
\end{split}
\end{equation}

\item
\label{item3:prop:risk:hess} 
it holds for all 
$ 
  \theta = ( \theta_1, \dots, \theta_{ \width + 1 } ) \in \R^{ \width + 1 } 
$, 
$ i \in \{ 1, 2, \dots, \width \} $ 
that
\begin{equation}
\textstyle 
  ( 
    \frac{ \partial^2 }{ \partial \theta_i \partial \theta_{ \width + 1 } } \loss 
  )( \theta ) 
  = 2 v_i 
  \int_{ \0 }^{ \1 } \dens(x) \indicator{ (0, 1) }( w_i x + \theta_i ) \, \d x
  = 2 v_i 
  \int_{ I_i^{ \theta } } 
%   \indicator{ (0, 1) }( w_i x + \theta_i ) 
  \dens(x) \, \d x,
\end{equation}
and

\item
\label{item4:prop:risk:hess} 
it holds for all 
$ 
  \theta = ( \theta_1, \dots, \theta_{ \width + 1 } ) \in \R^{ \width + 1 } 
$
that
\begin{equation}
  \textstyle 
  ( \frac{ \partial^2 }{ \partial \theta_{ \width + 1}^2 } \loss )( \theta ) 
  = 2 \int_{ \0 }^{ \1 } \dens(x) \, \d x .
\end{equation}
\end{enumerate}
\end{prop}
\begin{cproof}{prop:risk:hess}
\Nobs that the assumption that $ \dens^{ - 1 }( (0, \infty) ) = (\0, \1) $, 
\cref{prop:risk:gradient}, 
\cite[Corollary~2.3]{ShokhrukhLocalMinima}, 
and the Leibniz integral rule establish \cref{item1:prop:risk:hess,item2:prop:risk:hess,item3:prop:risk:hess,item4:prop:risk:hess}.
\end{cproof}

\subsubsection{Properties of critical points of the risk function}

In this part we analyze activity intervals of critical points of the risk function. The main result is \cref{cor:non:strict:saddle} below, which establishes a pointwise estimate for the slope of the realization function at non-descending critical points (cf.\ \cref{def:strict_saddle}). This will later be used in the convergence analysis, since convergence to descending critical points can be shown to occur with probability zero.

\begin{lemma}[Overlapping activity intervals of critical points must be nested]
\label{lem:crit:intervals:intersect}
Assume \cref{setting:ann:clip}, assume $\Lip(f) < \min_{i \in \{1, 2, \ldots, \width\}} v_i w_i$, let $\theta \in \R^{\width + 1}$ satisfy $(\nabla \loss) \allowbreak (\theta) \allowbreak = 0$, and let $i, j \in \{1, 2, \ldots, \width\}$ satisfy 
\begin{equation}
\label{eq:assumption_overlapping}
  ( I_i^{\theta} \cap I_j^{\theta} ) \neq \emptyset
  .
\end{equation}
Then 
$
%   ( 
    ( I_i^{ \theta } \subseteq I_j^{ \theta } ) 
    \vee 
    ( I_j^{ \theta } \subseteq I_i^{ \theta } )
%   )
$.
\end{lemma}
\begin{cproof}{lem:crit:intervals:intersect}
Throughout this proof let $q^i, r^i, q^j, r^j \in \R$ satisfy $I_i^{\theta} = (q^i, r^i)$ and $I_j^{\theta} = (q^j , r^j)$. 
\Nobs that for all $x, y, z, t \in \R$ with $(x - z)(y - t) \le 0$ and $(x, y) \cap (z, t) \neq \emptyset$ it holds that 
\begin{equation}
  [(x, y) \subseteq (z, t)] \vee [(z, t) \subseteq (x, y)]
  .
\end{equation}
Hence, for the sake of contradiction 
we assume without loss of generality that $q^i < q^j$ and $r^i < r^j$. 
\Nobs that \cref{prop:risk:gradient} and 
the assumption that $(\nabla \loss)(\theta) = 0$ imply that
\begin{equation} \label{lem:crit:intervals:intersect:eq1}
\textstyle \int_{q^i}^{r^i } ( \realization{\theta} ( x ) - f ( x ) ) \dens ( x ) \, \d x = 
\int_{q^j}^{r^j } ( \realization{\theta} ( x ) - f ( x ) ) \dens ( x ) \, \d x = 0 .
\end{equation}
\Nobs that \cref{eq:assumption_overlapping}
% the assumption that $I_i^\theta \cap I_j^\theta \not= \emptyset$ 
and the assumption that 
$
  (q^i < q^j) \wedge (r^i < r^j)
$ 
imply that $q^j < r^i$. 
In addition, \cref{lem:crit:intervals:intersect:eq1} 
and the mean value theorem ensure that 
there exists $u \in (q^j , r^j)$ which satisfies 
\begin{equation}
  \realization{\theta}(u) = f(u)
  .
\end{equation}
This, the assumption that $\Lip ( f ) < \min_{i \in \{1, 2, \ldots, \width\}} v_i w_i$, and \cref{prop:active:intervals} show for all $x \in (q^i, u)$ that $\realization{\theta}(x) < f(x)$. Combining this with \cref{lem:crit:intervals:intersect:eq1} ensures that $u < r^i$. Hence, we obtain for all $x \in (r^i , r^j )$ that $\realization{\theta}(x) > f(x)$. \cref{lem:crit:intervals:intersect:eq1} therefore demonstrates that
\begin{equation}
\begin{split}
0 & \textstyle < \int_{ r^i}^{r^j} (\realization{\theta}(x) - f(x)) \dens(x) \, \d x = - \int_{q^j}^{r^i} (\realization{\theta}(x) - f(x)) \dens(x) \, \d x \\
& \textstyle = \int_{q^i}^{q^j} (\realization{\theta}(x) - f(x)) \dens(x) \, \d x < 0,
\end{split}
\end{equation}
which is a contradiction.
\end{cproof}

\cfclear
\begin{definition}[Descending critical point]\label{def:strict_saddle}
Let $ \fh \in \N $, $ f \in C^2( \R^{ \fh }, \R ) $, $ x \in \R^{\fh} $. Then we say that $ x $ is a \descritic critical point of $ f $ if and only if it holds that 
\begin{enumerate}[label = (\roman*)]
\item we have that 
$
  ( \nabla f )( x ) = 0
$
and 
\item 
we have that there exists $ v \in \R^{ \fh } $ such that
\begin{equation}
\textstyle 
%   ( \nabla f )( x ) = 0 
%   \qqandqq 
%   \Exists v \in \R^{\fh} \colon 
  \scalar{v, ((\operatorname{Hess} f)(x))v} < 0.
\end{equation}
\end{enumerate}
\end{definition}

In the situation of \cref{def:strict_saddle} 
\nobs that 
the assumption that $ f \in C^2( \R^\width , \R) $ 
ensures that 
the matrix $ ( \operatorname{Hess} f)( x ) $ is symmetric and therefore diagonalizable with real eigenvalues. 
In the context of \cref{def:strict_saddle} we hence obtain 
that $ x $ is a \descritic critical point of $ f $ if and only if 
$ x $ is a critical point with the property that 
$ ( \operatorname{Hess} f)( x ) $ 
admits at least one negative eigenvalue.

In the scientific literature descending critical points are often called strict saddle points (cf., e.g., \cite[Definition 27]{BahRauhutTerstiege2021} and \cite[Definition 1]{LeePanageasRecht2019}).
However, note that a descending critical point in the sense of \cref{def:strict_saddle} is not necessarily a saddle point, but might as well be a local maximum.

\cfclear
\begin{lemma}[Upper bounds for the slopes of realization 
functions of suitable non-descending critical points 
under local regularity assumptions on the density]
\label{lem:crit:dens:bounded:below}
Assume \cref{setting:ann:clip}, 
assume $ \Lip( f ) < \min_{ i \in \{ 1, 2, \dots, \width \} } v_i w_i $, 
let $ \theta \in \R^{ \width + 1 } $ 
satisfy 
$
  ( \nabla \loss ) \allowbreak ( \theta ) \allowbreak = 0 
$, 
assume that $ \theta $ is not a \descritic critical point of $ \loss $ 
\cfadd{def:strict_saddle}, 
let 
$
  \scrx \in ( \0, \1 ) 
$, 
and assume that there exists 
$ 
  i \in \{ 1, 2, \dots, \width \} 
$ 
such that 
\begin{equation} 
\textstyle 
  \scrx \in I_i^{\theta} 
\qqandqq
  \inf_{ y \in I_i^{ \theta } } 
  \dens(y) 
  \ge (\min_{j \in \{1, 2, \ldots, \width\}} w_j)^{-1} \Lip(\dens)
\end{equation}
\cfload. 
Then
\begin{equation}
\textstyle 
\sum_{j \in \{ 1, 2, \dots, \width \}, \, \scrx \in I_j^{\theta}} v_j w_j \le 3 \bigl[\max_{j \in \{ 1, 2, \dots, \width \}, \, \scrx \in I_j^{\theta}} v_j w_j \bigr].
\end{equation}
\end{lemma}
\begin{cproof}{lem:crit:dens:bounded:below}
Throughout this proof assume that
\begin{equation}
\label{lem:crit:dens:bounded:below:eq1}
\textstyle 
\sum_{ j \in \{ 1, 2, \dots, \width \}, \, \scrx \in I_j^{ \theta }} v_j w_j > 3 \bigl[\max_{j \in \{ 1, 2, \dots, \width \}, \, \scrx \in I_j^{ \theta } } v_j w_j\bigr]. 
\end{equation}
\Nobs that \cref{lem:crit:intervals:intersect} ensures 
that there exists $ i \in \{ 1, 2, \dots, \width \} $ 
which satisfies $\scrx \in I_i^\theta$, $\inf_{y \in I_i^\theta} \dens(y) \ge [w_i]^{-1} \Lip(\dens) > 0$, and which satisfies for all $j \in \{1, 2, \ldots, \width\}$ with $\scrx \in I_j^{\theta}$ that $I_i^{\theta} \subseteq I_j^{\theta}$. Let $q, r \in \R$ satisfy $I_i^{\theta} = (q, r)$. \Nobs that \cref{prop:risk:hess} implies that
\begin{equation}\label{eqn:lem:crit:dens:bounded:below:parderi2_ii}
\textstyle \rbr[\big]{\frac{\partial^2}{\partial \theta_i^2} \loss}(\theta) = 2 v_i^2 \int_q^r \dens(x) \, \d x - 2 \frac{v_i}{w_i} \br{( \realization{\theta}(x) - f(x)) \dens(x)}_{x=q}^{x=r}.
\end{equation}
Furthermore, the assumption that 
$ \Lip(f) < \min_{ j \in \{ 1, 2, \dots, \width \} } v_j w_j $, 
\cref{prop:active:intervals}, and \cref{lem:crit:dens:bounded:below:eq1} show for all $x, y \in [q, r]$ with $x < y$ that $( \realization{\theta}(y) - f(y)) - (\realization{\theta}(x) - f(x)) > 2 v_i w_i (y - x) > 0 $. This and the fact that $\int_q^r ( \realization{\theta}(x) - f(x)) \dens(x) \, \d x  = 0$ assure that $\realization{\theta}(r) - f(r) > 0 > \realization{\theta}(q) - f(q)$. Therefore, we obtain that
\begin{equation}\label{eqn:lem:crit:dens:bounded:below:int_bound}
\begin{split}
\textstyle \br*{(\realization{\theta}(x) - f(x)) \dens(x)}_{x=q}^{x=r} & \textstyle \ge \rbr[\big]{(\realization{\theta}(r) - f(r)) - (\realization{\theta}(q) - f(q))} \inf\nolimits_{x \in [q,r] } \dens(x) \\
& > 2 v_i w_i (r - q) \inf\nolimits_{x \in [q,r]} \dens(x).
\end{split}
\end{equation}
Moreover, the fact that $0 < r - q \le (-[w_i]^{-1} \theta_i + [w_i]^{-1}) - (- [w_i]^{-1} \theta_i) = [w_i]^{-1}$ and the fact that $\inf_{x \in I_i^\theta} \dens(x) \ge [w_i]^{-1} \Lip(\dens)$
prove that
\begin{equation} 
\begin{split}
\textstyle \sup_{x \in [q,r] } \dens ( x ) 
& \textstyle \le  \br[\big]{\inf_{x \in [q, r]} \dens(x)} + (\Lip(\dens))(r - q) \le \br[\big]{\inf_{x \in [q, r]} \dens(x)} + (\Lip(\dens)) [w_i]^{-1} \\
& \textstyle \le 2 \inf_{x \in [q, r]} \dens(x).
\end{split}
\end{equation}
Combining this with \cref{eqn:lem:crit:dens:bounded:below:parderi2_ii} and \cref{eqn:lem:crit:dens:bounded:below:int_bound} ensures that
\begin{equation}
\rbr[\big]{ \tfrac{ \partial ^2 }{ \partial \theta_i ^2 } \loss } ( \theta )
< \tfrac{2v_i}{w_i} \rbr*{ v_i w_i ( r - q ) \br[\big]{ \sup\nolimits_{x \in [q,r]} \dens ( x ) } - 2 v_i w_i ( r - q ) \br[\big]{ \inf\nolimits_{x \in [q,r]} \dens ( x ) } } \le 0.
\end{equation}
This shows that $\theta$ is a \descritic critical point of $\loss$, which is a contradiction.
\end{cproof}

\begin{lemma} 
\label{lem:zero:mean:help}
Let $ a \in \R $, $ b \in (a, \infty) $, 
let $ \dens \in C ( [a, b] , [0, \infty) ) $ be strictly increasing, 
let $ g \in C( [a, b], \R ) $, $ L \in (0, \infty) $ satisfy 
for all $ x, y \in [a,b] $ with $ x < y $ that 
\begin{equation}\label{eqn:lem:zero:mean:help}
\textstyle 
%   \Forall x \in [a, b], y \in (x, b] \colon 
  g(y) - g(x) > 3 L (y - x),
\end{equation}
and assume $\int_a^b g(x) \dens(x) \, \d x = 0$.
Then 
\begin{equation}
\label{lem:zero:mean:help:eqclaim}
\textstyle L \int_a^b \dens(x) \, \d x < g(b) \dens(b).
\end{equation}
\end{lemma}
\begin{cproof}{lem:zero:mean:help}
\Nobs that \cref{eqn:lem:zero:mean:help} ensures that $g$ is strictly increasing. Combining this with the assumption that $\int_a^b g ( x ) \dens ( x ) \, \d x = 0$ shows that there exists $\xi \in (a,b)$ which satisfies $g(a) < 0 = g (\xi) < g(b)$. \Nobs that \cref{eqn:lem:zero:mean:help} shows that $g(b) > 3L(b - \xi)$.

In the following we distinguish between the case $g(b) \ge - g(a)$ and the case $g(b) < - g(a)$.
 First we prove \cref{lem:zero:mean:help:eqclaim} in the case 
 \begin{equation}
 \label{lem:zero:mean:help:case1}
 g(b) \ge - g(a).
 \end{equation}
  \Nobs that \cref{eqn:lem:zero:mean:help} and the fact that $g(a) < 0 = g (\xi) < g(b)$ assure that $g(a) < 3 L (a - \xi)$. Combining this with \cref{lem:zero:mean:help:case1} shows that $g(b) > 3 L (\xi - a)$. This, the fact that $g(b) > 3L(b - \xi)$, and the assumption that $\dens$ is strictly increasing demonstrate that
\begin{equation}
\begin{split}
\textstyle g(b) \dens(b) & \textstyle = 2 \frac{g(b) \dens(b)}{2} > \frac{3L}{2} (\xi - a) \dens(b) + \frac{3L}{2} (b - \xi) \dens(b) \\
& \textstyle \ge  L (b - a) \dens(b) > L \int_a^b \dens(x) \, \d x.
\end{split}
\end{equation}
This establishes \cref{lem:zero:mean:help:eqclaim} in the case $g(b) \ge - g(a)$.
 Next we prove \cref{lem:zero:mean:help:eqclaim} in the case 
 \begin{equation}
 \label{lem:zero:mean:help:case2}
 g(b) < - g(a).
 \end{equation} 
 \Nobs that \cref{lem:zero:mean:help:case2} and the intermediate value theorem assure that there exists $\zeta \in (a, \xi)$ which satisfies $g(\zeta) = - g(b)$. Moreover, \cref{eqn:lem:zero:mean:help} shows that $g(b) = g(\xi) - g(\zeta) > 3 L(\xi - \zeta)$. The fact that for all $x \in (a, \zeta )$,
$y \in (\zeta , \xi )$ it holds that $g(x) < g(\zeta ) = - g(b) = g(\zeta) < g(y) < 0$, the fact that $g$ is strictly increasing, and the assumption that $\int_a^b g(x) \dens(x)\, \d x = 0$ therefore demonstrate that
\begin{equation}
\textstyle g(b) \int_a^\zeta \dens(x) \, \d x \textstyle \le \int_a^\zeta(- g(x)) \dens(x) \, \d x = \int_\zeta^b g(x) \dens(x) \, \d x \le \int_\xi^b g(x) \dens(x) \, \d x \le g(b) \int_\xi^b \dens(x) \, \d x.
\end{equation}
Hence, we obtain that $\int_a^\zeta \dens(x) \, \d x \le \int_\xi ^b \dens ( x ) \, \d x$. This and the fact that $g ( b ) >3 L \max \cu{ \xi - \zeta , b - \xi }$ assure that
\begin{equation}
\begin{split}
\textstyle g(b) \dens(b) & \textstyle > 2 L \dens(b) (b - \xi) + L \dens(b) (\xi - \zeta) \ge 2 L \int_\xi^b \dens(x) \, \d x + L \int_\zeta^\xi \dens(x) \, \d x \\
& \textstyle \ge L \int_\xi^b \dens(x) \, \d x + L \int_a^\zeta \dens(x) \, \d x + L \int_\zeta^\xi \dens(x) \, \d x = L \int_a^b \dens(x) \, \d x .
\end{split}
\end{equation}
This establishes \cref{lem:zero:mean:help:eqclaim} in the case $g(b) < - g(a)$.
\end{cproof}

\begin{cor}
\label{cor:zero:mean:help}
Let $a \in \R$, $b \in (a, \infty )$, let $\dens \in C( [a,b] , [0, \infty ) )$ be strictly monotonous, 
let $g \in C([a, b] , \R )$, $L \in (0, \infty)$ satisfy for all 
$ x, y \in [a,b] $
with $ x < y $ that
\begin{equation}
  g(y) - \allowbreak g(x) \allowbreak > 3 L (y - x)
  ,
\end{equation}
and assume $\int_a^b g(x) \dens(x) \, \d x = 0$. Then 
\begin{equation}
\label{cor:zero:mean:help:eqclaim}
\textstyle L \int_a^b \dens(x) \, \d x < g(b) \dens(b) - g(a) \dens(a).
\end{equation}
\end{cor}
\begin{cproof}{cor:zero:mean:help}
Throughout this proof let 
$
  \overline{\dens}, \overline{g} \colon [a, b] \to \R
$ 
satisfy for all $x \in [a, b]$ that 
$\overline{\dens} ( x ) = \dens ( a + b - x )$ and $\overline{g} ( x ) = - g (a + b - x ) $.

In the following we distinguish between the case where $ \dens $ is strictly increasing 
and the case where $ \dens $ is strictly decreasing. 
First we prove \cref{cor:zero:mean:help:eqclaim} in the case where $ \dens $ is strictly increasing. 
\Nobs that \cref{lem:zero:mean:help} and the fact that $ g(a) \dens(a) \le 0 $ show that 
$ L \int_a^b \dens(x) \, \d x < g(b) \dens(b) - g(a) \dens(a) $, 
which establishes \cref{cor:zero:mean:help:eqclaim}. 
 
 Next we prove \cref{cor:zero:mean:help:eqclaim} in the case where $\dens$ is strictly decreasing. \Nobs that \cref{lem:zero:mean:help} (applied with $\dens \with \overline{\dens}$, $g \with \overline{g}$ in the notation of \cref{lem:zero:mean:help}) and the integral transformation theorem demonstrate that
\begin{equation}
\textstyle L \int_a^b \dens(x) \, \d x = L \int_a^b \overline{\dens} (x) \, \d x < \overline{g}(b) \overline{\dens}(b) = - g(a) \dens(a) \le g(b) \dens(b) - g(a) \dens(a).
\end{equation}
\end{cproof}

\cfclear
\begin{lemma}[Upper bounds for the slopes of realization 
functions of suitable non-descending critical points 
under local regularity assumptions on the density]
\label{lem:crit:dens:mono} 
Assume \cref{setting:ann:clip}, assume $\Lip(f) < \min_{i \in \{1, 2, \ldots, \width\}} v_i w_i$, let $\theta \in \R^{\width + 1}$ satisfy $(\nabla \loss) \allowbreak(\theta) \allowbreak = 0$, let $ \scrx \in (\0, \1) $, assume that there exists 
$ i \in \{ j \in \{ 1, 2, \dots, \width \} \colon \scrx \in I_j^\theta \} $ 
such that $\dens|_{I_i^\theta }$ is strictly monotonous, and assume that $\theta$ is not a \descritic critical point of $\loss$ \cfadd{def:strict_saddle}\cfload. Then
\begin{equation}
\textstyle \sum_{j \in \{1, 2, \ldots, \width\}, \, \scrx \in I_j^{\theta}} v_j w_j \le 4 \max_{j \in \{1, 2, \ldots, \width\}, \, \scrx \in I_j^{\theta}} v_j w_j.
\end{equation}
\end{lemma}
\begin{cproof}{lem:crit:dens:mono}
Throughout this proof assume
\begin{equation}
\label{lem:crit:dens:mono:eq1}
\textstyle \sum_{j \in \{1, 2, \ldots, \width\}, \, \scrx \in I_j^{\theta}} v_j w_j > 4 \max_{j \in \{1, 2, \ldots, \width\}, \, \scrx \in I_j^{\theta}} v_j w_j.
\end{equation}
\Nobs that \cref{lem:crit:intervals:intersect} ensures that there exists $i \in \{1, 2, \ldots, \width\}$ which satisfies for all $j \in \{1, \allowbreak 2, \allowbreak \ldots, \allowbreak \width\}$ with $\scrx \in I_j^{\theta}$ that $\scrx \in I_i^\theta \subseteq I_j^\theta$ and $\dens |_{I_i^\theta }$ 
is strictly monotonous. Let $q, r \in \R$ satisfy $I_i^\theta = (q, r)$. \Nobs that \cref{prop:risk:hess} implies that
\begin{equation}
\label{lem:crit:dens:mono:eq2}
\textstyle \rbr[\big]{ \tfrac{ \partial ^2 }{ \partial \theta_i ^2 } \loss } ( \theta ) = 2 v_i ^2 \int_q^r \dens ( x ) \, \d x - 2 \tfrac{v_i}{w_i} \br*{ ( \realization{\theta} ( x ) - f ( x ) ) \dens ( x ) } _{x=q}^{x=r}.
\end{equation}
Furthermore, \cref{lem:crit:dens:mono:eq1} ensures for all $x \in [q,r]$, $y \in (x, r]$ that $( \realization{\theta} ( y ) - f ( y ) ) - ( \realization{\theta} ( x ) - f ( x ) ) > 3 v_i w_i ( y - x ) $. This, the fact that $\int_q^r ( \realization{\theta} ( x ) - f ( x ) ) \dens ( x ) \, \d x = 0$, and \cref{cor:zero:mean:help} assure that
\begin{equation}
\textstyle v_i w_i \int_q^r \dens ( x ) \, \d x < \br*{ ( \realization{\theta} ( x ) - f ( x ) ) \dens ( x ) } _{x=q}^{x=r}.
\end{equation}
Combining this with \cref{lem:crit:dens:mono:eq2} demonstrates that 
$\rbr[\big]{ \tfrac{ \partial ^2 }{ \partial \theta_i ^2 } \loss } ( \theta ) < 0$.
This shows that $\theta$ is a \descritic \allowbreak critical point of $\loss$, which is a contradiction.
\end{cproof}

\cfclear
\begin{cor}[Upper bounds for the slopes of realization 
functions of suitable non-descending critical points 
under global regularity assumptions on the density]
\label{cor:non:strict:saddle}
Assume \cref{setting:ann:clip}, assume $\Lip ( f ) < \min_{i \in \{1, 2, \ldots, \width\}} v_i w_i$, let $\theta \in \R^{\width + 1 }$ satisfy $( \nabla \loss ) ( \theta ) \allowbreak = 0$,
 let $ \eps, \allowbreak \delta \in \allowbreak (0, \infty)$ satisfy 
\begin{equation}
\textstyle 
  \eps - \delta > 2 \allowbreak [\min_{j \in \{1, 2, \ldots, \width\}} w_j]^{-1} 
\qqandqq 
  \inf_{x \in [ \0 + \delta , \1 - \delta ] } \dens ( x ) \ge [\min_{j \in \{1, 2, \ldots, \width\}} w_j]^{-1} \allowbreak \Lip( \dens )
  ,
\end{equation}
assume that $\dens_{[\0 , \0 + \varepsilon ] }$ is strictly increasing, assume that $\dens_{[\1 -  \varepsilon  , \1 ] }$ is strictly decreasing,
and assume that $\theta$ is not a \descritic critical point of $\loss$ \cfadd{def:strict_saddle}\cfload. Then it holds for all $x \in ( \0 , \1 )$ that
\begin{equation}\label{eqn:cor:non:strict:saddle}
\textstyle 
\sum_{j \in \{1, 2, \ldots, \width\}, \, x \in I_j^{\theta}} v_j w_j \le 4 \max_{j \in \{1, 2, \ldots, \width\}, \, x \in I_j^{\theta}} v_j w_j.
\end{equation}
\end{cor}
\begin{cproof}{cor:non:strict:saddle}
Throughout this proof let $\scrx \in (\0 , \1 )$, $i \in \{1, 2, \ldots, \width\}$ satisfy $\scrx \in I_i^\theta$ and let $\scrw \in \R$ satisfy $\scrw = [\min_{j \in \{1, 2, \ldots, \width\}} w_j]^{-1}$. In the following we distinguish between the case $\scrx \in [ \0 + \delta + \scrw, \1 - \delta - \scrw]$ and the case $\scrx \in (\0 ,  \0 + \delta + \scrw) \cup (\1 - \delta - \scrw, \1 )$. First we prove \cref{eqn:cor:non:strict:saddle} in the case $\scrx \in [\0 + \delta + \scrw, \1 - \delta - \scrw]$. In this case, we have that $I_i^\theta \subseteq [\0 + \delta , \1 - \delta]$. This and the assumption that $\inf_{x \in [ \0 + \delta , \1 - \delta ] } \dens(x) \ge \scrw \Lip(\dens)$ ensure that $\inf_{x \in I_i^\theta } \dens(x) \ge \scrw \Lip(\dens)$. Combining this with \cref{lem:crit:dens:bounded:below} shows that 
\begin{equation}
\textstyle \sum_{j \in \{1, 2, \ldots, \width\}, \, \scrx \in I_j^{\theta}} v_j w_j \le 3 \max_{j \in \{1, 2, \ldots, \width\}, \, \scrx \in I_j^{\theta}} v_j w_j \le 4 \max_{j \in \{1, 2, \ldots, \width\}, \, \scrx \in I_j^{\theta}} v_j w_j.
\end{equation}
This establishes \cref{eqn:cor:non:strict:saddle} in the case $\scrx \in [\0 + \delta + \scrw, \1 - \delta - \scrw]$. In the next step we prove \cref{eqn:cor:non:strict:saddle} in the case $\scrx \in (\0 ,  \0 + \delta + \scrw) \cup (\1 - \delta - \scrw, \1 )$. \Nobs that the assumption that $\eps - \delta > 2 \scrw$ assures that $I_i^\theta \allowbreak \subseteq \allowbreak (\0, \allowbreak \0 + \varepsilon) \cup (\1 - \varepsilon, \1)$. \cref{lem:crit:dens:mono} therefore demonstrates that $\sum_{j \in \{1, 2, \ldots, \width\}, \, \scrx \in I_j^{\theta}} v_j w_j \le 4 \max_{j \in \{1, 2, \ldots, \width\}, \, \scrx \in I_j^{\theta}} v_j w_j$. This establishes \cref{eqn:cor:non:strict:saddle} in the case $\scrx \in (\0 ,  \0 + \delta + \scrw) \cup (\1 - \delta - \scrw, \1 )$.
\end{cproof}

\subsubsection{Pointwise estimates for realization functions at critical points}

We next derive in \cref{lem:realization:pointwise:bound} an estimate for the pointwise distance between the realization function $\cN^\vartheta$ and the target function $f$ at a critical point $\vartheta$. If the outer weights $v_i$ are chosen to be small this will result in a useful estimate for the risk value $\loss ( \vartheta )$.
Due to the fact that $\cN ^\vartheta$ is constant outside the activity intervals $I_i^\vartheta$ it suffices to prove this inside these intervals, the bound at the remaining points then follows from the monotonicity of $f$.

\begin{lemma}
\label{lem:realization:pointwise:bound}
Assume \cref{setting:ann:clip}, assume that $f$ is non-decreasing and satisfies $\Lip ( f ) < \min_{i \in \{1, 2, \ldots, \width\}} v_i w_i$, let $\vartheta \in \R^{\width + 1 }$ satisfy $(\nabla \loss)(\vartheta) = 0$, and let $V \in \R$ satisfy 
\begin{equation}
\textstyle V = \sup_{x \in (\0, \1)} \sum_{j \in \{1, 2, \ldots, \width\}, \, x \in I_j^{\theta}} v_j.
\end{equation}
Then
\begin{enumerate}[label = (\roman*)]
\item
\label{lem:realization:pointwise:bound:item1} it holds for all $j \in \{1, 2, \ldots, \width\}$, $x \in I_j^\vartheta$ that $\abs{\realization{\vartheta} ( x ) - f ( x ) } \le V$, 

\item
\label{lem:realization:pointwise:bound:item2} it holds for all $j , k \in \{1, 2, \ldots, \width\}$ with $I_j^\vartheta, I_k^\vartheta \not= \emptyset$ and all $x \in [ \sup I_j^\vartheta, \inf I_k^\vartheta ]$
that $\abs{ \realization{\vartheta} ( x ) - f ( x ) } \le V$,

\item
\label{lem:realization:pointwise:bound:item3} it holds for all $j \in \{1, 2, \ldots, \width\}$ with $I_j^\vartheta \not= \emptyset$ and all
$x \in [ \sup I_j^\vartheta , \1 ]$ that $\abs{\realization{\vartheta} ( x ) - f ( x ) } \le \max \cu{ V , \abs{ f ( \1 ) - \realization{\vartheta} ( \1 ) } }$, and

\item
\label{lem:realization:pointwise:bound:item4} it holds for all $j \in \{1, 2, \ldots, \width\}$ with $I_j^\vartheta \not= \emptyset$ and all
$x \in [ \0 ,\inf I_j^\vartheta ]$ that $\abs{\realization{\vartheta} ( x ) - f ( x ) } \le \max \cu{ V , \abs{ f ( \0 ) - \realization{\vartheta} ( \0 ) } }$.
\end{enumerate}
\end{lemma}
\begin{cproof}{lem:realization:pointwise:bound}
Throughout this proof let $k, l, m \in \{1, 2, \ldots, \width\}$, $q, s \in [\0, \1]$, $r \in [q, \1]$, $\scrx \allowbreak \in \allowbreak (\0, \1)$ satisfy $I_k^\vartheta \neq \emptyset$, $I_l^\vartheta \neq \emptyset$, $I_m^\vartheta \not= \emptyset$, $q = \sup I_k^\vartheta$, $s = \sup I_m^\vartheta$, $r = \inf I_l^\vartheta$, and $\cu{ j \in \{1, \allowbreak 2, \allowbreak \ldots, \allowbreak \width\} \allowbreak \colon \allowbreak \scrx \allowbreak \in \allowbreak I_j^\vartheta} \allowbreak \allowbreak \neq \emptyset$. \Nobs that \cref{lem:crit:intervals:intersect} assures that there exists $\scri \in \{1, 2, \ldots, \width\}$ which satisfies for all $j \in \allowbreak \{1, \allowbreak 2, \allowbreak \ldots, \width\}$ with $\scrx \in I_j^\vartheta$ that $\scrx \in I_{\scri}^\vartheta \subseteq I_j^\vartheta$. Moreover, the integral mean value theorem and the fact that 
\begin{equation}
\textstyle \int_{I_{\scri}^\vartheta} ( \realization{\vartheta} ( x ) - f ( x ) ) \allowbreak \dens ( x ) \, \d x \allowbreak = 0
\end{equation}
demonstrate that there exists $\scry \in I_{\scri}^\vartheta$ which satisfies $\realization{\vartheta} ( \scry ) = f ( \scry )$. In the following we assume without loss of generality that $\scry \le \scrx$ and let $\scru \in \R$ satisfy $\scru = \sup \cu{ z \in I_{\scri}^\vartheta \cap ( \0 , \scrx ] \colon \realization{\vartheta} ( z ) = f(z)}$. \Nobs that
\begin{equation} \label{lem:realization:pointwise:bound:eq1}
\begin{split}
& \textstyle \realization{\vartheta} ( \scrx ) \le \realization{\vartheta} ( \scru ) + \sum_{j \in \{1, 2, \ldots, \width\}, \, I_j^\vartheta \cap (\scru , \scrx ) \neq \emptyset} v_j \\
& \textstyle = f ( \scru ) + \sum_{j \in \{1, 2, \ldots, \width\}, \, I_j^\vartheta \cap (\scru , \scrx ) \neq \emptyset} v_j \le f ( \scrx ) + \sum_{j \in \{1, 2, \ldots, \width\}, \, I_j^\vartheta \cap (\scru , \scrx ) \neq \emptyset} v_j.
\end{split}
\end{equation}
In addition, the fact that $\Lip ( f ) < \min_{j \in \{1, 2, \ldots, \width\}} v_j w_j $, the fact that $\scrx, \scru \in I_{\scri}^\vartheta$, and \cref{prop:active:intervals} show that $\realization{\vartheta} ( \scrx ) \ge f ( \scrx )$. Next let $\scrj \in \{1, 2, \ldots, \width\}$ satisfy $I_{\scrj}^\vartheta \cap ( \scru , \scrx ) \neq \emptyset$. We now prove that $\scru \in I_{\scrj}^\vartheta$. \Nobs that the integral mean value theorem and the fact that
\begin{equation}
\textstyle \int_{I_{\scrj}^\vartheta} ( \realization{\vartheta} ( x ) - f ( x ) ) \allowbreak \dens ( x ) \, \d x \allowbreak = 0
\end{equation}
show that there exists $\scrz \in I_{\scrj}^\vartheta$ which satisfies $\realization{\vartheta} ( \scrz ) = f ( \scrz )$.

In the following we distinguish between the case $\scrz \le \scru$ and the case $\scrz > \scru$. First we prove $\scru \in I_{\scrj}^\vartheta$ in the case $\scrz \le \scru$. \Nobs that the assumption that $I_{\scrj}^\vartheta \cap ( \scru , \scrx ) \neq \emptyset$ and the fact that $I_{\scrj}^\vartheta$ is connected imply that $\scru \in I_{\scrj}^\vartheta$. This proves $\scru \in I_{\scrj}^\vartheta$ in the case $\scrz \le \scru$.

 In the next step we prove $\scru \in I_{\scrj}^\vartheta$ in the case $\scrz > \scru$. \Nobs that the fact that $\scru = \sup \cu{ z \in I_{\scri}^\vartheta \cap ( \0 , \scrx ] \colon \realization{\vartheta} ( z ) = f(z)}$ implies that $\scrz > \scrx$. Hence, we obtain that $\scrx \in I_{\scrj}^\vartheta$. Combining this with the definition of $\scri$ ensures that $I_{\scri}^\vartheta \subseteq I_{\scrj}^\vartheta$. Therefore, we obtain that $\scru \in I_{\scri}^\vartheta \subseteq I_{\scrj}^\vartheta$. This proves $\scru \in I_{\scrj}^\vartheta$ in the case $\scrz > \scru$. Hence, we obtain that
\begin{equation}
\textstyle \sum_{j \in \{1, 2, \ldots, \width\}, \, I_j^\vartheta \cap (\scru , \scrx ) \neq \emptyset} v_j \le \ssum_{j \in \{1, 2, \ldots, \width\}, \, \scru \in I_j^\vartheta} v_j \le V.
\end{equation}
Combining this with \cref{lem:realization:pointwise:bound:eq1} and the fact that $\realization{\vartheta} ( \scrx ) \ge f (\scrx)$ establishes \cref{lem:realization:pointwise:bound:item1}. Next we assume without loss of generality that 
\begin{equation}\label{eqn:lem:realization:pointwise:bound:k_l}
\textstyle \forall \, i \in \{1, 2, \ldots, \width\} \backslash \cu{k, l} \colon I_i^\vartheta \cap [q, r] = \emptyset.
\end{equation}
\Nobs that \cref{eqn:lem:realization:pointwise:bound:k_l} ensures for all $x \in [q, r]$ that $\realization{\vartheta} (x) = \realization{\vartheta}(q) = \realization{\vartheta}(r)$ and $f(q) \le f(x) \le f(r)$. Combining this with \cref{lem:realization:pointwise:bound:item1} assures for all
$x \in [q, r]$ that
\begin{equation}
\begin{split}
\textstyle \abs{\realization{\vartheta} ( x ) - f(x) } & \textstyle \le \max \cu[\big]{\abs{\realization{\vartheta} ( x ) - f(q) } , \abs{\realization{\vartheta} ( x ) - f(r) }} \\
& \textstyle = \max \cu[\big]{\abs{\realization{\vartheta} ( q ) - f(q) } , \abs{\realization{\vartheta} ( r ) - f(r) }}  \le V.
\end{split}
\end{equation}
This establishes \cref{lem:realization:pointwise:bound:item2}. In the next step we assume without loss of generality that $\forall \, i \in \allowbreak \{1, \allowbreak 2, \allowbreak \ldots, \width\} \backslash \cu{ m } \colon I_i^\vartheta \cap [s, \1] = \emptyset$ (otherwise the claim in \cref{lem:realization:pointwise:bound:item3} follows from \cref{lem:realization:pointwise:bound:item1,lem:realization:pointwise:bound:item2}). \Nobs that for all $x \in [s, \1]$ we have that $\realization{\vartheta} ( x ) = \realization{\vartheta} ( s ) = \realization{\vartheta} ( \1 )$ and $f( s ) \le f ( x ) \le f ( \1 )$. Combining this with \cref{lem:realization:pointwise:bound:item1} assures for all $x \in [s , \1 ]$ that 
\begin{equation}
\begin{split}
\textstyle \abs{\realization{\vartheta} ( x ) - f(x) } & \textstyle \le \max \cu[\big]{\abs{\realization{\vartheta} ( x ) - f(s) } , \abs{\realization{\vartheta} ( x ) - f(\1) }} \\
& \textstyle = \max \cu[\big]{\abs{\realization{\vartheta} ( s ) - f(s) } , \abs{\realization{\vartheta} (\1) - f(\1) }}  \le \max \cu{ V , \abs{f(\1) - \realization{\vartheta}(\1)}}.
\end{split}
\end{equation}
This establishes \cref{lem:realization:pointwise:bound:item3}. \Cref{lem:realization:pointwise:bound:item4} can be shown analogously as \cref{lem:realization:pointwise:bound:item3}.
\end{cproof}

\begin{cor} \label{cor:realization:endpoint:bound}
Assume \cref{setting:ann:clip}, assume that $f$ is non-decreasing and satisfies $\Lip ( f ) < \min_{i \in \{1, 2, \ldots, \width\}} v_i w_i$, let $\vartheta \in \R^{\width + 1 }$ satisfy $( \nabla \loss ) ( \vartheta ) = 0$, and let $V \in \R$ satisfy 
\begin{equation}
\textstyle V = \sup_{x \in ( \0 , \1 ) } \sum_{i \in \{1, 2, \ldots, \width\}, \, x \in I_i^\vartheta } v_i.
\end{equation}
Then
\begin{equation} \label{cor:realization:endpoint:bound:eq1}
\max \cu[\big]{ \realization{\vartheta} ( \1 ) - f ( \1 )  , f ( \0 ) - \realization{\vartheta} ( \0 ) } \le V.
\end{equation}
\end{cor} 

\begin{cproof}{cor:realization:endpoint:bound}
In the following we distinguish between the case $\forall \, i \in \{1, 2, \ldots, \width\} \colon \allowbreak I_i^\vartheta \allowbreak = \emptyset$ and the case $\exists \, i \in \{1, 2, \ldots, \width\} \colon I_i^\vartheta \neq \emptyset$. First we prove \cref{cor:realization:endpoint:bound:eq1} in the case $\forall \, i \in \{1, 2, \ldots, \allowbreak \width\} \allowbreak \colon \allowbreak I_i^\vartheta \allowbreak = \emptyset$. \Nobs that for all $x \in (\0 , \1 )$ we have that $\realization{\vartheta} ( x ) = \realization{\vartheta} ( \1 )$. This, \cref{prop:risk:gradient}, and the assumption that $( \nabla \loss ) ( \vartheta ) = 0$ demonstrate that 
\begin{equation}
\textstyle \int_\0^\1 ( \realization{\vartheta} ( \1 ) - f ( x ) ) \dens ( x ) \, \d x = 0.
\end{equation}
The integral mean value theorem therefore assures that there exists $x \in (\0 , \1)$ such that $f(x) = \realization{\vartheta} ( \1 ) = \realization{\vartheta} ( \0 )$. Combining this with the fact that $f$ is non-decreasing shows that $f(\0) \le \realization{\vartheta} ( \0 )$ and $f ( \1 ) \ge \realization{\vartheta} ( \1 )$. This proves \cref{cor:realization:endpoint:bound:eq1} in the case $\forall \, i \in \{1, 2, \ldots, \width\} \colon I_i^\vartheta = \emptyset$.

 Next we prove \cref{cor:realization:endpoint:bound:eq1} in the case $\exists \, i \in \{1, 2, \ldots, \width\} \colon I_i^\vartheta \neq \emptyset$. Let $\scrz \in \R$ satisfy $\scrz = \sup \rbr{\cup_{i \in \{1, 2, \ldots, \width\}} I_i^\vartheta}$. \Nobs that $\realization{\vartheta} ( \1 ) = \realization{\vartheta} (\scrz)$. Combining this with \cref{lem:realization:pointwise:bound} demonstrates that $\realization{\vartheta} ( \1 ) - f ( \1 ) \le \realization{\vartheta} ( \scrz ) - f ( \scrz ) \le V$. Analogously, we obtain that $\realization{\vartheta} ( \0 ) - f ( \0 ) \ge -  V$. This proves \cref{cor:realization:endpoint:bound:eq1} in the case $\exists \, i \in \{1, 2, \ldots, \width\} \colon I_i^\vartheta \neq \emptyset$.
\end{cproof}

\begin{cor}
\label{cor:endpoint:bound:2}
Assume \cref{setting:ann:clip}, assume that $f$ is non-decreasing and satisfies $\Lip ( f ) < \min_{i \in \{1, 2, \ldots, \width\}} v_i w_i$, let $\vartheta \in \R^{\width + 1 }$ satisfy $( \nabla \loss ) ( \vartheta ) = 0$, and let $V \in \R$ satisfy 
\begin{equation}
\textstyle V = \sup_{x \in ( \0 , \1 ) } \sum_{i \in \{1, 2, \ldots, \width\}, \, x \in I_i^\vartheta } v_i.
\end{equation}
Then
\begin{equation}
\max \cu[\big]{ \abs{ \realization{\vartheta} ( \0 ) - f ( \0 ) } , 
\abs{ \realization{\vartheta} ( \1 ) - f ( \1 ) } } \le \max \cu*{ f ( \1 ) - f ( \0 ) , V } .
\end{equation}
\end{cor}

\begin{cproof}{cor:endpoint:bound:2}
\Nobs that \cref{cor:realization:endpoint:bound} ensures that $\realization{\vartheta} ( \1 ) - f ( \1 ) \le V$ and $f(\0 ) - \realization{\vartheta} ( \0 ) \le V$. Next, \cref{prop:risk:gradient}, the fact that $f$ and $\realization{\vartheta}$ are non-decreasing, and the fact that $\rbr[\big]{ \frac{\partial}{\partial \theta_{\width + 1 } } \loss } ( \vartheta ) = 0$ demonstrate that
\begin{equation}
\begin{split}
\textstyle \realization{\vartheta} ( \1 ) \int_\0^\1  \dens ( x ) \, \d x & \textstyle = \int_\0^\1 \realization{\vartheta} ( \1 ) \dens ( x ) \, \d x \ge \int_\0^\1 \realization{\vartheta} ( x ) \dens ( x ) \, \d x = \int_\0^\1 f ( x ) \dens ( x ) \, \d x \\
& \textstyle \ge \int_\0^\1 f ( \0 ) \dens ( x ) \, \d x = f ( \0 ) \int_\0^\1  \dens ( x ) \, \d x.
\end{split}
\end{equation}
Hence, we obtain that $f( \1 ) - \realization{\vartheta} ( \1 ) \le f ( \1 ) - f ( \0 )$. In the same way we can show that $\realization{\vartheta} ( \0 ) - f ( \0 ) \le f ( \1 ) - f ( \0 )$.
\end{cproof}

\subsection{Convergence analysis for gradient flows (GFs)}
\label{subsec:gf_conv}
In this section we study the \GF\ dynamics for the \ANNs\ considered in \cref{setting:ann:clip}. In particular, we show that the \GF\ trajectory always converges to a critical point, and we establish estimates for the risk at this limit point.

\subsubsection{Subanalytic functions and Kurdyka-\L ojasiewicz inequalities}

\newcommand{\ssC}{\mathscr{A}}

In this section we show the Kurdyka-\L ojasiewicz (KL) inequality for the risk function $\loss$ under the assumption that both the target function $f$ and the density $\dens$ are piecewise analytic; see \cref{cor:KL_analytic} below.
For this we need some notions regarding subanalytic sets and functions which can also be found, e.g., in Bolte et al.~\cite[Definition~2.1]{BolteDaniilidis2006} and Shiota~\cite{Shiota1997}.

Our results generalize our previous work \cite{JentzenRiekert2021ExistenceGlobMin}, where a corresponding KL inequality has been established under the more restrictive assumption of a piecewise polynomial target function $f$ and density function $\dens$.

\cfclear
\begin{definition}[Semi-analytic and subanalytic sets]
	\label{def:sub_analytic}
	\begin{enumerate}[label=(\roman*)]
		\item For every $n \in \N$, a set $A \subseteq \R^n$ is called semi-analytic if
		for all $ v \in \R^n $ there exist 
		$M , N \in \N $, 
		an open $ U \subseteq \R^n $, 
		and real analytic functions
		$P_{ i, j, k } \colon U \to \R$,
		$( i, j, k ) \in \cu{1, 2, \ldots, M } \times \cu{1, 2, \ldots, N } \times \cu{ 0, 1 } $
		such that 
		$ v \in U $ and 
		\begin{equation} 
			\label{eq:def_semi_analytic_set}
			\textstyle
			A \cap U = 
			\bigcup_{ i = 1 }^M
			\bigcap_{ j = 1 }^N 
			\rbr[\big]{ \cu*{ 
				x \in U \colon P_{ i, j, 0 }( x ) = 0 } \cap \cu*{ x \in U \colon P_{ i, j, 1 }( x ) > 0 } } .
		\end{equation}
		\item For every $n \in \N$, a set $A \subseteq \R^n$ is called subanalytic if for all $ v \in \R^n $ there exist 
		$ m \in \N $, 
		an open $ U \subseteq \R^n $, 
		and 
		a bounded semi-analytic set
		$ B \subseteq \R^{ n + m } $
		such that 
		$ v \in U $ and 
		\begin{equation} 
			\label{eq:def_sub_analytic_set}
			\textstyle
			A 
			\cap 
			U
			= \{ 
			x \in \R^n \colon 
			(
			\exists \, y \in \R^m \colon
			(x,y) \in B
			) \} .
		\end{equation}
		\item For every $m, n \in \N$, a function $f \colon \R^m \to \R^n$ is called subanalytic if $
		\operatorname{Graph}( f ) = \cu{(x, f(x) ) \colon x \in \R^m } \subseteq \R^{m + n }$
		is a subanalytic set.
	\end{enumerate}
\end{definition}

We now define our notion of one-dimensional piecewise analytic functions, which generalizes the notion of piecewise polynomial functions introduced in \cite[Definition 5.1]{JentzenRiekert2021ExistenceGlobMin}.

\begin{definition}[Piecewise analytic functions]
	\label{def:piecewise_analytic}
	\begin{enumerate}[label=(\Roman*)]
		\item For every $A \subseteq \R$, we say that a function $g \colon A \to \R$ is analytic if there exist an open $U \subseteq \R$ and an analytic function $G \colon U \to \R$
		such that $A \subseteq U$ and $g = G | _A$.
		\item We say that a function $ f \colon \R\to \R$ is piecewise analytic if and only if 
		there exist
		$ K \in \N $, $ \tau_1, \ldots, \tau_K \in \R $
		such that 
		\begin{enumerate}[label=(\roman*)]	
			\item 
			it holds that 
			$  \tau_1 < \tau_2 <  \cdots < \tau_K$,	
			\item 
			it holds for all 
			$ i \in \cu{ 2, 3, \dots, K } $ that 
			$ f|_{ [ \tau_{ i - 1 }, \tau_i ] } $ is a analytic, and
			\item it holds that $f _{ (- \infty, \tau_1  ]}$ and $f _ { [\tau_K, \infty )} $ are analytic. 		
		\end{enumerate}
	\end{enumerate}
\end{definition}

We next define in \cref{def:analytic_function_class} a suitable subclass of subanalytic functions.
\cref{def:analytic_function_class} is inspired by Definition 4.1 in our previous article \cite{JentzenRiekert2021ExistenceGlobMin}, but is in a sense more general since we now allow analytic functions instead of polynomials and rational functions.

\cfclear
\begin{definition}
	\label{def:analytic_function_class}
	Let $n \in \N$.
	Then we denote by
	$\ssC_n \subseteq \cu{f \colon \R^n \times \R \to \R}$ the $\R$-vector space of functions given by
	\begin{equation}\label{eq:def_analytic_function_class}
    \begin{split}
    & \ssC_n = \operatorname{span}_\R \Bigl( f \colon
		\R^n \times \R \ni (\theta , x ) \mapsto f ( \theta , x ) \in \R \colon \\
    & \Bigl[ 
		\Exists N \in \N, s, \alpha_1, \ldots, \alpha_N \in \R, t \in [s , \infty ), \text{polynomials } P_0,  \ldots, P_N \colon \R^n \to \R, \\
		& \text{an analytic } Q \colon [s, t ] \to \R,
		A_1, \ldots, A_N \in \cu{\cu{0}, [0 , \infty ), (- \infty , 0 ] } \colon \Forall \theta \in \R^n, x \in \R \colon \\
		& f ( \theta ,x ) = P_0 ( \theta ) Q ( \min \cu{\max \cu{x, s } , t } ) \indicator{[s, t ] } ( x ) \textstyle\prod_{i=1}^N \indicator{A_i} ( \alpha_i x - P_i ( \theta) ) \Bigr] \Bigr).
	\end{split}
	\end{equation}
\end{definition}

\Nobs that the functions in \cref{def:analytic_function_class} depend on two variables $\theta \in \R^n$ and $x \in \R$.
In the following we think of $\theta$ as the parameter vector of the considered \ANNs\ and of $x$ as the input data.

\cfclear
\begin{lemma}
	\label{lem:risk_function_prop}
	Assume \cref{setting:ann:clip} and assume that $\dens$ and $f$ are piecewise analytic \cfadd{def:piecewise_analytic}\cfload.
	Then $[ \R^{ \width + 1 } \times \R \ni (\theta , x ) \mapsto (\realization{\theta} ( x ) - f ( x ) ) ^2 \dens ( x ) \in \R ] \in \ssC_{\width + 1 }$ \cfadd{def:analytic_function_class}\cfload.
\end{lemma}

\begin{cproof}{lem:risk_function_prop}
	\Nobs that \cref{eq:def_analytic_function_class} assures that $\ssC_{\width + 1 }$ is an algebra.
	Moreover, the assumption that $\dens^{-1} ( \R \backslash \cu{0} ) = ( \0 , \1 )$ ensures for all $\theta \in \R^{ \width + 1 }$, $x \in \R$ that
	\begin{equation}
		 (\realization{\theta} ( x ) - f ( x ) ) ^2 \dens ( x ) =  (\realization{\theta} ( x ) \indicator{[\0, \1]} ( x )  - f ( x ) \indicator{[\0, \1]} ( x )  ) ^2 \dens ( x ).
	\end{equation}
	Furthermore, the assumption that $\dens$ and $f$ are piecewise analytic implies that
	$ [\R^{\width + 1 } \ni (\theta , x ) \mapsto \dens ( x ) \in \R ] \in \ssC_{\width + 1 }$ and 
	$ [\R^{\width + 1 } \ni (\theta , x ) \mapsto f ( x ) \indicator{[\0, \1]} ( x ) \in \R ] \in \ssC_{\width + 1 }$.
	In addition, \nobs that for all $i \in \cu{1, 2, \dots, \width }$, $x \in \R$, $\theta \in \R^{ \width + 1 }$ we have that
	\begin{equation}
		\begin{split}
			 \fc ( w_i x + \theta _i ) 
			&= 
			 (w_i x + \theta_i ) \indicator{[0 , \infty ) } ( w_i x + \theta_i) \indicator{[0 , \infty ) } ( 1 - w_i x - \theta_i )
			+ \indicator{( 0 , \infty ) } (w_i x + \theta_i - 1 ) \\
			&= (w_i x + \theta_i ) \indicator{[0 , \infty ) } ( w_i x + \theta_i) \indicator{[0 , \infty ) } ( 1 - w_i x - \theta_i ) \\
			& \quad + ( 1 - \indicator{ ( - \infty , 0 ] } (w_i x + \theta_i - 1 ) ) .
		\end{split}
	\end{equation}
	This demonstrates for all $i \in \cu{1, 2, \dots, \width }$ that 
	$ [\R^{\width + 1 } \ni (\theta , x ) \mapsto v_i \fc ( w_i x + \theta _i ) \indicator { [\0 , \1 ] } \in \R ] \in \ssC_{\width + 1 }$.
	Finally, \nobs that $[\R^{ \width + 1 } \ni ( \theta , x ) \mapsto \theta_{\width + 1 } \indicator{[\0, \1]} ( x ) ] \in \ssC_{\width + 1 }$.
	Hence, we obtain that $[\R^{ \width + 1 } \ni ( \theta , x ) \mapsto \realization{\theta} ( x )  \indicator{[\0, \1]} ( x ) ] \in \ssC_{\width + 1 }$.
\end{cproof}

\cfclear
\begin{prop}
	\label{prop:integral_prop}
	Let $n \in \N$, $g \in \ssC_n$ \cfadd{def:analytic_function_class}\cfload.
	Then
	\begin{enumerate}[label=(\roman*)]
		\item \label{prop:integral_prop_i}
		 it holds for all $\theta \in \R^n$ that $\int_\R \abs{g(\theta , x)} \, \d x < \infty $ and
		\item \label{prop:integral_prop_ii}
		 the function $\R^n \ni \theta \mapsto \int_\R g(\theta , x) \, \d x \in \R$ is subanalytic
	\end{enumerate}
\cfadd{def:sub_analytic}\cfload. 
\end{prop}
\begin{cproof}{prop:integral_prop}
	We will use the following facts (cf.~Shiota~\cite[Chapter I.2]{Shiota1997}):
	\begin{itemize}
		\item If $f , g \colon \R^d \to \R$ are subanalytic and at least one of them is locally bounded, then $f+g$ is subanalytic.
		\item If $f , g \colon \R^d \to \R$ are subanalytic and locally bounded, then $fg$ is subanalytic.
		\item If $f \colon \R^d \to \R^k$, $g \colon \R^k \to \R$ are subanalytic and $f$ is locally bounded, then $g \circ f$ is subanalytic.
	\end{itemize}
	For every $s, t, x \in \R$ let $c_{s,t}(x) := \min \cu{ \max \cu{x, s}, t}$.
	By linearity of the integral it suffices to consider a function $f$ of the form
	\begin{equation}
		\label{integral_prop_eq:f}
		f ( \theta , x ) =  P_0 ( \theta ) Q ( c_{s, t } ( x ) ) \indicator{[s, t ] } ( x ) \textstyle\prod_{i=1}^N \indicator{A_i} ( \alpha_i x - P_i ( \theta) ) 
	\end{equation}
	for some $N \in \N$, $s , \alpha_1, \ldots, \alpha_N \in \R$, $t \in [s , \infty )$,
	polynomials $P_0, P_1, \ldots, P_N \colon \R^n \to \R$,
	an analytic function $Q \colon [s, t ] \to \R$,
	and sets $A_1, \ldots, A_N \in \cu{\cu{0}, [0 , \infty ), (- \infty , 0 ] }$.
	\Nobs that \cref{integral_prop_eq:f} establishes \cref{prop:integral_prop_i}.
	Next, by linearity and since $\forall \, x \in \R \colon \indicator{[0 , \infty)}(x) = 1 - \indicator{(- \infty , 0 ] } ( x ) $ we may assume for all $i \in \cu{1, 2, \ldots, N }$ that $A_i = [0 , \infty )$.
	Let $I, J, K \subseteq \cu{1, 2, \ldots, N }$ satisfy
	\begin{multline}
		I = \cu{i \in \cu{1, 2, \ldots, N } \colon \alpha_i > 0},
		\quad J = \cu{i \in \cu{1, 2, \ldots, N } \colon \alpha_i < 0}, \\
		\qandq K = \cu{i \in \cu{1, 2, \ldots, N } \colon \alpha_i = 0} .
	\end{multline} 
	\Nobs that for all $\theta \in \R^n$ it holds that
	\begin{equation}
		\begin{split}
		\int_\R f(\theta , x ) \, \d x 
		&= P_0 (\theta)  \rbr*{\prod_{i \in K } \indicator{[0 , \infty )} ( P_i ( \theta ) ) }
		\indicator{ [0 , \infty ) } \rbr*{ \textstyle  \br*{ \min_{j \in J } \frac{P_j ( \theta ) }{ \alpha_j } } -  \br*{ \max_{ i \in I } \frac{P_i ( \theta ) }{ \alpha_i } } } \\
		& \quad \times \displaystyle \int_{c_{s,t} ( \max_{ i \in I }  P_i ( \theta ) / \alpha_i ) }
		^{ c_{s, t } ( \min_{j \in J } P_j ( \theta ) / \alpha_j ) }
		Q(x) \, \d x .
		\end{split}
	\end{equation}
	Furthermore, it is well-known that $Q$ admits an analytic primitive $\tilde{\mathbf{Q}} \colon [s, t] \to \R$,
	which can be extended to a bounded subanalytic function $\mathbf{Q} \colon \R \to \R$ (by setting it to zero outside the interval $[s, t ]$).
	We therefore obtain for all $\theta \in \R^n$ that
	\begin{equation}
		\begin{split}
		\int_\R f(\theta , x ) \, \d x
		& = P_0 (\theta) \rbr*{\prod_{i \in K } \indicator{[0 , \infty )} ( P_i ( \theta ) ) }
		\indicator{ [0 , \infty ) } \rbr*{ \textstyle \br*{ \min_{j \in J } \frac{P_j ( \theta ) }{ \alpha_j } } - \br*{ \max_{ i \in I } \frac{P_i ( \theta ) }{ \alpha_i } } } \\
		& \quad \times \textstyle \br*{
			\mathbf{Q} \rbr*{ c_{s, t }\rbr*{ \min_{j \in J } \frac{P_j ( \theta ) }{ \alpha_j } } } -
			\mathbf{Q} \rbr*{c_{s,t} \rbr*{\max_{ i \in I } \frac{P_i ( \theta ) }{ \alpha_i } } } }
		\end{split}
	\end{equation}
	where all involved functions are locally bounded and subanalytic.
	This establishes \cref{prop:integral_prop_ii}.
\end{cproof}

\begin{cor}
	\label{cor:KL_analytic}
	Assume \cref{setting:ann:clip}, let $ \theta \in \R^{ \width + 1 } $, and assume that $ f $ and $ \dens $ are piecewise analytic \cfadd{def:piecewise_analytic}\cfload. Then there exist $\varepsilon, c \in (0, \infty)$, $\fA \in (0, 1)$ such that for all $\vartheta \in \{ \psi \in \R^{\width + 1} \colon \norm{ \psi - \theta} < \varepsilon\}$, $\fa \in (\fA, 1]$ it holds that
	\begin{equation} \label{cor:KL_eq}
		\textstyle \abs{\loss(\theta) - \loss(\vartheta)}^{\fa} \le c \norm{(\nabla \loss)(\vartheta)}.
	\end{equation}
\end{cor}

\cfclear
\begin{cproof}{cor:KL_analytic}
	\Nobs that \cref{lem:risk_function_prop,prop:integral_prop} ensure that $\loss$ is subanalytic \cfadd{def:sub_analytic}\cfload.
	Combining this with Bolte et al.~\cite[Theorem 3.1]{BolteDaniilidis2006} (see also Kurdyka~\cite[Theorem \L I]{Kurdyka1998})
	establishes \cref{cor:KL_eq}.
\end{cproof}

\begin{prop} \label{prop:gf:chain:rule}
	Assume \cref{setting:ann:clip} and let $\Theta \in C([0, \infty), \R^{\width + 1})$ satisfy for all $t \in [0, \infty)$ that $\Theta_t = \Theta_0 - \int_0^t ( \nabla \loss ) ( \Theta_s ) \, \d s$. Then it holds for all $t \in [0, \infty )$ that $\loss ( \Theta_t ) = \loss ( \Theta_0 ) - \int_0^t \norm{( \nabla \loss ) ( \Theta_s ) } ^2 \, \d s$.
\end{prop}
\begin{cproof}{prop:gf:chain:rule}
	This is a direct consequence of \cref{prop:risk:gradient}, the chain rule, and the fundamental theorem of calculus.
\end{cproof}

\subsubsection{Boundedness of GF trajectories}

We next verify in \cref{cor:a_priori_bound} that the \GF\ trajectory is bounded. This will be used to establish convergence in \cref{cor:limit_point}. 

\begin{lemma}[Upper bounds for realization functions]
\label{lem:realization_estimate}
Assume \cref{setting:ann:clip} and let $\theta = (\theta_1, \allowbreak \dots, \allowbreak \theta_{\width + 1}) \in \R^{\width + 1}$, $x \in [\0, \1]$. Then
\begin{equation}\label{eqn:lem:realization_estimate}
\textstyle \abs{\realization{\theta} ( x ) - \theta_{\width + 1 }} \leq \sum_{ i = 1 }^{\width} v_i.
\end{equation} 
\end{lemma}
\begin{cproof}{lem:realization_estimate}
\Nobs that the fact that for all $x \in \R$ it holds that $0 \le \fc(x) \le 1$ establishes \cref{eqn:lem:realization_estimate}.
\end{cproof}

\begin{prop}[Upper bounds for outer biases]
\label{prop:a_priori_bound0}
Assume \cref{setting:ann:clip} and let $\theta = (\theta_1, \allowbreak \dots, \allowbreak \theta_{\width + 1}) \in \R^{\width + 1}$. Then
\begin{equation} \label{prop:a_priori_bound:eqclaim}
\textstyle \abs{\theta_{\width + 1}} \le \br[\big]{\int_{\0}^{\1} \dens(x) \, \d x}^{- 1/2} \rbr[\big]{[\loss(\theta)]^{1/2} + \br[\big]{\int_{\0}^{\1} \abs{f(x)}^2 \dens(x) \, \d x}^{1/2}} + \sum_{i = 1}^{\width} v_i < \infty.
\end{equation}
\end{prop}
\begin{cproof}{prop:a_priori_bound0}
\Nobs that \cref{lem:realization_estimate} and Minkowski's inequality ensure that
\begin{equation}
\begin{split}
\textstyle [\loss( \theta )]^{1/2} & \textstyle = \br[\big]{\int_{\0}^{\1} \br[\big]{(\realization{\theta}(x) - \theta_{\width + 1}) - f(x) + \theta_{\width + 1}}^2 \dens(x) \, \d x}^{1/2} \\ 
& \textstyle \ge \br[\big]{\int_{\0}^{\1} (\theta_{\width + 1})^2 \dens(x) \, \d x}^{1/2} - \br[\big]{\int_{\0}^{\1} \br[\big]{(\realization{\theta}(x) - \theta_{\width + 1}) - f(x)}^2 \dens(x) \, \d x}^{1/2} \\ 
& \textstyle \ge \abs{\theta_{\width + 1}} \br[\big]{\int_{\0}^{\1} \dens(x) \, \d x}^{1/2} - \br[\big]{\int_{\0}^{\1} \br[\big]{\abs{f(x)} + (\sum_{i = 1}^{\width} v_i)}^2 \dens(x) \, \d x}^{1/2}.
\end{split}
\end{equation}
\Hence 
\begin{equation}
\begin{split}
\textstyle \abs{\theta_{\width + 1}} \br[\big]{\int_{\0}^{\1} \dens(x) \, \d x}^{1/2} & \textstyle \le [\loss(\theta)]^{1/2} + \br[\big]{\int_{\0}^{\1} \br[\big]{\abs{f(x)} + (\sum_{i = 1}^{\width} v_i)}^2 \dens(x) \, \d x}^{1/2} \\
& \textstyle \le [\loss(\theta)]^{1/2} + \br[\big]{\int_{\0}^{\1} \abs{f(x)}^2 \dens(x) \, \d x}^{1/2} + \br[\big]{\sum_{i = 1}^{\width} v_i} \br[\big]{\int_{\0}^{\1} \dens(x) \, \d x}^{1/2}.
\end{split}
\end{equation}
This establishes \cref{prop:a_priori_bound:eqclaim}.
\end{cproof}

\begin{cor}[Upper bounds for \GF\ trajectories]
\label{cor:a_priori_bound}
Assume \cref{setting:ann:clip} and let 
$\Theta = (\Theta^1, \dots, \allowbreak \Theta^{\width + 1}) \in C([0,\infty), \R^{\width + 1})$ 
satisfy for all $t \in [0, \infty)$ that $\Theta_t = \Theta_0 - \int_0^t (\nabla \loss)(\Theta_s) \, \d s$. Then
\begin{enumerate}[label=(\roman*)]
\item
\label{item1:cor:a_priori_bound} it holds for all $k \in \{1, 2, \ldots, \width\}$ that
\begin{equation}
\textstyle \sup\nolimits_{t \in [0, \infty)} \abs{\Theta^k_t} \le \max\{ w_k \abs{ \0 - [w_k ] ^{-1} }, w_k \abs{ \1 } , \abs{ \Theta^k_0 } \} < \infty,
\end{equation}

\item
\label{item2:cor:a_priori_bound} it holds that
\begin{equation}
\begin{split}
&
\textstyle 
  \sup_{ t \in [0,\infty) } \abs{\Theta^{\width + 1}_t} 
\\ & 
\textstyle 
  \le \sum_{i = 1}^{\width} v_i + \br[\big]{\int_{\0}^{\1} \dens(x) \, \d x}^{- 1/2} \rbr[\big]{[\loss( \Theta_0 )]^{1/2} + \br[\big]{\int_{\0}^{\1} \abs{f(x)}^2 \dens(x) \, \d x}^{1/2}} < \infty ,
\end{split}
\end{equation}

and 
\item
\label{item3:cor:a_priori_bound} it holds that
\begin{equation}
\begin{split}
\textstyle \sup\nolimits_{t \in [0,\infty)} \norm{\Theta_t} & \textstyle \le \sum_{i = 1}^{\width} v_i + \width \max\{(\abs{\0} + 1) \norm{w}, \abs{\1}\norm{w}, \norm{\Theta_0}\} \\
& \textstyle \quad + \br[\big]{\int_{\0}^{\1} \dens(x) \, \d x}^{- 1/2}\rbr[\big]{[\loss( \Theta_0 )]^{1/2} + \big[\int_{\0}^{\1}
\abs{f(x)}^2 \dens(x) \, \d x \big]^{1/2}} < \infty.
\end{split}
\end{equation}
\end{enumerate}
\end{cor}
\begin{cproof}{cor:a_priori_bound}
\Nobs that the fact that for all $k \in \cu{1, 2, \ldots, \width }$, $t \in [0 , \infty )$ with $\abs{\Theta_t^k } \allowbreak > \max \cu{w_k \abs{ \0 - [w_k ] ^{-1} }, w_k \abs{ \1 } }$ it holds that $\frac{\d}{\d t } \Theta_t ^k = - ( \frac{\partial}{\partial \theta_k} \loss ) ( \Theta_t ) = 0$ establishes \cref{item1:cor:a_priori_bound}. 
Combining \cref{item1:cor:a_priori_bound} with
 \cref{prop:gf:chain:rule} and \cref{prop:a_priori_bound0} 
 establishes \cref{item2:cor:a_priori_bound,item3:cor:a_priori_bound}.
\end{cproof}

\subsubsection{Existence of a limiting value}

\cfclear
\begin{cor}
\label{cor:limit_point}
Assume \cref{setting:ann:clip}, assume that $f$ and $\dens$ are piecewise analytic, 
and let $\Theta \in C( [0,\infty), \R^{ \width + 1 } ) $ satisfy for all $ t \in [0,\infty)$ that 
$\Theta_t = \Theta_0 - \int_0^t ( \nabla \loss ) ( \Theta_s ) \, \d s$ \cfadd{def:piecewise_analytic}\cfload. 
Then there exist $ \vartheta \in \R^{ \width + 1 } $, $ \alpha, \scrC \in (0,\infty) $
such that for all $ t \in [0,\infty) $ it holds that
\begin{equation}
\textstyle \norm{\Theta_t - \vartheta} \le \scrC (1 + t)^{- \alpha }, \; 0 \le \loss(\Theta_t) - \loss(\vartheta) \le \scrC (1 + t)^{- 1}, \; \text{and}\; \limsup_{s \to \infty} \norm{\Theta_s - \vartheta} = 0.
\end{equation}
\end{cor}
\begin{cproof}{cor:limit_point}
Combining \cref{cor:KL_analytic}, \cref{cor:a_priori_bound}, and the fact that $\loss \in C^1 ( \R^{ \width + 1 } , \allowbreak \R )$ with, e.g., \cite[Corollary~8.4]{JentzenRiekert2021ExistenceGlobMin} and \cref{prop:gf:chain:rule} ensures that there exist $\vartheta \in \R^{\width + 1}$, $\alpha , \scrC \in (0, \infty)$ such that for all $t \in [0,\infty)$ it holds that
\begin{equation}
\textstyle \norm{\Theta_t - \vartheta} \le \scrC ( 1 + t )^{- \alpha } \qqandqq 0 \le \loss(\Theta_t) - \loss(\vartheta) \le \scrC (1 + t)^{- 1}.
\end{equation}
\end{cproof}

\subsubsection{Activity of neurons during training}

In this section we analyze the behavior of the activity intervals $I_k^{ \Theta_t }$ during training.
We first prove in \cref{lem:never:leave:int} that a neuron can never become inactive (in the sense that its activity interval $I_k^{ \Theta_t }$ is empty)
in finite time.

\begin{lemma} 
\label{lem:never:leave:int}
Assume \cref{setting:ann:clip}, let $\Theta \in C ( [ 0 , \infty ) , \R ^{ \width + 1 } ) $ satisfy for all $t \in [0, \infty )$ that $\Theta_t = \Theta_0 - \int_0^t ( \nabla \loss ) ( \Theta_s ) \, \d s$,
and let $k \in \{1, 2, \ldots, \width\}$ satisfy $I_k^{ \Theta_0 } \neq \emptyset$. Then it holds for all $t \in [0, \infty )$ that $I_k^{ \Theta_t } \neq \emptyset$.
\end{lemma}

The proof of \cref{lem:never:leave:int} essentially relies on well-known arguments which are used to prove uniqueness for ODEs with locally Lipschitz continuous coefficients.

\begin{cproof}{lem:never:leave:int}
Throughout this proof let $\tau \in \R \cup \{\infty\}$ satisfy $\tau = \inf \cu{t \in [0, \infty ) \colon I_k^{ \Theta_t } = \emptyset}$. We assume for the sake of contradiction that $\tau  < \infty $. \Nobs that the fact that $[0, \infty ) \ni t \mapsto \psi_k ( \Theta_t ) \in \R $ is continuous and \cref{prop:active:intervals} ensure for all $t \in [0, \tau )$ that $\psi_k ( \Theta_t ) \in (\0 - \allowbreak [w_k]^{-1}, \allowbreak \1)$ and $\psi_k ( \Theta_\tau ) \in \cu{ \0 - [w_k]^{-1} , \1 }$. Moreover, we have that $\tau > 0$. In the following we assume without loss of generality that $\psi_k (\Theta_\tau ) = \1 $. Let $\fC \in ( 0, \infty )$ satisfy
\begin{equation} \label{lem:never:leave:int:eq1}
\textstyle \fC = 2 \rbr[\big]{\sum_{i=1}^\width v_i + \br[\big]{\sup\nolimits_{t \in [ 0 , \tau ] } \abs{\Theta_t^{ \width + 1 } } }
+ \br[\big]{ \sup\nolimits_{x \in [ \0 , \1 ] } \abs { f ( x ) } } }  v_k [w_k]^{-1} \sup\nolimits_{x \in [ \0 , \1 ] } \dens ( x ).
\end{equation}
\Nobs that \cref{prop:risk:gradient}, \cref{lem:realization_estimate}, and \cref{lem:never:leave:int:eq1} demonstrate for all $t \in [ 0 , \tau ]$ that
\begin{equation} \label{lem:never:leave:int:eq2}
\textstyle [w_k]^{-1} \abs[\big]{\rbr[\big]{ \tfrac{\partial}{\partial \theta_k } \loss } ( \Theta_t)} \le 2 v_k [w_k]^{-1} \int_{I_k^{\Theta_t }} \abs{\realization{\Theta_t} ( x ) - f ( x ) } \dens ( x ) \, \d x \le \fC \abs{ \psi_k ( \Theta_t ) - \1 }.
\end{equation}
Next let $\delta \in (0, \infty )$ satisfy $\delta = \min \cu{[2 \fC]^{-1}, \tau}$ and let $M \in \R$ satisfy $M = \sup_{t \in [\tau - \delta, \tau]} \abs{\psi_k (\Theta_t) - \1}$. \Nobs that the definition of $\tau$ proves that $M > 0$. Furthermore, \cref{lem:never:leave:int:eq2} shows for all $s \in [\tau - \delta, \tau]$ that
\begin{equation}
\begin{split}
\textstyle \abs{\psi_k ( \Theta_{s } ) - \1 } & \textstyle = \abs{\psi_k ( \Theta_{s } ) - \psi_k ( \Theta_{\tau  } ) }
= \abs[\big]{\int_{s}^\tau  [w_k]^{-1} \rbr[\big]{ \tfrac{\partial}{\partial \theta_k } \loss } ( \Theta_t ) \, \d t } \\
& \textstyle \le \int_{s}^\tau [w_k]^{-1} \abs[\big]{\rbr[\big]{\frac{\partial}{\partial \theta_k} \loss} (\Theta_t)} \, \d t \le \fC ( \tau - s ) M \le \fC \delta M \le \tfrac{M}{2}.
\end{split}
\end{equation}
Hence, we obtain that $M \le \tfrac{M}{2}$, which contradicts the fact that $M > 0$.
\end{cproof}

Next, we show that if a neuron becomes inactive in the limit (in the sense that $I_k^\vartheta = \emptyset$) then a certain endpoint bound for the realization $\realization{\vartheta}$ must be satisfied. Roughly speaking, the activity interval must leave the interval $( \0 , \1 )$ either on the right or on the left, which is why we formulate two separate results in \cref{lem:leave:interval:bound} and \cref{lem:leave:interval:bound:2}.

\begin{lemma}
\label{lem:leave:interval:bound}
Assume \cref{setting:ann:clip}, let $\Theta \in C([0, \infty ), \R^{\width + 1})$ satisfy for all $t \in [0, \infty)$ that $\Theta_t = \Theta_0 - \int_0^t ( \nabla \loss ) ( \Theta_s ) \, \d s$, let $\vartheta \in \R^{\width + 1 }$ satisfy $\limsup_{t \to \infty} \norm{\Theta_t - \vartheta } = 0$, and let $k \in \{1, 2, \ldots, \width\}$ satisfy $I_k^{ \Theta_0 } \neq \emptyset$ and $\psi_k (\vartheta) = \1$. Then $f(\1) \le \realization{\vartheta}(\1)$.
\end{lemma}
\begin{cproof}{lem:leave:interval:bound}
Throughout this proof let $\delta \in \R$ satisfy $\delta =  \frac{1}{2} (f(\1) - \realization{\vartheta}(\1))$ and assume for the sake of contradiction that $\delta > 0$. \Nobs that the assumption that $\limsup_{t \to \infty } \norm{\Theta_t - \vartheta } = 0$ implies that $\lim_{t \to \infty} \psi_k ( \Theta_t ) = \1$ and $\lim_{t \to \infty} \realization{\Theta_t} ( \psi_k ( \vartheta ) ) = \realization{\vartheta} ( \1 )$. In addition, the assumption that $\delta > 0$, \cref{lem:never:leave:int}, and the fact that $f$ is continuous demonstrate that there exists $T \in [0, \infty )$ which satisfies for all $t \in [T, \infty )$ that
\begin{equation}
\textstyle f ( \psi_{k} ( \Theta_t ) ) > \realization{\Theta_t} ( \psi_k ( \Theta_t ) ) + \delta \qandq \psi_k ( \Theta_t ) \in \rbr[\big]{ \1 - \min \cu[\big]{[w_k]^{-1}, \delta \br[\big]{\sum_{i=1}^\width v_i w_i}^{-1}}, \1}.
\end{equation}
Hence, we obtain for all $t \in [T, \infty)$, $x \in I_k^{\Theta_t}$ that $I_k^{\Theta_t} = (\psi_k (\Theta_t), \1)$ and
\begin{equation} \label{lem:leave:interval:bound:eq:mainest}
\begin{split}
\textstyle \realization{\Theta_t}(x) & \textstyle \le \realization{\Theta_t} (\psi_k ( \Theta_t ) ) + \rbr[\big]{\sum_{i = 1}^\width v_i w_i } (x - \psi_k (\Theta_t)) \\
& \textstyle \le \realization{\Theta_t} ( \psi_k ( \Theta_t ) ) + \rbr[\big]{\sum_{i=1}^\width v_i w_i } ( \1 - \psi_k ( \Theta_t ) ) \\
& \textstyle \le \realization{\Theta_t} ( \psi_k ( \Theta_t ) ) + \delta < f ( \psi_k ( \Theta_t ) ) \le f ( x ) .
\end{split}
\end{equation}
This and \cref{prop:risk:gradient} prove for all $t \in [T , \infty )$
that $\rbr{\frac{\partial}{\partial \theta_k } \loss } ( \Theta_t ) \le 0$. Combining this with the fact that for all $t \in [0, \infty )$
it holds that $\psi_k ( \Theta_t ) = \psi_k ( \Theta_0 ) + [w_k]^{-1} \int_0^t \rbr[\big]{ \frac{\partial}{\partial \theta_k } \loss } ( \Theta_s ) \, \d s$ ensures that $[T , \infty ) \ni t \mapsto \psi_{k } ( \Theta_t ) \in \R$ is non-increasing. However, this contradicts the fact that $\lim_{t \to \infty} \psi_{ k } ( \Theta_t ) = \1 > \psi_{k } ( \Theta_T)$.
\end{cproof}

Analogously we can prove the following result for the other endpoint.

\begin{lemma}
\label{lem:leave:interval:bound:2}
Assume \cref{setting:ann:clip}, let $\Theta \in C([0, \infty ), \R^{\width + 1 } )$ satisfy for all $t \in [0, \infty )$ that $\Theta_t = \Theta_0 - \int_0^t ( \nabla \loss ) ( \Theta_s ) \, \d s$, let $\vartheta \in \R^{\width + 1 }$ satisfy $\limsup_{t \to \infty } \norm{\Theta_t - \vartheta } = 0$, and let $k \in \{1, 2, \ldots, \width\}$ satisfy $I_k^{ \Theta_0 } \not= \emptyset$ and $\psi_k (\vartheta) = \0 - [w_k]^{-1}$. Then $f(\0) \ge \realization{\vartheta}(\0)$.
\end{lemma}

Next, we show in \cref{lem:neuron:leaves:interval} that under a suitable condition on the outer weights $v_i$, there necessarily exists a neuron which becomes inactive in the limit, and thus at least one of the previous endpoint bounds is applicable.

\begin{lemma} 
\label{lem:neuron:leaves:interval}
Assume \cref{setting:ann:clip}, assume that $f$ is non-decreasing
and satisfies $\Lip ( f ) < \min_{i \in \{1, 2, \ldots, \width\}} \allowbreak v_i w_i$, let $\Theta \in C([0, \infty ), \R^{\width + 1 } )$ satisfy for all $t \in [0, \infty )$ that $\Theta_t = \allowbreak \Theta_0 - \allowbreak \int_0^t ( \nabla \loss ) ( \Theta_s ) \allowbreak \, \d s$, let $\vartheta \in \R^{\width + 1 }$ satisfy $\limsup_{t \to \infty } \norm{\Theta_t - \vartheta } = 0$, let $V \in \R$ satisfy 
\begin{equation}
\textstyle V = \sup_{x \in ( \0 , \1 ) } \sum_{i \in \{1, 2, \ldots, \width\}, \, x \in I_i^\vartheta } v_i,
\end{equation}
and assume
\begin{equation}  \label{lem:neuron:leaves:interval:eqass}
\textstyle \sum_{i=1}^\width  v_i \indicator{(\0 - [w_i]^{-1} , \1)}(\psi_i (\Theta_0)) > f(\1) - f(\0) + 4 V.
\end{equation}
Then there exists $i \in \{1, 2, \ldots, \width\}$ such that $I_i^{\Theta_0} \neq \emptyset$ and $I_i^{\vartheta} = \emptyset$.
\end{lemma}
\begin{cproof}{lem:neuron:leaves:interval}
First, \cref{prop:risk:gradient} ensures that $( \nabla \loss ) ( \vartheta ) = 0$. Furthermore, \nobs that
\begin{equation}
\textstyle \sum_{i \in \{1, 2, \ldots, \width\}, \, \psi_i ( \vartheta ) \in (\1 - [w_i]^{-1}, \1 )} v_i \le V \qandq \sum_{i \in \{1, 2, \ldots, \width\}, \, \psi_i ( \vartheta ) \in (\0 - [w_i]^{-1}, \0 )} v_i \le V.
\end{equation}
Hence, we obtain that
\begin{equation}
\textstyle \realization{\vartheta} ( \1 ) - \realization{\vartheta} ( \0 ) \ge \sum_{i \in \{1, 2, \ldots, \width\}, \, \psi_i ( \vartheta ) \in [ \0 , \1 - [ w_i]^{-1} ] } v_i \ge \br[\big]{\sum_{i \in \{1, 2, \ldots, \width\}, \, I_i^\vartheta \neq \emptyset} v_i} - 2 V.
\end{equation}
This, \cref{prop:active:intervals}, \cref{lem:neuron:leaves:interval:eqass}, and \cref{cor:realization:endpoint:bound} imply that
\begin{equation}
\begin{split}
& \textstyle \br[\Big]{\sum_{i \in \{1, 2, \ldots, \width\}, \, I_i^{\Theta_0} \neq \emptyset} v_i} - (f(\1) - f(\0)) \\
& \textstyle = \br[\big]{\sum_{i \in \{1, 2, \ldots, \width\}, \, \psi_i (\Theta_0) \in (\0 - [w_i]^{-1}, \1)} v_i} - (f(\1) - f(\0)) \\
& \textstyle > 4 V \ge (\realization{\vartheta}(\1) - f(\1)) + (f(\0) - \realization{\vartheta}(\0)) + 2 V \\
& \textstyle = (\realization{\vartheta}(\1) - \realization{\vartheta}(\0) + 2 V) - (f(\1) - f(\0)) \\
& \textstyle \ge \br[\big]{\sum_{i \in \{1, 2, \ldots, \width\}, \, I_i^\vartheta \not= \emptyset} v_i} - ( f ( \1 ) - f ( \0 ) ).
\end{split}
\end{equation}
Therefore, we obtain that $\cu{i \in \{1, 2, \ldots, \width\} \colon I_i^\vartheta \neq \emptyset} \subsetneq \cu{i \in \{1, 2, \ldots, \width\} \colon I_i^{\Theta_0 } \neq \emptyset}$.
\end{cproof}

In the next step we verify in \cref{lem:limit:not:const} that there must exist at least one active neuron for the limit point $\vartheta$.

\begin{lemma}\label{lem:limit:not:const}
Assume \cref{setting:ann:clip}, assume that $f$ is non-decreasing and satisfies $\Lip ( f ) < \min_{i \in \{1, 2, \ldots, \width\}} v_i w_i$, let $\Theta \in C([0, \infty ), \R^{\width + 1 } )$ satisfy for all $t \in [0, \infty )$ that $\Theta_t = \allowbreak \Theta_0 - \allowbreak \int_0^t ( \nabla \loss ) ( \Theta_s ) \allowbreak \, \d s$,
let $\vartheta \in \R^{\width + 1 }$ satisfy $\limsup_{t \to \infty } \norm{\Theta_t - \vartheta } = 0$, assume
\begin{equation} \label{lem:limit:not:const:eq:ass1}
\textstyle \sup_{x \in (\0 , \1)} \sum_{i \in \{1, 2, \ldots, \width\}, \, x \in I_i^\vartheta} v_i < \frac{\min\cu{\int_\0^\1 ( f(\1 ) - f(x) ) \dens ( x ) \, \d x, \int_\0^\1 ( f( x ) - f ( \0 ) ) \dens ( x ) \, \d x} }{\int_\0^\1 \dens ( x ) \, \d x},
\end{equation}
and assume $\cu{i \in \{1, 2, \ldots, \width\} \colon \psi_i ( \Theta_0 ) \in ( \0 - [w_i]^{-1}, \1)} \neq \emptyset$. Then $\cu{i \in \{1, 2, \ldots, \width\} \colon \psi_i ( \vartheta ) \in ( \0 - [w_i]^{-1}, \1)} \neq \emptyset$.
\end{lemma}
\begin{cproof}{lem:limit:not:const}
Throughout this proof we assume for the sake of contradiction that
\begin{equation} \label{lem:limit:not:const:eq:cont}
\cu{i \in \{1, 2, \ldots, \width\} \colon \psi_i ( \vartheta ) \in ( \0 - [w_i]^{-1}, \1)} = \emptyset.
\end{equation}
\Nobs that \cref{lem:limit:not:const:eq:cont} demonstrates for all $x \in (\0, \1)$ that $\realization{\vartheta}(x) = \realization{\vartheta}(\1) = \realization{\vartheta}(\0)$. In addition, \cref{lem:limit:not:const:eq:cont} implies that there exists $k \in \{1, 2, \ldots, \width\}$ which satisfies $\psi_k ( \Theta_0 ) \in (\0 - [w_k]^{- 1}, \1)$ and $\psi_k ( \vartheta ) \allowbreak \in \allowbreak \cu{\0 \allowbreak - [w_k]^{-1}, \1}$. Assume without loss of generality that $\psi_k ( \vartheta ) = \1$. \Nobs that
the fact that $(\nabla \loss ) ( \vartheta ) = 0$,
 \cref{cor:realization:endpoint:bound}, and \cref{lem:leave:interval:bound} assure that
\begin{equation} \label{lem:limit:not:const:eq2}
\textstyle \abs{\realization{\vartheta} ( \1 ) - f ( \1 ) } = \realization{\vartheta} ( \1 ) - f ( \1 ) \le \sup\nolimits_{x \in (\0 , \1 ) } \sum_{i \in \{1, 2, \ldots, \width\}, \, x \in I_i^\vartheta} v_i.
\end{equation}
Furthermore, the fact that $( \frac{ \partial}{\partial \theta_{\width + 1 } } \loss ) ( \vartheta ) = 0$ implies that
\begin{equation}
\textstyle \realization{\vartheta} ( \1 ) = \br[\big]{\int_\0^\1 \dens ( x ) \, \d x }^{- 1} \rbr[\big]{\int_\0^\1 \realization{\vartheta} ( x ) \dens ( x ) \, \d x } = \br[\big]{\int_\0^\1 \dens ( x ) \, \d x}^{- 1} \rbr[\big]{\int_\0^\1 f(x) \dens(x) \, \d x}.
\end{equation}
Combining this with \cref{lem:limit:not:const:eq:ass1,lem:limit:not:const:eq2} yields a contradiction.
\end{cproof}

\subsubsection{Convergence with deterministic initial value}

Finally, we are in a position to establish in \cref{theo:conv:bias:deterministic} the main result about the convergence of \GF\ with deterministic initialization, which relies on the previous lemmas from this section.

\begin{theorem}
\label{theo:conv:bias:deterministic}
Assume \cref{setting:ann:clip}, assume that $f$ is non-decreasing and satisfies $\Lip ( f ) < \min_{i \in \{1, 2, \ldots, \width\}} v_i w_i$, let $\Theta \in C([0, \infty ), \R^{\width + 1 } )$ satisfy for all $t \in [0, \infty )$ that $\Theta_t = \allowbreak \Theta_0 - \allowbreak \int_0^t ( \nabla \loss ) ( \Theta_s ) \allowbreak \, \d s$, let $\vartheta \in \R^{\width + 1 }$ satisfy $\limsup_{t \to \infty } \norm{\Theta_t - \vartheta } = 0$, let $V \in \R$ satisfy 
\begin{equation}
\textstyle V = \sup_{x \in ( \0 , \1 ) } \sum_{i \in \{1, 2, \ldots, \width\}, \, x \in I_i^\vartheta } v_i,
\end{equation}
assume
\begin{equation}
\textstyle V < \br[\big]{\int_\0^\1 \dens ( x ) \, \d x } ^{-1} 
\min\cu[\big]{ \int_\0^\1 ( f(\1 ) - f(x) ) \dens ( x ) \, \d x, \int_\0^\1 ( f( x ) - f(\0)) \dens ( x ) \, \d x} ,
\end{equation}
and assume $\sum_{i = 1 }^ \width  v_i \indicator{ ( \0 - [w_i]^{-1} , \1 ) } ( \psi_i ( \Theta_0 ) ) > f ( \1 ) - f ( \0 ) + 4 V$. Then 
\begin{equation}
\label{theo:conv:bias:deterministic:eqclaim}
\textstyle \loss(\vartheta) \le \rbr{2 V ^2 + ( f(\1 ) - f(\0)) V} \int_\0^\1 \dens(x) \, \d x.
\end{equation}
\end{theorem}
\begin{cproof}{theo:conv:bias:deterministic}
First, \cref{prop:risk:gradient} implies that $(\nabla \loss ) ( \vartheta ) = 0$. Combining this with \cref{lem:limit:not:const}
ensures that $\cup_{i \in \{1, 2, \ldots, \width\}} I_i^\vartheta \not= \emptyset$. Let $q ,r \in [ \0 , \1 ]$ satisfy $q = \inf (\cup_{i \in \{1, 2, \ldots, \width\}} I_i^\vartheta)$ and $r = \sup (\cup_{i \in \{1, 2, \ldots, \width\}} I_i^\vartheta)$. \Nobs that \cref{lem:realization:pointwise:bound} implies that
\begin{equation} \label{theo:conv:bias:deterministic:eq1}
\textstyle \sup_{x \in (q, r)} \abs{\realization{\vartheta}(x) - f(x)} \le V,
\end{equation}
\begin{equation} \label{theo:conv:bias:deterministic:eq3}
\textstyle \sup_{x \in (\0, q)} \abs{\realization{\vartheta}(x) - f(x)} \le \max \cu{V, \abs{\realization{\vartheta}(\0) - f(\0)}},
\end{equation}
and
\begin{equation} \label{theo:conv:bias:deterministic:eq4}
\textstyle \sup_{x \in (r, \1)} \abs{\realization{\vartheta}(x) - f(x)} \le \max \cu{V, \abs{\realization{\vartheta}(\1) - f(\1)}}.
\end{equation}
We assume without loss of generality that there exists $k \in \{1, 2, \ldots, \width\}$ which satisfies $\psi_k (\Theta_0) \allowbreak \in (\0 - [w_k]^{-1}, \1)$ and $\psi_k (\vartheta) = \1$ (cf.\ \cref{lem:neuron:leaves:interval}). Combining this with \cref{cor:realization:endpoint:bound} and \cref{lem:leave:interval:bound} assures that $\abs{\realization{\vartheta}(\1) - f(\1)} \le V$.

 In the following we distinguish between the case $\realization{\vartheta}(\0) \le f(\0)$ and the case $\realization{\vartheta}(\0) > f(\0)$. First we prove \cref{theo:conv:bias:deterministic:eqclaim} in the case $\realization{\vartheta}(\0) \le f(\0)$. \Nobs that \cref{cor:realization:endpoint:bound} assures that $\abs{\realization{\vartheta} ( \0 ) - f ( \0 ) } \le V$. Combining this with \cref{theo:conv:bias:deterministic:eq1,theo:conv:bias:deterministic:eq3,theo:conv:bias:deterministic:eq4} demonstrates that
\begin{equation}
\textstyle \loss(\vartheta) = \int_\0^\1 (\realization{\vartheta}(x) - f(x))^2 \dens(x) \, \d x \le V^2 \int_\0^\1 \dens(x) \, \d x.
\end{equation}
This proves \cref{theo:conv:bias:deterministic:eqclaim} in the case $\realization{\vartheta}(\0) \le f(\0)$. 

In the next step we prove \cref{theo:conv:bias:deterministic:eqclaim} in the case $\realization{\vartheta}(\0) > f(\0)$. Let $i \in \{1, 2, \ldots, \width\}$ satisfy $q = \inf (I_i^\vartheta)$. \Nobs that \cref{prop:risk:gradient} and the integral mean value theorem imply that there exists $z \in I_i^\vartheta$ which satisfies $\realization{\vartheta}(z) = f(z)$. This and the assumption that $\Lip(f) < \min_{i \in \{1, 2, \ldots, \width\}} v_i w_i$ ensure that $\realization{\vartheta}(q) < f(q)$. Combining this with the intermediate value theorem demonstrates that $\cu{u \in (\0, q) \colon \realization{\vartheta}(u) = f(u)} \neq \emptyset$. Let $y = \inf \cu{u \in (\0, q) \colon \realization{\vartheta}(u) = f(u)}$. \Nobs that the fact that $\forall \, x \in (\0 , q ) \colon \realization{\vartheta}(x) = \realization{\vartheta}(q)$ and \cref{theo:conv:bias:deterministic:eq1} prove that $\sup_{x \in (y, q ) } \abs { \realization{\vartheta} ( x ) - f ( x ) } \le V$. The fact that $( \frac{\partial}{\partial \theta_{\width + 1 } } \loss ) (\vartheta) = 0$, \cref{theo:conv:bias:deterministic:eq1}, and \cref{theo:conv:bias:deterministic:eq4} therefore show that
\begin{equation}
\begin{split}
\textstyle \int_\0^y ( \realization{\vartheta} ( x ) - f ( x ) ) \dens ( x ) \, \d x & \textstyle = - \int_y^\1 ( \realization{\vartheta} ( x ) - f ( x ) ) \dens ( x ) \, \d x \\
& \textstyle \le \int_y^\1 \abs{\realization{\vartheta} ( x ) - f ( x ) } \dens ( x ) \, \d x \le V \int_y^\1 \dens ( x ) \, \d x \le V \int_\0^\1 \dens ( x ) \, \d x.
\end{split}
\end{equation}
Combining this with the fact that for all $x \in (\0 , y )$ it holds that $\realization{\vartheta} ( x ) \ge f ( x )$ and \cref{theo:conv:bias:deterministic:eq3} demonstrates that
\begin{equation}
\begin{split}
\textstyle \int_\0^y ( \realization{\vartheta} ( x ) - f ( x ) )^2 \dens ( x ) \, \d x & \textstyle \le \br[\big]{ \sup\nolimits_{x \in (\0 , y ) } \abs{\realization{\vartheta} ( x ) - f ( x ) } } \int_\0^y ( \realization{\vartheta} ( x ) - f ( x ) ) \dens ( x ) \, \d x \\
& \textstyle \le V \max \cu{f(\1 ) - f(\0), V} \int_\0 ^\1 \dens ( x ) \, \d x.
\end{split}
\end{equation}
This, \cref{theo:conv:bias:deterministic:eq3}, and \cref{theo:conv:bias:deterministic:eq4} show that
\begin{equation}
\textstyle \loss(\vartheta) \le \rbr{V^2 + V \max \cu{f(\1) - f(\0), V}} \int_\0^\1 \dens ( x ) \, \d x.
\end{equation}
This establishes \cref{theo:conv:bias:deterministic:eqclaim} in the case $\realization{\vartheta}(\0) > f(\0)$.
\end{cproof}

\subsection{Convergence results for GF with random initializations}
Our goal in this section is to verify that for our particular initialization the conditions of \cref{theo:conv:bias:deterministic} are satisfied with high probability and, thereby, prove our main convergence result for \GF\ with a single random initialization.

\subsubsection{Products of half-normally distributed random variables}

In our main convergence result, \cref{theo:main:conv}, the weights $v_i, w_i$ will be initialized randomly. Since they are required to be positive, we will use a half-normal distribution instead of the popular normal distribution.
The first goal is to prove \cref{lem:random:init} below. For this we require several technical results on products of two independent half-normal random variables and sums of random variables of this type.

\begin{definition}[Modified Bessel function]
\label{def:bessel}
We denote by $K_0 \colon (0, \infty ) \to [ 0 , \infty ]$ the function which satisfies for all $x \in (0, \infty )$ that
\begin{equation}
\textstyle K_0 (x) = \int_0^\infty \exp(-x \cosh(w)) \, \d w.
\end{equation}
\end{definition}

The modified Bessel function of the first kind appears when considering the product of two independent normal random variables.
The next result is well-known in the physics literature; cf., e.g., Abramowitz \& Stegun~\cite[Chapter 9]{AbramowitzStegun1992}.

\cfclear
\begin{lemma}[Properties of the Modified Bessel function]
\label{lem:bessel:properties} \cfadd{def:bessel}
We have that
\begin{enumerate} [label = (\roman*)]
\item it holds that $K_0$ is strictly decreasing,

\item it holds for all $x \in (0, \infty)$ that $0 < K_0(x) < \infty$,

\item it holds that $\liminf_{x \searrow 0} K_0(x) = \infty$, and

\item it holds that $\limsup_{x \searrow 0} \Abs{\frac{K_0(x)}{\abs{\ln(x)}} - 1} = 0$
\end{enumerate}
\cfload.
\end{lemma}

As a consequence,
we obtain the following well-known result on products of standard normal random variables. The detailed proof is only included for completeness.

\cfclear
\begin{prop}
\label{prop:product:normal}
Let $(\Omega, \mathcal{F}, \P)$ be a probability space, let $X \colon \Omega \to (0, \infty )$ and $Y \colon \Omega \to (0, \infty )$ be i.i.d.\ random variables, and assume for all $x \in ( 0, \infty )$ that $\P ( X \le x ) = [\frac{2}{\pi}]^{1/2} \int_0^x \exp(- \frac{y^2}{2}) \, \d y$. Then 
\begin{enumerate} [label = (\roman*)]
\item
\label{prop:product:normal:item1} it holds for all $z \in (0, \infty)$ that $\P (XY \le z) = \frac{2}{\pi} \int_0^z K_0 ( y ) \, \d y$ and \cfadd{def:bessel}

\item
\label{prop:product:normal:item2} it holds that there exists $c \in (0, \infty )$ such that for all $z \in (0, \infty)$ it holds that $\P ( X Y \le z ) \le c z (1 + \abs{\ln(z)})$
\end{enumerate}
\cfload.
\end{prop}

The distribution of $X$ and $Y$ is sometimes called half-normal distribution or folded normal distribution.

\begin{cproof}{prop:product:normal}
	First, the fact that $X,Y$ are independent, the integral transformation theorem, and Fubini's theorem ensure for all $z \in ( 0 , \infty )$ that
	\begin{equation} \label{prop:product:normal:eq1}
	\begin{split}
	\textstyle \P ( X Y \le z ) & \textstyle = \frac{2}{\pi } \int_0^\infty \int_0^\infty \exp ( - \frac{x^2 + y^2 }{ 2 } ) \indicator{[0 , z ] } ( x y ) \, \d y \, \d x \\
	& \textstyle = \frac{2}{\pi } \int_{0}^{\infty} \int_{0}^{z/x} \exp ( - \frac{x^2 + y^2 }{ 2 } ) \, \d y \, \d x \\
	& \textstyle = \frac{2}{\pi } \int_{0}^{\infty} \int_{0}^{z} \exp ( - \frac{x^2 + (u/x)^2 }{2} ) \frac{1}{x} \, \d u \, \d x \\
	& \textstyle = \frac{2}{\pi} \int_{0}^{z} \rbr[\big]{\int_{0}^{\infty} \frac{1}{x} \exp(- \frac{x^2 + (u/x)^2}{2}) \, \d x} \, \d u.
	\end{split}
	\end{equation}
	Furthermore, the integral transformation theorem with the diffeomorphism $\R \ni w \mapsto \allowbreak u^{1/2} \allowbreak \exp(\frac{w}{2}) \allowbreak \in (0, \infty)$ and the fact that for all $w \in \R$ it holds that $\cosh(w) = \cosh(- w)$ show that for all $u \in (0, \infty )$ it holds that
	\begin{equation}
	\begin{split}
	& \textstyle \int_{ ( 0 , \infty ) } \frac{1}{x} \exp \rbr[\big]{ - \tfrac{x^2 + (u/x)^2 }{2} } \, \d x = \int_\R  \frac{1}{\sqrt{u}} \exp ( - \frac{w}{2} ) \exp \rbr[\big]{ - \frac{u e^w + u e^{-w} }{2} } \frac{ \sqrt{u }}{2} \exp ( \frac{w}{2} ) \, \d w \\
	& \textstyle = \frac{1}{2} \int_\R  \exp ( - u \cosh ( w ) ) \, \d w 
	= \int_0^\infty \exp ( - u \cosh ( w ) ) \, \d w = K_0 ( u ).
	\end{split}
	\end{equation}
	Combining this with \cref{prop:product:normal:eq1} establishes \cref{prop:product:normal:item1}. Next, \cref{lem:bessel:properties}
	implies that there exists $\eps \in (0 , 1 )$ which satisfies for all $x \in (0, \eps )$ that $K_0 ( x ) \le 2 \abs{\ln(x)}$. Combining this with \cref{prop:product:normal:item1} proves for all $z \in (0, \eps]$ that
	\begin{equation}
	\textstyle \P ( X Y \le z ) \le \frac{4}{\pi } \int_0^z \abs{\ln(y)} \, \d y = \frac{4}{\pi } z ( 1 + \abs{\ln(z)}).
	\end{equation}
	This, the fact that $\forall \, z \in ( 0, \infty ) \colon \P ( X Y \le z ) \le 1$, and the fact that $\forall \, z \in [ \eps, \infty ) \colon z ( 1 + \abs{\ln(z)}) \allowbreak \ge \eps$ demonstrate for all $z \in (0 , \infty )$ that $\P ( X Y \le z ) \le \max \cu{\eps^{-1} , 4 \pi ^{-1} } z ( 1 + \abs{\ln(z)})$. This establishes \cref{prop:product:normal:item2}.
\end{cproof}

\subsubsection{Properties of order statistics}

\begin{lemma} \label{lem:log:estimate:help}
Let $c, \delta \in (0, \infty)$, let $g \colon (0,1) \to (0, \infty)$ be non-decreasing, and assume for all $u \in (0 , 1 )$ that $g ( u ) ( 1 + \abs{\ln(g(u))}) \ge c u$. Then $\inf_{u \in (0, 1)} u^{-1 - \delta} g(u) > 0$.
\end{lemma}
\begin{cproof}{lem:log:estimate:help}
Throughout this proof let $h \colon (0, 1) \to (0, \infty)$ satisfy $\forall \, u \in (0 , 1 ) \colon h ( u ) = u ^{-1 - \delta } g ( u )$. \Nobs that the assumption that $\forall \, u \in (0 , 1 ) \colon g ( u ) ( 1 + \abs{\ln(g(u))}) \ge c u$ demonstrates for all $u \in (0, 1)$ that
\begin{equation}\label{eqn:lem:log:estimate:help:lowerbound}
\textstyle u^{\delta} h ( u ) ( 1 + ( 1 + \delta ) \abs{\ln h ( u ) } + \abs{\ln u } ) \ge c > 0.
\end{equation}
Next we prove that there exists $\varepsilon \in (0, \infty)$ such that $\forall \, u \in (0 , \varepsilon ) \colon h ( u ) > \varepsilon$. For the sake of contradiction assume that there exists a sequence $(u_n) \downarrow 0$ which satisfies $\lim_{n \to \infty} h (u_n) = 0$. \Nobs that the fact that $\Forall \eta > 0 \colon \lim_{x \downarrow 0 } x^\eta \abs{\ln(x)} = 0$ ensures that $\lim_{n \to \infty} ( u_n^{\delta} h ( u_n ) ( 1 + ( 1 + \delta ) \abs{\ln(h(u_n))} + \abs{\ln(u_n)})) = 0$. This contradicts to \cref{eqn:lem:log:estimate:help:lowerbound}. Hence, we obtain that there exists $\varepsilon > 0$ which satisfies $\forall \, u \in (0 , \varepsilon ) \colon h ( u ) > \varepsilon$. This implies that $\inf_{u \in (0 , \varepsilon ) } h ( u ) > 0$. Combining this with the fact that $\inf_{u \in [ \varepsilon , 1 ) } h ( u ) \ge g (\eps) > 0$ shows that $\inf_{u \in (0, 1)} h(u) > 0$.
\end{cproof}

We next show some properties of the sums of the $k$ smallest among $N$ random variables with the distribution appearing in \cref{prop:product:normal}. For this we require a few known results about order statistics from \cite{Papadatos2001} and \cite{Bickel1967}.

\cfclear 
\begin{lemma}\label{lem:sum:order:stat}

Let $(\Omega, \cF, \P)$ be a probability space, for every $N \in \N$ let $X_1^N, \ldots, X_N^N \colon \Omega \to (0, \infty)$ be i.i.d.\ random variables, assume for all $N \in \N$, $x \in ( 0, \infty )$ that $\P (X_1^N \le x) = \frac{2}{\pi} \allowbreak \int_0^x K_0(y) \, \d y$, and for every $N \in \N$ let $X^N_{(1)}, \ldots, X^N_{(N)} \colon \Omega \to (0, \infty)$ be the order statistics\footnote{The random permutation defined by sorting the values of $X^N_i$, $i \in \{1, 2, \ldots, N\}$, in non-decreasing order, i.e., $\Forall N \in \N \colon X^N_{(1)} \le \cdots \le X^N_{(N)}$.} of $X^N_1, \ldots, X^N_N$ \cfadd{def:bessel}\cfload. Then
\begin{enumerate}[label = (\roman*)]
\item
\label{lem:sum:order:stat:item1} it holds for all $\delta \in (0, \infty)$ that there exists $c \in (0, \infty)$ such that for all $N \in \N$, $k \in \{1, 2, \ldots, N\}$ it holds that
\begin{equation}
\textstyle \E \br[\big]{\sum_{i=1}^k X_{(i)}^N } \ge \frac{c k^{2 + \delta}}{N^{1 + \delta}}
\end{equation}

and
\item
\label{lem:sum:order:stat:item2} it holds that for all $N \in \N$, $k \in \{1, 2, \ldots, N\}$ that
\begin{equation}
\textstyle \operatorname{Var} \rbr[\big]{\sum_{i=1}^k X^N_{(i)}} \le \operatorname{Var} \rbr[\big]{\sum_{i=1}^N X^N_i} = N \operatorname{Var} (X^N_1) \le N.
\end{equation}
\end{enumerate}
\end{lemma}
\begin{cproof}{lem:sum:order:stat}
Throughout this proof let $\delta \in (0, \infty)$ and let $g \colon (0, 1) \to (0, \infty)$ satisfy for all $u \in (0, 1)$ that $g(u) = \inf \cu{x \in (0, \infty) \colon \frac{2}{\pi} \int_0^x K_0(y) \, \d y \ge u}$. \Nobs that, e.g., Papadatos~\cite[Theorem~3.1]{Papadatos2001} implies for all $N \in \N$, $i \in \{1, 2, \ldots, N\}$ that 
\begin{equation} \label{lem:sum:order:stat:eq1}
\textstyle \E [X^N_{(i)}] \ge \frac{N}{i} \int_0^{i / N} g(u) \, \d u.
\end{equation}
Moreover, \cref{prop:product:normal} shows that there exists $c \in (0, \infty)$ which satisfies for all $N \in \N$, $x \in (0, \infty )$ that $\P(X_1^N \le x) \le c x (1 + \abs{\ln(x)})$. Hence, we obtain for all $u \in (0, 1)$ that $c g ( u ) ( 1 + \abs{ \ln(g(u))}) \ge u$. Combining this with \cref{lem:log:estimate:help} assures that there exists $\kappa \in (0, \infty)$ which satisfies $\forall \, u \in (0,1) \colon g ( u ) \ge \kappa u ^{ 1 + \delta }$. This and \cref{lem:sum:order:stat:eq1} imply for all $N \in \N$, $i \in \{1, 2, \ldots, N\}$ that
\begin{equation}
\textstyle \E [X^N_{(i)}] \ge \frac{\kappa N}{i} \int_0^{i/N} u^{1 + \delta} \, \d u = \frac{\kappa}{2 + \delta} \rbr[\big]{\frac{i}{N}}^{1 + \delta}.
\end{equation}
Combining this with the fact that for all $k \in \N$ it holds that $\sum_{i=1}^k i^{1 + \delta} \ge \int_0^k x ^{ 1 + \delta } \, \d x \allowbreak = \frac{1}{2 + \delta } k ^{ 2 + \delta }$ establishes \cref{lem:sum:order:stat:item1}. Next \nobs that for all $N \in \N$, $i, j \in \{1, 2, \ldots, N\}$ it holds that $\operatorname{Cov} ( X^N_{ ( i ) } , X^N_{ ( j ) } ) \ge 0$ (cf., e.g., Bickel~\cite[Theorem 2.1]{Bickel1967}). The assumption that $X^N_i$, $i \in \{1, \allowbreak 2, \allowbreak \ldots, \allowbreak N\}$, are i.i.d.\ therefore assures for all $N \in \N$, $k \in \{1, 2, \ldots, N\}$ that
\begin{equation}
\begin{split}
\textstyle N \operatorname{Var} ( X_1^N ) & \textstyle =
\operatorname{Var} \rbr[\big]{\sum_{i=1}^N X^N_i} = \operatorname{Var} \rbr[\big]{\sum_{i=1}^N X^N_{(i)}} = \sum_{i, j=1}^N \operatorname{Cov} (X^N_{(i)}, X^N_{(j)}) \\
& \textstyle \ge \sum_{i,j=1}^k \operatorname{Cov} (X^N_{(i)}, X^N_{(j)}) = \operatorname{Var} \rbr[\big]{\sum_{i=1}^k X^N_{(i)}}.
\end{split}
\end{equation}
This and the fact that $\forall \, N \in \N \colon \operatorname{Var} ( X^N_1 ) \le \E [ (X_1 ^N )^2 ] = 1$ establish \cref{lem:sum:order:stat:item2}.
\end{cproof}

\cfclear 
\begin{cor}
\label{cor:sum:order:stat}
Let $(\Omega, \cF, \P)$ be a probability space, for every $N \in \N$
let $X^N_1, \ldots, X^N_N \colon \Omega \to (0, \infty)$ be i.i.d.\ random variables, assume for all $N \in \N$, $x \in ( 0, \infty )$ that $\P(X^N_1 \le x) = \frac{2}{\pi} \allowbreak \int_0^x K_0(y) \, \d y $, for every $N \in \N$ let $X^N_{(1)}, \ldots, X^N_{(N)} \colon \Omega \to (0, \infty)$ be the order statistics of $X^N_1, \allowbreak \ldots, \allowbreak X^N_N$, and let $\fC \in \R$, $\eps \in (0, 1)$, $\gamma \in (\frac{1 + \varepsilon}{2}, 1)$, $\delta \in (0, \frac{2 \gamma - \varepsilon - 1 }{ 1 - \gamma } ) $ \cfadd{def:bessel}\cfload. Then there exists $c \in \R$ such that for all $N \in \N$ it holds that
\begin{equation}
\textstyle \P \rbr[\Big]{\sum_{i=1}^{ \lceil N ^\gamma \rceil } X_ { ( i ) } ^N \le \fC N ^{ \varepsilon } } \le c N ^{3 + 2 \delta - 4 \gamma - 2 \gamma \delta } .
	\end{equation}
\end{cor}
\begin{cproof}{cor:sum:order:stat}
First, \cref{lem:sum:order:stat} assures that there exists $\kappa \in (0, \infty )$ which satisfies for all $N \in \N$ that
\begin{equation}
\textstyle \E \br[\Big]{\sum_{i=1}^{\lceil N^\gamma \rceil } X_ { ( i ) } ^N } \ge 2 \kappa \frac{ (N ^\gamma ) ^{ 2 + \delta } }{ N ^{ 1 + \delta}} = 2 \kappa N^{\gamma (2 + \delta) - (1 + \delta)}.
\end{equation}
In addition, the fact that $\gamma (2 + \delta) - (1 + \delta) = (\gamma - 1) \delta + (2 \gamma - 1) > \varepsilon$ proves that there exists $\bfn \in \N$ which satisfies for all $N \in \N \cap [\bfn, \infty)$ that $\kappa N^{\gamma (2 + \delta) - (1 + \delta)} > \fC N^\varepsilon$. The Chebyshev's inequality and \cref{lem:sum:order:stat} hence prove that for all $N \in \N \cap [ \bfn , \infty )$ it holds that
\begin{equation}
\begin{split}
& \textstyle \P \rbr[\Big]{\sum_{i=1}^{\lceil N^\gamma \rceil } X_ { ( i ) } ^N \le \fC N^{\varepsilon}} \le \P \rbr[\Big]{\sum_{i=1}^{\lceil N^\gamma \rceil} X_{(i)}^N \le \E \br[\Big]{\sum_{i=1}^{\lceil N^\gamma \rceil} X_{(i)}^N} - \kappa N^{\gamma (2 + \delta) - (1 + \delta)}} \\
& \textstyle \le \P \rbr[\Big]{\abs[\Big]{\sum_{i=1}^{\lceil N^\gamma \rceil} X_{(i)}^N - \E \br[\Big]{\sum_{i=1}^{\lceil N^\gamma \rceil} X_{(i)}^N}} \ge \kappa N^{\gamma (2 + \delta) - (1 + \delta)}} \\
& \textstyle \le (\kappa N^{\gamma (2 + \delta) - (1 + \delta)})^{-2}
\operatorname{Var} \rbr[\big]{\sum_{i=1}^{\lceil N^\gamma \rceil} X_{(i)}^N} \le (\kappa N^{\gamma (2 + \delta) - (1 + \delta)})^{-2} N
= \kappa^{-2} N^{3 + 2 \delta - 4 \gamma - 2 \gamma \delta}.
\end{split}
\end{equation}
This, the fact that $\forall \, N \in \{1, 2, \ldots, \bfn\} \colon N ^{3 + 2 \delta - 4 \gamma - 2 \gamma \delta} \ge \min \cu{ 1 , \bfn^{ 3 + 2 \delta - 4 \gamma - 2 \gamma \delta}}$, and the fact that $\forall \, N \in \N \colon  \P \rbr{\sum_{i=1}^{\lceil N^\gamma \rceil } X_{(i)}^N \le \fC N^{\varepsilon}} \le 1$ demonstrate for all $N \in \N$ that
\begin{equation}
\textstyle \P \rbr[\Big]{\sum_{i=1}^{\lceil N^\gamma \rceil} X_{(i)}^N \le \fC N^{\eps}} \le \max \cu[\big]{1, \bfn^{- (3 + 2 \delta - 4 \gamma - 2 \gamma \delta)}, \kappa^{-2}} N^{3 + 2 \delta - 4 \gamma - 2 \gamma \delta}.
\end{equation}
\end{cproof}

\begin{lemma}
\label{lem:size:estimate}
Let $(\Omega, \cF, \P)$ be a probability space, let $N \in \N$, $\alpha \in \R$, $\varepsilon \in (0, 1)$, $\delta \in (0, \varepsilon)$, let $X_1, \ldots, X_N \colon \Omega \to ( 0 , \infty )$ be random variables, and assume for all $i \in \{1, 2, \ldots, N\}$, $x \in (0, \infty)$ that $\P(N^\alpha X_i \le x) = [\frac{2}{\pi}]^{1/2} \int_0^x \exp(- \frac{y^2}{2}) \, \d y$. Then
\begin{equation}
\textstyle \P \rbr[\big]{\max \cu{X_1, \ldots, X_N} \ge N^{-\alpha + \varepsilon}} \le (\varepsilon - \delta)^{-1} N^{- \delta}.
\end{equation}
\end{lemma}
\begin{cproof}{lem:size:estimate}
Throughout this proof let $\Gamma \colon (0, \infty) \to (0, \infty)$ satisfy for all $x \in (0, \infty)$ that $\Gamma(x) = \int_0^{\infty} t^{x - 1} e^{-t} \, \d t$. First, the well-known formula for absolute moments of the normal distribution (cf., e.g., (18) in Winkelbauer~\cite{Winkelbauer2012}) and the fact that $\forall \, x \in (0 , \infty ) \colon \Gamma (x + 1) \le \max \cu{1 , x^x}$ (cf., e.g., \cite[Proposition~5.3]{JentzenWelti2020}) assure for all $p \in [1, \infty)$ that
\begin{equation}
\begin{split}
\textstyle (\E[(X_1)^p])^{1/p} & \textstyle = N^{-\alpha} (\E[(N^\alpha X_1)^p])^{1/p} = N^{-\alpha} \rbr[\big]{\pi^{- 1/2} 2^{p/2} \Gamma(\frac{p+1}{2})}^{1/p} \\
& \textstyle \le N^{-\alpha} \pi^{- 1/(2 p)} 2^{1/2} \max \cu[\big]{1, (\tfrac{p-1}{2})^{(p-1)/(2 p)}} \\
& \textstyle \le N^{-\alpha} \rbr[\big]{2^{1/2} (\frac{p}{2})^{1/2} \indicator{[3, \infty)} (p) + 2^{1/2} \indicator{[2^{1/2}, 3)} (p) + \pi^{- 2^{- 3/2}} 2^{1/2} \indicator{[1, 2^{1/2})} (p)} \\
& \textstyle \le N^{-\alpha} p. 
\end{split}
\end{equation}
This, Jensen's inequality, and the fact that $(\varepsilon - \delta)^{-1} \in [1 , \infty )$
demonstrate that
\begin{equation}
\begin{split}
& \textstyle \P \rbr[\big]{\max \cu{X_1, \ldots, X_N} \ge N^{- \alpha + \varepsilon}} = \P \rbr[\big]{N^{\alpha - \varepsilon} \max \cu{X_1, \ldots, X_N} \ge 1} \\
& \textstyle \le \E \br[\big]{N^{\alpha - \varepsilon} \max \cu{X_1, \ldots, X_N}} \le N^{\alpha - \varepsilon} \rbr[\big]{\E \br[\big]{(\max \cu{X_1, \ldots, X_N})^{(\varepsilon - \delta)^{-1}}}}^{\varepsilon - \delta} \\
& \textstyle \le N^{\alpha - \varepsilon} \rbr[\big]{\E \br[\big]{(X_1)^{(\varepsilon - \delta)^{-1}} + \cdots + (X_N)^{(\varepsilon - \delta)^{-1}}}}^{\varepsilon - \delta} = N^{\alpha - \varepsilon} \rbr[\big]{N \E \br[\big]{( X_1)^{(\varepsilon - \delta)^{-1}}}}^{\varepsilon - \delta} \\
& \textstyle \le N^{\alpha - \varepsilon} N^{\varepsilon - \delta} N^{- \alpha}(\varepsilon - \delta)^{-1} = (\varepsilon - \delta)^{-1} N^{- \delta}.
\end{split}
\end{equation}
\end{cproof}

\begin{cor}
\label{cor:product:bias:sum}
Let $(\Omega, \cF, \P)$ be a probability space, let $\alpha, \beta \in \R$, for every $N \in \N$ let $V_1^N, \ldots, \allowbreak V_N^N, \allowbreak W_1^N, \allowbreak \ldots, W_N^N \colon \Omega \to (0, \infty)$ be random variables, assume for every $N \in \N$ that $V_i^N$, $i \in \{1, 2, \ldots, N\}$, are i.i.d.\ and that $W_i^N$, $i \in \{1, 2, \ldots, N\}$, are i.i.d., assume for all $N \in \N$, $x \in (0, \infty)$ that $\P (N^\alpha V_1^N \le x) = \P (N^\beta W_1^N \le x) = [\frac{2}{\pi}]^{1/2} \int_0^x \exp (- \frac{y^2}{2}) \, \d y$, and let $\fC \in \R$, $\gamma \in (\nicefrac{3}{4}, 1)$, $\eta \in (0 , \min \cu{4 \gamma - 3, \gamma - \frac{1}{2}})$. Then there exists $c \in \R$ which satisfies for all $N \in \N$ that
\begin{equation}
\label{cor:product:bias:sum:eqclaim}
\textstyle \P \rbr[\big]{\Exists J \subseteq \{1, 2, \ldots, N\} \colon \br[\big]{\# J \ge N^\gamma, \, \sum_{i \in J} V_i^N W_i^N \le \fC \max_{i \in \{1, 2, \ldots, N\}} V_i^N W_i^N}} \le c N^{- \eta}.
\end{equation}
\end{cor}

\cfclear
\begin{cproof}{cor:product:bias:sum}
Throughout this proof for every $N \in \N$, $i \in \{1, 2, \ldots, N\}$ let $X_i^N = V_i^N W_i^N \colon \Omega \to (0, \infty)$ be a random variable, and for every $N \in \N$ let $X_{(1)}^N, \ldots, X_{(N)}^N \colon \Omega \to (0, \infty)$ be the order statistics of $X_1^N, \ldots, X_N^N$. \Nobs that by homogeneity we may assume without loss of generality that $\alpha = \beta = 0$. Moreover, we have for every $N \in \N$ that $X_1^N, \ldots, X_N^N$ are i.i.d.\ random variables. Next, \cref{prop:product:normal} assures for all $N \in \N$, $i \in \{1, 2, \ldots, N\}$, $x \in [0, \infty)$ that $\P(X_i^N \le x) = \frac{2}{\pi} \int_0^x K_0(y) \, \d y$\cfadd{def:bessel} \cfload. Furthermore, the fact that $\gamma > \nicefrac{3}{4}$ and the fact that $\eta < \min \cu{4 \gamma - 3, \gamma - \frac{1}{2}}$ demonstrate that there exist $\varepsilon, \delta \in (0, \infty)$ which satisfy
\begin{equation}
\label{cor:product:bias:sum:eq:defdelta}
\textstyle \varepsilon < 2 \gamma - 1, \quad \varepsilon > 2 \eta, \quad \delta < \frac{2 \gamma - \varepsilon - 1}{1 - \gamma}, \qandq 3 + 2 \delta - 4 \gamma - 2 \gamma \delta < - \eta.
\end{equation}
In the next step \nobs that the fact that for all $N \in \N$, $J \subseteq \{1, 2, \ldots, N\}$ with $\# J \ge N^\gamma$ it holds that $\sum_{i \in J} V_i^N W_i^N = \sum_{i \in J } X_i^N \ge \sum_{i=1}^{ \lceil N^\gamma \rceil} X_{(i)}^N$ shows that
\begin{equation} \label{cor:product:bias:sum:eq1}
\begin{split}
& \textstyle \P \rbr[\big]{\exists \, J \subseteq \{1, 2, \ldots, N\} \colon \br*{\# J \ge N^\gamma, \, \sum_{i \in J} V_i^N W_i^N \le \fC \max_{i \in \{1, 2, \ldots, N\}} V_i^N W_i^N}} \\
& \textstyle \le \P \rbr*{\sum_{i=1}^{ \lceil N^\gamma \rceil} X_{(i)}^N \le \fC X_{(N)}^N} \le \P \rbr*{\sum_{i=1}^{\lceil N^\gamma \rceil} X_{(i)}^N \le \fC N^\eps} + \P \rbr[\big]{X_{(N)}^N \ge N^\eps}.
\end{split}
\end{equation}
Next, \cref{lem:size:estimate} (applied with $\alpha \with 0 $, $\eps \with \frac{\eps}{2}$, $\delta \with \eta$ in the notation of \cref{lem:size:estimate}) and \cref{cor:product:bias:sum:eq:defdelta} ensure that there exists $c \in \R$ which satisfies for all $N \in \N$ that
\begin{equation}
\begin{split}
& \textstyle \P \rbr[\big]{X_{(N)}^N \ge N^{\varepsilon}} \\
& \textstyle \le \P \rbr[\big]{\max_{i \in \{1, 2, \ldots, N\}} V_i^N \ge N^{\nicefrac{\varepsilon}{2}}} + \P \rbr[\big]{\max_{i \in \{1, 2, \ldots, N\}} W_i^N \ge N^{\nicefrac{\varepsilon}{2}}} \le c N^{- \eta}.
\end{split}
\end{equation}
Combining this with \cref{cor:sum:order:stat} assures that there exists $\kappa \in \R$ which satisfies for all $N \in \N$ that
\begin{equation}
\textstyle \P \rbr[\Big]{\ssum_{i = 1}^{\lceil N^\gamma \rceil} X_{(i)}^N \le \fC N^\varepsilon} + \P \rbr[\big]{X_{(N)}^N \ge N^\varepsilon} \le \kappa N^{3 + 2 \delta - 4 \gamma - 2 \gamma \delta} + c N^{- \eta}.
\end{equation}
This, \cref{cor:product:bias:sum:eq:defdelta}, and \cref{cor:product:bias:sum:eq1} establish \cref{cor:product:bias:sum:eqclaim}.
\end{cproof}

\subsubsection{Properties of random ANN initializations}

Now we are in a position to establish in \cref{lem:random:init} some important properties of the random initializations we consider. In the notation of \cref{lem:random:init}
the random variables $\Theta^{\width , i }$, $1 \le i \le \width$, describe the inner biases,
the random variables $\Theta^{\width , \width + i }$, $1 \le i \le \width$, describe the inner weights,
and $\Theta^{\width , 2 \width + i }$, $1 \le i \le \width$,
describe the outer weights.
While the biases are normally distributed, the weights are half-normally distributed with appropriate scaling in the width $\width$.

\begin{lemma}
\label{lem:random:init}
Let $\0 \in \R$, $\1 \in (\0, \infty)$, $\alpha \in (\nicefrac{3}{4}, 1)$, $\beta \in ( \alpha + 2 , \infty )$, $\gamma \in (\nicefrac{3}{4} , \alpha )$, for every $\width \in \N$, $i \in \{1, 2, \ldots, \width\}$, $\theta = ( \theta_1, \ldots, \theta_{3 \width } ) \in \R^{ 3 \width }$ let $I_i^\theta \subseteq \R$ satisfy $I_i^\theta = \cu{x \in ( \0 , \1 ) \colon \theta_{\width + i } x + \theta_{ i } \in (0, 1 )}$, let $\Theta^\width = (\Theta^{\width, 1}, \ldots, \Theta^{ \width , 3 \width } ) \colon \Omega \to \R^{ 3 \width }$, $\width \in \N$, be random variables, assume for all $\width \in \N$ that $\Theta^{ \width , 1 } , \ldots, \Theta^{ \width , 3 \width }$ are independent, assume for all $\width \in \N$, $k \in \{1, 2, \ldots, \width\}$ that $\width^{- \beta} \Theta^{ \width , k }$ is standard normal, and assume for all $\width \in \N$, $k \in \{1, 2, \ldots, \width\}$, $x \in (0, \infty)$ that $\P(\width^{- \beta} \Theta^{\width, \width + k} \leq x) = \P(\width^{\alpha} \Theta^{\width, 2 \width + k} \le x) = [\frac{2}{\pi}]^{1/2} \int_0^x \exp(- \frac{y^2}{2}) \, \d y$. Then
\begin{enumerate} [ label = (\roman*)]
\item
\label{lem:random:init:item1} there exist $c_1 \in (0, 1)$, $c_2 \in (c_1, 1)$ such that for all $\width \in \N$, $i \in \{1, 2, \ldots, \width\}$ it holds that $c_1 \le \P(I_i^{\Theta^\width} \neq \emptyset) \le c_2$,

\item
\label{lem:random:init:item2} for all $\delta \in (0, \frac{\alpha - \gamma}{2})$ there exists $c_3 \in \R$ such that for all $\width \in \N$ it holds that
\begin{equation}
\textstyle \P \rbr[\big]{\max_{i \in \{1, 2, \ldots, \width\}} \Theta^{ \width , 2 \width + i } \le \width ^{ - \gamma - \delta } } > 1 - c_3 \width ^{ - \delta },
\end{equation}

\item
\label{lem:random:init:itemnew} for all $\fC \in \R$, $\delta \in (0, \frac{\alpha - \gamma }{2} )$, $\eta \in (0 , \min \cu{ 4 \gamma - 3 , \gamma - \tfrac{1}{2} } )$ there exists $c_4 \in \R$ such that for all $\width \in \N$ it holds that
\begin{multline}
\textstyle \P\big(\Exists J \subseteq \{1, 2, \ldots, \width\} \colon
\big[\rbr[\big]{\sum_{i \in J} \Theta^{\width, \width + i} \Theta^{\width, 2 \width + i} \le \fC \max_{i \in \{1, 2, \ldots, \width\}} \Theta^{\width, \width + i} \Theta^{\width, 2 \width + i}} \\
\textstyle \wedge \rbr[\big]{\sum_{i \in J} \Theta^{\width, 2 \width + i} \ge \width^{- \delta}}\big]\big) \le c_4 \width^{- \min \cu{\delta, \eta}},
\end{multline}

\item
\label{lem:random:init:item3} for all $q \in [\alpha^{-1} , \infty)$ there exists $c_5 \in \R$ such that for all $\width \in \N$ it holds that
\begin{equation}
\textstyle \rbr[\Big]{\E \br[\Big]{\abs[\big]{\sum_{i=1}^\width \Theta^{\width, 2 \width + i}}^{2 q}}}^{\nicefrac{1}{q}} \le c_5 \width^{2 - 2 \alpha},
\end{equation}

\item
\label{lem:random:init:item4} for all $\kappa \in (0, \infty)$ there exists $c_6 \in (0, \infty )$ such that for all $\width \in \N$ it holds that
\begin{equation}
\textstyle \P \rbr*{\sum_{i \in \{1, 2, \ldots, \width\}, \, I_i^{\Theta ^\width } \neq \emptyset} \Theta^{\width, 2 \width + i} \le \kappa } \le \frac{c_6}{\width},
\end{equation}

and
\item 
\label{lem:random:init:item5} for all $\rho \in (0, \infty)$ there exists $c_7 \in \R$ such that for all $\width \in \N$ it holds that $\P \rbr{\min_{i \in \{1, 2, \ldots, \width\}} \allowbreak \Theta^{\width, \width + i} \allowbreak \ge \allowbreak \rho} \ge 1 - c_7 \width^{1 - \beta}$ and $\P \rbr{\min_{i \in \{1, 2, \ldots, \width\}} \Theta^{\width, \width + i} \Theta^{\width, 2 \width + i} \ge \rho} \ge 1 - c_7 \width^{(2 - \beta + \alpha)/2}$.
\end{enumerate}
\end{lemma}
\begin{cproof}{lem:random:init}
First, the fact that the distribution of $(\width^{- \beta} \Theta^{\width, i}, \width^{- \beta} \Theta^{\width, \width + i}) \colon \Omega \to \allowbreak [0, \infty) \allowbreak \times \R$, $\width \in \N$, $i \in \{1, 2, \ldots, \width\}$, is absolutely continuous with a positive density and independent of $\width \in \N$, $i \in \{1, 2, \ldots, \width\}$ and the fact that for all $\width \in \N$, $i \in \{1, 2, \ldots, \width\}$ it holds that
\begin{equation}
\begin{split}
 \textstyle 0 & < \P \rbr[\big]{ \width ^{ - \beta } ( \0 \Theta^{ \width , \width + i } + \Theta^{ \width , i } ) < 0 < \width ^{ - \beta } ( \1 \Theta^{ \width , \width + i } + \Theta^{ \width , i } ) } \\
 & \le \textstyle \P (I_i^{ \Theta^\width} \neq \emptyset) \le 1 -  \P \rbr[\big]{ \width ^{ - \beta } ( \1 \Theta^{ \width , \width + i } + \Theta^{ \width , i } < 0)} < 1
\end{split}
\end{equation}
establish \cref{lem:random:init:item1}. In the next step \nobs that \cref{lem:size:estimate} (applied for every $\width \in \N$, $\delta \in (0, \frac{\alpha - \gamma }{2} )$ with $N \with \width$, $\0 \with \alpha$, $\varepsilon \with \alpha - \gamma - \delta$, $\delta \with \delta$ in the notation of \cref{lem:size:estimate}) shows for all $\delta \in (0, \frac{\alpha - \gamma}{2})$, $\width \in \N$ that
\begin{equation}
\textstyle \P \rbr[\big]{ \max\nolimits_{i \in \{1, 2, \ldots, \width\}} \Theta^{ \width , 2 \width + i } \ge \width ^{- \gamma - \delta } } \le (\alpha - \gamma - 2 \delta ) ^{-1} \width^{- \delta}.
\end{equation}
This establishes \cref{lem:random:init:item2}. Next \nobs that for all $\width \in \N$, $J \subseteq \{1, 2, \ldots, \width\}$, $\omega \in \Omega$ it holds that $\# J \ge \width^\gamma$ or $\sum_{j \in J } \Theta^{ \width , 2 \width + j } ( \omega ) \le \width ^\gamma \max_{i \in \{1, 2, \ldots, \width\}} \Theta^{ \width , 2 \width + i } ( \omega)$. Hence, we obtain for all $\width \in \N$, $\fC \in \R$, $\delta \in (0, \frac{\alpha - \gamma}{2})$ that
\begin{equation}
\begin{split}
& \textstyle \P \rbr[\big]{\Exists J \subseteq \{1, 2, \ldots, \width\} \colon \textstyle \br[\big]{\sum_{i \in J } \Theta^{\width, \width + i} \Theta^{\width, 2 \width + i } \le \fC \max_{i \in \{1, 2, \ldots, \width\}} \Theta^{ \width , \width + i } \Theta^{\width , 2 \width + i } } \\
& \textstyle \qquad  \wedge \br[\big]{\sum_{i \in J} \Theta^{\width , 2 \width + i } \ge \width ^{ - \delta } } } \\
& \textstyle \le \P \rbr[\big]{\Exists J \subseteq \{1, 2, \ldots, \width\} \colon \br[\big]{\# J \ge \width^\gamma} \wedge \\
& \textstyle \qquad  \br[\big]{\sum_{i \in J} \Theta^{ \width , \width + i } \Theta^{ \width , 2 \width + i} \le \fC \max_{i \in \{1, 2, \ldots, \width\}} \Theta^{\width, \width + i} \Theta^{ \width , 2 \width + i } } }  \\
& \textstyle + \P \rbr[\big]{\max_{i \in \{1, 2, \ldots, \width\}} \Theta^{\width, 2 \width + i} \ge \width^{- \gamma - \delta }}.
\end{split}
\end{equation}
Combining this with \cref{lem:random:init:item2} and \cref{cor:product:bias:sum} demonstrates that for all $\fC \in \R$, $\delta \in (0 , \frac{ \alpha - \gamma}{2} )$, $\eta \in (0, \min \cu{4 \gamma - 3, \gamma - \tfrac{1}{2}, \delta})$ there exist $c_3, c_4 \in \R$ which satisfy for all $\width \in \N$ that
\begin{multline}
\textstyle \P \big(\Exists J \subseteq \{1, 2, \ldots, \width\} \colon \br[\big]{\sum_{i \in J} \Theta^{ \width , \width + i } \Theta^{ \width , 2 \width + i } \le \fC \max\nolimits_{i \in \{1, 2, \ldots, \width\}} \Theta^{ \width , \width + i } \Theta^{ \width , 2 \width + i }  } \\ 
\textstyle \wedge \br[\big]{\sum_{i \in J} \Theta^{\width, 2 \width + i} \ge \width ^{ - \delta}}\big) \le c_4 \width^{- \eta} + c_3 \width^{- \delta} \le (c_3 + c_4) \width^{- \min \cu{\delta , \eta}}.
\end{multline}
This establishes \cref{lem:random:init:itemnew}. Next, by the well-known formula for absolute moments of the normal distribution (cf., e.g., (18) in Winkelbauer~\cite{Winkelbauer2012}), for every $q \in [1, \infty )$ there exists $\fC_q \in (0, \infty )$ which satisfies for all $\width \in \N$ that $\E [ ( \Theta^{ \width , 2 \width + 1 } ) ^{ 2 q } ] = \fC_q \width^{- 2 q \alpha}$. Combining this with the fact that for all $q \in [1, \infty)$ it holds that $[0, \infty) \ni x \mapsto x^{2q} \in [0, \infty)$ is convex and Jensen's inequality assures that for all $q \in [1, \infty)$, $\width \in \N$ we have that
\begin{equation}
\textstyle \E\br*{\abs[\big]{\sum_{i=1}^\width \Theta^{\width, 2 \width + i}}^{2 q}} \le \E\br*{\width^{2 q - 1} \sum_{i = 1}^\width (\Theta^{\width, 2 \width + i})^{2 q}} = \width^{2 q} \E\br[\big]{(\Theta^{\width, 2 \width + 1})^{2 q}} = \fC_q \width^{2 q (1 - \alpha)}.
\end{equation}
This establishes \cref{lem:random:init:item3}. For every $\width \in \N$ let $X_i^\width \colon \Omega \to \R$, $i \in \{1, 2, \ldots, \width\}$, and $Z^\width \colon \Omega \to \R$ be random variables which satisfy for all $i \in \{1, 2, \ldots, \width\}$ that \begin{equation}
\textstyle X_i^\width (\omega) =
\begin{cases}
\Theta^{\width, 2 \width + i}(\omega) & \colon I_i^{\Theta^\width} \neq \emptyset  \\
0 & \colon \text{else}
\end{cases}
\end{equation}
and $Z^\width = \sum_{i=1}^\width X_i$. \Nobs that for all $\width \in \N$ it holds that $X_i^\width$, $i \in \{1, 2, \ldots, \width\}$, are i.i.d.\ random variables. Furthermore, \cref{lem:random:init:item1} and the fact that for all $\width \in \N$, $i \in \{1, 2, \ldots, \width\}$ it holds that there exists a standard normal random variable $U$ such that $\Theta^{\width, 2 \width + i} \sim \width^{- \alpha} \abs{U}$ imply that for all $\width \in \N$, $i \in \{1, 2, \ldots, \width\}$ we have that
\begin{equation}
\textstyle \E[X_i^\width] = \E[\Theta^{\width, 2 \width + i}] \, \P(I_i^{\Theta^\width} \neq \emptyset) \ge c_1 \width^{- \alpha} [\frac{2}{\pi}]^{1/2} \qqandqq \operatorname{Var}(X_i^\width) \le \E[(X_i^{\width})^2] \le \width^{-2 \alpha}.
\end{equation}
Hence, we obtain that there exists $\lambda \in (0, \infty)$ which satisfy for all $\width \in \N$ that $\E[Z^\width] \ge \lambda \width^{1 - \alpha}$ and $\operatorname{Var}(Z^\width) \le \width^{1 - 2 \alpha}$. Combining this with Chebyshev's inequality assures for all $\kappa \in (0, \infty)$, $\width \in \N \cap ((2 \kappa)^{(1 - \alpha)^{-1}} \lambda^{- (1 - \alpha)^{-1}}, \infty)$ that $\lambda \width^{1 - \alpha} - \kappa \ge \frac{\lambda \width^{- \alpha}}{2}$ and, therefore,
\begin{equation}
\begin{split}
\textstyle \P(Z^\width \le \kappa) & \textstyle = \P\rbr[\big]{\lambda \width^{1 - \alpha} - Z^\width \ge \lambda \width^{1 - \alpha} - \kappa} \le \P\rbr[\big]{\E[Z^\width] - Z^\width \ge \lambda \width^{1 - \alpha} - \kappa} \\
& \textstyle \le \P\rbr[\big]{\abs{Z^\width - \E[Z^\width]} \ge \frac{\lambda \width^{- \alpha}}{2}} \le \frac{4 \width^{1 - 2 \alpha}}{\lambda^2 \width^{- 2 \alpha}} = \frac{4}{\lambda^2 \width}.
\end{split}
\end{equation}
This and the fact that for all $\kappa \in (0, \infty)$, $\width \in \N \cap [1, (2 \kappa)^{(1 - \alpha)^{- 1}} \lambda^{- (1 - \alpha)^{-1}}]$ it holds that $\P(Z^\width \le \kappa) \le 1 \le (\frac{2 \kappa}{\lambda})^{(1 - \alpha)^{- 1}} \width^{- 1}$ demonstrate that there exists $c_5 \in \R$ which satisfies for all $\kappa \in (0, \infty )$, $\width \in \N$ that $\P(Z^\width \le \kappa ) \le c_5 \width^{- 1}$. This establishes \cref{lem:random:init:item4}.

Next, the fact that for all $\width \in \N$, $x \in [0, \infty)$ it holds that $\P(\Theta^{\width, \width + 1} \le x) = \P(\width^{- \beta} \Theta^{\width, 1} \le \width^{- \beta} x) \allowbreak = \allowbreak [\frac{2}{\pi}]^{1/2} \int_{0}^{\width^{- \beta} x} \exp(- \frac{y^2}{2}) \, \d y \le \width^{- \beta} x$ and $\P(\Theta^{\width, 2 \width + 1} \le x) \le \width^{\alpha} x$, the fact that for all $\rho \in (0, \allowbreak \infty), \allowbreak \width \allowbreak \in \allowbreak \N \allowbreak \cap [\rho^{(\beta - \alpha)^{-1}}, \infty)$ it holds that $\rho^{1/2} \width^{- \beta} \width^{(\alpha + \beta)/2} \le 1$, and Bernoulli's inequality establish that for all $\rho \allowbreak \in \allowbreak (0, \allowbreak \infty)$, $\width \in \N \cap [\rho^{(\beta - \alpha)^{- 1}}, \infty)$ we have that
\begin{equation}
\begin{split}
& \textstyle \P \rbr[\big]{\min_{i \in \{1, 2, \ldots, \width\}} \Theta^{\width, \width + i} \Theta^{\width, 2 \width + i} \ge \rho} \\
& \textstyle = \P \rbr[\big]{\Forall i \in \{1, 2, \ldots, \width\} \colon \Theta^{\width, \width + i} \Theta^{\width, 2 \width + i} \ge \rho} = \br[\big]{\P\rbr[\big]{\Theta^{\width, \width +  1} \Theta^{\width, 2 \width + 1} \ge \rho}}^\width \\
& \textstyle \ge \br[\big]{\P \rbr[\big]{\Theta^{\width, \width + 1} \ge \rho^{1/2} \width^{(\alpha + \beta)/2}} \P \rbr[\big]{\Theta^{\width, 2 \width + 1} \ge \rho^{1/2} \width^{- (\alpha + \beta)/2}}}^\width \\
& \textstyle \ge \rbr[\big]{1 - \rho^{1/2} \width^{- \beta} \width^{(\alpha + \beta)/2}}^\width \rbr[\big]{1 - \rho^{1/2} \width^{\alpha} \width^{- (\alpha + \beta)/2}}^\width \\
& \textstyle = \rbr[\big]{1 - \rho^{1/2} \width^{(\alpha - \beta)/2}}^{2 \width} \ge \rbr[\big]{1 - \rho^{1/2} \width^{(2 + \alpha - \beta)/2}}^2 \ge 1 - 2 \rho^{1/2} \width^{(2 + \alpha - \beta)/2}.
\end{split}
\end{equation}
This and the fact that for all $\rho \in (0, \infty)$ there exists $c_\rho \in \R$ which satisfies for all $\width \in \N \cap (0, \allowbreak \rho^{(\beta - \alpha)^{- 1}}]$ that $c_\rho \width^{(2 + \alpha - \beta)/2} \ge 1$ demonstrate for all $\rho \in (0, \infty)$, $\width \in \N$ that $\P \rbr[\big]{\min_{i \in \{1, 2, \ldots, \width\}} \allowbreak \Theta^{\width, i} \allowbreak \Theta^{\width, 2 \width + i} \allowbreak \ge \rho} \ge 1 - \max \cu{2 \rho^{1/2}, c_\rho} \width^{(2 + \alpha - \beta)/2}$. In the same way, we can use Bernoulli's inequality to show that for every $\rho \in (0, \infty)$ there exists $\const \in \R$ which satisfies for all $\width \in \N$ that $\P \rbr{\min_{i \in \{1, 2, \ldots, \width\}} \Theta^{\width, \width + i} \ge \rho} \ge 1 - \const \width^{ 1 - \beta}$. This establishes \cref{lem:random:init:item5}.
\end{cproof}

%------------------------------------------------------------------------------%
%----------------------------------Subsection----------------------------------%
%------------------------------------------------------------------------------%
\subsubsection{Almost sure non-convergence of GF to descending critical points}

Our final building block is the non-convergence of randomly initialized gradient flows to descending critical points in the sense of \cref{def:strict_saddle}.
For this we use \cref{prop:non-convergence_to_strict_saddle} below, which is well-known in the scientific literature.

\cfclear
\begin{proposition}\label{prop:non-convergence_to_strict_saddle}
	Let $\fh \in \N$, $f \in C^2 (\R^{\fh}, [0, \infty))$, for every $x \in \R^{\fh}$ let $X^x = (X_t^x)_{t \in [0, \infty)} \in C([0, \infty), \R^{\fh})$ satisfy for all $t \in [0, \infty)$ that 
	\begin{equation}\label{eqn:prop:non-convergence_to_strict_saddle:gf}
	\textstyle X_0^x = x \qqandqq X_t^x = X_0^x - \int_0^{t} (\nabla f)(X_s^x) \, \d s,
	\end{equation}
	and let $\lambda \colon \cB(\R^{\fh}) \to [0, \infty]$ be the Lebesgue--Borel measure on $\R^{\fh}$. Then
	\begin{equation}\label{eqn:prop:non-convergence_to_strict_sadd:zero_set}
	\textstyle \lambda \bigl(\bigl\{x \in \R^{\fh} \colon \bigl[\Exists y, v \in \R^{\fh} \colon (\limsup_{t \to \infty} \norm{X_t^x - y} = 0  > \scalar{v, ((\operatorname{Hess} f)(y)) v})\bigr]\bigr\}\bigr) = 0.
	\end{equation}
\end{proposition}
\begin{cproof}{prop:non-convergence_to_strict_saddle}
	\Nobs that Bah et al.~\cite[Theorem~28]{BahRauhutTerstiege2021} (applied with $\cM \with \R^{\fh}$, $L \with f$ in the notation of \cite[Theorem~28]{BahRauhutTerstiege2021}) establishes \cref{eqn:prop:non-convergence_to_strict_sadd:zero_set}.
\end{cproof}

We also refer to Lee et al.~\cite{MR3960812} for a similar result with discrete-time \GD. We now formulate a corollary for the case where only some parameters are trained and the remaining ones are initialized randomly and independently, which is suitable for our \ANN\ setting with untrained weights. In combination with \cref{cor:non:strict:saddle} this will allow us to estimate the risk of the critical limit point.

\cfclear 
\begin{cor}
\label{cor:strict:saddle:conv:random}
Let $\fn , \fh \in \N$,
$f \in C ( \R^{ \fn } \times \R^{ \fh } , [ 0 , \infty ) )$,
let $\cG \colon \R^\fn \times \R^\fh \to \R^\fh$ be measurable and satisfy for all $x \in \R^\fn$,
$y \in \cu{z \in \R^\fh \colon f ( x , \cdot ) \text{ is differentiable at } z}$
that $\cG ( x , y ) = (\nabla _y f) ( x , y )$,
let $\lambda \colon \cB(\R^{\fh}) \to [0, \infty]$ be the Lebesgue--Borel measure on $\R^{\fn}$,
let $A \subseteq \R^\fn$ be measurable and satisfy $\lambda ( \R^\fn \backslash A ) = 0$,
assume for all $x \in A$ that $f ( x , \cdot ) \in C^2 ( \R^\fh , \R )$,
let $\cS_x \subseteq \R^\fh$, $x \in \R^\fn$,
satisfy
for every $x \in A$
that $\cS_x = \cu{y \in \R^\fh \colon y \text{ is a \descritic critical point of } f(x, \cdot ) }$\cfadd{def:strict_saddle}, let $(\Omega , \cF , \P)$ be a probability space,
let $X \colon [0 , \infty ) \times \Omega \to \R^\fn \times \R^\fh$ and $Y \colon [0 , \infty ) \times \Omega \to \R^\fn \times \R^\fh$ be stochastic processes with continuous sample paths,
assume that $X_0 \colon \Omega \to \R^\fn$ and $Y_0 \colon \Omega \to \R^\fh$ are independent random variables,
assume that the distribution of $Y_0$ is absolutely continuous with respect to the Lebesgue--Borel measure on $\R^\fh$,
assume that $\P ( X_0 \in A  ) = 1$,
and assume for all $t \in [0 , \infty )$, $\omega \in \Omega $
that
\begin{equation}
\label{cor:strict:saddle:conv:random:eq1}
\textstyle
X_t ( \omega ) = X_0 ( \omega ) \qqandqq Y_t ( \omega ) = Y_0 ( \omega ) - \int_0^t \cG ( X_s ( \omega ) , Y_s ( \omega ) ) \, \d s\ifnocf.
\end{equation}
\cfload. Then
\begin{equation}
\label{cor:strict:saddle:conv:random:eqclaim}
\textstyle
\P \rbr*{ \exists \, ( x , y ) \in \R^\fn \times \R^\fh \colon \br*{ \lim_{t \to \infty} X_t = x }
\wedge \br*{ \lim_{t \to \infty} Y_t = y }
\wedge \br*{ y \in \cS_x } } = 0 .
\end{equation}
\end{cor}
\begin{cproof}{cor:strict:saddle:conv:random}
First \nobs that \cref{cor:strict:saddle:conv:random:eq1} and the fact that $\P ( X_0 \in A ) = 1$ ensure that
\begin{equation}
\begin{split}
& \textstyle \P \rbr*{ \exists \, ( x , y ) \in \R^\fn \times \R^\fh \colon \br*{ \lim_{t \to \infty} X_t = x }
\wedge \br*{ \lim_{t \to \infty} Y_t = y }
\wedge \br*{ y \in \cS_x } } \\
& \textstyle = \P \rbr*{ \br*{ X_0 \in A } \wedge \br*{\exists \, y \in \cS_{ X_0 } \colon \lim_{t \to \infty} Y_t = y  } } .
\end{split}
\end{equation}
Furthermore, the fact that for all $x \in A$ it holds that $f ( x , \cdot ) \in C^2 ( \R^\width , \R )$,
the fact that for all $\omega \in \Omega$,
$t \in [0 , \infty )$ with $X_0 ( \omega ) \in A$
it holds that
$Y_t ( \omega ) = Y_0 ( \omega ) - \int_0^t \nabla_y f ( X_0 ( \omega ) , Y_s ( \omega ) ) \, \d s$,
and the Picard-Lindelöf uniqueness theorem for ODEs demonstrate that 
there exists a measurable function $\Phi \colon \R^\fn \times \R^\fh \to C( [ 0 , \infty ) , \R^\fh )$
which satisfies for all $\omega \in \Omega$, $t \in [0 , \infty )$ with $X_0 ( \omega ) \in A$
that $Y_t ( \omega ) = \Phi_t ( X_0 ( \omega ) , Y_0 ( \omega ) )$.
In particular, there exists a measurable function $F \colon \R^\fn \times \R^\fh \to \R$ which satisfies
\begin{equation}
F ( X_0, Y_0 ) = \begin{cases}
1 & \colon \br*{ X_0 \in A } \wedge \br*{\exists \, y \in \cS_{ X_0 } \colon \lim_{t \to \infty} Y_t = y  } \\
0 & \colon \text{else.}
\end{cases}
\end{equation}
In addition, \cref{prop:non-convergence_to_strict_saddle} and the assumption that the distribution of $Y_0$ is absolutely continuous with respect to the Lebesgue measure on $\R^\fh$ show for all $x \in A$ that
\begin{equation}
\textstyle
\E [ F ( x , Y_0 ) ] = \P \rbr*{ \exists \, y \in \cS_x \colon \lim_{t \to \infty} \Phi_t ( x , Y_0 ) = y } = 0.
\end{equation}
Combining this with the fact that $X_0$ and $Y_0$ are independent and Lemma 2.8 in the \href{https://arxiv.org/abs/1801.09324}{arXiv-version} of \cite{JentzenKuckuckNeufeldVonWurstemberger2021}
implies that
\begin{equation}
\textstyle \P \rbr*{ \br*{ X_0 \in A } \wedge \br*{\exists \, y \in \cS_{ X_0 } \colon \lim_{t \to \infty} Y_t = y  } } = \E [ F ( X_0 , Y_0 ) ] = \E \br*{ \indicator{A} ( X_0 ) F ( X_0 , Y_0 ) } = 0.
\end{equation}
\end{cproof}

\subsubsection{Convergence of GF with random initialization}
\label{sec:convergence_of_GF_with_random_init_clipping}

In \cref{theo:main:conv} we establish the convergence of \GF\ with random initializations, which is one of the main results of this article. \cref{theorem_intro} in the introduction is a direct consequence of \cref{theo:main:conv}.
 For convenience of the reader, \cref{theo:main:conv} is formulated in a self-contained way.

\cfclear
\begin{theorem} \label{theo:main:conv}
Let $\0 \in \R$, $\1 \in (\0 , \infty)$, let $f, \dens \in C(\R , \R)$ be piecewise analytic\cfadd{def:piecewise_analytic}, assume for all $x \in \R$ that $\dens(x) \ge 0$ and $\dens^{- 1}(\R \backslash \{0\}) = (\0, \1)$, assume that $f$ is strictly increasing, let $(\Omega, \cF, \P)$ be a probability space, for every $\width \in \N$, $\theta = ( \theta_1, \ldots, \theta_{\fd_\width } ) \in \R^{\fd_\width}$ let $\fd _\width = 3 \width + 1$, let $\fN^{ \theta}_\width \colon \R \to \R$
satisfy for all $x \in \R$ that
\begin{equation}
\textstyle \fN^{ \theta }_\width (x) = \theta_{\fd_\width} + \sum_{i = 1}^{\width} \theta_{2 \width + i} \max \cu{\min\cu{\theta_{\width + i } x + \theta_{ i }, 1}, 0},
\end{equation}
let $\loss_\width \colon \R^{ \fd_\width } \to \R$ satisfy
\begin{equation} \label{theo:main:conv:eq:risk}
\textstyle \loss_{\width}(\theta) = \int_{\R}
\rbr[\big]{\fN^{\theta}_\width (x) - f(x)}^2 \dens(x) \, \d x,
\end{equation}
let $\cG_\width = (\cG_\width^1 , \ldots, \cG_\width^{\fd_\width}) \colon \R^{\fd_\width} \to \R^{\fd_\width}$ satisfy for all $k \in \{1, 2, \ldots, \fd_\width\}$ that $\cG_\width^k (\theta) \allowbreak = \allowbreak (\frac{\partial^-}{\partial \theta_k} \loss_\width) (\theta) \allowbreak \indicator{\N \backslash (\width, 3 \width]} (k)$, and let $\Theta^\width = (\Theta^{\width, 1}, \ldots, \Theta^{\width, \fd_\width}) \colon [0, \infty) \times \Omega \to \R^{\fd_\width}$
be a stochastic process with continuous sample paths, assume for all $\width \in \N$, $t \in [0 , \infty )$ that 
\begin{equation} \label{theo:main:conv:eq:assgf}
\textstyle \P \rbr[\big]{\Theta_t^\width = \Theta_0^\width - \int_0^t \cG_\width (\Theta_s^\width) \, \d s} = 1,
\end{equation}
assume for all $\width \in \N$ that $\Theta_0^{\width, 1}, \ldots, \Theta_0^{\width, \fd_\width}$ are independent, let $\alpha \in (\nicefrac{3}{4}, 1)$, $\beta \in (\alpha + 2, \infty)$, assume for all $\width \in \N$, $k \in \{1, 2, \ldots, \width\}$ that $\width^{- \beta} \Theta_0^{\width, k}$ is standard normal,
assume for all $\width \in \N $, $k \in \{1, 2, \ldots, \width\}$, $x \in (0, \infty)$ that
\begin{equation}\label{eqn:theo:main:conv:weights}
\textstyle  \P\rbr[\big]{\width^{- \beta } \Theta_0^{ \width, \width + k } \leq x} = \P\rbr[\big]{\width^{\alpha} \Theta_0^{\width, 2 \width + k} \le x} = [\frac{2}{\pi}]^{1/2} \int_0^x \exp (- \frac{y^2}{2}) \, \d y,
\end{equation}
assume for all $\width \in \N$
that $\E [ \abs{ \Theta_0^{ \width , \fd_\width } } ^2 ] < \infty $,
and let $\delta \in (0, \min \cu{\frac{\beta - \alpha - 2 }{2} , \frac{ 4 \alpha - 3 }{ 9 } } ) $ \cfload. Then
\begin{enumerate} [label = (\roman*)]
\item
\label{theo:main:conv:item1}
it holds for all $\width \in \N$ that $ \sup_{t \in [0, \infty ) } \E [ \loss_\width ( \Theta_t^\width ) ]
= \E [ 
\sup_{t \in [0, \infty ) } \loss_\width ( \Theta_t^\width ) ] < \infty $,
\item
\label{theo:main:conv:item2}
there exists $c \in (0, \infty )$ such that for all $\width \in \N$ it holds that
\begin{equation}
\textstyle \P \rbr[\big]{ \limsup\nolimits_{t \to \infty} \loss_\width ( \Theta_t^\width ) \ge c}
\ge 
\P \rbr[\big]{ \inf\nolimits_{t \in [0 , \infty ) } \loss_\width ( \Theta_t^\width ) \ge c}
\ge c^\width > 0,
\end{equation}
\item
\label{theo:main:conv:item3}
there exist $\const \in (0, \infty )$
and random variables $\fC_\width \colon \Omega \to \R$,
$\width \in \N$,
which satisfy for all
$\width \in \N$
that
\begin{equation}
\textstyle \P \rbr[\big]{ \forall \, t \in (0, \infty ) \colon \loss_\width ( \Theta_t ^\width ) \le \const \width ^{ - \delta } + \fC_\width t^{-1} } \ge 1 - \const \width ^{ - \delta },
\end{equation}
and
\item
\label{theo:main:conv:item4}
for every $\varepsilon \in ( - \infty , \delta - 2 ( 1 - \alpha ) )$ there exists $\fC \in (0 , \infty )$ such that for all
$\width \in \N$
we have that
\begin{equation}
\textstyle \limsup\nolimits_{ t \to \infty }
\rbr[\big]{ \E[ 
\loss_{ \width }( 
\Theta^{ \width }_t ) ] }
\le \E \br[\big]{ \limsup\nolimits_{t \to \infty} \loss_{ \width }( 
\Theta^{ \width }_t) }
\le \fC \width ^{ - \varepsilon }.
\end{equation}
\end{enumerate}
\end{theorem}

Note that \cref{theo:main:conv:item2} shows for every width $\width \in \N$ that, with a positive probability, the risk $\loss _ \width ( \Theta_t ^\width )$ does not converge to zero.
Nevertheless, \cref{theo:main:conv:item3} reveals that the risk converges with probability close to $1$ at rate $\cO ( t^{-1} )$ to a limiting value of at most $\bfc \width ^{ - \delta }$.
Additionally, \cref{theo:main:conv:item4} establishes convergence of the risk in the $L^1$-norm with respect to the underlying probability distribution.

\begin{cproof}{theo:main:conv}
Throughout this proof let $\Lip \colon C(\R, \R) \to [0, \infty]$ satisfy for all $F \in C(\R, \allowbreak \R)$ that $\Lip(F) = \sup_{x, y \in [\0, \1], x \neq y}\frac{\abs{F(x) - F(y)}}{\abs{x - y}}$, let $\fD , \nu \in (0 , \infty )$ satisfy $\fD = \sup_{x \in [ \0 , \1 ]} \abs{ f(x) }$ and $\nu = \int_\R \dens ( x ) \, \d x$,
for every $\width \in \N$, $i \in \{1, 2, \ldots, \width\}$, $\theta = (\theta_1, \ldots, \theta_{\fd_\width } ) \in \R^{ \fd_\width}$ let $I_i^\theta \subseteq \R$ satisfy $I_i^\theta = \cu{ x \in ( \0 , \1 ) \colon \theta_{\width + i} x + \theta_{ i } \in (0 , 1 ) } $, and let $G \in \cF$ satisfy
\begin{equation} \label{theo:main:conv:eq:def:g}
\textstyle G = \cu[\big]{\omega \in \Omega \colon \forall \, t \in [0, \infty ) \colon \Theta_t^\width  (\omega)
= \Theta_0^\width ( \omega )  - \tint_0^t \cG_\width ( \Theta_s ^\width ( \omega ) ) \, \d s } .
\end{equation}
First, \cref{theo:main:conv:eq:def:g} and \cref{prop:gf:chain:rule} assure for all $\width \in \N$, $t \in [0, \infty )$, $\omega \in G$
that $\loss_\width ( \Theta_t^\width ( \omega )) \le \loss_\width ( \Theta_0^\width ( \omega ))$. This and the fact that $\P ( G ) = 1$ (cf.~\cref{theo:main:conv:eq:assgf}) show for all $\width \in \N$ that 
\begin{equation} \label{theo:main:conv:eq:supt}
\sup\nolimits_{t \in [0, \infty )} \E \br[\big]{ \loss_\width ( \Theta_t^\width ) } = \E \br[\big]{ \sup\nolimits_{t \in [0, \infty) } \loss_\width ( \Theta_t^\width ) } = \E \br[\big]{ \loss_\width ( \Theta_0 ^\width )  }.
\end{equation}
Furthermore, the assumption that $\dens^{-1}(\R \backslash \cu{0}) = (\0, \1)$ and the Cauchy-Schwarz inequality imply for all $\width \in \N$, $\theta = ( \theta_1, \ldots, \theta_{\fd_\width } ) \in \R^{ \fd_\width }$ that
\begin{equation}
\begin{split}
\textstyle \loss_\width ( \theta )
& \textstyle \le ( \width + 2 ) \rbr[\big]{
\int_\0^\1 f(x) ^2 \dens ( x ) \, \d x + \int_\0^\1 \abs{\theta_{\fd_\width } } ^2 \dens ( x ) \, \d x 
+ \sum_{i=1}^\width \int_\0 ^\1 \abs{ \theta_{2 \width + i } } ^2 \dens ( x ) \, \d x
}
\\
& \textstyle \le (\width + 2 ) \nu
\rbr[\big]{ \fD ^2 + \abs{\theta_{\fd_\width } } ^2 + \ssum_{i=1}^\width  \abs{ \theta_{ 2 \width + i } } ^2 }.
\end{split}
\end{equation}
Combining this with the fact that for all $\width \in \N$, $i \in \{1, 2, \ldots, \width\}$ it holds that $\E [ \abs{ \Theta^{\width , \fd_\width }_0 } ^2 ] < \infty$ and $\E [ \abs{ \Theta^{\width , 2 \width + i } _0 } ^2 ] = \width^{-2 \alpha }$ demonstrates that for all $\width \in \N$ we have that
\begin{equation}
\E [ \loss_\width ( \Theta^\width_0 ) ]
\le ( \width + 2 ) \nu \rbr[\big]{ \fD ^2 + \width ^{ 1 - 2 \alpha } + \E [ \abs{ \Theta^{\width , \fd_\width }_0 } ^2 ] } < \infty .
\end{equation}
This and \cref{theo:main:conv:eq:supt} establish \cref{theo:main:conv:item1}. Next, let $\xi \in (0, \infty )$ satisfy
\begin{equation}
\textstyle \xi = \inf_{\mu \in \R } \int_\0^\1 (f(x) - \mu )^2 \dens ( x ) \, \d x
\end{equation}
and for every $\width \in \N$ let $B_\width \in \cF$ satisfy 
\begin{equation}
\textstyle B_\width = \cu[\big]{\omega \in G \colon \big[\forall \, i \in \{1, 2, \ldots, \width\} \colon I_i^{\Theta_0^\width} = \emptyset\big]}.
\end{equation}
\Nobs that \cref{lem:random:init}, \cref{eqn:theo:main:conv:weights}, and the fact that for all $\width \in \N$ it holds that $\Theta_0^{ \width , i }$, $i \in \{1, 2, \ldots, \fd_\width \}$, are independent demonstrate that there exists $c \in (0,1)$ which satisfies for all $\width \in \N$ that $\P ( B_\width ) \ge c^\width $. Moreover, the fact that $\dens^{-1 } ( \R \backslash \cu{0} ) = ( \0 , \1 )$, \cref{theo:main:conv:eq:def:g}, and \cref{prop:risk:gradient} prove for all $\width \in \N$, $\omega \in B_\width$, $t \in [0, \infty )$ that $ \forall \, i \in \{1, 2, \ldots, \width\} \colon I_i^{ \Theta_t^\width } = \emptyset $. Hence, we obtain for all $\width \in \N$, $\omega \in B_\width$, $t \in [0, \infty )$ that $(\0 , \1 ) \ni x \mapsto \fN^ {\Theta_t^\width  }_\width ( x ) \in \R $ is constant. This implies for all $\width \in \N$, $\omega \in B_\width$ that
\begin{equation}
\textstyle \Forall t \in [0, \infty ) \colon \loss_\width ( \Theta_t^\width ( \omega ) ) \ge \xi > 0.
\end{equation}
Therefore, we obtain for all $\width \in \N$ that
\begin{equation}
\textstyle \P \rbr[\big]{\limsup\nolimits_{t \to \infty} \loss_\width ( \Theta_t^\width ) \ge \xi} \ge \P \rbr[\big]{ \inf\nolimits_{t \in [0 , \infty ) } \loss_\width ( \Theta_t^\width ) \ge \xi} \ge \P(B_\width)
\ge c^\width > 0.
\end{equation}
This establishes \cref{theo:main:conv:item2}. Next, the assumption that $\dens^{-1} (\R \backslash \{0\}) = (\0, \1)$ and the assumption that $\dens$ is piecewise analytic imply that there exist $\varepsilon_1, \varepsilon_2  \in (0 , \infty )$, $\varepsilon_3 \in (0 , \varepsilon_1 )$ which satisfy that $\dens _ { [ \0 , \0 + \varepsilon_1 ] }$ is strictly increasing, $\dens_{[\1 - \varepsilon_1, \1 ]}$ is strictly decreasing, and $\inf_{x \in [ \0 + \varepsilon_3 , \1 - \varepsilon_3 ] } > \varepsilon_2 \Lip ( \dens )$. In addition, the fact that $0 < \delta < (4 \alpha - 3) / 9$ and the fact that $\nicefrac{3}{4} < \alpha < 1$ assure that there exists $\gamma \in (\nicefrac{3}{4} , \alpha )$ which satisfies
\begin{equation} \label{theo:main:conv:eq:gamma}
\textstyle \delta < \min \cu*{ 4 \gamma - 3 , \gamma - \tfrac{1}{2} , \tfrac{ \alpha - \gamma}{2} }.
\end{equation}
Moreover, \cref{cor:limit_point} proves that for all $\width \in \N$, $\omega \in G$ there exist $\vartheta^\width ( \omega ) \in \R^{\fd_\width }$, $\fC_\width ( \omega ) \in \R$ which satisfy for all $t \in [0, \infty )$ that $\limsup_{s \to \infty} \norm{\Theta^\width_s ( \omega ) - \vartheta^\width ( \omega ) } = 0$ and $0 \le \loss_\width ( \Theta_t^\width ( \omega ) ) - \loss_\width ( \vartheta^\width ( \omega ) ) \le \fC_\width ( \omega ) ( 1+t)^{-1}$. 
In the following let $\bfm \in (0 , \infty )$
satisfy
\begin{equation}
\textstyle \bfm =  \nu ^{ - 1 }
\min \cu[\big]{ \int_\0^\1 ( f(\1 ) - f(x) ) \dens ( x ) \, \d x ,
	\int_\0^\1 ( f( x ) - f ( \0 ) ) \dens ( x ) \, \d x} 
\end{equation}
and for every $\width \in \N$ let $A_i^\width \in \cF$, $i \in \{1, 2, \ldots, 5\}$, satisfy
\begin{equation}
\textstyle A_1^\width = \cu*{ \omega \in G \colon \sup\nolimits_{x \in ( \0 , \1 ) } \sum_{i \in \{1, 2, \ldots, \fh\}, \, x \in  I_i^{\vartheta^\width ( \omega ) } } \vartheta^{ \width , 2 \width + i } ( \omega ) < \min \cu{ \bfm , \width ^{ - \delta } } },
\end{equation}
\begin{equation}
\textstyle A_2^\width = \cu*{ \omega \in \Omega \colon
\sum_{i \in \{1, 2, \ldots, \fh\}, \, I_i^{ \Theta_0 ( \omega ) } \not= \emptyset }
\Theta_0^{ \width , 2 \width + i } ( \omega ) 
> f ( \1 ) - f ( \0 ) + 4 \bfm  },
\end{equation}
\begin{multline}
\textstyle A_3^\width = \Big\{ \omega \in G \colon 
\sup\nolimits_{x \in ( \0 , \1 ) }  \sum_{i \in \{1, 2, \ldots, \fh\}, \, x \in I_i^{ \vartheta^\width ( \omega ) } } \vartheta^{ \width , \width + i } ( \omega ) \vartheta^{ \width , 2 \width + i } ( \omega ) \\
 > 4 \max \nolimits_{i \in \{1, 2, \ldots, \fh\} }  \vartheta^{ \width , \width + i } ( \omega ) \vartheta^{ \width , 2 \width + i } ( \omega ) \Big\},
\end{multline}
\begin{equation}
A_4^\width = \cu[\big]{ \omega \in G \colon \min \nolimits_{i \in \{1, 2, \ldots, \fh\} } \Theta_0^{ \width , \width + i } ( \omega ) \Theta_0^{ \width , 2 \width + i } ( \omega ) > \Lip ( f ) },
\end{equation}
and
\begin{equation}
\textstyle A_5^\width = \cu[\big]{ \omega \in G \colon \min\nolimits_{i \in \{1, 2, \ldots, \fh\} } \Theta_0^{ \width , \width + i } ( \omega ) > \max \cu{ ( \varepsilon_2 ) ^{-1} , 2 ( \varepsilon_1 - \varepsilon_3 ) ^{-1} } }.
\end{equation}
\Nobs that the fact that $\forall \, \width \in \N , \, \omega \in  G \colon \cG_\width ( \vartheta^\width ( \omega ) ) = 0 $ and \cref{cor:non:strict:saddle} demonstrate for all $\width \in \N$, $\omega \in A_3^\width \cap A_4^\width \cap A_5^\width$ that the Hessian of $\loss_ \width ( \vartheta^\width ( \omega ) )$ restricted to the coordinates $\vartheta_i^\width$, $i \in ( \N \cap [1  , \width ] ) \cup \cu{ \fd_\width }$, has a negative eigenvalue. Combining this with \cref{cor:strict:saddle:conv:random} shows for all $\width \in \N $ that $ \P ( A_3^\width \cap A_4^\width \cap A_5^\width ) = 0$. \Nobs that \cref{theo:conv:bias:deterministic} proves that there exists $\scrC \in \R$ which satisfies for all $\width \in \N$, $\omega \in A_1^\width \cap A_2^\width \cap A_4^\width$ that
\begin{equation} \label{theo:main:conv:eq:limit:risk}
\textstyle \loss_\width ( \vartheta^\width ( \omega ) )
\le \scrC ( \width ^{ - \delta } + \width ^{ - 2 \delta } )
\le 2 \scrC \width ^{ - \delta }.
\end{equation}
Hence, we obtain for all $\width \in \N$ that
\begin{equation} \label{theo:conv:random:eq:risk:prob}
\begin{split}
\textstyle \P \rbr[\big]{ \forall \, t \in [0, \infty ) \colon \loss_\width ( \Theta_t^\width ) \le 2 \scrC \width ^{ - \delta } + \fC_\width t^{-1}} & \textstyle \ge \P \rbr[\big]{ \loss_\width ( \vartheta^\width ( \omega ) ) \le 2 \scrC \width ^{ - \delta } } \\
& \textstyle \ge  \P \rbr[\big]{ A_1^\width \cap A_2^\width \cap A_4^\width}.
\end{split}
\end{equation}
Furthermore, \cref{lem:random:init} proves that there exist $c_1, c_2 \in (0 , \infty)$ which satisfy for all $\width \in \N$ that $\P ( A_2^\width ) \ge 1 - c_1 \width ^{-1} \ge 1 - c_2 \width ^{ - \delta }$, $\P ( A_4^\width ) \ge 1 - c_2 \width ^{1 - \nicefrac{( \beta - \alpha )}{2} } \ge 1 - c_2 \width ^{ - \delta }$, and $\P ( A_5^\width ) \ge 1 - c_2 \width^{ 1 - \nicefrac{ \beta}{2} } \ge 1 - c_2 \width ^{ - \delta }$. Hence, we obtain for all $\width \in \N$ that
\begin{equation} \label{theo:conv:random:eq:prop:lower1}
\begin{split}
\textstyle \P \rbr[\big]{ A_1^\width \cap A_2^\width \cap A_4^\width } & \textstyle \ge 1 - \P ( \Omega \backslash A_1^\width ) - \P (\Omega \backslash A_2^\width) - \P ( \Omega \backslash A_4^\width ) \\
& \textstyle \ge 1 - \P \rbr[\big]{ ( \Omega \backslash A_1^\width ) \cap A_3^\width } - \P \rbr[\big]{ ( \Omega \backslash A_1^\width ) \cap (\Omega \backslash A_3^\width) } - (c_1 + c_2) \width^{- \delta}.
\end{split}
\end{equation}
Moreover, the fact that $\forall \, \width \in \N$, $i \in \N \cap ( \width , 3 \width ] \colon \vartheta^{ \width , i } = \Theta_0^{ \width , i }$, \cref{theo:main:conv:eq:gamma}, and \cref{lem:random:init:itemnew} in \cref{lem:random:init} (applied with $\eta \with \delta$ in the notation of \cref{lem:random:init:itemnew} in \cref{lem:random:init}) establish that there exists $c_3 \in (0 , \infty )$ which satisfies for all $\width \in \N \cap ( \bfm^{ - \nicefrac{1}{\delta } } , \infty )$ that 
\begin{equation} 
\begin{split}
& \textstyle \P \rbr[\big]{ ( \Omega \backslash A_1^\width ) \cap ( \Omega \backslash A_3^\width ) } \le \P \rbr[\big]{ \exists \,  J \subseteq \{1, 2, \ldots, \width\} \colon \\
& \textstyle \br[\big]{ \sum_{i \in J } \Theta_0^{ \width , \width + i } \Theta_0^{ \width , 2 \width + i } \le 4 \max\nolimits_{i \in \{1, 2, \ldots, \width\} } \Theta^{ \width , \width + i } \Theta_0^{ \width , 2 \width + i } \wedge  \ssum_{j \in J } \Theta_0^{ \width , 2 \width + j } \ge  \width ^{ - \delta } } } \le c_3 \width ^{ - \delta }.
\end{split}
\end{equation}
Hence, we obtain for all $\width \in \N$ that
\begin{equation} \label{theo:conv:random:eq:prop:lower2}
\textstyle \P \rbr[\big]{ ( \Omega \backslash A_1^\width ) \cap ( \Omega \backslash A_3^\width ) }
\le \max \cu{ c_3 , \bfm^{-1} } \width ^{ - \delta }.
\end{equation}
In addition, the fact that $ \P ( A_3^\width \cap A_4^\width \cap A_5^\width ) = 0$ proves for all $\width \in \N$ that $\P ( ( \Omega \backslash A_1^\width ) \cap A_3^\width ) \le \P ( A_3^\width ) \le \P ( \Omega \backslash A_4^\width ) + \P ( \Omega \backslash A_5^\width ) \le 2 c_2 \width ^{ - \delta }$. This, \cref{theo:conv:random:eq:prop:lower1}, and \cref{theo:conv:random:eq:prop:lower2} ensure for all $\width \in \N$ that
\begin{equation} \label{theo:conv:random:eq:prop:lower3}
\textstyle \P \rbr[\big]{ A_1^\width \cap A_2^\width \cap A_4^\width } \ge 1 - ( c_1 + 3c_2 + \max \cu{ c_3 , \bfm^{-1} } ) \width ^{ - \delta }.
\end{equation}
Combining this with \cref{theo:conv:random:eq:risk:prob} establishes \cref{theo:main:conv:item3}. Next \nobs that for all $\width \in \N$,
$\omega \in G$
it holds that $\limsup_{t \to \infty} \loss_\width ( \Theta_t^\width ( \omega ) ) = \loss_\width ( \vartheta^\width ( \omega ) )$. Combining this with \cref{theo:main:conv:eq:limit:risk},
the fact that $\P ( G ) = 1$, and Hölder's inequality establishes for all $\width \in \N$, $p , q\in (1 , \infty )$ with $\frac{1}{p} + \frac{1}{q} = 1$ that
\begin{equation} \label{theo:main:conv:eq:limsup:est}
\begin{split}
\textstyle \E \br[\big]{ \limsup\nolimits_{t \to \infty} \loss_\width ( \Theta_t^\width)} & \textstyle = \E \br[\big]{\loss_\width (\vartheta^\width)} = \E \br[\big]{\loss_\width ( \vartheta^\width ) \indicator{ A_1^\width \cap A_2^\width \cap A_4^\width } } + 
\E \br[\big]{\loss_\width ( \vartheta^\width ) \indicator{\Omega \backslash( A_1^\width \cap A_2^\width \cap A_4^\width ) } } \\
& \textstyle \le 2 \scrC \width^{- \delta} + \rbr*{ \E \br[\big]{\indicator{\Omega  \backslash( A_1^\width \cap A_2^\width \cap A_4^\width ) } } }^{ \nicefrac{1}{p}} \rbr[\big]{ \E \br[\big]{ ( \loss _ \width ( \vartheta ^\width ) )^{ q } } }^{ \nicefrac{1}{q} } \\
& \textstyle =  2 \scrC \width^{- \delta} + \rbr[\big]{1 - \P  \rbr[\big]{ A_1^\width \cap A_2^\width \cap A_4^\width}}^{\nicefrac{1}{p}} \rbr[\big]{ \E \br[\big]{ ( \loss _ \width ( \vartheta ^\width ) )^{q}}}^{\nicefrac{1}{q}}.
\end{split}
\end{equation}
In addition, the fact that $\forall \, \width \in \N , \omega \in G \colon \cG_\width ( \vartheta^\width ( \omega ) ) = 0 $ and \cref{lem:realization:pointwise:bound} assure for all $\width \in \N$,
$\omega \in G$ that
\begin{equation}
\begin{split}
\textstyle \loss_\width ( \vartheta^{ \width } ( \omega ))
& \textstyle \le \nu \rbr[\big]{ (f(\1 ) - f(\0) ) ^2 + \abs[\big]{ \sum_{i=1}^\width \vartheta^{ \width , 2 \width + i } ( \omega) }^2 } \\
& \textstyle =  \nu \rbr[\big]{ ( f ( \1 ) - f ( \0 ) )^2 + \abs[\big]{\sum_{i=1}^\width  \Theta_0^{ \width , 2 \width + i}(\omega)}^2}.
\end{split}
\end{equation}
This, the fact that $\P ( G ) = 1$, and the fact that $ \forall \, x,y \in [0, \infty ) , q \in [1 , \infty ) \colon   ( x + y )^q \le 2^{ q - 1 } ( x^q + y^q )$ demonstrate for all $\width \in \N$, $q \in [1 , \infty )$ that
\begin{equation}
\textstyle \E \br[\big]{ ( \loss _ \width ( \vartheta ^\width ) )^{ q } } \le 2^{q-1} \nu^q \rbr*{ (f(\1 ) - f(\0 ) ) ^{ 2 q } + \abs[\big]{ \sum_{i = 1}^\width \Theta_0^{\width, 2 \width + i}}^{2 q}}.
\end{equation}
Combining this with \cref{lem:random:init:item3} in \cref{lem:random:init} implies that for all $q \in [  1 , \infty )$
there exists $\fC \in \R$ such that for all $\width \in \N$ it holds
that
\begin{equation}
\textstyle \rbr[\big]{\E \br[\big]{ ( \loss _ \width ( \vartheta ^\width ))^{q}}}^{\nicefrac{1}{q}} \le \fC ( 1 + \width ^{ 2 - 2 \alpha } ) \le 2 \fC \width^{2 - 2 \alpha}.
\end{equation}
This, \cref{theo:conv:random:eq:prop:lower3}, and \cref{theo:main:conv:eq:limsup:est} show that for all $p \in (1, \infty )$ there exists $\fC \in \R$ which satisfies for all $\width \in \N$ that $ \E \br{ \limsup\nolimits_{t \to \infty} \loss_\width ( \Theta_t^\width ) } \le 2 \scrC \width ^{ - \delta } + \fC  \width ^{ 2 - 2 \alpha - \nicefrac{\delta }{p}} $. The reverse Fatou's lemma and the fact that $\forall \, p \in (1, \infty ), \varepsilon \in ( - \infty ,  \nicefrac{\delta }{p} - 2 ( 1 - \alpha ) ] , \width \in \N \colon  \width ^{ 2 - 2 \alpha - \nicefrac{\delta }{p}} \le \width^{ - \varepsilon}$ therefore establish \cref{theo:main:conv:item4}.
\end{cproof}

As a consequence, we obtain in \cref{cor:main:conv} a convergence rate of $\width ^{- \varepsilon}$ for $\varepsilon \in (0, \infty)$ arbitrarily close to $\frac{1}{9}$, if the scaling exponents $\alpha, \beta$ are chosen appropriately. Thus the rate is independent of any specific features of the target function $f$ or the density $\dens$. For the proof, we simply need to choose the scaling parameters $\alpha, \beta$ in \cref{theo:main:conv} in such a way that the conditions on $\varepsilon$ and $\delta$ are satisfied, which is possible iff $\varepsilon < \frac{1}{9} $.

\cfclear
\begin{cor} \label{cor:main:conv}
Let $\0 \in \R$, $\1 \in (\0 , \infty )$, let $f, \dens \in C ( \R , \R)$ be piecewise analytic\cfadd{def:piecewise_analytic}, assume for all $x \in \R$ that $\dens ( x ) \ge 0$ and $\dens^{- 1} (\R \backslash \{0\}) = (\0, \1)$, assume that $f$ is strictly increasing, let $(\Omega, \cF, \P)$ be a probability space, for every $\width \in \N$, $\theta = ( \theta_1, \ldots, \theta_{\fd_\width } ) \in \R^{ \fd_\width }$ let $\fd _\width = 3 \width + 1 $, let $\fN^{ \theta}_\width \colon \R \to \R$ satisfy for all $x \in \R$ that
\begin{equation}
\textstyle \fN^{ \theta }_\width ( x ) = \theta_{ \fd_\width } +
\sum_{ i = 1 }^{ \width } \theta_{2 \width + i } \max \cu{ \min \cu{ \theta_{\width + i } x + \theta_{ i } , 1 } , 0 } ,
\end{equation}
let $\loss_\width \colon \R^{ \fd_\width } \to \R$ satisfy
\begin{equation} 
\textstyle \loss_{ \width }( \theta ) 
= \int_{ \R }
\rbr[\big]{  \fN^{\theta}_\width ( x ) - f ( x ) }^2 \dens( x ) \, \d x,
\end{equation}
and let $\cG_\width = ( \cG_\width ^1 , \ldots, \cG_\width^{\fd_\width}) \colon \R^{\fd_\width} \to \R^{\fd_\width}$ satisfy for all $k \in \{1, 2, \ldots, \fd_\width\}$ that $\cG_\width^k (\theta) = (\frac{\partial^-}{\partial \theta_k} \loss_\width) (\theta) \allowbreak \indicator{\N \backslash ( \width , 3 \width ] } ( k )$, for every $\alpha, \beta \in \R$ let $\Theta^{\width, \alpha, \beta} = (\Theta^{\width, \alpha, \beta, 1}, \ldots, \Theta^{ \width , \alpha , \beta , \fd_\width } ) \colon [0, \infty ) \times \Omega \to \R^{ \fd_\width }$, $\width \in \N$, be stochastic processes with continuous sample paths, assume for all $\width \in \N$, $\alpha , \beta \in \R$, $t \in [0 , \infty ) $ that 
\begin{equation}
\P \bigl(\Theta_t^{ \width , \alpha , \beta }  
= \Theta_0^{ \width , \alpha , \beta }   - \tint_0^t \cG_\width ( \Theta_s ^{\width , \alpha , \beta }  ) \, \d s \bigr) = 1 ,
\end{equation}
assume for all $\width \in \N$, $\alpha , \beta \in \R$ that $\Theta_0^{ \width , \alpha , \beta , 1 } , \ldots, \Theta_0^{ \width , \alpha , \beta , \fd_\width }$ are independent, assume for all $\width \in \N$, $\alpha , \beta \in \R$, $k \in \{1, 2, \ldots, \width\}$ that $\width^{- \beta } \Theta_0^{ \width , \alpha , \beta , k }$ is standard normal, assume for all $\width \in \N$, $\alpha, \beta \in \R$, $ k \in \{1, 2, \ldots, \width\}$, $x \in (0, \infty)$ that
\begin{equation}
\textstyle \P(\width^{- \beta} \Theta_0^{\width, \alpha, \beta, \width + k} \leq x) = \P (\width^{\alpha} \Theta_0^{\width, \alpha, \beta, 2 \width + k} \le x) = [\frac{2}{\pi}]^{1/2} \int_0^x \exp(-  \frac{y^2}{2} ) \, \d y,
\end{equation}
assume for all $\width \in \N$, $\alpha , \beta \in \R$ that $\E [ \abs{ \Theta_0^{ \width , \alpha , \beta , \fd_\width } } ^2 ] < \infty $, and let $\varepsilon \in (0, \nicefrac{1}{9})$ \cfload. Then there exist $\alpha \in (\nicefrac{3}{4}, 1)$, $\beta \in (\alpha + 2, \infty)$ such that
\begin{enumerate} [label = (\roman*)]
\item
\label{cor:main:conv:item1} it holds for all $\width \in \N$ that $ \sup_{t \in [0, \infty ) } \E [ \loss_\width (\Theta_t^{\width, \alpha, \beta})] = \E [\sup_{t \in [0, \infty ) } \loss_\width ( \Theta_t^{ \width , \alpha , \beta }  ) ] < \infty $,

\item
\label{cor:main:conv:item2} there exist $\const \in (0, \infty )$
and random variables $\fC_\width \colon \Omega \to \R$, $\width \in \N$, such that for all $\width \in \N$ it holds that
\begin{equation}
\textstyle \P \rbr[\big]{ \forall \, t \in (0, \infty ) \colon \loss_\width ( \Theta_t ^{ \width , \alpha , \beta }  ) \le \const \width ^{ - \varepsilon } + \fC_\width t^{-1} } \ge 1 - \const \width ^{ - \varepsilon },
\end{equation}
and

\item
\label{cor:main:conv:item3} there exists $\scrC \in (0 , \infty )$ such that for all $\width \in \N$ it holds that
\begin{equation}
\textstyle \limsup\nolimits_{ t \to \infty } \rbr[\big]{\E[\loss_{\width}(\Theta^{\width, \alpha, \beta}_t)]} \le \E \br[\big]{ \limsup\nolimits_{t \to \infty} \loss_{ \width }(\Theta^{\width , \alpha , \beta } _t) } \le \scrC \width ^{ - \varepsilon }.
\end{equation}
\end{enumerate}
\end{cor}
\begin{cproof}{cor:main:conv}
Throughout this proof let $\delta \in (\varepsilon, \nicefrac{1}{9})$. \Nobs that $\frac{9 \delta + 3 }{4} < 1$. This ensures that there exist $\alpha \in (\nicefrac{3}{4}, 1)$, $\beta \in (\alpha + 2, \infty )$ which satisfy
\begin{equation}
\delta < \min \cu*{ \tfrac{ 4 \alpha - 3 }{9} , \tfrac{ \beta - \alpha - 2 }{2} } \qqandqq \varepsilon <  \delta - 2 ( 1 - \alpha ).
\end{equation}
\cref{theo:main:conv} and the fact that for all $\width \in \N$ it holds that $\width^{- \delta} \le \width^{- \varepsilon}$ hence establish \cref{cor:main:conv:item1,cor:main:conv:item2,cor:main:conv:item3}.
\end{cproof}

%==============================================================================%
%------------------------------------------------------------------------------%
%=================================-----Section-----============================%
%------------------------------------------------------------------------------%
%==============================================================================%
\section{ANNs with rectified linear unit (ReLU) activation}
\label{sec:anns_with_relu}

The main result of this section is \cref{cor:main} below, which in particular implies \cref{theorem_intro} from the introduction.

%------------------------------------------------------------------------------%
%----------------------------------Subsection----------------------------------%
%------------------------------------------------------------------------------%
\subsection{Differentiability of the risk function}

\subsubsection{ANNs with density}

We first introduce in \cref{setting1} the framework for shallow \ANNs\ with \ReLU\ activation which will be employed throughout this section.

\begin{setting}\label{setting1}
Let $\0, c \in \R$, $\1 \in (\max\{\0, 0\}, \infty)$, $\fh \in \N$, $v = (v_1, \ldots, v_{\fh}) \in \R^{\fh}$, $f \in C^1(\R, \allowbreak \R)$, $\fp, \fc \in C(\R, \R)$ satisfy for all $x \in \R$ that 
\begin{equation}
\textstyle f'(x) > 0 = f(\0), \qquad \fp^{-1}(\R \backslash \{0\}) = (\0, \1), \qqandqq \fc(x) = \min\{\max\{x, \0\}, \1\},
\end{equation}
for every $\theta = (\theta_1, \ldots, \theta_{\fh}) \in \R^{\fh}$ let $\cN^{\theta} \in C(\R, \R)$ satisfy for all $x \in \R$ that $\cN^{\theta}(x) = c + \sum_{i = 1}^{\fh} v_j \max\{x - \theta_{j}, 0\}$, let $\loss \colon \R^{\fh} \to \R$ satisfy for all $\theta \in \R^{\fh}$ that
\begin{equation}
\textstyle \loss(\theta) = \int_{\R} (\cN^{\theta}(x) - f(x))^2 \fp(x) \, \d x,
\end{equation}
and let $\cG = (\cG_1, \ldots, \cG_{\fh}) \colon \R^{\fh} \to \R^{\fh}$ satisfy for all $\theta \in \{\vartheta \in \R^{\fh} \colon \loss \text{ is differentiable at } \vartheta\}$ that
\begin{equation}
\textstyle \cG(\theta) = - (\nabla \loss)(\theta). 
\end{equation}
\end{setting}

In \cref{setting1} the outer weight parameters $v_i \in \R$ and the outer bias $c \in \R$ are fixed real numbers, and the inner biases $\theta_i \in \R$ are the variable parameters.
The inner weights are for the moment assumed to be $1$; however, by homogeneity of the \ReLU\ activation function $\R \ni x \mapsto \max \cu{x , 0 } \in \R$ this is not a restrictive assumption.
For every parameter vector $\theta \in \R^\width$ we again denote by $\cN ^\theta \in C ( \R , \R )$ the realization function of the considered \ANN\ with \ReLU\ activation and by $\loss ( \theta ) \in [0 , \infty )$ the corresponding $L^2$-risk.

%------------------------------------------------------------------------------%
%----------------------------------Subsection----------------------------------%
%------------------------------------------------------------------------------%
\subsubsection{Derivatives of the risk function}

As before, the assumption that $\dens^{-1} ( \R \backslash \cu{0}) = (\0, \1)$ can be used to show in \cref{lemma:derivatives_of_the_risk} that the risk function $\loss$ is twice continuously differentiable.

\cfclear
\begin{lemma}\label{lemma:derivatives_of_the_risk}
Assume \cref{setting1}. Then
\begin{enumerate}[label=(\roman*)]
\item
\label{item1:lemma:derivatives_of_the_risk} it holds that $\loss \in C^2 (\R^{\fh}, \R)$,

\item
\label{item2:lemma:derivatives_of_the_risk} it holds for all $\theta = (\theta_1, \ldots, \theta_{\fh}) \in \R^{\fh}$, $j \in \{1, 2, \ldots, \fh\}$ that
\begin{equation}\label{eqn:1st_derivative}
\textstyle \bigl(\frac{\partial}{\partial \theta_{j}} \loss\bigr)(\theta) = - 2 v_j \int_{\fc(\theta_j)}^{\1} (\cN^{\theta}(x) - f(x)) \fp(x) \d x,
\end{equation}
and
\item
\label{item3:lemma:derivatives_of_the_risk} it holds for all $\theta = (\theta_1, \ldots, \theta_{\fh}) \in \R^{\fh}$, $i, j \in \{1, 2, \ldots, \fh\}$ that
\begin{equation}\label{eqn:2nd_derivative}
\begin{split}
\textstyle \bigl(\frac{\partial^2}{\partial \theta_i \partial \theta_j} \loss\bigr)(\theta) & \textstyle = 2 v_i v_j \int_{\fc(\max\{\theta_i, \theta_j\})}^{\1} \fp(x) \, \d x + 2 v_j (\cN^{\theta}(\theta_j) - f(\theta_j)) \fp(\theta_j) \mathbbm{1}_{\{j\}}(i) \mathbbm{1}_{[\0, \1]}(\theta_j).
\end{split}
\end{equation}
\end{enumerate}
\end{lemma}
\begin{cproof}{lemma:derivatives_of_the_risk}
\Nobs that \cite[Lemma~14]{JentzenRiekert2022Piecewise}, the assumption that $\fp^{-1}(\R \backslash \{0\}) \allowbreak = (\0, \allowbreak \1)$, and the Leibniz integral rule establish \cref{item1:lemma:derivatives_of_the_risk,item2:lemma:derivatives_of_the_risk,item3:lemma:derivatives_of_the_risk}.
\end{cproof}

\cfclear
\begin{lemma}[Critical points]\label{lemma:critical_points}
Assume \cref{setting1}, let $\theta = (\theta_1, \ldots, \theta_{\fh}) \in \cG^{-1} (\{0\})$, $k \in \{1, \allowbreak 2, \allowbreak \ldots, \allowbreak \fh\}$, and let $s \colon \{1, 2, \ldots, k\} \to \{1, 2, \ldots, k\}$ satisfy for all $j \in \{1, 2, \ldots, k\}$ that $\abs{v_{s(j)}} > 0$ and $\0 \allowbreak \le \allowbreak \theta_{s(1)} \allowbreak < \theta_{s(2)} < \ldots < \theta_{s(k)} \le \1$. Then
\begin{equation}\label{eqn:lemma:critical_points}
\begin{split}
0 & \textstyle = \bigl[\int_{\0}^{\theta_{s(1)}} (\cN^{\theta}(x) - f(x))\fp(x) \, \d x \bigr] \mathbbm{1}_{(-\infty, \0]}(\min_{j \in \{1, 2, \ldots, \fh\}}(\theta_j)) \\
& \textstyle = \int_{\theta_{s(k)}}^{\1} (\cN^{\theta}(x) - f(x)) \fp(x) \, \d x = \sum_{j = 1}^{k - 1} \Abs{\int_{\theta_{s(j)}}^{\theta_{s(j + 1)}} (\cN^{\theta}(x) - f(x)) \fp(x) \, \d x}.
\end{split}
\end{equation}
\end{lemma}
\begin{cproof}{lemma:critical_points}
\Nobs that \cref{item2:lemma:derivatives_of_the_risk} in \cref{lemma:derivatives_of_the_risk} and the assumption that $\theta = (\theta_1, \allowbreak \ldots, \allowbreak \theta_{\fh}) \allowbreak \in \allowbreak \cG^{-1}(\{0\})$ establish \cref{eqn:lemma:critical_points}.
\end{cproof}

%------------------------------------------------------------------------------%
%----------------------------------Subsection----------------------------------%
%------------------------------------------------------------------------------%
\subsection{Convergence of GF to suitable critical points}

In this section we establish convergence of the considered \GF\ process to a critical point, and we also analyze the critical points of the risk. The main upper bound for the risk of critical points is given in \cref{cor:suitable_critical_upper_bound} below.

\subsubsection{Distinctness property of the kinks of active neurons of critical points}

\cfclear
\begin{lemma}[Distinctness of the active kinks]\label{lemma:distinctness_active_kinks}
Assume \cref{setting1}, assume for all $x \in \R$ that $\fp(x) \ge 0$, let $\fH \subseteq \N$ satisfy $\fH = \{1, 2, \ldots, \fh\}$, assume $\min_{j \in \fH} v_j > 0$, let $\theta = (\theta_1, \ldots, \theta_{\fh}) \in \cG^{-1}(\{0\})$, and assume that $\theta$ is not a \descritic critical point of $\loss$ \cfload. Then it holds for all $j \in \fH$ with $\0 < \theta_j < \1$ and $\cN^{\theta}(\theta_j) < f(\theta_j)$ that
\begin{equation}\label{eqn:lemma:distinctness_active_kinks}
\textstyle \{i \in \fH \backslash \{j\} \colon \theta_i = \theta_j\} = \emptyset.
\end{equation}
\end{lemma}
\begin{cproof}{lemma:distinctness_active_kinks}
\Nobs that the assumption that $\theta \in \cG^{-1}(\{0\})$ and the assumption that $\theta$ is not a \descritic critical point of $\loss$ ensure that for all $v \in \R^{\fh}$ it holds that
\begin{equation}
\textstyle \scalar{v, ((\Hs \loss)(\theta))v} = v ((\Hs \loss)(\theta)) v^T \ge 0\ifnocf.
\end{equation}
\cfload[.]This implies that $(\Hs \loss)(\theta)$ is positive-semidefinite. The Sylvester's criterion therefore shows that\footnote{\Nobs that the Sylvester's criterion ensures that a Hermitian matrix is positive-semidefinite if and only if all of its principal minors are non-negative.} for all $i, j \in \fH$ with $i \neq j$ it holds that $\frac{\partial^2}{\partial \theta_j^2} \loss(\theta) \ge 0 $ and
\begin{equation}\label{eqn:lemma:distinctness_active_kinks:minors}
\textstyle \det\!
\begin{pmatrix}
\bigl(\frac{\partial^2}{\partial \theta_i^2} \loss\bigr)(\theta) & \bigl(\frac{\partial^2}{\partial \theta_i \partial \theta_j} \loss\bigr)(\theta) \\[1ex]
\bigl(\frac{\partial^2}{\partial \theta_j \partial \theta_i} \loss\bigr)(\theta) & \bigl(\frac{\partial^2}{\partial \theta_j^2} \loss\bigr)(\theta)
\end{pmatrix}\! = \bigl[\bigl(\frac{\partial^2}{\partial \theta_i^2} \loss\bigr)(\theta)\bigr] \bigl[\bigl(\frac{\partial^2}{\partial \theta_j^2} \loss\bigr)(\theta)\bigr] - \bigl[\bigl(\frac{\partial^2}{\partial \theta_i \partial \theta_j} \loss\bigr)(\theta)\bigr]^2 \ge 0.
\end{equation}
Assume that there exist $i ,j \in \fH$ which satisfy
\begin{equation}
\textstyle i \neq j, \qquad \0 < \theta_i = \theta_j < \1, \qqandqq \cN^{\theta}(\theta_j) < f(\theta_j).
\end{equation}
\Nobs that the assumption that $\min_{\fj \in \fH} v_{\fj} > 0$, \cref{eqn:lemma:distinctness_active_kinks:minors}, and \cref{lemma:derivatives_of_the_risk} demonstrate that
\begin{equation}
\begin{split}
\textstyle \frac{\partial^2}{\partial \theta_j^2} \loss(\theta) & \textstyle = 2 [v_j]^2 \int_{\theta_j}^{\1} \fp(x) \, \d x + 2 v_j (\cN^{\theta}(\theta_j) - f(\theta_j)) \fp(\theta_j) \\
& \textstyle = 2 v_j \bigl[v_j \int_{\theta_j}^{\1} \fp(x) \, \d x + (\cN^{\theta}(\theta_j) - f(\theta_j)) \fp(\theta_j)\bigr] \ge 0
\end{split}
\end{equation}
and
\begin{equation}
\begin{split}
& \textstyle \bigl[\frac{\partial^2}{\partial \theta_i^2} \loss(\theta)\bigr] \bigl[\frac{\partial^2}{\partial \theta_j^2} \loss(\theta)\bigr] - \bigl[\frac{\partial^2}{\partial \theta_i \partial \theta_j} \loss(\theta)\bigr] \bigl[\frac{\partial^2}{\partial \theta_j \partial \theta_i} \loss(\theta)\bigr] \\
& \textstyle = 4 v_i v_j \bigl[v_j \int_{\theta_j}^{\1} \fp(x) \, \d x + (\cN^{\theta}(\theta_j) - f(\theta_j)) \fp(\theta_j)\bigr] \bigl[v_i \int_{\theta_j}^{\1} \fp(x) \, \d x + (\cN^{\theta}(\theta_j) - f(\theta_j)) \fp(\theta_j)\bigr] \\
& \textstyle \quad - 4 [v_i v_j]^2 \bigl[\int_{\theta_j}^{\1} \fp(x) \, \d x\bigr]^2 \\
& \textstyle = 4 v_i v_j (\cN^{\theta}(\theta_j) - f(\theta_j)) \fp(\theta_j) \bigl[(v_i + v_j) \int_{\theta_j}^{\1} \fp(x) \, \d x + (\cN^{\theta}(\theta_j) - f(\theta_j)) \fp(\theta_j)\bigr] \ge 0.
\end{split}
\end{equation}
This contradicts to the assumption that $\cN^{\theta}(\theta_j) < f(\theta_j)$. Hence, we obtain that for all $l \in \fH$ with $\0 < \theta_l < \1$ and $\cN^{\theta}(\theta_l) < f(\theta_l)$ it holds that
\begin{equation}
\textstyle \{k \in \fH \backslash \{l\} \colon \theta_k = \theta_l\} = \emptyset.
\end{equation}
\end{cproof}

%------------------------------------------------------------------------------%
%----------------------------------Subsection----------------------------------%
%------------------------------------------------------------------------------%
\subsubsection{Upper bound for the risk of suitable critical points}

Since the weight parameters $v_i$ are assumed to be positive, the realization function $\cN^\theta$ is always convex. Thus we need the target function $f$ to be convex as well in \cref{lemma:general_L2_risk_bound} and subsequent results in order for the approximation error to converge to zero.

\cfclear
\begin{lemma}\label{lemma:nonemptyset_endpoints}
Let $\scra \in \R$, $\scrb \in (\scra, \infty)$, $f \in C^1([\scra, \scrb], \R)$, $\fp \in C([\scra, \scrb], \R)$ satisfy for all $x \in [\scra, \scrb]$ that $f'(x) > 0 \le \fp(x)$ and $\int_{\0}^{\1} \abs{\fp(y)} \, \d y > 0$, assume that $f$ is strictly convex\footnote{\Nobs that for all $\0 \in \R$, $\1 \in (\0, \infty)$ and all functions $f \colon [\0, \1] \to \R$ it holds that $f$ is strictly convex if and only if it holds for all $x, y \in [\0, \1]$, $\lambda \in (0, 1)$ with $x \neq y$ that $f(\lambda x + (1 - \lambda) y) < \lambda f(x) + (1 - \lambda) f(y)$.}, let $g \colon [\scra, \scrb] \to \R$ be affine linear, and assume $\int_{\scra}^{\scrb} (f(x) - g(x)) \fp(x) \, \d x = 0$. Then
\begin{equation}
\textstyle \{x \in \{\scra, \scrb\} \colon g(x) < f(x)\} \neq \emptyset.
\end{equation}
\end{lemma}
\begin{cproof}{lemma:nonemptyset_endpoints}
For the sake of contradiction we assume that 
\begin{equation}\label{eqn:lemma:nonemptyset_endpoints:contradict}
\textstyle g(\0) \ge f(\0) \qqandqq g(\1) \ge f(\1).
\end{equation}
\Nobs that \cref{eqn:lemma:nonemptyset_endpoints:contradict}, the assumption that $f$ is strictly convex, and the assumption that $g$ is affine linear ensure that for all $x \in (\0, \1)$ it holds that $g(x) > f(x)$. Combining this with the assumption that $\Forall x \in [\0, \1] \colon \fp(x) \ge 0$ and the assumption that $\int_{\0}^{\1} \abs{\fp(y)} \, \d y > 0$ shows that $\int_{\scra}^{\scrb} (f(x) - g(x)) \fp(x) \, \d x < 0$. This contradicts to the assumption that $\int_{\scra}^{\scrb} (f(x) - g(x)) \fp(x) \, \d x = 0$. Therefore, we obtain that $\{x \in \{\scra, \scrb\} \colon g(x) < f(x)\} \neq \emptyset$.
\end{cproof}

\cfclear
\begin{lemma}\label{lemma:general_L2_risk_bound}
Let $\scra \in \R$, $\scrb \in (\scra, \infty)$, $f \in C^1([\scra, \scrb], \R)$, $\fp \in C([\scra, \scrb], \R)$ satisfy for all $x \in [\scra, \scrb]$ that $f'(x) > 0 \le \fp(x)$ and $\int_{\0}^{\1} \abs{\fp(y)} \, \d y > 0$, assume that $f$ is strictly convex, let $g \colon [\scra, \scrb] \to \R$ be affine linear, assume that $\int_{\scra}^{\scrb} (f(x) - g(x)) \fp(x) \, \d x = 0$, and let $\lambda \in \R$ satisfy $\lambda = \max_{x \in \{y \in \{\scra, \scrb\} \colon g(y) < f(y)\}} \allowbreak \abs{f(x) - g(x)}$ (cf.\ \cref{lemma:nonemptyset_endpoints}). Then
\begin{equation}\label{eqn:lemma:general_L2_risk_bound}
\textstyle \int_{\scra}^{\scrb} (f(x) - g(x))^2 \fp(x) \, \d x \le \lambda [\lambda + \sup_{x \in [\scra, \scrb]} f'(x)][\sup_{x \in [\scra, \scrb]} \fp(x)] (\scrb - \scra).
\end{equation}
\end{lemma}
\begin{cproof}{lemma:general_L2_risk_bound}
Throughout this proof let $L \in (0, \infty)$ satisfy $L = \sup_{x \in [\scra, \scrb]} f'(x)$, let $\sigma \in \R$ satisfy $\sigma = \sup_{x \in [\scra, \scrb]} \fp(x)$, and let $A, B \in \R$ satisfy for all $x \in [\scra, \scrb]$ that $g(x) = A x + B$. \Nobs that the assumption that $f$ is strictly convex ensures that for all $x, y \in [\0, \1]$ it holds that
\begin{equation}\label{eqn:convex_f-n}
\textstyle f(y) \ge f'(x) (y - x) + f(x).
\end{equation}
In the following we distinguish between the case $(g(\scra) < f(\scra)) \wedge (g(\scrb) \ge f(\scrb))$, the case $(g(\scra) \ge f(\scra)) \wedge (g(\scrb) < f(\scrb))$, and the case $(g(\scra) < f(\scra)) \wedge (g(\scrb) < f(\scrb))$. First we prove \cref{eqn:lemma:general_L2_risk_bound} in the case 
\begin{equation}\label{eqn:lemma:general_L2_risk_bound:case_1}
\textstyle (g(\scra) < f(\scra)) \wedge (g(\scrb) \ge f(\scrb)).
\end{equation}
Let $q \in [\scra, \scrb)$ satisfy $g(q) = f(q)$. \Nobs that \cref{eqn:lemma:general_L2_risk_bound:case_1} and the assumption that $f$ is strictly convex prove that
\begin{equation}
\textstyle \Forall x \in [\scra, q) \colon g(x) < f(x) \qqandqq \Forall x \in (q, \scrb) \colon g(x) > f(x).
\end{equation}
Combining this with \cref{eqn:convex_f-n} shows that
\begin{equation}
\begin{split}
& \textstyle \sup_{x \in [\scra, q]} \abs{g(x) - f(x)} \\
& \textstyle = f(\scra) - g(\scra) \ge  f'(q)(\scra - q) + f(q) - g(\scra) \\
& \textstyle = f'(q)(\scra - q) + g(q) - g(\scra) = (A - f'(q)) (q - \scra) \ge 0
\end{split}
\end{equation}
and
\begin{equation}
\begin{split}
& \textstyle \sup_{x \in [q, \scrb]} \abs{g(x) - f(x)} \\
& \textstyle \le g(\scrb) - [f'(q)(\scrb - q) + f(q)] = g(\scrb) - g(q) - f'(q)(\scrb - q) \\
& \textstyle = (A - f'(q)) (\scrb - q).
\end{split}
\end{equation}
The assumption that $\int_{\scra}^{\scrb} (f(x) - g(x)) \fp(x) \, \d x = 0$ therefore ensures that
\begin{equation}
\textstyle \int_{\scra}^{q} (g(x) - f(x))^2 \fp(x) \, \d x \le \abs{g(\scra) - f(\scra)}^2 (q - \scra) \sigma \le \lambda^2 \sigma (q - \scra)
\end{equation}
and
\begin{equation}
\begin{split}
& \textstyle \int_{q}^{\scrb} (g(x) - f(x))^2 \fp(x) \, \d x \\
& \textstyle \le (A - f'(q)) (\scrb - q) \int_{q}^{\scrb} (g(x) - f(x)) \fp(x) \, \d x \\
& \textstyle = (A - f'(q)) (\scrb - q) \int_{\scra}^{q} (f(x) - g(x)) \fp(x) \, \d x \\
& \textstyle \le (A - f'(q)) (\scrb - q) (f(\scra) - g(\scra)) (q - \scra) \sigma \\
& \textstyle \le \abs{f(\scra) - g(\scra)}^2 (\scrb - q) \sigma \le \lambda^2 \sigma (\scrb - q).
\end{split}
\end{equation}
Hence, we obtain that
\begin{equation}
\begin{split}
& \textstyle \int_{\scra}^{\scrb} (g(x) - f(x))^2 \fp(x) \, \d x \\
& \textstyle = \int_{\scra}^{q} (g(x) - f(x))^2 \fp(x) \, \d x + \int_{q}^{\scrb} (g(x) - f(x))^2 \fp(x) \, \d x \\
& \le \lambda^2 \sigma (q - \scra) + \lambda^2 \sigma (\scrb - q) = \lambda^2 \sigma (\scrb - \scra).
\end{split}
\end{equation}
This proves \cref{eqn:lemma:general_L2_risk_bound} in the case $(g(\scra) < f(\scra)) \wedge (g(\scrb) \ge f(\scrb))$. Next we prove \cref{eqn:lemma:general_L2_risk_bound} in the case 
\begin{equation}\label{eqn:lemma:general_L2_risk_bound:case_2}
\textstyle (g(\scra) \ge f(\scra)) \wedge (g(\scrb) < f(\scrb)).
\end{equation}
Let $q \in (\scra, \scrb]$ satisfy $g(q) = f(q)$. \Nobs that \cref{eqn:lemma:general_L2_risk_bound:case_2} and the assumption that $f$ is strictly convex assure that
\begin{equation}
\textstyle \Forall x \in (\scra, q) \colon g(x) > f(x) \qqandqq \Forall x \in (q, \scrb] \colon g(x) < f(x).
\end{equation}
Combining this with \cref{eqn:convex_f-n} proves that
\begin{equation}
\begin{split}
& \textstyle \sup_{x \in [\scra, q]} \abs{g(x) - f(x)} \\
& \textstyle \le g(\scra) - f'(q)(\scra - q) - f(q) = g(\scra) - g(q) - f'(q)(\scra - q) \\
& \textstyle = (f'(q) - A) (q - \scra)
\end{split}
\end{equation}
and
\begin{equation}
\begin{split}
& \textstyle \sup_{x \in [q, \scrb]} \abs{g(x) - f(x)} \\
& \textstyle = f(\scrb) - g(\scrb) \ge f'(q)(\scrb - q) + f(q) - g(\scrb) \\
& = f'(q)(\scrb - q) + g(q) - g(\scrb) \\
& \textstyle = (f'(q) - A) (\scrb - q) \ge 0.
\end{split}
\end{equation}
The assumption that $\int_{\scra}^{\scrb} (f(x) - g(x)) \fp(x) \, \d x = 0$ hence demonstrates that
\begin{equation}
\begin{split}
& \textstyle \int_{\scra}^{q} (g(x) - f(x))^2 \fp(x) \, \d x \\
& \textstyle \le (f'(q) - A) (q - \scra) \int_{\scra}^{q} (g(x) - f(x)) \fp(x) \, \d x \\
& \textstyle = (f'(q) - A) (q - \scra) \int_{q}^{\scrb} (f(x) - g(x)) \fp(x) \, \d x \\
& \textstyle \le (f'(q) - A) (q - \scra) (f(\scrb) - g(\scrb)) (\scrb - q) \sigma \\
& \textstyle \le \abs{f(\scrb) - g(\scrb)}^2 (q - \scra) \sigma \le \lambda^2 \sigma (q - \scra)
\end{split}
\end{equation}
and
\begin{equation}
\textstyle \int_{q}^{\scrb} (g(x) - f(x))^2 \, \d x \le \abs{g(\scrb) - f(\scrb)}^2 (\scrb - q) \sigma \le \lambda^2 \sigma (\scrb - q).
\end{equation}
Therefore, we obtain that
\begin{equation}
\begin{split}
& \textstyle \int_{\scra}^{\scrb} (g(x) - f(x))^2 \fp(x) \, \d x \\
& \textstyle = \int_{\scra}^{q} (g(x) - f(x))^2 \fp(x) \, \d x + \int_{q}^{\scrb} (g(x) - f(x))^2 \fp(x) \, \d x \\
& \le \lambda^2 \sigma (q - \scra) + \lambda^2 \sigma (\scrb - q) = \lambda^2 \sigma (\scrb - \scra).
\end{split}
\end{equation}
This proves \cref{eqn:lemma:general_L2_risk_bound} in the case $(g(\scra) \ge f(\scra)) \wedge (g(\scrb) < f(\scrb))$. In the next step we prove \cref{eqn:lemma:general_L2_risk_bound} in the case
\begin{equation}\label{eqn:lemma:general_L2_risk_bound:case_3}
\textstyle (g(\scra) < f(\scra)) \wedge (g(\scrb) < f(\scrb)).
\end{equation}

Let $q_1, q_2 \in (\scra, \scrb)$ satisfy $q_1 < q_2$ and $\abs{g(q_1) - f(q_1)} + \abs{g(q_2) - f(q_2)} = 0$. \Nobs that \cref{eqn:lemma:general_L2_risk_bound:case_3} and the assumption that $f$ is strictly convex show that
\begin{equation}
\textstyle \Forall x \in [\scra, q_1) \cup (q_2, \scrb] \colon g(x) < f(x) \qqandqq \Forall x \in (q_1, q_2) \colon g(x) > f(x).
\end{equation}
Combining this with \cref{eqn:convex_f-n} ensures that
\begin{equation}
\begin{split}
& \textstyle \sup_{x \in [\scra, q_1]} \abs{g(x) - f(x)} \\
& \textstyle = f(\scra) - g(\scra) \ge f'(q_1)(\scra - q_1) + f(q_1) - g(\scra) \\
& \textstyle = f'(q_1)(\scra - q_1) + g(q_1) - g(\scra) = (A - f'(q_1)) (q_1 - \scra) \ge 0,
\end{split}
\end{equation}
\begin{equation}
\begin{split}
& \textstyle \sup_{x \in [q_1, q_2]} \abs{g(x) - f(x)} \\
& \textstyle \le g(q_1) - f'(q_2) (q_1 - q_2) - f(q_2) = g(q_1) - g(q_2) - f'(q_2) (q_1 - q_2) \\
& \textstyle = (f'(q_2) - A)(q_2 - q_1),
\end{split}
\end{equation}
\begin{equation}
\begin{split}
& \textstyle \sup_{x \in [q_1, q_2]} \abs{g(x) - f(x)} \\
& \textstyle \le g(q_2) - f'(q_1) (q_2 - q_1) - f(q_1) = g(q_2) - g(q_1) - f'(q_1) (q_2 - q_1) \\
& \textstyle = (A - f'(q_1))(q_2 - q_1),
\end{split}
\end{equation}
and
\begin{equation}
\begin{split}
& \textstyle \sup_{x \in [q_2, \scrb]} \abs{g(x) - f(x)} \\
& \textstyle = f(\scrb) - g(\scrb) \ge f'(q_2)(\scrb - q_2) + f(q_2) - g(\scrb) \\
& \textstyle = f'(q_2)(\scrb - q_2) + g(q_2) - g(\scrb) = (f'(q_2) - A) (\scrb - q_2) \ge 0.
\end{split}
\end{equation}
The assumption that $\int_{\scra}^{\scrb} (f(x) - g(x)) \fp(x) \, \d x = 0$ therefore proves that
\begin{equation}
\textstyle \int_{\scra}^{q_1} (g(x) - f(x))^2 \fp(x) \, \d x \le (f(\scra) - g(\scra))^2 (q_1 - \scra) \sigma \le \lambda^2 \sigma (q_1 - \scra),
\end{equation}
\begin{equation}
\begin{split}
& \textstyle \int_{q_1}^{q_2} (g(x) - f(x))^2 \fp(x) \, \d x \le (f'(q_2) - A) (q_2 - q_1) \int_{q_1}^{q_2} (g(x) - f(x)) \fp(x) \, \d x \\
& \textstyle = (f'(q_2) - A) (q_2 - q_1) \bigl[\int_{\scra}^{q_1} (f(x) - g(x)) \fp(x) \, \d x + \int_{q_2}^{\scrb} (f(x) - g(x)) \fp(x) \, \d x\bigr] \\
& \textstyle \le (f'(q_2) - A) (q_2 - q_1) [(f(\scra) - g(\scra)) (q_1 - \scra) \sigma \\
& \textstyle \quad + (f(\scrb) - g(\scrb)) (\scrb - q_2) \sigma] \\
& \textstyle \le L (q_2 - q_1) (\lambda (q_1 - \scra) \sigma + \lambda (\scrb - q_2) \sigma) \le L \lambda \sigma (q_2 - q_1),
\end{split}
\end{equation}
and
\begin{equation}
\textstyle \int_{q_2}^{\scrb} (g(x) - f(x))^2 \fp(x) \, \d x \le (f(\scrb) - g(\scrb))^2 (\scrb - q_2) \sigma \le \lambda^2 \sigma (\scrb - q_2).
\end{equation}
Hence, we obtain that
\begin{equation}
\begin{split}
& \textstyle \int_{\scra}^{\scrb} (g(x) - f(x))^2 \fp(x) \, \d x \\
& \textstyle = \int_{\scra}^{q_1} (g(x) - f(x))^2 \fp(x) \, \d x  + \int_{q_1}^{q_2} (g(x) - f(x))^2 \fp(x) \, \d x + \int_{q_2}^{\scrb} (g(x) - f(x))^2 \fp(x) \, \d x \\
& \textstyle \le \lambda^2 \sigma (q_1 - \scra) + L \lambda \sigma (q_2 - q_1) + \lambda^2 \sigma (\scrb - q_2) \le \lambda \sigma (L + \lambda) (\scrb - \scra).
\end{split}
\end{equation}
This proves \cref{eqn:lemma:general_L2_risk_bound} in the case $(g(\scra) < f(\scra)) \wedge (g(\scrb) < f(\scrb))$.
\end{cproof}

In \cref{lemma:critical_upper_bound} we use the fact that $f$ is strictly convex and the realization is piecewise linear, and thus it is essentially sufficient to estimate the error at the breakpoints $\theta_i$.

\cfclear
\begin{lemma}[Upper bound for the risk of critical points]\label{lemma:critical_upper_bound}
Assume \cref{setting1}, assume for all $x \in \R$ that $\fp(x) \ge 0$, assume that $f$ is strictly convex, let $L \in (0, \infty)$ satisfy $L = \allowbreak \sup_{x \in [\0, \1]} f'(x)$, let $\fH \subseteq \N$, $\theta = (\theta_1, \ldots, \theta_{\fh}) \in \cG^{-1}(\{0\})$, $\lambda \in \R$ satisfy $\fH = \{1, 2, \ldots, \fh\}$ and 
\begin{equation}
\textstyle \lambda = \max_{j \in \{i \in \fH \colon [\0 < \theta_{i} < \1 \text{ and } \cN^{\theta}(\theta_i) < f(\theta_i)]\}} \abs{\cN^{\theta}(\theta_j) - f(\theta_j)},
\end{equation}
assume $\min_{j \in \fH} v_j > 0$, and assume $\cN^{\theta}(\1) \ge f(\1)$. Then it holds that
\begin{equation}\label{eqn:lemma:critical_upper_bound}
\textstyle \loss(\theta) \le \max\{\lambda, \abs{c}\} (L + \max\{\lambda, \abs{c}\}) [\sup_{x \in [\0, \1]} \fp(x)] (\1 - \0).
\end{equation}
\end{lemma}
\begin{cproof}{lemma:critical_upper_bound}
Throughout this proof let $k \in \N$, $s \colon \{1, 2, \ldots, k\} \to \{1, 2, \ldots, k\}$ satisfy
\begin{equation}
\textstyle k = \#((\cup_{j \in \fH}\{\theta_j\}) \cap (\0, \1)) \qqandqq \0 < \theta_{s(1)} < \theta_{s(2)} < \ldots < \theta_{s(k)} < \1
\end{equation}
and let $\sigma \in \R$ satisfy $\sigma = \sup_{x \in [\0, \1]} \fp(x)$. \Nobs that \cref{lemma:critical_points} and the assumption that $\cN^{\theta}(\1) \allowbreak \ge f(\1)$ ensure that $\cN^{\theta}(\theta_{s(k)}) \allowbreak < f(\theta_{s(k)})$ and
\begin{equation}
\textstyle \Abs{\int_{\theta_{s(k)}}^{\1} (\cN^{\theta}(x) - f(x)) \fp(x) \, \d x} + \sum_{j = 1}^{k - 1} \Abs{\int_{\theta_{s(j)}}^{\theta_{s(j + 1)}} (\cN^{\theta}(x) - f(x)) \fp(x) \, \d x} = 0.
\end{equation}
\cref{lemma:general_L2_risk_bound} and the fact that for all $j \in \{1, 2, \ldots, k - 1\}$ it holds that $[\theta_{s(j)}, \theta_{s(j + 1)}] \ni x \mapsto \cN^{\theta}(x) \allowbreak \in \allowbreak \R$ and $[\theta_{s(k)}, \1] \ni x \mapsto \cN^{\theta}(x) \in \R$ are affine linear therefore show that
\begin{equation}\label{eqn:lemma:critical_upper_bound:partial_bound_j}
\textstyle \Forall j \in \{1, 2, \ldots, k - 1\} \colon \int_{\theta_{s(j)}}^{\theta_{s(j + 1)}} (\cN^{\theta}(x) - f(x))^2 \fp(x) \, \d x \le \lambda \sigma (L + \lambda) (\theta_{s(j + 1)} - \theta_{s(j)}),
\end{equation}
and
\begin{equation}\label{eqn:lemma:critical_upper_bound:partial_bound_k}
\textstyle \int_{\theta_{s(k)}}^{\1} (\cN^{\theta}(x) - f(x))^2 \fp(x) \, \d x \le \lambda^2 \sigma (\1 - \theta_{s(k)}).
\end{equation}
In the next step we prove that
\begin{equation}\label{eqn:lemma:critical_upper_bound:partial_bound_1}
\textstyle \int_{\0}^{\theta_{s(1)}} (\cN^{\theta}(x) - f(x))^2 \fp(x) \, \d x \le \max\{\lambda, \abs{c}\} (L + \max\{\lambda, \abs{c}\}) (\theta_{s(1)} - \0) \sigma.
\end{equation}
In the following we distinguish between the case $\#\{j \in \fH \colon \theta_j \le \0\} > 0$ and the case $\#\{j \in \fH \colon \theta_j \le \0\} = 0$. First we prove \cref{eqn:lemma:critical_upper_bound:partial_bound_1} in the case $\#\{j \in \fH \colon \theta_j \le \0\} > 0$. \Nobs that the assumption that $\#\{j \in \fH \colon \theta_j \le \0\} > 0$ and \cref{lemma:critical_points} ensure that $\int_{\0}^{\theta_{s(1)}} (\cN^{\theta}(x) - f(x)) \fp(x) \, \d x = 0$. Combining this with \cref{lemma:general_L2_risk_bound} shows that
\begin{equation}
\textstyle \int_{\0}^{\theta_{s(1)}} (\cN^{\theta}(x) - f(x))^2 \fp(x) \, \d x \le
\begin{cases}
\lambda (L + \lambda) (\theta_{s(1)} - \0) \sigma & \colon \cN^{\theta}(\0) \ge 0 \\
\max\{\lambda, \abs{c}\} (L + \max\{\lambda, \abs{c}\}) (\theta_{s(1)} - \0) \sigma & \colon \cN^{\theta}(\0) < 0.
\end{cases}
\end{equation}
This proves \cref{eqn:lemma:critical_upper_bound:partial_bound_1} in the case $\#\{j \in \fH \colon \theta_j \le \0\} > 0$. Next we prove \cref{eqn:lemma:critical_upper_bound:partial_bound_1} in the case $\#\{j \in \fH \colon \theta_j \le \0\} = 0$. \Nobs that the assumption that $\#\{j \in \fH \colon \theta_j \le \0\} = 0$ shows that $\Forall x \in [\0, \theta_{s(1)}] \colon \cN^{\theta}(x) = c$. Combining this with the assumption that $f$ is strictly convex ensures that
\begin{equation}
\textstyle \sup_{x \in [\0, \theta_{s(1)}]} \abs{\cN^{\theta}(x) - f(x)} \le
\begin{cases}
\lambda & \colon [\cN^{\theta}(\theta_{s(1)}) < f(\theta_{s(1)})] \wedge [c < 0] \\
\max\{\lambda, \abs{c}\} & \colon [\cN^{\theta}(\theta_{s(1)}) < f(\theta_{s(1)})] \wedge [c \ge 0] \\
\abs{c} & \colon \cN^{\theta}(\theta_{s(1)}) \ge f(\theta_{s(1)}).
\end{cases}
\end{equation}
Hence, we obtain that
\begin{equation}
\textstyle \int_{\0}^{\theta_{s(1)}} (\cN^{\theta}(x) - f(x))^2 \fp(x) \, \d x \le [\max\{\lambda, \abs{c}\}]^2 (\theta_{s(1)} - \0) \sigma.
\end{equation}
This proves \cref{eqn:lemma:critical_upper_bound:partial_bound_1} in the case $\#\{j \in \fH \colon \theta_j \le \0\} = 0$. \Nobs that \cref{eqn:lemma:critical_upper_bound:partial_bound_j,eqn:lemma:critical_upper_bound:partial_bound_1,eqn:lemma:critical_upper_bound:partial_bound_k} demonstrate that
\begin{equation}
\begin{split}
& \textstyle \int_{\0}^{\1} (\cN^{\theta}(x) - f(x))^2 \fp(x) \, \d x \\
& \textstyle = \int_{\0}^{\theta_{s(1)}} (\cN^{\theta}(x) - f(x))^2 \fp(x) \, \d x + \bigl(\sum_{j = 1}^{k - 1} \bigl[\int_{\theta_{s(j)}}^{\theta_{s(j + 1)}} (\cN^{\theta}(x) - f(x))^2 \fp(x) \, \d x\bigr]\bigr) \\
& \textstyle \quad + \int_{\theta_{s(k)}}^{\1} (\cN^{\theta}(x) - f(x))^2 \fp(x) \, \d x \le \max\{\lambda, \abs{c}\} \sigma (L + \max\{\lambda, \abs{c}\}) (\theta_{s(1)} - \0) \\
& \textstyle \quad + \sum_{j = 1}^{k - 1} [\lambda \sigma (L + \lambda) (\theta_{s(j + 1)} - \theta_{s(j)})] + \lambda^2 \sigma (\1 - \theta_{s(k)}) \\
& \textstyle = \max\{\lambda, \abs{c}\} \sigma (L + \max\{\lambda, \abs{c}\}) (\theta_{s(1)} - \0) + \lambda \sigma (L + \lambda) (\theta_{s(k)} - \theta_{s(1)}) + \lambda^2 \sigma (\1 - \theta_{s(k)}) \\
& \textstyle \le \max\{\lambda, \abs{c}\} \sigma (L + \max\{\lambda, \abs{c}\}) (\1 - \0).
\end{split}
\end{equation}
This establishes \cref{eqn:lemma:critical_upper_bound}.
\end{cproof}

\cfclear
\begin{corollary}[Upper bound for the risk of suitable critical points]\label{cor:suitable_critical_upper_bound}
Assume \cref{setting1}, assume for all $x \in \R$ that $\fp(x) \ge 0$, assume that $f$ is strictly convex, let $L \in (0, \infty)$ satisfy $L = \allowbreak \sup_{x \in [\0, \1]} f'(x)$, let $\fH \subseteq \N$, $\theta = (\theta_1, \ldots, \theta_{\fh}) \in \cG^{-1}(\{0\})$ satisfy $\fH = \{1, 2, \ldots, \fh\}$, assume $\min_{j \in \fH} v_j > 0$, assume $\cN^{\theta}(\1) \ge f(\1)$, and assume that $\theta$ is not a \descritic critical point of $\loss$ \cfload. Then
\begin{enumerate}[label=(\roman*)]
\item
\label{item1:cor:suitable_critical_upper_bound} it holds for all $j \in \fH$ with $\0 < \theta_j < \1$ and $\cN^{\theta}(\theta_j) < f(\theta_j)$ that
\begin{equation}
\textstyle 0 < f(\theta_j) - \cN^{\theta}(\theta_j) \le v_j,
\end{equation}

\item
\label{item2:cor:suitable_critical_upper_bound} it holds that
\begin{equation}
\textstyle \max_{j \in \{i \in \fH \colon [\0 < \theta_i < \1 \text{ and } \cN^{\theta}(\theta_i) < f(\theta_i)]\}} \abs{\cN^{\theta}(\theta_j) - f(\theta_j)} \le \max_{j \in \fH} v_j,
\end{equation}

and
\item
\label{item3:cor:suitable_critical_upper_bound} it holds that
\begin{equation}\label{eqn:cor:suitable_critical_upper_bound}
\textstyle \loss(\theta) \le [\max\{\max_{j \in \fH} v_j, \abs{c}\}] (L + [\max\{\max_{j \in \fH} v_j, \abs{c}\}]) [\sup_{x \in [\0, \1]} \fp(x)] (\1 - \0).
\end{equation}
\end{enumerate}
\end{corollary}
\begin{cproof}{cor:suitable_critical_upper_bound}
Throughout this proof let $k \in \N$, $s \colon \{1, 2, \ldots, k\} \to \{1, 2, \ldots, k\}$, $\theta_{s(k + 1)} \allowbreak \in \allowbreak \R$ satisfy
\begin{equation}
\textstyle k = \# ((\cup_{j \in \fH} \{\theta_j\}) \cap (\0, \1)) \qqandqq \0 < \theta_{s(1)} < \theta_{s(2)} < \ldots < \theta_{s(k)} < \theta_{s(k + 1)} = \1
\end{equation}
and let $A_1, A_2, \ldots, A_{k + 1}, B_1, B_2, \ldots, B_{k + 1} \in \R$ satisfy
\begin{equation}
\begin{split}
& \textstyle (\Forall x \in [\0, \theta_{s(1)}] \colon \cN^{\theta}(x) = A_1 x + B_1), \quad (\Forall x \in [\theta_{s(k)}, \1] \colon \cN^{\theta}(x) = A_{k + 1} x + B_{k + 1}),  \\
& \textstyle \andq (\Forall j \in \{1, 2, \ldots, k - 1\}, x \in [\theta_{s(j)}, \theta_{s(j + 1)}] \colon \cN^{\theta}(x) = A_{j + 1} x + B_{j + 1}).
\end{split}
\end{equation}
\Nobs that for all $j \in \{1, 2, \ldots, k\}$ it holds that
\begin{equation}\label{eqn:cor:suitable_critical_upper_bound:slope_change}
\textstyle A_{j + 1} = A_j + \sum_{i \in \{l \in \fH \colon \theta_l = \theta_{s(j)}\}} v_i \qqandqq B_{j + 1} = B_j - \sum_{i \in \{l \in \fH \colon \theta_l = \theta_{s(j)}\}} v_i \theta_i. 
\end{equation}
Next \nobs that \cref{lemma:distinctness_active_kinks} and the assumption that $\theta$ is not a \descritic critical point of $\loss$ ensure that for all $j \in \{1, 2, \ldots, k\}$ with $\cN^{\theta}(\theta_{s(j)}) < f(\theta_{s(j)})$ it holds that
\begin{equation}
\textstyle \{i \in \fH \backslash \{j\} \colon \theta_i = \theta_{s(j)}\} = \emptyset.
\end{equation}
Combining this with \cref{eqn:cor:suitable_critical_upper_bound:slope_change} shows that for all $j \in \{1, 2, \ldots, k\}$ with $\cN^{\theta}(\theta_{s(j)}) < f(\theta_{s(j)})$ it holds that
\begin{equation}\label{eqn:cor:suitable_critical_upper_bound:slope_change_special}
\textstyle A_{j + 1} = A_j + v_{s(j)} \qqandqq B_{j + 1} = B_j - v_{s(j)} \theta_{s(j)}.
\end{equation}
Let $\fj \in \{1, 2, \ldots, k\}$, $q \in (\theta_{s(\fj)}, \theta_{s(\fj + 1)})$ satisfy $\cN^{\theta}(q) = f(q)$ and $\Forall x \in [\theta_{s(\fj)}, q) \colon \cN^{\theta}(x) < f(x)$. \Nobs that the assumption that $f$ is strictly convex, the Lagrange MVT, and \cref{eqn:cor:suitable_critical_upper_bound:slope_change_special} assure that there exists $c \in (\theta_{s(\fj)}, q)$ such that $[\Forall x \in (\theta_{s(\fj)}, q) \colon A_{\fj} < f'(x)]$ and
\begin{equation}
\begin{split}
0 & \textstyle < f(\theta_{s(\fj)}) - \cN^{\theta}(\theta_{s(\fj)}) = f(\theta_{s(\fj)}) - f(q) + \cN^{\theta}(q) - \cN^{\theta}(\theta_{s(\fj)}) \\
& \textstyle = f'(c) (\theta_{s(\fj)} - q) + A_{\fj + 1} (q - \theta_{s(\fj)}) = (A_{\fj + 1} - f'(c)) (q - \theta_{s(\fj)}) \\
& \textstyle = (A_{\fj} + v_{s(\fj)} - f'(c)) (q - \theta_{s(\fj)}) \le v_{s(\fj)} (q - \theta_{s(\fj)}) \le v_{s(\fj)}.
\end{split}
\end{equation}
Therefore, we obtain that for all $j \in \fH$ with $\0 < \theta_j < \1$ and $\cN^{\theta}(\theta_j) < f(\theta_j)$ it holds that
\begin{equation}
\textstyle 0 < f(\theta_j) - \cN^{\theta}(\theta_j) \le v_j.
\end{equation}
This establishes \cref{item1:cor:suitable_critical_upper_bound} and \cref{item2:cor:suitable_critical_upper_bound}. Furthermore, \nobs that \cref{item2:cor:suitable_critical_upper_bound} and \cref{lemma:critical_upper_bound} establish \cref{item3:cor:suitable_critical_upper_bound}.
\end{cproof}

%------------------------------------------------------------------------------%
%----------------------------------Subsection----------------------------------%
%------------------------------------------------------------------------------%
\subsubsection{Convergence of GF}

We next show convergence of the \GF\ trajectory with a suitable convergence rate. Similarly to the previous section, we use the KL inequality for the risk, after having established boundedness of the trajectory.

\cfclear
\begin{lemma}[Boundedness of \GF\ trajectory]\label{lemma:GF_boundedness}
Assume \cref{setting1}, assume $\min_{j \in \{1, 2, \ldots, \fh\}} v_j > 0$, assume for all $x \in \R$ that $\fp(x) \ge 0 < \int_{\R} \fp(y) \, \d y$, and let $\Theta \in C([0, \infty), \R^{\fh})$ satisfy for all $t \in [0, \allowbreak \infty)$ that
\begin{equation}\label{eqn:lemma:GF_boundedness:GF}
\textstyle \Theta_0 \in \R^{\fh} \qquad \text{and} \qquad \Theta_t = \Theta_0 + \int_0^t \cG(\Theta_s) \d s.
\end{equation}
Then
\begin{enumerate}[label=(\roman*)]
\item
\label{item1:lemma:GF_boundedness} it holds for all $j \in \{1, 2, \ldots, \fh\}$ that $(\Exists T \in [0, \infty) \colon \Theta_{j, T} > \1)$ implies that $(\Forall t \in [0, \allowbreak \infty) \allowbreak \colon \allowbreak \Theta_{j, t} \allowbreak = \Theta_{j, 0} > \1)$,

\item
\label{item2:lemma:GF_boundedness} it holds for all $j \in \{1, 2, \ldots, \fh\}$ that $(\Exists T \in [0, \infty) \colon \Theta_{j, T} \le \1)$ implies that for all $t \in [0, \infty)$ it holds that
\begin{equation}
\textstyle \1 \ge \Theta_{j, t} \ge \0 - [v_j]^{-1} \bigl([\loss(\Theta_0)]^{1/2} + \bigl[\int_{\R} (f(x) - c)^2 \fp(x) \, \d x\bigr]^{1/2}\bigr) \bigl[\int_{\R} \fp(x) \, \d x\bigr]^{-1/2},
\end{equation}

\item
\label{item3:lemma:GF_boundedness} it holds for all $j \in \{1, 2, \ldots, \fh\}$, $t \in [0, \infty)$ that
\begin{equation}
\textstyle \abs{\Theta_{j, t}} \le \max\Bigl\{\1, \abs{\Theta_{j, 0}}, \abs{\0} + \frac{[v_j]^{-1}}{[\int_{\R} \fp(x) \, \d x]^{1/2}} \bigl([\loss(\Theta_0)]^{1/2} + \bigl[\int_{\R} (f(x) - c)^2 \fp(x) \, \d x\bigr]^{1/2}\bigr)\Bigr\},
\end{equation}

and
\item
\label{item4:lemma:GF_boundedness} it holds for all $t \in [0, \infty)$ that
\begin{equation}
\textstyle \norm{\Theta_t} \le \fh (\1 + \abs{\0}) + \fh^{1/2} \norm{\Theta_0} + \frac{\sum_{j = 1}^{\fh} [v_{j}]^{-1}}{[\int_{\R} \fp(x) \, \d x]^{1/2}} \bigl([\loss(\Theta_0)]^{1/2} + \bigl[\int_{\R} (f(x) - c)^2 \fp(x) \, \d x\bigr]^{1/2}\bigr).
\end{equation}
\end{enumerate}
\end{lemma}
\begin{cproof}{lemma:GF_boundedness}
Throughout this proof let $j \in \{1, 2, \ldots, \fh\}$. \Nobs that \cref{eqn:lemma:GF_boundedness:GF} assures that for all $t, \ft \in [0, \infty)$ with $\Theta_{j, t} \ge \1$ it holds that $[0, \infty) \ni \tau \mapsto \Theta_{j, \tau} \in \R$ is differentiable at $\ft$ and
\begin{equation}\label{eqn:lemma:GF_boundedness:Theta_j>1}
\textstyle \frac{\partial}{\partial t} \Theta_{j, t} = \cG_j(\Theta_t) = - \frac{\partial}{\partial \Theta_{j, t}} \loss(\Theta_t) = 2 v_j \int_{\Theta_{j, t}}^{\1} (\cN^{\Theta_t}(x) - f(x)) \fp(x) \, \d x = 0.
\end{equation}
Assume that there exist $T \in [0, \infty)$, $\scrc \in \R$ which satisfy $\scrc = \Theta_{j, T} > \1$ and
\begin{multline}\label{eqn:lemma:GF_boundedness:assumption_inf_sup}
\textstyle \bigl[\inf\{t \in [0, T] \colon (\Forall x \in [t, T] \colon \Theta_{j, x} > \1)\} > 0\bigr] \\
\textstyle \vee \bigl[\sup\{t \in [T, \infty) \colon (\Forall x \in [T, t] \colon \Theta_{j, x} > \1)\} < \infty\bigr] = 1.
\end{multline}
Let $\bft, \bfT \in [0, \infty)$ satisfy $\bft = \inf\{t \in [0, T] \colon (\Forall x \in [t, T] \colon \theta_{j, x} > \1)\}$ and $\bfT = \sup\{t \in [T, \infty) \colon (\Forall x \in [T, t] \colon \Theta_{j, x} > \1)\}$. \Nobs that the fact that $[0, \infty) \ni t \mapsto \Theta_{j, t} \in \R$ is continuous shows that for all $t \in (\bft, \bfT)$ it holds that $0 < \bft < T < \bfT < \infty$ and $\Theta_{j, t} > \1$. Combining this with \cref{eqn:lemma:GF_boundedness:Theta_j>1} and Lagrange MVT demonstrates that for all $t \in (\bft, \bfT) \backslash \{T\}$ there exists $\scrt \in (\min\{t, T\}, \max\{t, T\})$ such that
\begin{equation}
\textstyle \Theta_{j, t} = \Theta_{j, t} - \Theta_{j, T} + \scrc = (t - T) \cG_j(\Theta_{\scrt}) + \scrc = \scrc.
\end{equation}
This implies that for all $t \in (\bft, \bfT)$ it holds that $\Theta_{j, t} = \scrc$. The fact that $[0, \infty) \ni t \mapsto \Theta_{j, t} \in \R$ is continuous hence proves that for all $t \in [\bft, \bfT]$ it holds that $\Theta_{j, t} = \scrc$. Combining this with the fact that $[0, \infty) \ni t \mapsto \Theta_{j, t} \in \R$ is continuous shows that there exists $\eps \in (0, \infty)$ which satisfies for all $t \in [\bft - \eps, \bfT + \eps]$ that $0 \le \bft - \eps$ and $\Theta_{j, t} > \1$. This, \cref{eqn:lemma:GF_boundedness:Theta_j>1}, and Lagrange MVT demonstrate that for all $t \in [\bft - \eps, \bfT + \eps] \backslash \{T\}$ there exists $\scrt \in (\min\{t, T\}, \max\{t, T\})$ such that
\begin{equation}
\textstyle \Theta_{j, t} = \Theta_{j, t} - \Theta_{j, T} + \scrc = (t - T) \cG_j(\Theta_{\scrt}) + \scrc = \scrc.
\end{equation}
Hence, we obtain that for all $t \in [\bft - \eps, \bfT + \eps]$ it holds that $\Theta_{j, t} = \Theta_{j, T} = \scrc$. This contradicts to the assumption in \cref{eqn:lemma:GF_boundedness:assumption_inf_sup}. Therefore, we obtain that for all $t \in (0, \infty)$ it holds that $\Theta_{j, t} = \Theta_{j, 0} = \scrc > \1$. Combining this and the fact that $[0, \infty) \ni t \mapsto \Theta_{j, t} \in \R$ is continuous ensures that for all $t \in [0, \infty)$ it holds that $\Theta_{j, t} = \Theta_{j, 0} = \scrc > \1$. This establishes \cref{item1:lemma:GF_boundedness}. Next \nobs that \cite[Lemma~3.1]{Riekert2021ConvergenceConstTF} and the fact that for all $\scrR \in [0, \infty)$ it holds that 
\begin{equation}
\textstyle \sup_{\theta \in \R^{\fh}, \norm{\theta} \le \scrR} \norm{(\nabla \loss)(\theta)} < \infty
\end{equation}
ensure that for all $t \in [0, \infty)$ it holds that
\begin{equation}\label{eqn:lemma:GF_boundedness:Risk_bound}
\textstyle 0 \le \loss(\Theta_t) = \loss(\Theta_0) - \int_0^t \norm{\cG(\Theta_s)}^2 \, \d s \le \loss(\Theta_0).
\end{equation}
H\"{o}lder inequality and the assumption that $\min_{i \in \{1, 2, \ldots, \fh\}} v_i > 0$ therefore demonstrate that for all $t \in [0, \infty)$ with $\Theta_{j, t} \le \min\{\0, 0\}$ it holds that
\begin{equation}\label{eqn:lemma:GF_boundedness:Theta_j<1}
\begin{split}
& \textstyle \bigl([\loss(\Theta_0)]^{1/2} + \bigl[\int_{\R} (f(x) - c)^2 \fp(x) \, \d x\bigr]^{1/2}\bigr)^2 \\
& \textstyle \ge \bigl([\loss(\Theta_t)]^{1/2} + \bigl[\int_{\R} (f(x) - c)^2 \fp(x) \, \d x\bigr]^{1/2}\bigr)^2 \ge \int_{\R} (\cN^{\Theta_t}(x) - c)^2 \fp(x) \, \d x \\
& \textstyle \ge \int_{\max\{\Theta_{j, t}, \0\}}^{\1} (v_j [x - \Theta_{j, t}])^2 \fp(x) \, \d x = \int_{\0}^{\1} (v_j [x - \Theta_{j, t}])^2 \fp(x) \, \d x \\
& \textstyle \ge (v_j [\0 - \Theta_{j, t}])^2 \int_{\0}^{\1} \fp(x) \, \d x.
\end{split}
\end{equation}
This shows that for all $t \in [0, \infty)$ with $\Theta_{j, t} \le \1$ it holds that
\begin{equation}
\textstyle \1 \ge \Theta_{j, t} \ge \0 - [v_j]^{-1} \bigl([\loss(\Theta_0)]^{1/2} + \bigl[\int_{\R} (f(x) - c)^2 \fp(x) \, \d x\bigr]^{1/2}\bigr) \bigl[\int_{\R} \fp(x) \, \d x\bigr]^{-1/2}.
\end{equation}
Combining this with \cref{item1:lemma:GF_boundedness} establishes \cref{item2:lemma:GF_boundedness}. Furthermore, \nobs that \cref{item1:lemma:GF_boundedness} and \cref{item2:lemma:GF_boundedness} assure that for all $t \in [0, \infty)$ it holds that
\begin{equation}
\textstyle \abs{\Theta_{j, t}} \le \max\bigl\{\1, \abs{\Theta_{j, 0}}, \abs{\0} + [v_j]^{-1} \bigl([\loss(\Theta_0)]^{1/2} + \bigl[\int_{\R} (f(x) - c)^2 \fp(x) \, \d x\bigr]^{1/2}\bigr) \bigl[\int_{\R} \fp(x) \, \d x\bigr]^{-1/2}\bigr\}.
\end{equation}
This establishes \cref{item3:lemma:GF_boundedness}. Next \nobs that \cref{item3:lemma:GF_boundedness} shows that for all $t \in [0, \infty)$ it holds that
\begin{equation}
\begin{split}
\textstyle \norm{\Theta_t} & \textstyle \le \fh \1 + \bigl[\sum_{\fj = 1}^{\fh} \abs{\Theta_{\fj, 0}}\bigr] + \fh \abs{\0} + \frac{\sum_{\fj = 1}^{\fh} [v_{\fj}]^{-1}}{[\int_{\R} \fp(x) \, \d x]^{1/2}} \bigl([\loss(\Theta_0)]^{1/2} + \bigl[\int_{\R} (f(x) - c)^2 \fp(x) \, \d x\bigr]^{1/2}\bigr) \\
& \textstyle \le \fh (\1 + \abs{\0}) + \fh^{1/2} \norm{\Theta_0} + \frac{\sum_{\fj = 1}^{\fh} [v_{\fj}]^{-1}}{[\int_{\R} \fp(x) \, \d x]^{1/2}} \bigl([\loss(\Theta_0)]^{1/2} + \bigl[\int_{\R} (f(x) - c)^2 \fp(x) \, \d x\bigr]^{1/2}\bigr).
\end{split}
\end{equation}
This establishes \cref{item4:lemma:GF_boundedness}.
\end{cproof}

\cfclear
\begin{lemma}[GF convergence]\label{lemma:GF_convergence}
Assume \cref{setting1}, assume that $f$ and $\fp$ are piecewise analytic, assume $\min_{j \in \{1, 2, \ldots, \fh\}} v_j > 0$, assume for all $x \in \R$ that $\fp(x) \ge 0$, and let $\Theta \in C([0, \infty), \allowbreak \R^{\fh})$ satisfy for all $t \in [0, \infty)$ that
\begin{equation}\label{eqn:lemma:GF_convergence:GF}
\textstyle \Theta_0 \in \R^{\fh} \qqandqq \Theta_t = \Theta_0 + \int_0^t \cG(\Theta_s) \, \d s
\end{equation}
(cf. \cref{def:piecewise_analytic}). Then there exist $\vartheta \in \R^{\fh}$, $\scrC, \alpha \in (0, \infty)$ such that for all $t \in [0, \infty)$ it holds that
\begin{equation}\label{eqn:lemma:GF_convergence:limits}
\textstyle \cG(\vartheta) = 0, \quad \norm{\Theta_t - \vartheta} \le \scrC (1 + t)^{- \alpha}, \qandq 0 \le \loss(\Theta_t) - \loss(\vartheta) \le \scrC (1 + t)^{-1}.
\end{equation}
\end{lemma}
\begin{cproof}{lemma:GF_convergence}
In the same way as in the proof of \cref{lem:risk_function_prop}, the assumption that $f$ and $\fp$ are piecewise analytic implies that 
\begin{equation}
	[\R^\width \times \R \ni ( \theta , x ) \mapsto (\realization{\theta} ( x ) - f ( x ) ) ^2 \fp ( x ) \in \R ] \in \ssC_ \fh . 
\end{equation}
Combining this with \cref{prop:integral_prop} proves that $\loss$ is a subanalytic function.
This, \cite[Corollary~7.11]{JentzenRiekert2021ExistenceGlobMin}, and \cite[Theorem~3.1 and $(4)$]{BolteDaniilidis2006} demonstrate that for all $\theta \in \R^{\fh}$ there exist $\eps, \fC \in (0, \infty)$, $\alpha \in (0, 1)$ such that for all $\vartheta \in \R^{\fh}$ with $\norm{\vartheta - \theta} < \eps$ it holds that $\abs{\loss(\vartheta) - \loss(\theta)}^{\alpha} \allowbreak \le \fC \norm{\cG(\vartheta)}$. Combining this and \cref{lemma:GF_boundedness} with \cite[Corollary~8.4]{JentzenRiekert2021ExistenceGlobMin} establishes \cref{eqn:lemma:GF_convergence:limits}.
\end{cproof}

%------------------------------------------------------------------------------%
%----------------------------------Subsection----------------------------------%
%------------------------------------------------------------------------------%
\subsubsection{Convergence of GF to suitable critical points}

Next we establish in \cref{prop:convergence_gf_suitable_cp} an upper bound for the risk of the limit point $\vartheta$ of the \GF\ trajectory, under the assumption that $\cN^{\vartheta}(\1) \ge f(\1)$.
Afterwards, we show in \cref{lemma:N>f_at_1} that this assumption is always satisfied under suitable assumptions on the initial value $\Theta_0$ of the \GF\ and the weight parameters $v_j$.
 
\cfclear
\begin{proposition}\label{prop:convergence_gf_suitable_cp}
Assume \cref{setting1}, assume for all $x \in \R$ that $\fp(x) \ge 0$, assume that $f$ is strictly convex, assume that $f$ and $\fp$ are piecewise analytic, let $\fH \subseteq \N$, $L \in \R$ satisfy $\fH = \{1, 2, \ldots, \fh\}$ and $L = \sup_{x \in [\0, \1]} f'(x)$, assume $\min_{j \in \fH} v_j > 0$, let $\Theta \in C([0, \infty), \R^{\fh})$ satisfy for all $t \in [0, \infty)$ that
\begin{equation}
\textstyle \Theta_0 \in \R^{\fh} \qandq \Theta_t = \Theta_0 + \int_0^t \cG(\Theta_s) \, \d s,
\end{equation}
let $\vartheta \in \R^{\fh}$ satisfy
\begin{equation}\label{eqn:prop:convergence_gf_suitable_cp:vartheta}
\textstyle \limsup_{t \to \infty} (\norm{\Theta_t - \vartheta} + \abs{\loss(\Theta_t) - \loss(\vartheta)} + \norm{\cG(\Theta_t)} + \norm{\cG(\vartheta)}) = 0
\end{equation}
(cf.\ \cref{lemma:GF_convergence}), assume $\cN^{\vartheta}(\1) \ge f(\1)$, and assume that $\vartheta$ is not a \descritic critical point of $\loss$ (cf. \cref{def:strict_saddle,def:piecewise_analytic}). Then
\begin{equation}\label{eqn:prop:convergence_gf_suitable_cp:risk_bound}
\textstyle \loss(\vartheta) \le [\max\{\max_{j \in \fH} v_j, \abs{c}\}] (L + [\max\{\max_{j \in \fH} v_j, \abs{c}\}]) [\sup_{x \in [\0, \1]} \fp(x)] (\1 - \0).
\end{equation}
\end{proposition}
\begin{cproof}{prop:convergence_gf_suitable_cp}
\Nobs that \cref{eqn:prop:convergence_gf_suitable_cp:vartheta} and \cref{item3:cor:suitable_critical_upper_bound} in \cref{cor:suitable_critical_upper_bound} establish \cref{eqn:prop:convergence_gf_suitable_cp:risk_bound}.
\end{cproof}

\begin{lemma}\label{lemma:N>f_at_1}
Assume \cref{setting1}, assume for all $x \in \R$ that $\fp(x) \ge 0$, assume that $f$ is strictly convex, let $L \in (0, \infty)$ satisfy $L = \sup_{x \in [\0, \1]} f'(x)$, assume $\min_{j \in \{1, 2, \ldots, \fh\}} v_j > 0$, let $\Theta \in \allowbreak C([0, \infty), \allowbreak \R^{\fh})$ satisfy for all $t \in [0, \infty)$ that
\begin{equation}
\textstyle \Theta_0 \in \R^{\fh} \qandq \Theta_t = \Theta_0 + \int_0^t \cG(\Theta_s) \, \d s,
\end{equation}
assume $\sum_{j \in \{i \in \{1, 2, \ldots, \fh\} \colon \Theta_{i, 0} < \1\}} v_j \ge L$, and let $\vartheta = (\vartheta_1, \ldots, \vartheta_{\fh}) \in \R^{\fh}$ satisfy
\begin{equation}\label{eqn:lemma:N>f_at_1:vartheta}
\textstyle \limsup_{t \to \infty} (\norm{\Theta_t - \vartheta} + \abs{\loss(\Theta_t) - \loss(\vartheta)} + \norm{\cG(\Theta_t)} + \norm{\cG(\vartheta)}) = 0.
\end{equation}
Then
\begin{equation}
\textstyle \cN^{\vartheta}(\1) \ge f(\1).
\end{equation}
\end{lemma}
\begin{cproof}{lemma:N>f_at_1}
Throughout this proof let $\scrv \in (\0, \1)$ satisfy 
\begin{equation}
\textstyle \scrv = \max_{j \in \{i \in \{1, 2, \ldots, \fh\} \colon \0 < \vartheta_i < \1\}} \vartheta_j.
\end{equation}
For the sake of contradiction we assume that $\cN^{\vartheta}(\1) < f(\1)$. \Nobs that \cref{eqn:lemma:N>f_at_1:vartheta}, \cref{lemma:critical_points}, and the assumption that $\cN^{\vartheta}(\1) < f(\1)$ ensure that $\cN^{\vartheta}(\scrv) \ge f(\scrv)$. Combining this with \cref{lemma:GF_boundedness}, the assumption that $\sum_{j \in \{i \in \{1, 2, \ldots, \fh\} \colon \Theta_{i, 0} < \1\}} v_j \ge L$, and the assumption that $f$ is strictly convex assures that there exists $\fj \in \{i \in \{1, 2, \ldots, \fh\} \colon \Theta_{i, 0} < \1\}$ which satisfies $\vartheta_{\fj} = \1$. This and \cref{eqn:lemma:N>f_at_1:vartheta} demonstrate that
\begin{equation}\label{eqn:lemma:N>f_at_1:i}
\textstyle \limsup_{t \to \infty} \abs{\Theta_{\fj, t} - \vartheta_{\fj}} = 0.
\end{equation}
The assumption that $\cN^{\vartheta}(\1) < f(\1)$ hence assures that there exists $T \in [0, \infty)$ which satisfies for all $t \in [T, \infty)$ that $\Theta_{\fj, T} < \1$, $\0 < \Theta_{\fj, t} \le \1$, and $\cN^{\Theta_t}(\Theta_{\fj, t}) < f(\Theta_{\fj, t})$. This and \cref{lemma:derivatives_of_the_risk} show that for all $t \in (T, \infty)$ it holds that
\begin{equation}
\textstyle \frac{\partial}{\partial t} \Theta_{\fj, t} = \cG_{\fj}(\Theta_t) = - \bigl(\frac{\partial}{\partial \Theta_{\fj, t}} \loss\bigr)(\Theta_t) = 2 v_{\fj} \int_{\Theta_{\fj, t}}^{\1} (\cN^{\Theta_t}(x) - f(x)) \fp(x) \, \d x < 0.
\end{equation}
This implies that $(T, \infty) \ni t \mapsto \Theta_{\fj, t} \in \R$ is decreasing. Therefore, we obtain that for all $t \in (T, \infty)$ it holds that $\Theta_{\fj, t} < \Theta_{\fj, T} < \1 = \vartheta_{\fj}$. This and \cref{eqn:lemma:N>f_at_1:i} contradict to the assumption that $\cN^{\vartheta}(\1) < f(\1)$. Hence, we obtain that $\cN^{\vartheta}(\1) \ge f(\1)$.
\end{cproof}

%------------------------------------------------------------------------------%
%----------------------------------Subsection----------------------------------%
%------------------------------------------------------------------------------%
\subsection{Convergence of GF with random initializations}

In this section we
establish in \cref{cor:main}, one of the main results of this article, convergence of \GF\ processes with suitable random initializations. For this we first derive several properties of these random initializations.
Again, the inner and outer weights are assumed to be half-normally distributed with a scaling factor depending on the number of neurons $\width$, and the inner biases are normally distributed.

\subsubsection{Upper bounds for the Gaussian tails}
\begin{lemma}\label{lemma:Gaussian_tail_bound}
It holds for all $y \in [0, \infty)$ that
\begin{equation}
\textstyle \int_{y}^{\infty} \bigl[\frac{2}{\pi}\bigr]^{1/2} \exp(- \frac{x^2}{2}) \, \d x \le \exp(- \frac{y^2}{2}).
\end{equation}
\end{lemma}
\begin{cproof}{lemma:Gaussian_tail_bound}
\Nobs that the integral transformation theorem and the fact that for all $x \in (0, 1)$ it holds that $\Gamma(x) \Gamma(1 - x) = \frac{\pi}{\sin(\pi x)}$ ensure that
\begin{equation}
\textstyle \int_0^{\infty} \bigl[\frac{2}{\pi}\bigr]^{1/2} \exp(- \frac{x^2}{2}) \, \d x = \frac{1}{\pi^{1/2}} \int_0^{\infty} t^{-1/2} \exp(-t) \, \d t = \frac{1}{\pi^{1/2}} \Gamma(\frac{1}{2}) = 1.
\end{equation}
Combining this with the integral transformation theorem shows that for all $y \in [0, \infty)$ it holds that
\begin{equation}
\begin{split}
\textstyle \int_{y}^{\infty} \bigl[\frac{2}{\pi}\bigr]^{1/2} \exp(- \frac{x^2}{2}) \, \d x & \textstyle = \int_0^{\infty} \bigl[\frac{2}{\pi}\bigr]^{1/2} \exp(- \frac{(x + y)^2}{2}) \, \d x = \int_0^{\infty} \bigl[\frac{2}{\pi}\bigr]^{1/2} \exp(- \frac{x^2 + 2 x y + y^2}{2}) \, \d x \\
& \textstyle \le \exp(- \frac{y^2}{2}) \int_0^{\infty} \bigl[\frac{2}{\pi}\bigr]^{1/2} \exp(- \frac{x^2}{2}) \, \d x = \exp(- \frac{y^2}{2}).
\end{split}
\end{equation}
\end{cproof}

%------------------------------------------------------------------------------%
%----------------------------------Subsection----------------------------------%
%------------------------------------------------------------------------------%
\subsubsection{Size estimates for the slope changes of the realizations of ANNs}

\begin{lemma}[Lower and upper bounds for the slope changes]\label{lemma:bound_slope_changes}
Let $\alpha, \beta \in \R$, let $(\Omega, \cF, \P)$ be a probability space, for every $\fh \in \N$ let $V_n^{\fh} \colon \Omega \to \R$, $n \in \{1, 2, \ldots, \fh\}$, and $W_n^{\fh} \colon \Omega \to \R$, $n \in \{1, 2, \ldots, \fh\}$, be random variables, and assume for all $\fh \in \N$, $n \in \{1, 2, \ldots, \fh\}$, $x \in \R$ that $V_1^{\fh}, \ldots, V_{\fh}^{\fh}, W_1^{\fh}, \ldots, W_{\fh}^{\fh}$ are independent and $\P(\fh^{\alpha} W_n^{\fh} \ge x) = \P(\fh^{\beta} V_n^{\fh} \ge x) = \int_{\max\{x, 0\}}^{\infty} [\frac{2}{\pi}]^{1/2} \exp(-\frac{y^2}{2}) \, \d y$. Then
\begin{enumerate}[label=(\roman*)]
\item
\label{item1:lemma:bound_slope_changes} it holds for all $\fh \in \N$, $r \in \R$ that
\begin{equation}
\textstyle \P(W_1^{\fh} \le \fh^{-r}) \ge 1 - \exp(- \frac{1}{2} \fh^{2 \alpha - 2 r}) \qandq \P(V_1^{\fh} \le \fh^{-r}) \ge 1 - \exp(- \frac{1}{2} \fh^{2 \beta - 2 r}),
\end{equation}

\item
\label{item2:lemma:bound_slope_changes} it holds for all $\fh \in \N$, $r, s \in \R$ that
\begin{equation}
\textstyle \P\bigl(\max_{j \in \{1, 2, \ldots, \fh\}}(V_j^{\fh} W_j^{\fh}) \le \fh^{- r}\bigr) \ge \bigl[1 - \exp(- \frac{1}{2} \fh^{2 \beta - 2 s})\bigr]^{\fh} \bigl[1 - \exp(- \frac{1}{2} \fh^{2 \alpha + 2 s - 2r})\bigr]^{\fh},
\end{equation}

and
\item
\label{item3:lemma:bound_slope_changes} it holds for all $\fh \in \N$, $\kappa \in \{1, 2, \ldots, \fh\}$, $r \in \R$ that
\begin{equation}
\textstyle \P\bigl(\sum_{j = 1}^{\kappa} V_j^{\fh} W_j^{\fh} \le \bigl[\frac{2}{\pi}\bigr] \kappa \fh^{- \alpha - \beta} - \fh^{r} \bigr) \le \kappa \fh^{- 2 \alpha - 2 \beta - 2 r} \bigl(1 - \frac{4}{\pi^2}\bigr).
\end{equation}
\end{enumerate}
\end{lemma}
\begin{cproof}{lemma:bound_slope_changes}
Throughout this proof for every $\fh \in \N$, $\kappa \in \{1, 2, \ldots, \fh\}$ let $S_{\kappa}^{\fh} \colon \Omega \to \R$ satisfy $S_{\kappa}^{\fh} = \sum_{j = 1}^{\kappa} V_j^{\fh} W_j^{\fh}$. \Nobs that \cref{lemma:Gaussian_tail_bound} ensures that for all $\fh \in \N$, $r \in \R$ it holds that
\begin{equation}
\begin{split}
& \textstyle \P(W_1^{\fh} \le \fh^{- r}) = \P(\fh^{\alpha} W_1^{\fh} \le \fh^{\alpha - r}) = 1 - \P(\fh^{\alpha} W_1^{\fh} > \fh^{\alpha - r}) \\
& \textstyle = 1 - \int_{\fh^{\alpha - r}}^{\infty} \bigl[\frac{2}{\pi}\bigr]^{1/2} \exp(-\frac{x^2}{2}) \, \d x \ge 1 - \exp(- \frac{1}{2} \fh^{2 \alpha - 2 r}) \ge 0
\end{split}
\end{equation}
and
\begin{equation}
\begin{split}
& \textstyle \P(V_1^{\fh} \le \fh^{- r}) = \P(\fh^{\beta} V_1^{\fh} \le \fh^{\beta - r}) = 1 - \P(\fh^{\beta} V_1^{\fh} > \fh^{\beta - r}) \\
& \textstyle = 1 - \int_{\fh^{\beta - r}}^{\infty} \bigl[\frac{2}{\pi}\bigr]^{1/2} \exp(-\frac{x^2}{2}) \, \d x \ge 1 - \exp(- \frac{1}{2} \fh^{2 \beta - 2 r}) \ge 0.
\end{split}
\end{equation}
This establishes \cref{item1:lemma:bound_slope_changes}. Next \nobs that \cref{item1:lemma:bound_slope_changes} and the assumption that for all $\fh \in \N$ it holds that $V_1^{\fh}, \ldots, V_{\fh}^{\fh}, W_1^{\fh}, \ldots, W_{\fh}^{\fh}$ are independent assure that for all $\fh \in \N$, $r, s \in \R$ it holds that
\begin{equation}
\begin{split}
& \textstyle \P(\max\{V_1^{\fh} W_1^{\fh}, V_2^{\fh} W_2^{\fh}, \ldots, V_{\fh}^{\fh} W_{\fh}^{\fh}\} \le \fh^{- r}) \\
& \textstyle = \bigl[\P(V_1^{\fh} W_1^{\fh} \le \fh^{- r})\bigr]^{\fh} \ge \bigl[\P([V_1^{\fh} \le \fh^{- s}] \cap [W_1^{\fh} \le \fh^{s - r}])\bigr]^{\fh} \\ 
& \textstyle = \bigl[\P([V_1^{\fh} \le \fh^{- s}]) \, \P([W_1^{\fh} \le \fh^{s - r}])\bigr]^{\fh} \\
& \textstyle \ge \bigl[1 - \exp(- \frac{1}{2} \fh^{2 \beta - 2 s})\bigr]^{\fh} \bigl[1 - \exp(- \frac{1}{2} \fh^{2 \alpha + 2 s - 2r})\bigr]^{\fh}.
\end{split}
\end{equation}
This establishes \cref{item2:lemma:bound_slope_changes}. In the next step \nobs that the assumption that for all $\fh \in \N$ it holds that $V_1^{\fh}$ and $W_1^{\fh}$ are independent and the integral transformation theorem prove that for all $\fh \in \N$ it holds that
\begin{equation}\label{eqn:lemma:bound_slope_changes:exp_d2}
\begin{split}
\textstyle \E(V_1^{\fh} W_1^{\fh}) & \textstyle = \E(V_1^{\fh}) \E(W_1^{\fh}) = \fh^{- \alpha - \beta} \bigl[\int_{- \infty}^{\infty} \frac{1}{(2 \pi)^{1/2}} \abs{x} \exp(- \frac{x^2}{2}) \, \d x\bigr]^2 \\
& \textstyle = \frac{2}{\pi} \fh^{- \alpha - \beta} \bigl[\int_0^{\infty} x \exp(- \frac{x^2}{2}) \, \d x\bigr]^2 = \frac{2}{\pi} \fh^{- \alpha - \beta} ([-\exp(- \frac{x^2}{2})]_{x=0}^{x=\infty})^2 = \frac{2}{\pi} \fh^{- \alpha - \beta}.
\end{split}
\end{equation}
and
\begin{equation}
\begin{split}
\textstyle \E([V_1^{\fh} W_1^{\fh}]^2) & \textstyle = \E([V_1^{\fh}]^2) \E([W_1^{\fh}]^2) = \fh^{- 2 \alpha - 2 \beta} \bigl[\int_{-\infty}^{\infty} \frac{1}{(2 \pi)^{1/2}} \abs{x}^2 \exp(- \frac{x^2}{2}) \, \d x\bigr]^2 \\
& \textstyle = \frac{2}{\pi} \fh^{- 2 \alpha - 2 \beta} \bigl[\int_0^{\infty} x^2 \exp(- \frac{x^2}{2}) \, \d x\bigr]^2 = \frac{4}{\pi} \fh^{- 2 \alpha - 2 \beta} \bigl[\int_0^{\infty} t^{1/2} \exp(-t) \, \d t\bigr]^2 \\
& \textstyle = \frac{4}{\pi} \fh^{- 2 \alpha - 2 \beta} \bigl[\Gamma(\frac{3}{2})\bigr]^2 = \frac{4}{\pi} \fh^{- 2 \alpha - 2 \beta} \bigl[\frac{1}{2}\Gamma(\frac{1}{2})\bigr]^2 = \fh^{- 2 \alpha - 2 \beta}.
\end{split}
\end{equation}
Hence, we obtain that for all $\fh \in \N$, $\kappa \in \{1, 2, \ldots, \fh\}$ it holds that
\begin{equation}
\begin{split}
& \textstyle \E(\abs{S_{\kappa}^{\fh} - \E(S_{\kappa}^{\fh})}^2) = \E([S_{\kappa}^{\fh}]^2) - [\E(S_{\kappa}^{\fh})]^2 \\
& \textstyle = \E\bigl(\sum_{j = 1}^{\kappa}[V_{j}^{\fh} W_{j}^{\fh}]^2 + 2 \sum_{1 \le i < j \le \kappa} V_i^{\fh} W_i^{\fh} V_j^{\fh} W_j^{\fh}\bigr) - \kappa^2 [\E(V_1^{\fh} W_1^{\fh})]^2 \\
& \textstyle = \kappa \E([V_1^{\fh} W_1^{\fh}]^2) + (\kappa^2 - \kappa) [\E(V_1^{\fh} W_1^{\fh})]^2 - \kappa^2 [\E(V_1^{\fh} W_1^{\fh})]^2 \\
& \textstyle = \kappa (\E([V_1^{\fh} W_1^{\fh}]^2) - [\E(V_1^{\fh} W_1^{\fh})]^2) = \kappa \fh^{- 2 \alpha - 2 \beta} \bigl(1 - \frac{4}{\pi^2}\bigr).
\end{split}
\end{equation}
Combining this with Markov's inequality and \cref{eqn:lemma:bound_slope_changes:exp_d2} demonstrates that for all $\fh \in \N$, $\kappa \in \{1, 2, \ldots, \fh\}$, $r \in \R$ it holds that
\begin{equation}
\begin{split}
& \textstyle \P(S_{\kappa}^{\fh} \le \bigl[\frac{2}{\pi}\bigr] \kappa \fh^{- \alpha - \beta} - \fh^{r}) \\
& \textstyle = \P(S_{\kappa}^{\fh} - \E(S_{\kappa}^{\fh}) \le - \fh^{r}) \le \P(\abs{S_{\kappa}^{\fh} - \E(S_{\kappa}^{\fh})} \ge \fh^{r}) = \P(\abs{S_{\kappa}^{\fh} - \E(S_{\kappa}^{\fh})}^2 \ge \fh^{2 r}) \\
& \textstyle \le \frac{\E(\abs{S_{\kappa}^{\fh} - \E(S_{\kappa}^{\fh})}^2)}{\fh^{2 r}} = \kappa \fh^{- 2 \alpha - 2 \beta -2 r} \bigl(1 - \frac{4}{\pi^2}\bigr).
\end{split}
\end{equation}
This establishes \cref{item3:lemma:bound_slope_changes}.
\end{cproof}

\begin{corollary}[Lower and upper bounds for the slope changes]\label{cor:bound_outer_weights}
Let $\alpha, \beta \in \R$ satisfy $0 < \alpha + \beta < 1$, let $(\Omega, \cF, \P)$ be a probability space, for every $\fh \in \N$ let $V_n^{\fh} \colon \Omega \to \R$, $n \in \{1, 2, \ldots, \fh\}$, and $W_n^{\fh} \colon \Omega \to \R$, $n \in \{1, 2, \ldots, \fh\}$, be random variables, and assume for all $\fh \in \N$, $n \in \{1, 2, \ldots, \fh\}$, $x \in \R$ that $V_1^{\fh}, \ldots, V_{\fh}^{\fh}, W_1^{\fh}, \ldots, W_{\fh}^{\fh}$ are independent and $\P(\fh^{\alpha} W_n^{\fh} \ge x) = \P(\fh^{\beta} V_n^{\fh} \ge x) = \int_{\max\{x, 0\}}^{\infty} [\frac{2}{\pi}]^{1/2} \exp(-\frac{y^2}{2}) \, \d y$. Then
\begin{enumerate}[label=(\roman*)]
\item
\label{item1:cor:bound_outer_weights} it holds for all $\fh \in \N$, $r \in (0, \alpha + \beta)$ that
\begin{equation}
\textstyle \P(\max\{V_1^{\fh} W_1^{\fh}, V_2^{\fh} W_2^{\fh}, \ldots, V_{\fh}^{\fh} W_{\fh}^{\fh}\} \le \fh^{r - \alpha - \beta}) \ge 1 - 2 \fh \exp(-\frac{1}{2} \fh^{r})
\end{equation}

and
\item
\label{item2:cor:bound_outer_weights} it holds for all $\sigma \in (0, 1]$, $\eps \in (\alpha + \beta, 1]$, $r \in (\max\{(\eps / 2) - \alpha - \beta, 0\}, \eps - \alpha - \beta)$ that there exists $\fC_{\sigma, \eps, r} \in (0, \infty)$ such that for all $\fh \in \N$ it holds that
\begin{equation}
\textstyle \P\Bigl(\sum_{j = 1}^{\ceil{\sigma \fh^{\eps}}} V_j^{\fh} W_j^{\fh} \ge h^{r} \Bigr) \ge 1 - \fC_{\sigma, \eps, r} \fh^{\eps - 2 \alpha - 2 \beta - 2r}.
\end{equation}
\end{enumerate}
\end{corollary}
\begin{cproof}{cor:bound_outer_weights}
Throughout this proof let $\sigma \in (0, 1]$, $\eps \in (\alpha + \beta, 1]$, $r \in (\max\{(\eps / 2) - \alpha - \beta, 0\}, \eps - \alpha - \beta)$ and let $\fC_{\sigma, \eps, r} \in \R$ satisfy 
\begin{equation}\label{eqn:cor:bound_outer_weights:fC_gamma}
\textstyle \fC_{\sigma, \eps, r} = \max\Bigl\{2 \sigma, \bigl[\frac{\pi}{\sigma}\bigr]^{\frac{2 r + 2 \alpha + 2 \beta - \eps}{\eps - \alpha - \beta - r}} \max_{\fh \in \N \cap (0, [\pi / \sigma]^{1/(\eps - \alpha - \beta - r)})} \P\Bigl(\sum_{j = 1}^{\ceil{\sigma \fh^{\eps}}} V_j^{\fh} W_j^{\fh} < \fh^{r} \Bigr)\Bigr\}.
\end{equation}
\Nobs that \cref{item2:lemma:bound_slope_changes} in \cref{lemma:bound_slope_changes} ensures that for all $r \in (0, \alpha + \beta)$, $\fh \in \N$ it holds that
\begin{equation}
\textstyle \P(\max\{V_1^{\fh} W_1^{\fh}, V_2^{\fh} W_2^{\fh}, \ldots, V_{\fh}^{\fh} W_{\fh}^{\fh}\} \le \fh^{r - \alpha - \beta}) \ge [1 - \exp(-\frac{1}{2} \fh^{r})]^{2 \fh}.
\end{equation}
Combining this and the fact that for all $\fh \in \N$, $r \in (0, \infty)$ it holds that $\exp(-\frac{1}{2} \fh^{r}) \allowbreak < \allowbreak 1$ with Bernoulli's inequality demonstrates that for all $\fh \in \N$, $r \in (0, \alpha + \beta)$ it holds that
\begin{equation}
\textstyle \P(\max\{V_1^{\fh} W_1^{\fh}, V_2^{\fh} W_2^{\fh}, \ldots, V_{\fh}^{\fh} W_{\fh}^{\fh}\} \le \fh^{r - \alpha - \beta}) \ge [1 - \exp(-\frac{1}{2} \fh^{r})]^{2 \fh} \ge 1 - 2 \fh \exp(-\frac{1}{2} \fh^{r}).
\end{equation}
This establishes \cref{item1:cor:bound_outer_weights}. Next \nobs that \cref{eqn:cor:bound_outer_weights:fC_gamma} ensures for all $\fh \in \N \cap (0, [\pi / \sigma]^{1/(\eps - \alpha - \beta - r)})$ that
\begin{equation}\label{eqn:cor:bound_outer_weights:prob_bound_small_fh}
\textstyle \P\Bigl(\sum_{j = 1}^{\ceil{\sigma \fh^{\eps}}} V_j^{\fh} W_j^{\fh} < \fh^{r} \Bigr) \le \bigl[\frac{\pi}{\sigma}\bigr]^{\frac{\eps - 2 \alpha - 2 \beta -2 r}{\eps - \alpha - \beta - r}} \fC_{\sigma, \eps, r} \le \fC_{\sigma, \eps, r} \fh^{\eps - 2 \alpha - 2 \beta -2 r}.
\end{equation}
In the next step \nobs that \cref{item3:lemma:bound_slope_changes} in \cref{lemma:bound_slope_changes} (applied with $\kappa \with \ceil{\sigma \fh^{\eps}}$ in the notation of \cref{item3:lemma:bound_slope_changes} in \cref{lemma:bound_slope_changes}) shows that for all $\fh \in \N$ it holds that
\begin{equation}
\begin{split}
& \textstyle \P\Bigl(\sum_{j = 1}^{\ceil{\sigma \fh^{\eps}}} V_j^{\fh} W_j^{\fh} < \frac{2}{\pi} \sigma h^{\eps - \alpha - \beta} - \fh^{r}\Bigr) \le \P\Bigl(\sum_{j = 1}^{\ceil{\sigma \fh^{\eps}}} V_j^{\fh} W_j^{\fh} \le \frac{2}{\pi} \ceil{\sigma \fh^{\eps}} h^{- \alpha - \beta} - \fh^{r}\Bigr) \\
& \textstyle \le \ceil{\sigma \fh^{\eps}} \fh^{-2 \beta - 2 \alpha - 2 r} (1 - \frac{4}{\pi^2}) \le 2 \sigma \fh^{\eps - 2 \alpha - 2 \beta - 2 r} (1 - \frac{4}{\pi^2}) \le 2 \sigma \fh^{\eps - 2 \alpha - 2 \beta - 2 r}.
\end{split}
\end{equation}
Combining this with \cref{eqn:cor:bound_outer_weights:fC_gamma} demonstrates that for all $\fh \in \N \cap [[\pi / \sigma]^{1/(\eps - \alpha - \beta - r)}, \infty)$ it holds that
\begin{equation}
\begin{split}
& \textstyle \P\Bigl(\sum_{j = 1}^{\ceil{\sigma \fh^{\eps}}} V_j^{\fh} W_j^{\fh} < \fh^{r} \Bigr) \le \P\Bigl(\sum_{j = 1}^{\ceil{\sigma \fh^{\eps}}} V_j^{\fh} W_j^{\fh} < \frac{2}{\pi} \sigma h^{\eps - \alpha - \beta} - \fh^{r}\Bigr) \\
& \textstyle \le 2 \sigma \fh^{\eps - 2 \alpha - 2 \beta - 2 r} \le \fC_{\sigma, \eps, r} \fh^{\eps - 2 \alpha - 2 \beta - 2 r}.
\end{split}
\end{equation}
This and \cref{eqn:cor:bound_outer_weights:prob_bound_small_fh} prove that for all $\fh \in \N$ it holds that
\begin{equation}
\textstyle \P\Bigl(\sum_{j = 1}^{\ceil{\sigma \fh^{\eps}}} V_j^{\fh} W_j^{\fh} < \fh^{r} \Bigr) \le \fC_{\sigma, \eps, r} \fh^{\eps - 2 \alpha - 2 \beta - 2 r}.
\end{equation}
Therefore, we obtain that for all $\fh \in \N$ it holds that
\begin{equation}
\textstyle \P\Bigl(\sum_{j = 1}^{\ceil{\sigma \fh^{\eps}}} V_j^{\fh} W_j^{\fh} \ge \fh^{r} \Bigr) = 1 - \P\Bigl(\sum_{j = 1}^{\ceil{\sigma \fh^{\eps}}} V_j^{\fh} W_j^{\fh} < \fh^{r} \Bigr) \ge 1 - \fC_{\sigma, \eps, r} \fh^{\eps - 2 \alpha - 2 \beta - 2 r}.
\end{equation}
This establishes \cref{item2:cor:bound_outer_weights}.
\end{cproof}

%------------------------------------------------------------------------------%
%----------------------------------Subsection----------------------------------%
%------------------------------------------------------------------------------%
\subsubsection{Location of the kinks of the realizations of ANNs}

\begin{lemma}[Location of the kinks]\label{lemma:kinks_location}
Let $\alpha, \gamma \in \R$, $\1 \in (0, \infty)$ satisfy $\alpha < \gamma$, let $(\Omega, \cF, \P)$ be a probability space, for every $\fh \in \N$ let $B_n^{\fh} \colon \Omega \to \R$, $n \in \{1, 2, \ldots, \fh\}$, and $W_n^{\fh} \colon \Omega \to \R$, $n \in \{1, 2, \ldots, \fh\}$, be random variables, and assume for all $\fh \in \N$, $n \in \{1, 2, \ldots, \fh\}$, $x \in \R$ that $B_1^{\fh}, \ldots, B_{\fh}^{\fh}, W_1^{\fh}, \ldots, W_{\fh}^{\fh}$ are independent, $\fh^{\gamma} B_n^{\fh}$ is standard normal, and $\P(\fh^{\alpha} W_n^{\fh} \ge x) = \int_{\max\{x, 0\}}^{\infty} [\frac{2}{\pi}]^{1/2} \exp(-\frac{y^2}{2}) \, \d y$. Then it holds for all $\fh \in \N \cap [(2 / (\pi \1))^{1 / (\gamma - \alpha)}, \infty)$, $\kappa \in \{1, 2, \ldots, \fh\}$ that
\begin{equation}
\textstyle \P(\min_{n \in \{1, 2, \ldots, \kappa\}}([B_n^{\fh} / W_n^{\fh}]) > - \1) \ge 1 - \frac{2 \kappa}{\pi \1} \fh^{\alpha - \gamma}. 
\end{equation}
\end{lemma}
\begin{cproof}{lemma:kinks_location}
Throughout this proof let $\xi_{1, n}^{\fh} \colon \Omega \to \R$, $\fh \in \N$, $n \in \{1, 2, \ldots, \fh\}$, and $\xi_{2, n}^{\fh} \colon \Omega \to \R$, $\fh \in \N$, $n \in \{1, 2, \ldots, \fh\}$, be standard normal distributions and assume for all $\fh \in \N$, $n \in \{1, 2, \ldots, \fh\}$ that $B_n^{\fh} = \fh^{- \gamma} \xi_{1, n}^{\fh}$ and $W_n^{\fh} = \fh^{- \alpha} \abs{\xi_{2, n}^{\fh}}$. \Nobs that for all $\fh \in \N$, $n \in \{1, 2, \ldots, \fh\}$ it holds that
\begin{equation}
\textstyle \Forall x \in \R \colon \P([\xi_{1, n}^{\fh} / \xi_{2, n}^{\fh}] \ge x) = \int_x^{\infty} \frac{1}{\pi (x^2 + 1)} \, \d x.
\end{equation}
The fact that for all $x \in [0, \infty]$ it holds that $\lim_{y \to x, y \in (0, \infty)} [\arctan(y) + \arctan(1/y)] = \pi/2$, the fact that for all $x \in [0, \infty)$ it holds that $\arctan(x) \le x$, and the integral transformation theorem therefore imply that for all $\fh \in \N$ it holds that
\begin{equation}
\begin{split}
& \textstyle \P([B_1^{\fh}/W_1^{\fh}] > - \1) = \P([\xi_{1, 1}^{\fh} / \abs{\xi_{2, 1}^{\fh}}] > - \1 \fh^{\gamma - \alpha}) \ge \P(\1 \fh^{\gamma - \alpha} > [\xi_{1, 1}^{\fh} / \abs{\xi_{2, 1}^{\fh}}] > - \1 \fh^{\gamma - \alpha}) \\
& \textstyle \ge \P(\abs{\xi_{1, 1}^{\fh} / \xi_{2, 1}^{\fh}} < \1 \fh^{\gamma - \alpha}) = \int_{- \1 \fh^{\gamma - \alpha}}^{\1 \fh^{\gamma - \alpha}} \frac{1}{\pi (x^2 + 1)} \, \d x = 1 - \int_{- \infty}^{- \1 \fh^{\gamma - \alpha}} \frac{1}{\pi (x^2 + 1)} \, \d x - \int_{\1 \fh^{\gamma - \alpha}}^{\infty} \frac{1}{\pi (x^2 + 1)} \, \d x \\
& \textstyle = 1 - 2 \int_{\1 \fh^{\gamma - \alpha}}^{\infty} \frac{1}{\pi (x^2 + 1)} \, \d x = 1 - \frac{2}{\pi}([\lim_{x \to \infty} \arctan(x)] - \arctan(\1 \fh^{\gamma - \alpha})) \\
& \textstyle = 1 - \frac{2}{\pi}(\frac{\pi}{2} - \arctan(\1 \fh^{\gamma - \alpha})) = 1 - \frac{2}{\pi} \arctan(\1^{-1} \fh^{\alpha - \gamma}) \ge 1 - \frac{2}{\pi \1} \fh^{\alpha - \gamma}.
\end{split}
\end{equation}
Combining this, the assumption that $\alpha < \gamma$, the assumption that for all $\fh \in \N$ it holds that $B_1^{\fh}, \ldots, B_{\fh}^{\fh}, W_1^{\fh}, \ldots, W_{\fh}^{\fh}$ are independent, and the fact that for all $\fh \in \N \cap [(2 / (\pi \1))^{1 / (\gamma - \alpha)}, \allowbreak \infty)$ it holds that $\frac{2}{\pi \1} \fh^{\alpha - \gamma} \le 1$ with Bernoulli's inequality ensures that for all $\fh \in \N \cap [(2 / (\pi \1))^{1 / (\gamma - \alpha)}, \allowbreak \infty)$, $\kappa \in \{1, 2, \ldots, \fh\}$ it holds that
\begin{equation}
\textstyle \P(\min_{n \in \{1, 2, \ldots, \kappa\}}([B_n^{\fh} / W_n^{\fh}]) > - \1) = [\P([B_1^{\fh}/W_1^{\fh}] > - \1)]^{\kappa} \ge [1 - \frac{2}{\pi \1} \fh^{\alpha - \gamma}]^{\kappa} \ge 1 - \frac{2 \kappa}{\pi \1} \fh^{\alpha - \gamma}.
\end{equation}
\end{cproof}

\begin{corollary}[Location of the kinks]\label{cor:kinks_location}
Let $\alpha, \gamma \in \R$, $\1 \in (0, \infty)$ satisfy $\alpha < \gamma$, let $(\Omega, \cF, \P)$ be a probability space, for every $\fh \in \N$ let $B_n^{\fh} \colon \Omega \to \R$, $n \in \{1, 2, \ldots, \fh\}$, and $W_n^{\fh} \colon \Omega \to \R$, $n \in \{1, 2, \ldots, \fh\}$, be random variables, and assume for all $\fh \in \N$, $n \in \{1, 2, \ldots, \fh\}$, $x \in \R$ that $B_1^{\fh}, \ldots, B_{\fh}^{\fh}, W_1^{\fh}, \ldots, W_{\fh}^{\fh}$ are independent, $\fh^{\gamma} B_n^{\fh}$ is standard normal, and $\P(\fh^{\alpha} W_n^{\fh} \ge x) = \int_{\max\{x, 0\}}^{\infty} [\frac{2}{\pi}]^{1/2} \exp(-\frac{y^2}{2}) \, \d y$. Then it holds for all $\fh \in \N \cap [(2 / (\pi \1))^{1 / (\gamma - \alpha)}, \allowbreak \infty)$, $\sigma, \eps \in (0, 1]$ that
\begin{equation}
\textstyle \P\bigl(\min_{n \in \{1, 2, \ldots, \ceil{\sigma \fh^{\eps}}\}}([B_n^{\fh} / W_n^{\fh}]) > - \1\bigr) \ge 1 - \frac{4 \sigma}{\pi \1} \fh^{\eps + \alpha - \gamma}. 
\end{equation}
\end{corollary}
\begin{cproof}{cor:kinks_location}
\Nobs that \cref{lemma:kinks_location} ensures that for all $\fh \in \N \cap [(2 / (\pi \1))^{1 / (\gamma - \alpha)}, \allowbreak \infty)$, $\sigma, \eps \in (0, 1]$ it holds that
\begin{equation}
\textstyle \P\bigl(\min_{n \in \{1, 2, \ldots, \ceil{\sigma \fh^{\eps}}\}}([B_n^{\fh} / W_n^{\fh}]) > - \1\bigr) \ge 1 - \frac{2 \ceil{\sigma \fh^{\eps}}}{\pi \1} \fh^{\alpha - \gamma} \ge 1 - \frac{4 \sigma \fh^{\eps}}{\pi \1} \fh^{\alpha - \gamma} = 1 - \frac{4 \sigma}{\pi \1} \fh^{\eps + \alpha - \gamma}.
\end{equation}
\end{cproof}

%------------------------------------------------------------------------------%
%----------------------------------Subsection----------------------------------%
%------------------------------------------------------------------------------%
\subsubsection{Size estimates for the outer biases of ANNs}

\begin{lemma}[Size estimates for the outer biases]\label{lemma:outer_bias_size}
Let $\delta \in \R$, $\fh \in \N$, let $(\Omega, \cF, \P)$ be a probability space, let $C \colon \Omega \to \R$ be random variable, and assume that $\fh^{\delta} C$ is standard normal. Then
\begin{equation}
\textstyle \Forall \kappa \in \R \colon \P(- \fh^{- \kappa} \le C \le \fh^{- \kappa}) \ge 1 - \exp(-\frac{1}{2} \fh^{2 (\delta - \kappa)}). 
\end{equation}
\end{lemma}
\begin{cproof}{lemma:outer_bias_size}
\Nobs that \cref{lemma:Gaussian_tail_bound}, the assumption that $\fh^{\delta} C$ is standard normal, and the integral transformation theorem ensure that for all $\kappa \in \R$ it holds that
\begin{equation}
\begin{split}
& \textstyle \P(- \fh^{- \kappa} \le C \le \fh^{- \kappa}) = \P(- \fh^{\delta - \kappa} \le \fh^{\delta} C \le \fh^{\delta - \kappa}) = \int_{- \fh^{\delta - \kappa}}^{\fh^{\delta - \kappa}} \frac{1}{(2 \pi)^{1/2}} \exp(- \frac{x^2}{2}) \, \d x \\
& \textstyle = \int_{0}^{\fh^{\delta - \kappa}} \bigl[\frac{2}{\pi}\bigr]^{1/2} \exp(- \frac{x^2}{2}) \, \d x = 1 - \int_{\fh^{\delta - \kappa}}^{\infty} \bigl[\frac{2}{\pi}\bigr]^{1/2} \exp(- \frac{x^2}{2}) \, \d x \ge 1 - \exp(-\frac{1}{2} \fh^{2 (\delta - \kappa)}).
\end{split}
\end{equation}
\end{cproof}

%------------------------------------------------------------------------------%
%----------------------------------Subsection----------------------------------%
%------------------------------------------------------------------------------%
\subsubsection{GF convergence with random initializations}
\label{sec:ANNs_with_ReLU_final_result}

We next combine in \cref{prop:main} the properties of the random initialization established above. Afterwards, we employ \cref{prop:convergence_gf_suitable_cp} to show in \cref{cor:main} convergence of the \GF\ to a good critical point with high probability.

\cfclear
\begin{proposition}\label{prop:main}
Let $\0 \in \R$, $\1 \in (\max\{\0, 0\}, \infty)$, let $f \in C^1(\R, \R)$ and $\fp \in C(\R, \R)$
be piecewise analytic,
assume for all $x \in (\0, \1)$ that $f'(x) > f(\0) = 0 \le \fp(x)$, $\fp^{-1}(\R \backslash \{0\}) = (\0, \1)$, and $f$ is strictly convex, let $(\Omega, \cF, \P)$ be a probability space, for every $\fh \in \N$ let $\dimension_{\fh} = 3 \fh + 1$, let $\loss_{\fh} \allowbreak \colon \allowbreak \R^{\dimension_{\fh}} \allowbreak \to \allowbreak \R$ satisfy for all $\theta = (\theta_1, \ldots, \theta_{\dimension_{\fh}}) \in \R^{\dimension_{\fh}}$ that
\begin{equation}
\textstyle \loss_{\fh}(\theta) = \int_{\R} (f(x) - \theta_{\dimension_{\fh}} - \sum_{j = 1}^{\fh} \theta_{2 \fh + j} \max\{\theta_{\fh + j} x + \theta_{j}\})^2 \fp(x) \, \d x,
\end{equation}
let $\cG_{\fh} = (\cG_{\fh}^1, \ldots, \cG_{\fh}^{\dimension_{\fh}}) \colon \R^{\dimension_{\fh}} \to \R^{\dimension_{\fh}}$ satisfy for all $\theta = (\theta_1, \ldots, \theta_{\dimension_{\fh}}) \in \R^{\dimension_{\fh}}$, $k \in \{1, 2, \ldots, \dimension_{\fh}\}$ that
\begin{equation}\label{eqn:prop:main:gradient}
\textstyle \cG_{\fh}^{k}(\theta) = (\frac{\partial^{-}}{\partial \theta_k} \loss_{\fh})(\theta) \mathbbm{1}_{\N \cap (0, \fh]}(k),
\end{equation}
and let $\Theta^{\fh} = (\Theta^{\fh, 1}, \ldots, \Theta^{\fh, \dimension_{\fh}}) \colon [0, \infty) \times \Omega \to \R^{\dimension_{\fh}}$ be a stoch. proc. with c.s.p., assume for all $\fh \in \N$, $t \in [0, \infty)$ that $\P(\Theta_t^{\fh} = \Theta_0^{\fh} - \int_0^{t} \cG_{\fh}(\Theta_s^{\fh}) \, \d s) = 1$, assume for all $\fh \in \N$ that $\Theta_0^{\fh, 1}, \ldots, \Theta_0^{\fh, \dimension_{\fh}}$ are independent, let $\alpha, \beta, \gamma \in \R$, $\delta \in (0, \infty)$ satisfy $0 < \alpha + \beta < 1$ and $2 \alpha + \beta < \gamma$, assume for all $\fh \in \N$, $k \in \N \cap (0, \fh]$ that $\fh^{\gamma} \Theta_0^{\fh, k}$ and $\fh^{\delta} \Theta_0^{\fh, \dimension_{\fh}}$ are standard normal, assume for all $\fh \in \N$, $k \in \N \cap (\fh, 2 \fh]$, $l \in \N \cap (2 \fh, 3 \fh]$, $x \in \R$ that 
\begin{equation}\label{eqn:prop:main:vjwj}
\textstyle \P\bigl(\fh^{\alpha} \Theta_0^{\fh, k} \ge x\bigr) = \P\bigl(\fh^{\beta} \Theta_0^{\fh, l} \ge x\bigr) = \int_{\max\{x, 0\}}^{\infty} [\frac{2}{\pi}]^{1/2} \exp(-\frac{y^2}{2}) \, \d y,
\end{equation}
and let $\eps_1, \eps_2 \in (0, \infty)$ satisfy $\eps_1 < \min\{\alpha + \beta, \delta\}$ and $\eps_2 \allowbreak < \min\{1, \gamma - 2 \alpha - \beta, (\gamma - \alpha)/2\}$. Then there exist $\Psi^{\fh} = (\Psi^{\fh, 1}, \ldots, \Psi^{\fh, \dimension_{\fh}}) \colon \Omega \to \R^{\dimension_{\fh}}$, $\fh \in \N$, and $\fc \in (0, \infty)$ such that
\begin{enumerate}[label=(\roman*)]
\item
\label{item1:prop:main} it holds for all $\fh \in \N$ that $\P(\limsup_{t \to \infty} \norm{\Theta_t^{\fh} - \Psi^{\fh}} = 0) = 1$,

\item
\label{item2:prop:main} it holds for all $\fh \in \N$ that $\P\bigl(\min_{j \in \{1, 2, \ldots, 2\fh\}} \Psi^{\fh, \fh + j} > 0\bigr) = \P\bigl(\min_{j \in \{1, 2, \ldots, 2\fh\}} \Theta_0^{\fh, \fh + j} > 0\bigr) = 1$,

\item
\label{item3:prop:main} it holds for all $\fh \in \N$ that $\P\bigl(\max_{j \in \{1, 2, \ldots, \fh\}}(\Psi^{\fh, \fh + j} \Psi^{\fh, 2 \fh + j} ) \le \fh^{- \eps_1}\bigr) \ge 1 - \fc \fh ^{- \eps_2}$,

\item
\label{item4:prop:main} it holds for all $\fh \in \N$ that $\P\bigl(\abs{\Psi^{\fh, \dimension_{\fh}}} \le \fh^{- \eps_1}\bigr) \ge 1 - \fc \fh^{- \eps_2}$,

\item
\label{item5:prop:main} it holds for all $\fh \in \N$ that $\P\bigl(\Psi^{\fh} \text{ is not a \descritic critical point of } \loss_{\fh}\bigr) = 1$, and

\item
\label{item6:prop:main} it holds for all $\fh \in \N$ that $\P\bigl(\cN^{\Psi^{\fh}}(\1) \ge f(\1)\bigr) \ge 1 - \fc \fh^{- \eps_2}$\ifnocf.
\end{enumerate}
(cf.\ \cref{def:strict_saddle}).
\end{proposition}
\begin{cproof}{prop:main}
Throughout this proof let $L, r_1, r_2, r_3 \in (0, \infty)$ satisfy
\begin{equation}\label{eqn:prop:main:consts}
\begin{split}
& \textstyle L = \sup_{x \in [\0, \1]} f'(x), \qquad \max\{\alpha + \beta, \eps_2\} < r_1 < \min\{\gamma - \alpha - \eps_2, 1\}, \\
& \textstyle \frac{r_1 + \eps_2}{2} - \alpha - \beta < r_2 < r_1 - \alpha - \beta, \qqandqq r_3 < \alpha + \beta - \eps_1,
\end{split}
\end{equation}
let $c_1, c_2, c_3 \in (0, \infty)$ satisfy\footnote{\Nobs that \cref{item2:cor:bound_outer_weights} in \cref{cor:bound_outer_weights} provides the existence of $c_2$.} for all $\fh \in \N$ that
\begin{equation}\label{eqn:prop:main:c1c2c3}
\begin{split}
& \textstyle c_1 \ge \sup_{\fk \in \N} \bigl(2 \fk^{1 + \eps_2} \exp(-\frac{1}{2} \fk^{\min\{r_3, 2 (\delta - \eps_1)\}})\bigr), \\
& \textstyle \P\bigl(\sum_{j = 1}^{\ceil{\fh^{r_1}}} \Theta_0^{\fh, \fh + j} \Theta_0^{\fh, \fh + 2 j} \ge \fh^{r_2}\bigr) \ge 1 - c_2 \fh^{r_1 - 2 \beta - 2 \alpha - 2 r_2}, \\
& \textstyle \andq c_3 \ge \max_{\fh \in \N \cap (0, \max\{L^{1/r_2}, (2 / (\pi \1))^{1 / (\gamma - \alpha)}\})} \fh^{\eps_2}\bigl[1 - \P\bigl(\cN^{\Psi^{\fh}}(\1) \ge f(\1)\bigr)\bigr],
\end{split}
\end{equation}
and let $\fc \in (0, \infty)$ satisfy $\fc \ge \max\{c_1, \frac{4}{\pi \1} + c_2, c_3\}$. \Nobs that \cref{eqn:prop:main:vjwj} ensures that for all $\fh \in \N$ it holds that
\begin{equation}\label{eqn:prop:main:vjwj>0}
\textstyle \P\bigl(\min_{j \in \{1, 2, \ldots, 2\fh\}} \Theta_0^{\fh, \fh + j} > 0\bigr) = \bigl[\P\bigl(\Theta_0^{\fh, \fh + 1} > 0\bigr)\bigr]^{2 \fh} = \bigl[\int_0^{\infty} \bigl[\frac{2}{\pi}\bigr]^{1/2} \exp(-\frac{y^2}{2}) \, \d y\bigr]^{2 \fh} = 1.
\end{equation}
Combining this with \cref{lemma:GF_convergence} shows that there exists $\Psi^{\fh} = (\Psi^{\fh, 1}, \ldots, \Psi^{\fh, \dimension_{\fh}}) \colon \Omega \to \R^{\dimension_{\fh}}$, $\fh \in \N$, which satisfies for all $\fh \in \N$ that 
\begin{equation}
\textstyle \P(\limsup_{t \to \infty} \norm{\Theta_t^{\fh} - \Psi^{\fh}} = 0) \ge \P\bigl(\min_{j \in \{1, 2, \ldots, 2\fh\}} \Theta_0^{\fh, \fh + j} > 0\bigr) = 1.
\end{equation}
This establishes \cref{item1:prop:main}. Next \nobs that \cref{eqn:prop:main:gradient} ensures that for all $\fh \in \N$, $k \in \N \cap (\fh, 3 \fh + 1]$, $\theta \in \R^{3 \fh}$ it holds that $\cG_{\fh}^{k}(\theta) = 0$. This shows that
\begin{equation}\label{eqn:prop:main:fixed_components_of_gf}
\textstyle \Forall \fh \in \N, k \in \N \cap (\fh, 3 \fh + 1] \colon \P(\Psi^{\fh, k} = \Theta_0^{\fh, k}) = 1.
\end{equation}
Combining this with \cref{eqn:prop:main:vjwj>0} implies that for all $\fh \in \N$ it holds that
\begin{equation}\textstyle \P\bigl(\min_{j \in \{1, 2, \ldots, 2\fh\}}\Psi^{\fh, \fh + j} > 0\bigr) = \P\bigl(\min_{j \in \{1, 2, \ldots, 2\fh\}}\Theta_0^{\fh, \fh + j} > 0\bigr) = 1.
\end{equation}
This establishes \cref{item2:prop:main}. Next \nobs that \cref{item1:cor:bound_outer_weights} in \cref{cor:bound_outer_weights}, \cref{eqn:prop:main:vjwj}, \cref{eqn:prop:main:consts}, \cref{eqn:prop:main:c1c2c3}, and \cref{eqn:prop:main:fixed_components_of_gf} assure that for all $\fh \in \N$ it holds that
\begin{equation}
\begin{split}
& \textstyle \P\bigl(\max_{j \in \{1, 2, \ldots, \fh\}}(\Psi^{\fh, \fh + j} \Psi^{\fh, 2 \fh + j}) \le \fh^{- \eps_1}\bigr) \ge \P\bigl(\max_{j \in \{1, 2, \ldots, \fh\}}(\Psi^{\fh, \fh + j} \Psi^{\fh, 2 \fh + j}) \le \fh^{r_3 - \alpha - \beta}\bigr) \\
& \textstyle = \P\bigl(\max_{j \in \{1, 2, \ldots, \fh\}}(\Theta_0^{\fh, \fh + j} \Theta_0^{\fh, 2 \fh + j}) \le \fh^{r_3 - \alpha - \beta}\bigr) \ge 1 - 2 \fh \exp(-\frac{1}{2} \fh^{r_3}) \\
& \textstyle \ge 1 - c_1 \fh^{- \eps_2} \ge 1 - \fc \fh^{- \eps_2}.
\end{split}
\end{equation}
This establishes \cref{item3:prop:main}. Furthermore, \nobs that \cref{lemma:outer_bias_size}, \cref{eqn:prop:main:c1c2c3}, and \cref{eqn:prop:main:fixed_components_of_gf} demonstrate that for all $\fh \in \N$ it holds that
\begin{equation}
\textstyle \P\bigl(\abs{\Psi^{\fh, \dimension_{\fh}}} \le \fh^{- \eps_1}\bigr) = \P\bigl(\abs{\Theta_0^{\fh, \dimension_{\fh}}} \le \fh^{- \eps_1}\bigr) \ge 1 - \exp(-\frac{1}{2} \fh^{2 (\delta - \eps_1)}) \ge 1 - c_1 \fh^{- \eps_2} \ge 1 - \fc \fh^{- \eps_2}.
\end{equation}
This establishes \cref{item4:prop:main}. In the next step \nobs that \cref{cor:strict:saddle:conv:random} proves that the set of initial values of the \GF\ such that \GF\ converges to a \descritic critical point of the corresponding risk function has Lebesgue measure zero. This establishes \cref{item5:prop:main}. In addition, \nobs that \cref{item2:cor:bound_outer_weights} in \cref{cor:bound_outer_weights}, \cref{eqn:prop:main:vjwj}, \cref{eqn:prop:main:consts}, \cref{eqn:prop:main:c1c2c3}, and \cref{eqn:prop:main:fixed_components_of_gf} demonstrate that for all $\fh \in \N$ it holds that
\begin{equation}\label{eqn:prop:main:item4}
\begin{split}
& \textstyle \P\bigl(\sum_{j = 1}^{\ceil{\fh^{r_1}}} \Psi^{\fh, \fh + j} \Psi^{\fh, \fh + 2 j} \ge \fh^{r_2}\bigr) \\
& \textstyle = \P\bigl(\sum_{j = 1}^{\ceil{\fh^{r_1}}} \Theta_0^{\fh, \fh + j} \Theta_0^{\fh, \fh + 2 j} \ge \fh^{r_2}\bigr) \ge 1 - c_2 \fh^{r_1 - 2 \alpha - 2 \beta - 2 r_2} \ge 1 - c_2 \fh^{- \eps_2}.
\end{split}
\end{equation}
Next \nobs that \cref{eqn:prop:main:c1c2c3} proves that for all $\fh \in \N \cap (0, \max\{L^{1/r_2}, (2 / (\pi \1))^{1 / (\gamma - \alpha)}\})$ it holds that
\begin{equation}\label{eqn:prop:main:N>1_at_1}
\textstyle \P\bigl(\cN^{\Psi^{\fh}}(\1) \ge f(\1)\bigr) = 1 - \fh^{- \eps_2} \bigl(\fh^{\eps_2}\bigl[1 - \P\bigl(\cN^{\Psi^{\fh}}(\1) \ge f(\1)\bigr)\bigr]\bigr) \ge 1 - c_3 \fh^{- \eps_2} \ge 1 - \fc \fh^{- \eps_2}.
\end{equation}
In the next step \nobs that \cref{lemma:N>f_at_1}, \cref{cor:kinks_location} (applied with $\sigma \with 1$, $\eps \with r_1$ in the notation of \cref{cor:kinks_location}), \cref{eqn:prop:main:consts}, \cref{eqn:prop:main:c1c2c3}, \cref{eqn:prop:main:fixed_components_of_gf}, and \cref{eqn:prop:main:item4} show that for all $\fh \in \N \cap [\max\{L^{1/r_2}, (2 / (\pi \1))^{1 / (\gamma - \alpha)}\}, \infty)$ it holds that
\begin{equation}
\begin{split}
& \textstyle \P\bigl(\cN^{\Psi^{\fh}}(\1) \ge f(\1)\bigr) \\
& \textstyle \ge \P\bigl(\bigl[\sum_{j = 1}^{\ceil{\fh^{r_1}}} \Theta_0^{\fh, \fh + j} \Theta_0^{\fh, \fh + 2 j} \ge L\bigr] \cap \bigl[\min_{j \in \{1, 2, \ldots, \ceil{\fh^{r_1}}\}}(\Theta_0^{\fh, j} / \Theta_0^{\fh, \fh + j}) > - \1\bigr]\bigr) \\
& \textstyle = 1 - \P\bigl(\bigl[\sum_{j = 1}^{\ceil{\fh^{r_1}}} \Theta_0^{\fh, \fh + j} \Theta_0^{\fh, \fh + 2 j} < L\bigr] \cup \bigl[\min_{j \in \{1, 2, \ldots, \ceil{\fh^{r_1}}\}}(\Theta_0^{\fh, j} / \Theta_0^{\fh, \fh + j}) \le - \1\bigr]\bigr) \\
& \textstyle \ge 1 - \P\bigl(\sum_{j = 1}^{\ceil{\fh^{r_1}}} \Theta_0^{\fh, \fh + j} \Theta_0^{\fh, \fh + 2 j} < L\bigr) - \P\bigl(\min_{j \in \{1, 2, \ldots, \ceil{\fh^{r_1}}\}}(\Theta_0^{\fh, j} / \Theta_0^{\fh, \fh + j}) \le - \1\bigr) \\
& \textstyle = \P\bigl(\sum_{j = 1}^{\ceil{\fh^{r_1}}} \Theta_0^{\fh, \fh + j} \Theta_0^{\fh, \fh + 2 j} \ge L\bigr) + \P\bigl(\min_{j \in \{1, 2, \ldots, \ceil{\fh^{r_1}}\}}(\Theta_0^{\fh, j} / \Theta_0^{\fh, \fh + j}) > - \1\bigr) - 1 \\
& \textstyle \ge \P\bigl(\sum_{j = 1}^{\ceil{\fh^{r_1}}} \Theta_0^{\fh, \fh + j} \Theta_0^{\fh, \fh + 2 j} \ge \fh^{r_2}\bigr) + \P\bigl(\min_{j \in \{1, 2, \ldots, \ceil{\fh^{r_1}}\}}(\Theta_0^{\fh, j} / \Theta_0^{\fh, \fh + j}) > - \1\bigr) - 1 \\
& \ge \textstyle 1 - c_2 \fh^{- \eps_2} + 1 - \frac{4}{\pi \1} \fh^{r_1 + \alpha - \gamma} - 1 \ge 1 - (c_2 + \frac{4}{\pi \1}) \fh^{- \eps_2} \ge 1 - \fc \fh^{- \eps_2}.
\end{split}
\end{equation}
Combining this with \cref{eqn:prop:main:N>1_at_1} proves that for all $\fh \in \N$ it holds that
\begin{equation}
\textstyle \P\bigl(\cN^{\Psi^{\fh}}(\1) \ge f(\1)\bigr) \ge 1 - \fc \fh^{- \eps_2}.
\end{equation}
This establishes \cref{item6:prop:main}.
\end{cproof}

\begin{theorem}
\label{cor:main}
Let $\0 \in \R$, $\1 \in (\max\{\0, 0\}, \infty)$, let $f \in C^1(\R, \R)$, $\fp \in C(\R, \R)$ be piecewise analytic, assume for all $x \in (\0, \1)$ that $f'(x) > f(\0) = 0 \le \fp(x)$, $\fp^{-1}(\R \backslash \{0\}) = (\0, \1)$, and $f$ is strictly convex, let $(\Omega, \cF, \P)$ be a probability space, for every $\fh \in \N$ let $\dimension_{\fh} = 3 \fh + 1$, let $\loss_{\fh} \allowbreak \colon \allowbreak \R^{\dimension_{\fh}} \allowbreak \to \allowbreak \R$ satisfy for all $\theta = (\theta_1, \ldots, \theta_{\dimension_{\fh}}) \in \R^{\dimension_{\fh}}$ that
\begin{equation}
\textstyle \loss_{\fh}(\theta) = \int_{\R} (f(x) - \theta_{\dimension_{\fh}} - \sum_{j = 1}^{\fh} \theta_{2 \fh + j} \max\{\theta_{\fh + j} x + \theta_{j}\})^2 \fp(x) \, \d x,
\end{equation}
let $\cG_{\fh} = (\cG_{\fh}^1, \ldots, \cG_{\fh}^{\dimension_{\fh}}) \colon \R^{\dimension_{\fh}} \to \R^{\dimension_{\fh}}$ satisfy for all $\theta = (\theta_1, \ldots, \theta_{\dimension_{\fh}}) \in \R^{\dimension_{\fh}}$, $k \in \{1, 2, \ldots, \dimension_{\fh}\}$ that
\begin{equation}\label{eqn:cor:main:gradient}
\textstyle \cG_{\fh}^{k}(\theta) = (\frac{\partial^{-}}{\partial \theta_k} \loss_{\fh})(\theta) \mathbbm{1}_{\N \cap (0, \fh]}(k),
\end{equation}
and let $\Theta^{\fh} = (\Theta^{\fh, 1}, \ldots, \Theta^{\fh, \dimension_{\fh}}) \colon [0, \infty) \times \Omega \to \R^{\dimension_{\fh}}$ be a stoch. proc. with c.s.p., assume for all $\fh \in \N$, $t \in [0, \infty)$ that $\P(\Theta_t^{\fh} = \Theta_0^{\fh} - \int_0^{t} \cG_{\fh}(\Theta_s^{\fh}) \, \d s) = 1$, assume for all $\fh \in \N$ that $\Theta_0^{\fh, 1}, \ldots, \Theta_0^{\fh, \dimension_{\fh}}$ are independent, let $\alpha, \beta, \gamma \in \R$, $\delta \in (0, \infty)$ satisfy $0 < \alpha + \beta < 1$ and $2 \alpha + \beta < \gamma$, assume for all $\fh \in \N$, $k \in \N \cap (0, \fh]$ that $\fh^{\gamma} \Theta_0^{\fh, k}$ and $\fh^{\delta} \Theta_0^{\fh, \dimension_{\fh}}$ are standard normal, and assume for all $\fh \in \N$, $k \in \N \cap (\fh, 2 \fh]$, $l \in \N \cap (2 \fh, 3 \fh]$, $x \in \R$ that 
\begin{equation}\label{eqn:cor:main:vjwj}
\textstyle \P\bigl(\fh^{\alpha} \Theta_0^{\fh, k} \ge x\bigr) = \P\bigl(\fh^{\beta} \Theta_0^{\fh, l} \ge x\bigr) = \int_{\max\{x, 0\}}^{\infty} [\frac{2}{\pi}]^{1/2} \exp(-\frac{y^2}{2}) \, \d y
\end{equation}
(cf.\ \cref{def:piecewise_analytic}). Then there exist $\fc, \degerror \in (0, \infty)$ and random variables $\fC_{\fh} \colon \Omega \to \R$, $\fh \in \N$, such that for all $\fh \in \N$ it holds that $\limsup_{t \to \infty} \E[\loss_{\fh}(\Theta_t^{\fh})] \le \fc \fh^{- \degerror}$ and
\begin{equation}
\textstyle \P(\Forall t \in (0, \infty) \colon \loss_{\fh}(\Theta_t^{\fh}) \le \fc \fh^{- \degerror} + \fC_{\fh} t^{- 1}) \ge 1 - \fc \fh^{- \degerror}.
\end{equation}
\end{theorem}
\begin{cproof}{cor:main}
Throughout this proof let $\eps_1, \eps_2 \in (0, \infty)$ satisfy $\eps_1 < \min\{\alpha + \beta, \delta\}$ and $\eps_2 < \min\{1, \gamma - 2 \alpha - \beta, (\gamma - \alpha) / 2\}$, let $\Psi^{\fh} = (\Psi^{\fh, 1}, \ldots, \Psi^{\fh, 2 \fh}) \colon \Omega \to \R^{2 \fh}$, $\fh \in \N$, and $c \in (0, \infty)$ satisfy that
\begin{enumerate}[label=(\Roman*)]
\item
\label{item1:cor:main} it holds for all $\fh \in \N$ that $\P(\limsup_{t \to \infty} \norm{\Theta_t^{\fh} - \Psi^{\fh}} = 0) = 1$,

\item
\label{item2:cor:main} it holds for all $\fh \in \N$ that $\P\bigl(\min_{j \in \{1, 2, \ldots, 2\fh\}} \Psi^{\fh, \fh + j} > 0\bigr) = \P\bigl(\min_{j \in \{1, 2, \ldots, 2\fh\}} \Theta_0^{\fh, \fh + j} > 0\bigr) = 1$,

\item
\label{item3:cor:main} it holds for all $\fh \in \N$ that $\P\bigl(\max_{j \in \{1, 2, \ldots, \fh\}}(\Psi^{\fh, \fh + j} \Psi^{\fh, 2 \fh + j} ) \le \fh^{- \eps_1}\bigr) \ge 1 - c \fh ^{- \eps_2}$,

\item
\label{item4:cor:main} it holds for all $\fh \in \N$ that $\P\bigl(\abs{\Psi^{\fh, \dimension_{\fh}}} \le \fh^{- \eps_1}\bigr) \ge 1 - c \fh^{- \eps_2}$,

\item
\label{item5:cor:main} it holds for all $\fh \in \N$ that $\P\bigl(\Psi^{\fh} \text{ is not a \descritic critical point of } \loss_{\fh}\bigr) = 1$, and

\item
\label{item6:cor:main} it holds for all $\fh \in \N$ that $\P\bigl(\cN^{\Psi^{\fh}}(\1) \ge f(\1)\bigr) \ge 1 - c \fh^{- \eps_2}$
\end{enumerate}
(cf.\ \cref{prop:main}), and let $C, \eps, L \in (0, \infty)$ satisfy
\begin{equation}\label{eqn:cor:main:Landfc}
\textstyle C \ge \max\{3 c, (L + 1) (\1 - \0) \sup_{x \in [\0, \1]} \fp(x)\}, \,\, \eps \le \min\{\eps_1, \eps_2\}, \,\, \text{and} \,\, L = \sup_{x \in [\0, \1]} f'(x)
\end{equation}
(cf.\ \cref{def:strict_saddle}). \Nobs that \cref{eqn:cor:main:gradient}, \cref{eqn:cor:main:vjwj}, \cref{item2:cor:main}, and \cref{lemma:GF_convergence} ensure that there exist random variables $\fC_{\fh} \colon \Omega \to \R$, $\fh \in \N$, which satisfy for all $\fh \in \N$, $T \in (0, \infty)$ that
\begin{equation}\label{eqn:cor:main:Lojasiewicz}
\begin{split}
& \textstyle \min\bigl\{\P\bigl(\Forall t \in (0, \infty) \colon 0 \le \loss_{\fh}(\Theta_t^{\fh}) - \loss_{\fh}(\Psi^{\fh}) \le \fC_{\fh} t^{-1}\bigr), \P\bigl(0 \le \loss_{\fh}(\Theta_T^{\fh}) - \loss_{\fh}(\Psi^{\fh}) \le \fC_{\fh} T^{-1}\bigr)\bigr\} \\
& \textstyle \ge \P\bigl(\min_{j \in \{1, 2, \ldots, 2\fh\}} \Theta_0^{\fh, \fh + j} > 0\bigr) = 1.
\end{split}
\end{equation}
Combining this with \cref{prop:convergence_gf_suitable_cp}, \cref{item1:cor:main}, \cref{item2:cor:main}, \cref{item3:cor:main}, \cref{item4:cor:main}, \cref{item5:cor:main}, \cref{item6:cor:main}, and \cref{eqn:cor:main:Landfc} demonstrates that for all $\fh \in \N$ it holds that
\begin{equation}\label{eqn:cor:main:ineqq}
\begin{split}
& \textstyle \P\bigl(\Forall t \in (0, \infty) \colon \loss_{\fh}(\Theta_t^{\fh}) \le C \fh^{- \eps} + \fC_{\fh} t ^{-1}\bigr) \\
& \textstyle \ge \P\bigl([\Forall t \in (0, \infty) \colon 0 \le \loss_{\fh}(\Theta_t^{\fh}) - \loss_{\fh}(\Psi^{\fh}) \le \fC_{\fh} t^{-1}] \cap [\loss_{\fh}(\Psi^{\fh}) \le C \fh^{- \eps}]\bigr) \\
& \textstyle \ge \P\bigl(\Forall t \in (0, \infty) \colon 0 \le \loss_{\fh}(\Theta_t^{\fh}) - \loss_{\fh}(\Psi^{\fh}) \le \fC_{\fh} t^{-1}\bigr) + \P\bigl(\loss_{\fh}(\Psi^{\fh}) \le C \fh^{- \eps}\bigr) - 1 \\
& \textstyle = \P\bigl(\loss_{\fh}(\Psi^{\fh}) \le C \fh^{- \eps}\bigr) \ge \P\bigl(\loss_{\fh}(\Psi^{\fh}) \le C \fh^{- \eps_1}\bigr) \\
& \textstyle \ge \P\bigl(\loss_{\fh}(\Psi^{\fh}) \le \fh^{- \eps_1} (L + \fh^{- \eps_1}) (\1 - \0) \sup_{x \in [\0, \1]} \fp(x)\bigr) \\
& \textstyle \ge \P\bigl([\min_{j \in \{1, 2, \ldots, 2\fh\}} \Psi^{\fh, \fh + j} > 0] \cap [\max_{j \in \{1, 2, \ldots, \fh\}}(\Psi^{\fh, 2 \fh + j} \Psi^{\fh, \fh + j}) \le \fh^{- \eps_1}] \\
& \textstyle \quad \cap [\Psi^{\fh} \text{ is not a \descritic critical point of } \loss_{\fh}] \cap [\cN^{\Psi^{\fh}}(\1) \ge f(\1)] \cap [\abs{\Psi^{\fh, \dimension_{\fh}}} \le \fh^{- \eps_1}]\bigr) \\
& \textstyle \ge \P\bigl(\min_{j \in \{1, 2, \ldots, 2\fh\}} \Psi^{\fh, \fh + j} > 0\bigr) + \P\bigl(\max_{j \in \{1, 2, \ldots, \fh\}}(\Psi^{\fh, 2 \fh + j} \Psi^{\fh, \fh + j}) \le \fh^{- \eps_1}\bigr) \\
& \textstyle \quad + \P\bigl(\Psi^{\fh} \text{ is not a \descritic critical point of } \loss_{\fh}\bigr) + \P\bigl(\cN^{\Psi^{\fh}}(\1) \ge f(\1)\bigr) \\
& \textstyle \quad + \P\bigl(\abs{\Psi^{\fh, \dimension_{\fh}}} \le \fh^{- \eps_1}\bigr) - 4 \ge 1 + (1 - c \fh^{- \eps_2}) + 1 + (1 - c \fh^{- \eps_2}) + (1 - c \fh^{- \eps_2}) - 4 \\
& \textstyle = 1 - 3 c \fh^{- \eps_2} \ge 1 - C \fh^{- \eps_2} \ge 1 - C \fh^{- \eps}.
\end{split}
\end{equation}
Next \nobs that \cref{eqn:cor:main:ineqq} implies that
\begin{equation}
\textstyle \P\bigl(\limsup_{\fh \to \infty} \loss_{\fh}(\Psi^{\fh}) = 0\bigr) = 1.
\end{equation}
This shows that there exists $\scrc \in (0, \infty)$ which satisfies
\begin{equation}
\textstyle \Forall \fh \in \N \colon \E[(\loss_{\fh}(\Psi^{\fh}))^2] < \scrc.
\end{equation}
In addition, \nobs that \cref{eqn:cor:main:ineqq} and the H\"{o}lder inequality show that for all $\fh \in \N$ it holds that 
\begin{equation}
\begin{split}
& \textstyle \E[\loss_{\fh}(\Psi^{\fh})] = \E[\loss_{\fh}(\Psi^{\fh}) \mathbbm{1}_{\{\loss_{\fh}(\Psi^{\fh}) \le C \fh^{- \eps}\}}] + \E[\loss_{\fh}(\Psi^{\fh}) \mathbbm{1}_{\{\loss_{\fh}(\Psi^{\fh}) > C \fh^{- \eps}\}}] \\
& \textstyle \le C \fh^{- \eps} + \bigl[\E[(\loss_{\fh}(\Psi^{\fh}))^2]\bigr]^{1/2} \bigl[\P(\loss_{\fh}(\Psi^{\fh}) > C \fh^{- \eps})\bigr]^{1/2} \\
& \textstyle \le C \fh^{- \eps} + \bigl[\E[(\loss_{\fh}(\Psi^{\fh}))^2]\bigr]^{1/2} [C \fh^{- \eps}]^{1/2} \le C \fh^{- \eps} + \scrc^{1/2} [C \fh^{- \eps}]^{1/2} \\
& \textstyle \le (C + [\scrc C]^{1 / 2}) \fh^{-\eps / 2}.
\end{split}
\end{equation}
Combining this with \cref{eqn:cor:main:Lojasiewicz} and the Dominated convergence theorem ensures that for all $\fh \in \N$ it holds that
\begin{equation}\label{eqn:cor:main:expectation_bound}
\textstyle \limsup_{t \to \infty} \E[\loss_{\fh}(\Theta_t^{\fh})] = \E[\limsup_{t \to \infty} \loss_{\fh}(\Theta_t^{\fh})] = \E[\loss_{\fh}(\Psi^{\fh})] \le (C + [\scrc C]^{1 / 2}) \fh^{-\eps / 2}.
\end{equation}
Let $\fc, \degerror \in (0, \infty)$ satisfy $\fc \ge C + [\scrc C]^{1/2}$ and $\degerror \le \eps / 2$. \Nobs that \cref{eqn:cor:main:ineqq} and \cref{eqn:cor:main:expectation_bound} demonstrate that for all $\fh \in \N$ it holds that $\limsup_{t \to \infty} \E[\loss_{\fh}(\Theta_t^{\fh})] \le (C + [\scrc C]^{1 / 2}) \fh^{-\eps / 2} \le \fc \fh^{- \degerror}$ and 
\begin{equation}
\begin{split}
& \textstyle \P\bigl(\Forall t \in (0, \infty) \colon \loss_{\fh}(\Theta_t^{\fh}) \le \fc \fh^{- \degerror} + \fC_{\fh} t ^{-1}\bigr) \\
& \textstyle \ge \P\bigl(\Forall t \in (0, \infty) \colon \loss_{\fh}(\Theta_t^{\fh}) \le C \fh^{- \eps} + \fC_{\fh} t ^{-1}\bigr) \ge 1 - C \fh^{- \eps} \ge 1 - \fc \fh^{- \degerror}.
\end{split}
\end{equation}
\end{cproof}

%==============================================================================%
%------------------------------------------------------------------------------%
%=================================-----Section-----============================%
%------------------------------------------------------------------------------%
%==============================================================================%
\section{Numerical simulations}
\label{sec:numerical_simulations}

In this section we complement the analytical findings of this work by means of several numerical simulation results for shallow \ANNs\ with just one hidden layer (see \cref{fig:all_biases_clipping_1_hidden_layer} in \cref{subsec_num_sim_clipping_biases_1_hidden_layer}, \cref{fig:inner_biases_ReLU_1_hidden_layer} in \cref{subsec_num_sim_ReLU_inner_bias_1_hidden_layer}, \cref{fig:all_parameters_ReLU_1_hidden_layer_random_normal} in \cref{subsec_num_sim_ReLU_all_params_1_hidden_layer_random_normal}, and \cref{fig:all_parameters_ReLU_1_hidden_layer_He} in \cref{subsec_num_sim_ReLU_all_params_1_hidden_layer_He}) as well as for deep \ANNs\ with two (see \cref{fig:all_parameters_ReLU_2_hidden_layers_random_normal} in \cref{subsec_num_sim_ReLU_all_params_2_hidden_layers_random_normal}, \cref{fig:all_parameters_ReLU_2_hidden_layers_Xavier} in \cref{subsec_num_sim_ReLU_all_params_2_hidden_layers_Xavier}, \cref{fig:all_parameters_ReLU_2_hidden_layers_He} in \cref{subsec_num_sim_ReLU_all_params_2_hidden_layers_He}, \cref{fig:all_parameters_ReLU_2_hidden_layers_Xavier_Adam} in \cref{subsec_num_sim_ReLU_all_params_2_hidden_layers_Xavier_Adam}, and \cref{fig:all_parameters_ReLU_2_hidden_layers_He_Adam} in \cref{subsec_num_sim_ReLU_all_params_2_hidden_layers_He_Adam}) or three (see \cref{fig:all_parameters_ReLU_3_hidden_layers_Xavier} in \cref{subsec_num_sim_ReLU_all_params_3_hidden_layers_Xavier} and \cref{fig:all_parameters_ReLU_3_hidden_layers_He} in \cref{subsec_num_sim_ReLU_all_params_3_hidden_layers_He}) hidden layers. These numerical simulation results all suggest that with high probability the considered \GD\ optimization method overcomes all bad non-global local minima, does not converge to a global minimum, but does converge in probability to good non-optimal generalized critical points (cf.\ \cref{def:limit:subdiff,cor:gradient:limit:zero:abstract}) whose risk values are very close to the optimal risk value.

\begin{remark}[Source codes]\label{remark:python_codes}
The source codes used to create these numerical simulation results are all presented and can be downloaded as separate \underline{\emph{.py}} files on {\sc GitHub} from \url{https://github.com/deeplearningmethods/overcome-bad-local-minima}. All numerical simulations have been performed in {\sc Python} using {\sc TensorFlow~2.9.2}.
\end{remark}

\subsection{Training ANNs via stochastic gradient descent (SGD) optimization methods}
\label{subsection:dnn:framework}

In \cref{setting:dnn} below we describe the mathematical framework which we employ in all of our numerical simulations in \cref{subsec_num_sim_clipping_biases_1_hidden_layer,subsec_num_sim_ReLU_inner_bias_1_hidden_layer,subsec_num_sim_ReLU_all_params_1_hidden_layer_random_normal,subsec_num_sim_ReLU_all_params_1_hidden_layer_He,subsec_num_sim_ReLU_all_params_2_hidden_layers_random_normal,subsec_num_sim_ReLU_all_params_2_hidden_layers_Xavier,subsec_num_sim_ReLU_all_params_2_hidden_layers_He,subsec_num_sim_ReLU_all_params_2_hidden_layers_Xavier_Adam,subsec_num_sim_ReLU_all_params_2_hidden_layers_He_Adam,subsec_num_sim_ReLU_all_params_3_hidden_layers_Xavier,subsec_num_sim_ReLU_all_params_3_hidden_layers_He} below. In \cref{eqn:setting:dnn:wb,eqn:setting:dnn:affine,eqn:setting:dnn:actmollifier,eqn:setting:dnn:multdim,eqn:setting:dnn:realization} in \cref{setting:dnn} below we introduce the fully-connected feedforward \ANNs\ which we employ to describe our numerical simulation results. We also refer to Figure~\ref{figure_dnn_illustration} for a graphical illustration of the \ANN\ architecture used in \cref{eqn:setting:dnn:wb,eqn:setting:dnn:affine,eqn:setting:dnn:actmollifier,eqn:setting:dnn:multdim,eqn:setting:dnn:realization} in \cref{setting:dnn}.

\begin{setting}
\label{setting:dnn}
Let $(\ell_k)_{k \in \N_0} \subseteq \N$, $L, \fd \in \N$ satisfy $\fd = \sum_{k = 1}^L \ell_k (\ell_{k - 1} + 1)$, for every $\theta = (\theta_1, \dots, \theta_{\fd}) \in \R^{\fd}$ let $\fw^{k, \theta} = (\fw^{k, \theta}_{i, j})_{(i,j) \in \{1, \ldots, \ell_k\} \times \{1, \ldots, \ell_{k - 1}\}} \in \R^{\ell_k \times \ell_{k - 1}}$, $k \in \N$, and $\fb^{k, \theta} = (\fb^{k, \theta}_1, \dots, \allowbreak \fb^{k, \theta}_{\ell_k}) \in \R^{\ell_k}$, $k \in \N$, satisfy for all $k \in \{1, \dots, L\}$, $i \in \{1, \ldots, \ell_k\}$, $j \in \{1, \ldots, \ell_{k - 1}\}$ that
\begin{equation}\label{eqn:setting:dnn:wb}
\fw^{k, \theta}_{i, j} = \theta_{(i - 1) \ell_{k - 1} + j + \sum_{h = 1}^{k - 1} \ell_h (\ell_{h - 1} + 1)} \qqandqq \fb^{k, \theta}_i = \theta_{\ell_k \ell_{k - 1} + i + \sum_{h = 1}^{k - 1} \ell_h (\ell_{h - 1} + 1)},
\end{equation}
for every $k \in \N$, $\theta \in \R^{\fd}$ let $\cA_k^{\theta} = (\cA_{k, 1}^{\theta}, \ldots, \cA_{k, \ell_k}^{\theta}) \colon \R^{\ell_{k - 1}} \to \R^{\ell_k}$ satisfy for all $x \in \R^{\ell_{k - 1}}$ that
\begin{equation}\label{eqn:setting:dnn:affine}
\cA_k^{\theta}(x) = \fb^{k, \theta} + \fw^{k, \theta} x,
\end{equation}
let $\act_r \colon \R \to \R$, $r \in \N \cup\cu{\infty}$, satisfy for all $x \in \R$ that $(\cup_{r \in \N} \cu{\act_r}) \subseteq C^1(\R, \R)$, $\sup_{r \in \N} \allowbreak \sup_{y \in \allowbreak [-\abs{x}, \abs{x}]} \allowbreak \abs{(\act_{r})'(y)} < \infty$, and
\begin{equation}\label{eqn:setting:dnn:actmollifier}
\textstyle \limsup_{R \to \infty} \bigl[\sum_{r = R}^{\infty} \indicator{(0,\infty)} \big(\abs {\act_r(x) - \act_\infty(x)} + \limsup_{h \nearrow 0} \bigl| (\act_r)' (x) - \frac{\act_\infty (x + h) - \act_\infty(x)}{h} \bigr|\big) \bigr] = 0,
\end{equation}
for every $r \in \N \cup\{\infty\}$, $k \in \N$ let $\fM_{r, k} \colon \R^{\ell_k} \to \R^{\ell_k}$ satisfy for all $x = (x_1, \ldots, x_{\ell_k}) \in \R^{\ell_k}$ that
\begin{equation}\label{eqn:setting:dnn:multdim}
\textstyle
\fM_{r, k}(x) = (\act_r(x_1), \ldots, \act_r(x_{\ell_k})),
\end{equation}
for every $r \in \N \cup \{\infty\}$, $\theta \in \R^{\fd}$ let $\cN^{k, \theta}_r = (\cN^{k, \theta}_{r, 1}, \ldots, \cN^{k, \theta}_{r, \ell_k}) \colon \R^{\ell_0} \to \R^{\ell_k}$, $k \in \N$, satisfy for all $k \in \N$ that
\begin{equation}\label{eqn:setting:dnn:realization}
\cN^{1, \theta}_r = \cA^{\theta}_1 \qqandqq \cN^{k + 1, \theta}_r = \cA_{k + 1}^{\theta} \circ \fM_{r, k} \circ \cN^{k, \theta}_r,
\end{equation}
let $f = (f_1, \ldots, f_{\ell_L}) \colon \R^{\ell_0} \to \R^{\ell_L}$ be measurable, let $(\Omega, \cF, \mathbbm{P})$ be a probability space, let $X^n_k \colon \Omega \to \R^{\ell_0}$, $k, n \in \N$, be random variables, let $\batchsize \in \N$, for every $n \in \N$, $r \in \N \cup \{\infty\}$ let $\lossapp^n_{r} \colon \R^{\fd} \times \Omega \to \R$ satisfy for all $\theta \in \R^{\fd}$, $\omega \in \Omega$ that
\begin{equation}\label{eqn:setting:dnn:loss}
\textstyle \lossapp^{n}_{r}(\theta, \omega) = \frac{1}{\batchsize} \sum_{k = 1}^{\batchsize} \norm{\cN^{L, \theta}_r(X^{n}_k (\omega)) - f(X_k^{n}(\omega))}^2,
\end{equation}
let $\indexset \subseteq \N$, for every $n\in \N$ let $\fG_n = (\fG_{n}^1, \ldots, \fG_{n}^{\fd}) \colon \R^{\fd} \times \Omega \to \R^{\fd}$ satisfy for all $\theta \allowbreak \in \allowbreak \R^{\fd}$, $\omega \allowbreak \in \allowbreak \{w \in \Omega \colon ((\nabla_\theta \lossapp^{n}_r)(\theta, w))_{r \in \N} \text{ is } \text{convergent}\}$, $\fj \in \{1, 2, \ldots, \fd\}$ that
\begin{equation}
\label{eq:generalized_stochastic_gradients}
\textstyle \fG_{n}^{\fj}(\theta, \omega) = \bigl[\lim_{r \to \infty}(\frac{\partial}{\partial \theta_{\fj}} \lossapp^{n}_r)(\theta, \omega)\bigr] \indicator{\indexset}(\fj),
\end{equation}
for every $n \in \N$ let $\psi_n = (\psi_n^1, \ldots, \psi_n^\fd) \colon \R^{n \fd} \to \R^\fd$ be a function, and let $\Theta = (\Theta^1, \ldots, \Theta^{\dimension}) \allowbreak \colon \allowbreak \N_0 \allowbreak \times \allowbreak \Omega \to \R^{\fd}$ be a stochastic process which satisfies for all $n \in \N$ that
\begin{equation}\label{eqn:setting:dnn:stochastic_pr}
\textstyle \Theta_n = \Theta_{n-1} - \psi_n(\fG_1(\Theta_0), \fG_2(\Theta_1), \ldots, \fG_n(\Theta_{n-1})).
\end{equation}
\end{setting}

In \cref{eqn:setting:dnn:stochastic_pr} above we make use of the setup in, e.g., (81) in Framework~6.1 in Becker et al.~\cite{Wurstemberger2022LRV} to formulate general stochastic gradient descent optimization methods covering, for instance, the standard stochastic gradient descent optimization method (cf.\ Becker et al.~\cite[Lemma~6.2]{Wurstemberger2022LRV} and \cref{lemma:StandardSGD} below) as well as the \Adam\ (cf.\ Kingma \& Ba~\cite{DBLP:journals/corr/KingmaB14}, Becker et al.~\cite[Lemma~6.8]{Wurstemberger2022LRV}, and \cref{lemma:AdamSGD} below).

\def\layersep{2.5cm}
\begin{figure}[H]
\centering
\begin{adjustbox}{width=\textwidth}
\begin{tikzpicture}[shorten >=1pt,->,draw=black!50, node distance=\layersep]
\tikzstyle{every pin edge}=[<-,shorten <=1pt]
\tikzstyle{input neuron}=[very thick, circle,draw=red, fill=red!20, minimum size=30pt,inner sep=0pt]
\tikzstyle{output neuron}=[very thick, circle, draw=green,fill=green!20,minimum size=60pt,inner sep=0pt]
\tikzstyle{hidden neuron}=[very thick, circle,draw=blue,fill=blue!20,minimum size=52pt,inner sep=0pt]
\tikzstyle{annot} = [text width=9em, text centered]
\tikzstyle{annot2} = [text width=4em, text centered]

%----------Neuron(s) of input layer----------

\node[input neuron] (I-1) at (-3,-1.4) {$x_1$};
\node[input neuron] (I-2) at (-3,-2.6) {$x_2$};
\node(I-dots) at (-3,-3.36) {\vdots};
\node[input neuron] (I-3) at (-3,-4.3) {$x_{\ell_0}$};

%----------Neuron(s) of 1st hidden layer----------
\path[yshift = 1.5cm]
node[hidden neuron](H0-1) at (0*\layersep, -1 cm) {$\cN_{\infty, 1}^{1,\theta}(x)$};
\path[yshift = 1.5cm]
node[hidden neuron](H0-2) at (0*\layersep, -3 cm) {$\cN_{\infty, 2}^{1,\theta}(x)$};
\path[yshift = 1.5cm]
node[hidden neuron](H0-3) at (0*\layersep, -5 cm) {$\cN_{\infty, 3}^{1,\theta}(x)$};
\path[yshift = 1.5cm]
node(H0-dots) at (0*\layersep, -6.2 cm) {\vdots};
\path[yshift = 1.5cm]
node[hidden neuron](H0-4) at (0*\layersep, -7.6 cm) {$\cN_{\infty, \ell_1}^{1,\theta}(x)$};

%----------Neuron(s) of 2nd hidden layer----------
\path[yshift = 1.5cm]
node[hidden neuron](H1-1) at (1.5*\layersep, -1 cm) {$\cN_{\infty, 1}^{2,\theta}(x)$};
\path[yshift = 1.5cm]
node[hidden neuron](H1-2) at (1.5*\layersep, -3 cm) {$\cN_{\infty, 2}^{2,\theta}(x)$};
\path[yshift = 1.5cm]
node[hidden neuron](H1-3) at (1.5*\layersep, -5 cm) {$\cN_{\infty, 3}^{2,\theta}(x)$};
\path[yshift = 1.5cm]
node(H1-dots) at (1.5*\layersep, -6.2 cm) {\vdots};
\path[yshift = 1.5cm]
node[hidden neuron](H1-4) at (1.5*\layersep, -7.6 cm) {$\cN_{\infty, \ell_2}^{2,\theta}(x)$};

%----------Neuron(s) of ... hidden layer----------
\path[yshift = 0.5cm]
node(Hdot-1) at (2.7*\layersep, -1.5 cm) {$\cdots$};
\path[yshift = 0.5cm]
node(Hdot-2) at (2.7*\layersep, -2.9 cm) {$\cdots$};
\path[yshift = 0.5cm]
node(Hdot-dots) at (2.7*\layersep, -4.1 cm) {$\ddots$};
\path[yshift = 0.5cm]
node(Hdot-3) at (2.7*\layersep, -5.5 cm) {$\cdots$};

%----------Neuron(s) of last hidden layer----------
\path[yshift = 1.5cm]
node[hidden neuron](H2-1) at (3.9*\layersep, -1 cm) {$\cN_{\infty, 1}^{L - 1,\theta}(x)$};
\path[yshift = 1.5cm]
node[hidden neuron](H2-2) at (3.9*\layersep, -3 cm) {$\cN_{\infty, 2}^{L - 1,\theta}(x)$};
\path[yshift = 1.5cm]
node[hidden neuron](H2-3) at (3.9*\layersep, -5 cm) {$\cN_{\infty, 3}^{L - 1,\theta}(x)$};
\path[yshift = 1.5cm]
node(H2-dots) at (3.9*\layersep, -6.2 cm) {\vdots};
\path[yshift = 1.5cm]
node[hidden neuron](H2-4) at (3.9*\layersep, -7.6 cm) {$\cN_{\infty, \ell_{L - 1}}^{L - 1,\theta}(x)$};

%----------Neuron(s) of output layer (last layer)----------
\path[yshift = 1.5cm]
node[output neuron](O-1) at (5.4*\layersep,-1.9 cm) {$\cN_{\infty, 1}^{L,\theta}(x)$}; 
\path[yshift = 1.5cm]
node[output neuron](O-2) at (5.4*\layersep,-4.2 cm) {$\cN_{\infty, 2}^{L,\theta}(x)$}; 
\path[yshift = 1.5cm]
node(O-dots) at (5.4*\layersep, -5.5 cm) {\vdots};
\path[yshift = 1.5cm]
node[output neuron](O-3) at (5.4*\layersep,-7 cm) {$\cN_{\infty, \ell_L}^{L,\theta}(x)$};
%----------Arrow(s) from 1st to 2nd layer----------
\foreach \source in {1,2,3}
\foreach \dest in {1,2,3,4}
\path[-{>[length=2mm, width=4mm]}, line width = 0.8] (I-\source) edge (H0-\dest);

%----------Arrow(s) from 2nd to 3rd layer----------
\foreach \source in {1,2,3,4}
\foreach \dest in {1,2,3,4}
\path[-{>[length=2mm, width=4mm]}, line width = 0.8] (H0-\source) edge (H1-\dest);

%----------Arrow(s) from 3rd to ... layer----------
\foreach \source in {1,2,3,4}
\foreach \dest in {1,2,3}
\draw[-{>[length=2mm, width=4mm]}, path fading=east] (H1-\source) -- (Hdot-\dest);

%----------Arrow(s) from ... to last hidden layer----------
\foreach \source in {1,2,3}
\foreach \dest in {1,2,3,4}
\draw[-{>[length=2mm, width=4mm]}, path fading=west] (Hdot-\source) -- (H2-\dest);

%----------Arrow(s) from last hidden layer to last layer----------
\foreach \source in {1,2,3,4}
\foreach \dest in {1,2,3}
\path[-{>[length=2mm, width=4mm]}, line width = 0.8] (H2-\source) edge (O-\dest);

%----------Annotate the layers----------
\node[annot,above of=H0-1, node distance=1.5cm, align=center] (hl) {$1^{\text{st}}$ hidden layer\\($2^{\text{nd}}$ layer)};
\node[annot,above of=H1-1, node distance=1.5cm, align=center] (hl) {$2^{\text{nd}}$ hidden layer\\($3^{\text{rd}}$ layer)};
\node[annot,above of=H2-1, node distance=1.5cm, align=center] (hl2) {${(L - 1)}^{\text{th}}$ hidden layer\\($L^{\text{th}}$ layer)};
\node[annot,above of=I-1, node distance=1.1cm, align=center] {Input layer\\ ($1^{\text{st}}$ layer)};
\node[annot,above of=O-1, node distance=1.62cm, align=center] {Output layer\\(${(L + 1)}^{\text{th}}$ layer)};

\node[annot2,below of=H0-4, node distance=1.5cm, align=center] (sl) {$\ell_1$ \\ neurons};
\node[annot2,below of=H1-4, node distance=1.5cm, align=center] (sl) {$\ell_2$ \\ neurons};
\node[annot2,below of=H2-4, node distance=1.5cm, align=center] (sl2) {$\ell_{L - 1}$ neurons};
\node[annot2,below of=I-3, node distance=1.1cm, align=center] {$\ell_0$ \\ neurons};
\node[annot2,below of=O-3, node distance=1.62cm, align=center] {$\ell_L$ \\ neurons};
\end{tikzpicture}
\end{adjustbox}
\caption{Graphical illustration of the \ANN\ architecture used in \cref{eqn:setting:dnn:wb,eqn:setting:dnn:affine,eqn:setting:dnn:actmollifier,eqn:setting:dnn:multdim,eqn:setting:dnn:realization} in \cref{setting:dnn}: The realization function $\cN_{\infty}^{L, \theta} \colon \R^{\ell_0} \to \R^{\ell_L}$ of the \ANN\ associated to the parameter vector $\theta = (\theta_1, \ldots, \theta_{\dimension}) \in \R^{\dimension}$ maps the $\ell_0$-dimensional input vector $x = (x_1, \ldots, x_{\ell_0}) \in \R^{\ell_0}$ to the $\ell_{L}$-dimensional output vector $\cN_{\infty}^{L, \theta}(x) = (\cN_{\infty, 1}^{L, \theta}(x), \ldots, \cN_{\infty, \ell_L}^{L, \theta}(x)) \in \R^{\ell_L}$.}
\label{figure_dnn_illustration} 
\end{figure}

%\begin{remark}
%Assume \cref{setting:dnn}
%\end{remark}
%
%
%
\subsection{Examples of popular SGD optimization methods}

\subsubsection{Standard SGD optimization methods}
\label{subsubsec_Standard_SGD}

The next elementary lemma, \cref{lemma:StandardSGD}, illustrates that the framework in \cref{eqn:setting:dnn:stochastic_pr} in \cref{setting:dnn} covers plain vanilla standard \SGD\ optimization processes. \cref{lemma:StandardSGD} is a slightly adapted variant of, e.g., Becker et al.~\cite[Lemma~6.2]{Wurstemberger2022LRV}.

\begin{lemma}\label{lemma:StandardSGD}
Assume \cref{setting:dnn}, let $(\gamma_n)_{n\in \N}\subseteq (0, \infty)$, and assume for all $n \in \N$, $g_1, g_2, \ldots, \allowbreak g_n \in \R^{\fd}$ that $\psi_n(g_1, g_2, \ldots, g_n) = \gamma_n g_n$. Then it holds for all $n \in \N$ that 
\begin{equation}\label{eqn:standard_SGD}
\textstyle \Theta_{n} = \Theta_{n-1} - \gamma_n \fG_n(\Theta_{n-1}).
\end{equation}
\end{lemma}
\begin{cproof}{lemma:StandardSGD}
\Nobs that \cref{eqn:setting:dnn:stochastic_pr} establishes \cref{eqn:standard_SGD}.
\end{cproof}

\subsubsection{Adaptive moment estimation (Adam) SGD optimization methods}
\label{subsubsec_Adam_SGD}

The next elementary lemma, \cref{lemma:AdamSGD}, illustrates that the framework in \cref{eqn:setting:dnn:stochastic_pr} in \cref{setting:dnn} covers also the popular \Adam\ optimization processes in Kingma \& Ba~\cite{DBLP:journals/corr/KingmaB14}. \cref{lemma:AdamSGD} is a slightly adapted variant of, e.g., Becker et al.~\cite[Lemma~6.8]{Wurstemberger2022LRV}.

\begin{lemma}\label{lemma:AdamSGD}
Assume \cref{setting:dnn}, let $\alpha, \beta \in (0, 1)$, $\eps \in (0,\infty)$, let $(\gamma_n)_{n\in \N}\subseteq (0, \infty)$, assume $\indexset = \{1, 2, \ldots, \fd\}$, and assume for all $n\in \N$, $i \in \indexset$, $g_1 = (g_1^1, \ldots, g_1^{\fd})$, $g_2 = (g_2^1, \ldots, g_2^{\fd})$, $\dots$, $g_n = (g_n^1, \ldots, g_n^{\fd}) \in \R^{\fd}$ that
\begin{equation}\label{eqn:lemma:AdamSGD:psi}
\textstyle \psi_n^{i} (g_1, g_2, \ldots, g_n) \textstyle = \gamma_n \left[\frac{\sum_{k = 1}^n \alpha^{n-k} (1-\alpha) g_k^{i}}{1 - \alpha^n} \right] \left[\varepsilon + \Bigl[\frac{\sum_{k = 1}^n \beta^{n-k}(1-\beta) \abs{g_k^{i}}^2}{1 - \beta^n}\Bigr]^{\nicefrac{1}{2}} \right]^{-1}.
\end{equation}
Then there exist $\mathbf{m}  = (\mathbf{m}^{1}, \ldots, \mathbf{m}^{\fd}) \colon \N_0 \times \Omega \to \R^{\fd}$ and $\mathbb{M} = (\mathbb{M}^{1}, \ldots, \mathbb{M}^{\fd}) \colon \N_0 \times \Omega \to \R^{\fd}$ such that for all $n \in \N$, $i \in \indexset$ it holds that
\begin{equation}\label{eqn:lemma:AdamSGD:m}
\textstyle \mathbf{m}_0 = 0, \qquad \mathbf{m}_n = \alpha \mathbf{m}_{n-1} + (1 - \alpha) \fG_n(\Theta_{n - 1}),
\end{equation}
\begin{equation}\label{eqn:lemma:AdamSGD:M}
\textstyle \mathbb{M}_0 = 0, \qquad \mathbb{M}_n^{i} = \beta\,\mathbb{M}_{n-1}^{i} + (1 - \beta) \abs{\fG_{n}^i(\Theta_{n - 1})}^2,
\end{equation}
\begin{equation}\label{eqn:lemma:AdamSGD:Theta}
\textstyle \text{and} \qquad \Theta_n^i = \Theta_{n - 1}^i - \gamma_n \left[\frac{\mathbf{m}_n^{i}}{1 - \alpha^n}\right] \left[\varepsilon + \Bigl[\frac{\mathbb{M}_n^{i}}{1 - \beta^n}\Bigr]^{\nicefrac{1}{2}}\right]^{-1}.
\end{equation}
\end{lemma}
\begin{cproof}{lemma:AdamSGD}
Throughout this proof let $\mathbf{m}  = (\mathbf{m}^{1}, \ldots, \mathbf{m}^{\fd}) \colon \N_0 \times \Omega \to \R^{\fd}$ and $\mathbb{M} = (\mathbb{M}^{1}, \ldots, \mathbb{M}^{\fd}) \colon \N_0 \times \Omega \to \R^{\fd}$ satisfy for all $n \in \N_0$, $i \in \indexset$ that
\begin{equation}\label{eqn:lemma:AdamSGD:m_i}
\textstyle \mathbf{m}_n^{i} = \sum_{k = 1}^n \alpha^{n - k} (1 - \alpha) \fG_k^i(\Theta_{k - 1})
\end{equation}
and
\begin{equation}\label{eqn:lemma:AdamSGD:M_i}
\textstyle \mathbb{M}_n^{i} = \sum_{k = 1}^n \beta^{n - k} (1 - \beta) \abs{\fG_k^i(\Theta_{k - 1})}^2.
\end{equation}
\Nobs that \cref{eqn:lemma:AdamSGD:m_i,eqn:lemma:AdamSGD:M_i} ensure that for all $n \in \N$, $i \in \indexset$ it holds that
\begin{equation}\label{eqn:lemma:AdamSGD:m_i:induct}
\begin{split}
& \textstyle \mathbf{m}_n^{i} - \alpha \mathbf{m}_{n - 1}^{i} \\
& \textstyle = \sum_{k = 1}^n \alpha^{n - k} (1 - \alpha) \fG_k^i(\Theta_{k - 1}) - \alpha \sum_{k = 1}^{n - 1} \alpha^{n - 1 - k} (1 - \alpha) \fG_k^i(\Theta_{k - 1}) \\
& \textstyle = (1 - \alpha) \bigl[\sum_{k = 1}^n \alpha^{n - k} \fG_k^i(\Theta_{k - 1}) - \sum_{k = 1}^{n - 1} \alpha^{n - k} \fG_k^i(\Theta_{k - 1})\bigr] = (1 - \alpha) \fG_n^i(\Theta_{n - 1}),
\end{split}
\end{equation}
\begin{equation}\label{eqn:lemma:AdamSGD:M_i:induct}
\begin{split}
& \textstyle \mathbb{M}_n^{i} - \beta \mathbb{M}_{n - 1}^{i} \\
& \textstyle = \sum_{k = 1}^n \beta^{n - k} (1 - \beta) \abs{\fG_k^i(\Theta_{k - 1})}^2 - \beta \sum_{k = 1}^{n - 1} \beta^{n - 1 - k} (1 - \beta) \abs{\fG_k^i(\Theta_{k - 1})}^2 \\
& \textstyle = (1 - \beta) \bigl[\sum_{k = 1}^n \beta^{n - k} \abs{\fG_k^i(\Theta_{k - 1})}^2 - \sum_{k = 1}^{n - 1} \beta^{n - k} \abs{\fG_k^i(\Theta_{k - 1})}^2\bigr] = (1 - \beta) \abs{\fG_n^i(\Theta_{n - 1})}^2,
\end{split}
\end{equation}
and $\mathbf{m}_0 = \mathbb{M}_0 = 0$. Next \nobs that \cref{eqn:setting:dnn:stochastic_pr}, \cref{eqn:lemma:AdamSGD:psi}, \cref{eqn:lemma:AdamSGD:m_i}, and \cref{eqn:lemma:AdamSGD:M_i} demonstrate that for all $n \in \N$, $i \in \indexset$ it holds that
\begin{equation}
\begin{split}
& \textstyle \Theta_n^i - \Theta_{n - 1}^i = - \psi_n^i (\fG_1(\Theta_0), \fG_2(\Theta_1), \ldots, \fG_n(\Theta_{n-1})) \\
& \textstyle = - \gamma_n \left[\frac{\sum_{k = 1}^n \alpha^{n - k} (1 - \alpha) \fG_k^i(\Theta_{k - 1})}{1 - \alpha^n}\right] \left[\varepsilon + \Bigl[\frac{\sum_{k = 1}^n \beta^{n - k} (1 - \beta) \abs{\fG_k^i(\Theta_{k - 1})}^2}{1 - \beta^n}\Bigr]^{\nicefrac{1}{2}}\right]^{-1} \\
& \textstyle = - \gamma_n \left[\frac{\mathbf{m}_n^{i}}{1 - \alpha^n}\right] \left[\varepsilon + \Bigl[\frac{\mathbb{M}_n^{i}}{1 - \beta^n}\Bigr]^{\nicefrac{1}{2}}\right]^{-1}\!.
\end{split}
\end{equation}
Combining this with \cref{eqn:lemma:AdamSGD:m_i:induct,eqn:lemma:AdamSGD:M_i:induct} establishes \cref{eqn:lemma:AdamSGD:m}, \cref{eqn:lemma:AdamSGD:M}, and \cref{eqn:lemma:AdamSGD:Theta}.
\end{cproof}

\subsection{Properties of subdifferentials of convergent stochastic processes}

In the subsequent subsections of this section (\cref{subsec_num_sim_clipping_biases_1_hidden_layer,subsec_num_sim_ReLU_inner_bias_1_hidden_layer,subsec_num_sim_ReLU_all_params_1_hidden_layer_random_normal,subsec_num_sim_ReLU_all_params_1_hidden_layer_He,subsec_num_sim_ReLU_all_params_2_hidden_layers_random_normal,subsec_num_sim_ReLU_all_params_2_hidden_layers_Xavier,subsec_num_sim_ReLU_all_params_2_hidden_layers_He,subsec_num_sim_ReLU_all_params_2_hidden_layers_Xavier_Adam,subsec_num_sim_ReLU_all_params_2_hidden_layers_He_Adam,subsec_num_sim_ReLU_all_params_3_hidden_layers_Xavier,subsec_num_sim_ReLU_all_params_3_hidden_layers_He}) we complement the analytical findings of this article through several numerical simulation results for shallow and deep \ANNs. These numerical simulation results all suggest that the sum of the $L^2(\P; \R^{\dimension})$-norm of the generalized gradient function (cf., \cref{eq:generalized_stochastic_gradients} and \cite[Items~$($v$)$ and $($vi$)$ in Theorem~2.9]{DNNReLUarXiv}) of the risk function composed with the considered \SGD\ type optimization process and of the $L^2(\P; \R^{\dimension})$-distance between the \SGD\ type optimization process and a random limit point converges to zero. In this subsection we show in the elementary results in \cref{lem:gradient:limit:zero:abstract,cor:gradient:limit:zero:abstract} below that one can conclude from this convergence property (see \cref{eqn:Frechet_convergence} below) that the random limit point must be a generalized critical point in the sense of limiting Fr\'{e}chet subdifferentiability (see \cref{def:limit:subdiff}). In our formulation of the conclusion of \cref{cor:gradient:limit:zero:abstract} we thus employ the notion of the set of limiting Fr\'{e}chet subdifferentials which we briefly recall in \cref{def:limit:subdiff} below (cf., e.g., \cite[Chapter~10]{RockafellarWets1998}, \cite[Definition~2.9]{EberleJentzenRiekertWeiss2023}, and \cite[Definition~3.1]{JentzenRiekert2021ExistenceGlobMin}).

\begin{lemma}
\label{lem:gradient:limit:zero:abstract}
Let $\fd \in \N$, $p \in ( 0 , \infty )$, let $(\Omega , \cF , \P)$ be a probability space, let $\Theta = (\Theta_n )_{n \in \N_0 } \colon \N_0 \times \Omega \to \R^\fd$ be a stochastic process, let $\vartheta \colon \Omega \to \R^\fd$ be a random variable, let $\fg \colon \R^\fd \to [0 , \infty ]$ be lower semicontinuous, assume for all $n \in \N_0$ that $\E [ \norm{\Theta_n } ^p + \norm{\vartheta } ^p + \abs{ \fg ( \Theta_n ) }  ^p ] < \infty$, and assume
\begin{equation}\label{eqn:lem:gradient:limit:zero:abstract}
\textstyle \limsup_{n \to \infty} \E [ \norm{ \Theta_n - \vartheta } ^p + \abs{\fg ( \Theta_n ) } ^p ] = 0.
\end{equation}
Then $\P ( \fg ( \vartheta ) = 0 ) = 1$.
\end{lemma}

\begin{cproof}{lem:gradient:limit:zero:abstract}
\Nobs that \cref{eqn:lem:gradient:limit:zero:abstract} implies that there exist $E \in \cF$ and a strictly increasing $k \colon \N_0 \to \N_0$ which satisfy $\P ( E ) = 1$ and
\begin{equation}\label{eqn:subseq_k_n}
\textstyle \Forall \omega \in E \colon \limsup_{n \to \infty} ( \norm{ \Theta_{ k ( n ) } ( \omega ) - \vartheta ( \omega ) } ^p + \abs{ \fg  ( \Theta _ { k ( n ) } ( \omega ) ) } ^p ) = 0.
\end{equation}
\Nobs that \cref{eqn:subseq_k_n} and the fact that $\fg$ is lower semicontinuous show for all $\omega \in E$ that
\begin{equation}
\textstyle
\fg ( \vartheta ( \omega ) ) = \fg \rbr*{ \lim_{n \to \infty } \Theta_{ k ( n ) } ( \omega  ) }
\le \liminf_{n \to \infty } \fg  ( \Theta _ { k ( n ) } ( \omega ) ) = 0 .
\end{equation}
Hence, we obtain for all $\omega \in E$ that $\fg ( \vartheta ( \omega ) ) = 0$.
\end{cproof}

\begin{definition}[Fr\'{e}chet subdifferentials and limiting Fr\'{e}chet subdifferentials]
\label{def:limit:subdiff}
Let $n \in \N$, $\loss \in C(\R^n, \R)$, $x \in \R^n$. Then we denote by 
$(\cD \loss)(x) \subseteq \R^n$ the set given by
\begin{equation}
\textstyle (\cD \loss)(x) = \cu*{y \in \R^n \colon 
\br*{ \liminf_{\R^n \backslash \cu{  0 } \ni h \to 0 } 
\rbr*{ \frac{\loss(x + h) - \loss(x) - \spro{y, h}}{\norm{h}}} \geq 0}}
\end{equation}
and we denote by $(\bD \loss)(x) \subseteq \R^n$ the set given by
\begin{equation}
\textstyle (\bD \loss)(x) = \bigcap_{\varepsilon \in (0, \infty) } 
\overline{\br*{\cup_{y \in \cu{ z \in \R^n \colon \norm{ x - z } < \varepsilon}}(\cD \loss)(y)}}.
\end{equation}
\end{definition}

\cfclear
\begin{cor}
\label{cor:gradient:limit:zero:abstract}
Let $\fd \in \N$, $\loss \in C ( \R^\fd , \R )$, $p \in ( 0 , \infty )$, let $\cG \colon \R^\fd \to \R^\fd$ be measurable, assume for all $\theta \in \R^\fd$ that $\cG ( \theta ) \in ( \bD \loss ) (\theta)$, let $(\Omega , \cF , \P)$ be a probability space, let $\Theta = (\Theta_n )_{n \in \N_0 } \colon \N_0 \times \Omega \to \R^\fd$ be a stochastic process, let $\vartheta \colon \Omega \to \R^\fd$ be a random variable, assume for all $n \in \N_0$ that $\E[\norm{\Theta_n}^p + \norm{\vartheta}^p + \norm{\cG(\Theta_n)}^p] < \infty$, and assume
\begin{equation}\label{eqn:Frechet_convergence}
\textstyle \limsup_{n \to \infty} \E [ \norm{ \Theta_n - \vartheta } ^p + \norm{\cG ( \Theta_n ) } ^p ] = 0
\end{equation}
\cfadd{def:limit:subdiff}\cfload. Then $\P(0 \in (\bD \loss) (\vartheta)) = 1$.
\end{cor}
\begin{cproof}{cor:gradient:limit:zero:abstract}
Throughout this proof let $\fg \colon \R^\fd \to [0 , \infty ]$ satisfy for all $\theta \in \R^\fd$ that
\begin{equation}
\textstyle \fg ( \theta ) = \inf \rbr{ \cu{ \norm{v } \colon v \in ( \bD \loss ) ( \theta ) } \cup \cu{\infty } }.
\end{equation}
(cf., e.g., (4) in Bolte et al.~\cite{BolteDaniilidis2006}). \Nobs that, e.g., \cite[Lemma 6.1]{JentzenRiekert2021ExistenceGlobMin} ensures that $\fg$ is lower semicontinuous. In addition, \nobs that the assumption that $\forall \, \theta \in \R^\fd \colon \cG ( \theta ) \in ( \bD \loss ) ( \theta ) $ assures for all $\theta \in \R^\fd$ that $\fg ( \theta ) \le \norm{\cG ( \theta ) }$. Hence, we obtain that 
\begin{equation}
\textstyle \limsup_{n \to \infty} \E [\norm{\Theta_n - \vartheta}^p + \abs{\fg (\Theta_n)}^p] = 0.
\end{equation}
Combining this with \cref{lem:gradient:limit:zero:abstract} demonstrates that $\P ( \fg ( \vartheta ) = 0 ) =1$. This and the fact that for all $\theta \in \R^\fd$ we have that $(\bD \loss ) (\theta ) \subseteq \R^\fd$
is closed (cf., e.g., Rockafellar~\cite[Theorem~8.6]{RockafellarWets1998} and Jentzen \& Riekert~\cite[Item~$($v$)$ in Lemma~3.3]{JentzenRiekert2021ExistenceGlobMin}) demonstrate that $\P ( 0 \in ( \bD \loss ) ( \vartheta ) ) = 1$.
\end{cproof}

\subsection{SGD training of all biases of clipping ANNs (3 layers)}
\label{subsec_num_sim_clipping_biases_1_hidden_layer}

In this subsection we present numerical simulation results (see \cref{fig:all_biases_clipping_1_hidden_layer} below) in a supervised learning framework for the training of shallow clipping \ANNs\ with 3 layers (corresponding to $L = 2$ in \cref{setting:dnn}) with a 1-dimensional input layer (corresponding to $\ell_0 = 1$ in \cref{setting:dnn}), an $\ell_1$-dimensional hidden layer, and a 1-dimensional output layer (corresponding to $\ell_L = 1$ in \cref{setting:dnn}) where $\ell_1 \in \{10, 100, 1000\}$ and where, as in the statement of \cref{theorem_intro_2} above, not all \ANN\ parameters but only the bias parameters are modified during the training process. We also refer to \cref{figure_deep_3_layers} below for a graphical illustration of the \ANN\ architectures considered in this subsection.

Assume \cref{setting:dnn}, assume $L = 2$, $\ell_0 = \ell_2 = 1$, $\ell_1 \in \{10, 100, 1000\}$, and $\batchsize = 1024$, assume for all $x \in \R$ that $\act_{\infty} (x) = \allowbreak \min\{\max\{x, \allowbreak 0\}, 1\}$ and $f(x) = \abs{x}^{1/4}$, assume for all $k \in \{1, \allowbreak 2, \allowbreak \ldots, \allowbreak \batchsize\}$, $n \in \N$, $x \in [0, 1]$ that $\P(X_k^n < x) = x$, assume
\begin{equation}\label{eqn:index_all_biases_clipping}
\textstyle \indexset = (\N \cap (\ell_1, 2 \ell_1]) \cup \{3 \ell_1 + 1\},
\end{equation}
assume\footnote{\Nobs that the function $\floor{\cdot} \colon \R \to \R$ satisfies that for all $x \in \R$ it holds that $\floor{x} = \max((- \infty, x] \cap \Z)$.} for all $n \in \N$, $g_1, g_2, \ldots, g_n \in \R^{\dimension}$ that $\psi_n(g_1, \allowbreak g_2, \allowbreak \ldots, \allowbreak g_n) = 10^{-2} \, 2^{- \floor{\frac{n}{150}}} g_n$, assume for all $j \in \allowbreak \{1, \allowbreak 2, \allowbreak \ldots, \allowbreak \ell_1\}$, $x \in \R$ that
\begin{equation}\label{eqn:clipping_distribution}
\textstyle \P\big(\ell_1^{-3} \Theta_{0}^j \le x\big) = \P\big(\ell_1^{7/8} \Theta_{0}^{2 \ell_1 + j} \le x\big) = \int_0^x [\frac{2}{\pi}]^{1/2} \exp(-\frac{y^2}{2}) \, \d y,
\end{equation}
and assume for all $j \in \{1, 2, \ldots, \ell_1\}$ that $\ell_1^{-3} \Theta_0^{\ell_1 + j}$ and $\Theta_0^{3 \ell_1 + 1}$ are standard normal random variables.

\def\layersep{4cm}
\begin{figure}[H]
\centering
\begin{adjustbox}{width=\textwidth}
\begin{tikzpicture}[shorten >=1pt,->,draw=black!50, node distance=\layersep]
\tikzstyle{every pin edge}=[<-,shorten <=1pt]
\tikzstyle{input neuron}=[very thick, circle,draw=red, fill=red!20, minimum size=15pt,inner sep=0pt]
\tikzstyle{output neuron}=[very thick, circle, draw=green,fill=green!20,minimum size=30pt,inner sep=0pt]
\tikzstyle{hidden neuron}=[very thick, circle,draw=blue,fill=blue!20,minimum size=20pt,inner sep=0pt]
\tikzstyle{annot} = [text width=9em, text centered]
\tikzstyle{annot2} = [text width=4em, text centered]

%----------Neuron(s) of input layer----------
\node[input neuron] (I) at (0,-1.6) {$x$};

%----------Neuron(s) of 1st hidden layer----------
%\foreach \name / \y in {1, ..., 8}
\path[yshift = 1.5cm]
%node[hidden neuron] (H0-\name) at (\layersep, -\y cm) {};
node[hidden neuron](H0-1) at (\layersep, -1 cm) {};
\path[yshift = 1.5cm]
node[hidden neuron](H0-2) at (\layersep, -2 cm) {};
\path[yshift = 1.5cm]
node(H0-dots) at (\layersep, -2.9 cm) {\vdots};
\path[yshift = 1.5cm]
node[hidden neuron](H0-999) at (\layersep, -4 cm) {};
\path[yshift = 1.5cm]
node[hidden neuron](H0-1000) at (\layersep, -5 cm) {};

%----------Neuron(s) of output layer----------
\path[yshift = 1.5cm]
node[output neuron](O) at (2*\layersep,-3.1 cm) {$\cN^{\theta}(x)$};

\foreach \dest in {1,2,999,1000}
\path[line width = 0.8 ] (I) edge (H0-\dest);

\foreach \source in {1,2,999,1000}
\path[line width = 0.8 ] (H0-\source) edge (O);

% Annotate the layers
\node[annot,above of=H0-1, node distance=1cm, align=center] (hl) {Hidden layer\\($2^{\text{nd}}$ layer)};
\node[annot,left of=hl, align=center] {Input layer\\ ($1^{\text{st}}$ layer)};
\node[annot,right of=hl, align=center] {Output layer\\($3^{\text{rd}}$ layer)};

\node[annot2,below of=H0-1000, node distance=1cm, align=center] (sl) {$\ell_1$ \\ neurons};
\node[annot2,left of=sl, align=center] {$\ell_0=1$ neurons};
\node[annot2,right of=sl, align=center] {$\ell_2=1$ neurons};
\end{tikzpicture}
\end{adjustbox}
\caption{Graphical illustration of the considered shallow \ANN\ architectures used in Subsections~\ref{subsec_num_sim_clipping_biases_1_hidden_layer}, \ref{subsec_num_sim_ReLU_inner_bias_1_hidden_layer}, \ref{subsec_num_sim_ReLU_all_params_1_hidden_layer_random_normal}, and \ref{subsec_num_sim_ReLU_all_params_1_hidden_layer_He}.}
\label{figure_deep_3_layers}
\end{figure}
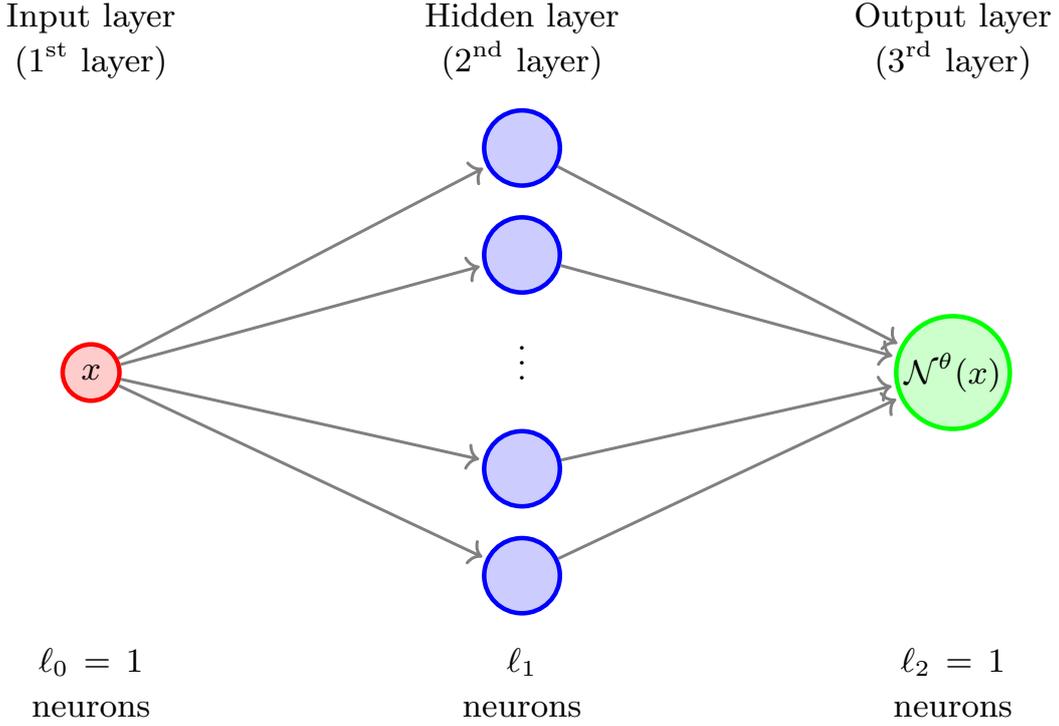

Let us add a few comments regarding the above presented setup. \Nobs that the set $\indexset$ in \cref{eqn:index_all_biases_clipping} represents the index set labeling the bias parameters of the considered shallow \ANNs. Within this subsection only those bias parameters of the \ANNs\ but not the weight parameters of the \ANNs\ are modified during the training procedure. Furthermore, \nobs that the distribution of $\Theta_0$ is nothing else but a special case of the initialization distributions appearing in \cref{theorem_intro_2} above. In addition, \nobs that the above described setup (see \cref{eqn:index_all_biases_clipping,eqn:clipping_distribution} above) ensures that \cref{theorem_intro_2} is applicable with $\alpha = \nicefrac{7}{8}$ and $\beta = 3$. Moreover, \nobs that the assumption that for all $n \in \N$, $g_1, g_2, \ldots, g_n \in \R^{\dimension}$ it holds that $\psi_n(g_1, g_2, \ldots, g_n) = 10^{-2} \, 2^{- \floor{\frac{n}{150}}} g_n$ assures that the stochastic process $\Theta = (\Theta^1, \ldots, \Theta^{\dimension}) \colon \N_0 \times \Omega \to \R^{\dimension}$ in \cref{eqn:setting:dnn:stochastic_pr} satisfies for all $n \in \N$ that
\begin{equation}
\textstyle \Theta_n = \Theta_{n - 1} - 10^{-2} \, 2^{- \floor{\frac{n}{150}}} \fG_n(\Theta_{n - 1})
\end{equation}
(cf.\ Subsection~\ref{subsubsec_Standard_SGD}). Next \nobs that the assumption that $L = 2$, the assumption that $\ell_0 = \ell_2 = 1$, and the assumption that $\ell_1 \in \{10, \allowbreak 100, \allowbreak 1000\}$ demonstrate that the number $\dimension \in \N$ of the employed real parameters to describe the considered shallow \ANNs\ (see \cref{setting:dnn}) satisfies
\begin{equation}
\textstyle \dimension = \ell_1(\ell_0 + 1) + \ell_2 (\ell_1 + 1) = 3 \ell_1 + 1 = 
\begin{cases}
31 & \colon \ell_1 = 10 \\
301 & \colon \ell_1 = 100 \\
3001 & \colon \ell_1 = 1000.
\end{cases}
\end{equation}
In addition, \nobs that \cref{eqn:setting:dnn:loss}, the assumption that $L = 2$, and the assumption that $\batchsize = 1024$ assure that for all $n \in \N$, $\theta \in \R^{\dimension}$, $\omega \in \Omega$ it holds that $\lossapp_{\infty}^{n} (\theta, \omega) = \frac{1}{1024} \sum_{k = 1}^{1024} \abs{\cN_{\infty}^{2, \theta}(X_k^n (\omega)) - f(X_k^n (\omega))}^2$.

In the left image of Figure~\ref{fig:all_biases_clipping_1_hidden_layer} we approximately plot 300 samples of the random variable $\lossapp_{\infty}^{10001}(\Theta_{10000})$ in a histogram with 40 distinct subintervals (with 40 bins) on the $x$-axis. In the upper-right image of Figure~\ref{fig:all_biases_clipping_1_hidden_layer} we approximately plot the $L^2 (\P; \R^{\dimension})$-norm $(\E[\norm{\fG_{n + 1}(\Theta_{n})}]^2)^{1/2}$ of the gradient against the number $n \in (\cup_{k = 0}^{498} \{20 k + 1\})$ of \SGD\ training steps. In the lower-right image of Figure~\ref{fig:all_biases_clipping_1_hidden_layer} we approximately plot the $L^2 (\P; \R^{\dimension})$-distance $(\E[\norm{\Theta_{n} - \Theta_{9981}}]^2)^{1/2}$ between $\Theta_n$ and $\Theta_{9981}$ against the number $n \in (\cup_{k = 0}^{498} \{20 k + 1\})$ of \SGD\ training steps. In, both, the upper-right image of Figure~\ref{fig:all_biases_clipping_1_hidden_layer} and the lower-right image of Figure~\ref{fig:all_biases_clipping_1_hidden_layer} we approximate the expectations in the $L^2 (\P; \R^{\dimension})$-norms by means of 300 independent Monte Carlo samples. The source code used to create Figure~\ref{fig:all_biases_clipping_1_hidden_layer} can be found at \url{https://github.com/deeplearningmethods/overcome-bad-local-minima} (cf.\ \cref{remark:python_codes} above).

\begin{figure}[H]
\centering
\subfloat{\includegraphics[width=16cm]{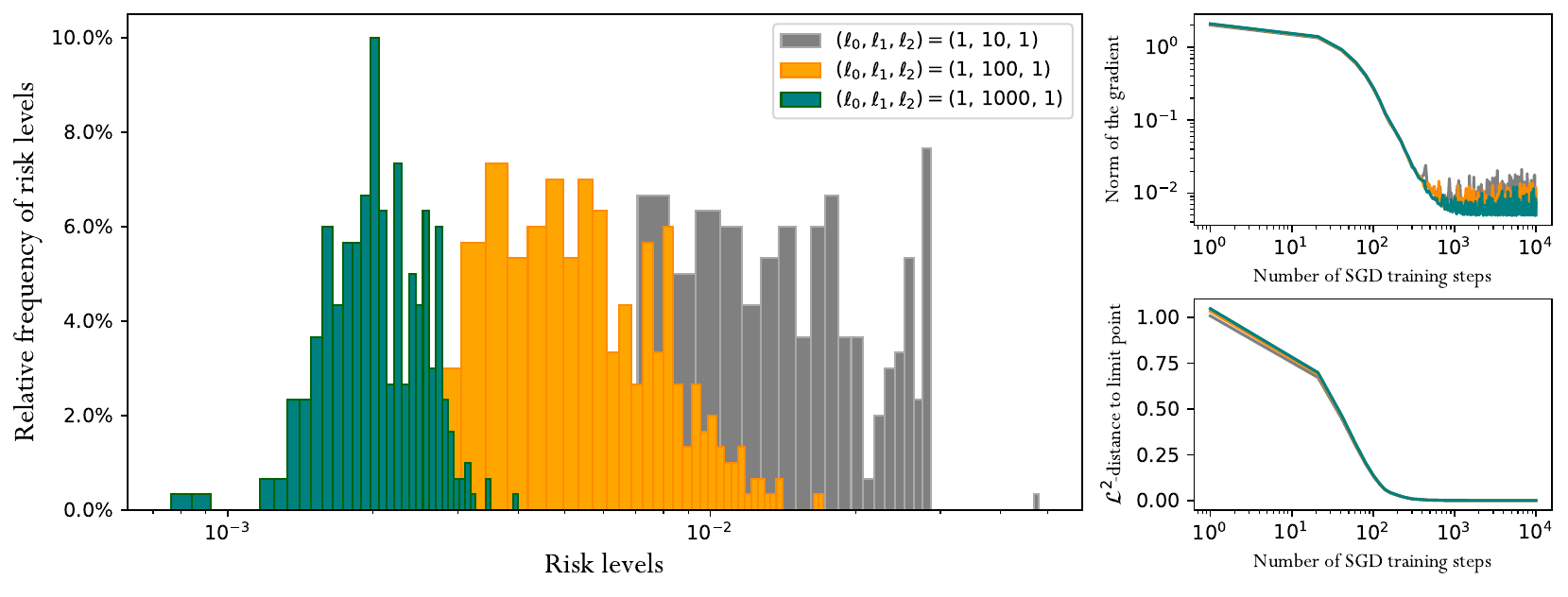}}
\caption{Result of the numerical simulations described in  Subsection~\ref{subsec_num_sim_clipping_biases_1_hidden_layer}.}
\label{fig:all_biases_clipping_1_hidden_layer}
\end{figure}

\subsection{SGD training of inner biases of ReLU ANNs (3 layers)}
\label{subsec_num_sim_ReLU_inner_bias_1_hidden_layer}

In this subsection we present numerical simulation results (see \cref{fig:inner_biases_ReLU_1_hidden_layer} below) in a supervised learning framework for the training of shallow \ReLU\ \ANNs\ with 3 layers (corresponding to $L = 2$ in \cref{setting:dnn}) with a 1-dimensional input layer (corresponding to $\ell_0 = 1$ in \cref{setting:dnn}), an $\ell_1$-dimensional hidden layer, and a 1-dimensional output layer (corresponding to $\ell_L = 1$ in \cref{setting:dnn}) where $\ell_1 \in \{5, 20, 100\}$ and where, as in the statement of \cref{theorem_intro} above, not all \ANN\ parameters but only the inner bias parameters are modified during the training process. We also refer to \cref{figure_deep_3_layers} above for a graphical illustration of the \ANN\ architectures considered in this subsection.

Assume \cref{setting:dnn}, assume $L = 2$, $\ell_0 = \ell_2 = 1$, $\ell_1 \in \{5, 20, 100\}$, and $\batchsize = 1024$, assume for all $x \in \R$ that $\act_{\infty} (x) = \allowbreak \max\{x, \allowbreak 0\}$ and $f(x) = x^4 \sin(x)$, assume for all $k \in \{1, \allowbreak 2, \allowbreak \ldots, \allowbreak \batchsize\}$, $n \in \N$, $x \in [0, 1]$ that $\P(X_k^n < x) = x$, assume
\begin{equation}\label{eqn:index_inner_biases_ReLU}
\textstyle \indexset = \N \cap (\ell_1, 2 \ell_1],
\end{equation}
assume for all $n \in \N$, $g_1, g_2, \ldots, g_n \in \R^{\dimension}$ that $\psi_n(g_1, \allowbreak g_2, \allowbreak \ldots, \allowbreak g_n) = 10^{-1} \, 2^{- \floor{\frac{n}{500}}} g_n$, assume for all $j \in \allowbreak \{1, \allowbreak 2, \allowbreak \ldots, \allowbreak \dimension - 1\} \backslash \indexset$, $x \in \R$ that
\begin{equation}\label{eqn:ReLU_inner_bias_distribution}
\textstyle \P\big(\ell_1^{4/15} \Theta_{0}^j \le x\big) = \int_0^x [\frac{2}{\pi}]^{1/2} \exp(-\frac{y^2}{2}) \, \d y,
\end{equation}
and assume for all $j \in \indexset \cup \{\dimension\}$ that $\ell_1^{9/10} \Theta_{0}^j$ is a standard normal random variable.

Let us add a few comments regarding the above presented setup. \Nobs that the set $\indexset$ in \cref{eqn:index_inner_biases_ReLU} represents the index set labeling the inner bias parameters of the considered shallow \ANNs. Within this subsection only those inner bias parameters of the \ANNs\ but not the outer bias parameter and the weight parameters of the \ANNs\ are modified during the training procedure. Furthermore, \nobs that the distribution of $\Theta_0$ is nothing else but a special case of the initialization distributions appearing in \cref{theorem_intro} above. In addition, \nobs that the above described setup (see \cref{eqn:index_inner_biases_ReLU,eqn:ReLU_inner_bias_distribution} above) ensures that \cref{theorem_intro} is applicable with $\alpha = \nicefrac{9}{10}$ and $\beta = \nicefrac{4}{15}$. Moreover, \nobs that the assumption that for all $n \in \N$, $g_1, \allowbreak g_2, \allowbreak \ldots, \allowbreak g_n \allowbreak \in \R^{\dimension}$ it holds that $\psi_n(g_1, g_2, \ldots, g_n) = 10^{-1} \, 2^{- \floor{\frac{n}{500}}} g_n$ assures that the stochastic process $\Theta = (\Theta^1, \ldots, \Theta^{\dimension}) \colon \N_0 \times \Omega \to \R^{\dimension}$ in \cref{eqn:setting:dnn:stochastic_pr} satisfies for all $n \in \N$ that
\begin{equation}
\textstyle \Theta_n = \Theta_{n - 1} - 10^{-1} \, 2^{- \floor{\frac{n}{500}}} \fG_n(\Theta_{n - 1})
\end{equation}
(cf.\ Subsection~\ref{subsubsec_Standard_SGD}). Next \nobs that the assumption that $L = 2$, the assumption that $\ell_0 = \ell_2 = 1$, and the assumption that $\ell_1 \in \{5, 20, 100\}$ demonstrate that the number $\dimension \in \N$ of the employed real parameters to describe the considered shallow \ANNs\ (see \cref{setting:dnn}) satisfies
\begin{equation}
\textstyle \dimension = \ell_1(\ell_0 + 1) + \ell_2 (\ell_1 + 1) = 3 \ell_1 + 1 = 
\begin{cases}
16 & \colon \ell_1 = 5 \\
61 & \colon \ell_1 = 20 \\
301 & \colon \ell_1 = 100.
\end{cases}
\end{equation}
In addition, \nobs that \cref{eqn:setting:dnn:loss}, the assumption that $L = 2$, and the assumption that $\batchsize = 1024$ assure that for all $n \in \N$, $\theta \in \R^{\dimension}$, $\omega \in \Omega$ it holds that $\lossapp_{\infty}^{n} (\theta, \omega) = \frac{1}{1024} \sum_{k = 1}^{1024} \abs{\cN_{\infty}^{2, \theta}(X_k^n (\omega)) - f(X_k^n (\omega))}^2$.

In the left image of Figure~\ref{fig:inner_biases_ReLU_1_hidden_layer} we approximately plot 300 samples of the random variable $\lossapp_{\infty}^{10001}(\Theta_{10000})$ in a histogram with 200 distinct subintervals (with 200 bins) on the $x$-axis. In the upper-right image of Figure~\ref{fig:inner_biases_ReLU_1_hidden_layer} we approximately plot the $L^2 (\P; \R^{\dimension})$-norm $(\E[\norm{\fG_{n + 1}(\Theta_{n})}]^2)^{1/2}$ of the gradient against the number $n \in (\cup_{k = 0}^{498} \{20 k + 1\})$ of \SGD\ training steps. In the lower-right image of Figure~\ref{fig:inner_biases_ReLU_1_hidden_layer} we approximately plot the $L^2 (\P; \R^{\dimension})$-distance $(\E[\norm{\Theta_{n} - \Theta_{9981}}]^2)^{1/2}$ between $\Theta_n$ and $\Theta_{9981}$ against the number $n \in (\cup_{k = 0}^{498} \{20 k + 1\})$ of \SGD\ training steps. In, both, the upper-right image of Figure~\ref{fig:inner_biases_ReLU_1_hidden_layer} and the lower-right image of Figure~\ref{fig:inner_biases_ReLU_1_hidden_layer} we approximate the expectations in the $L^2 (\P; \R^{\dimension})$-norms by means of 300 independent Monte Carlo samples. The source code used to create Figure~\ref{fig:inner_biases_ReLU_1_hidden_layer} can be found at \url{https://github.com/deeplearningmethods/overcome-bad-local-minima} (cf.\ \cref{remark:python_codes} above).

\begin{figure}[H]
\centering
\subfloat{\includegraphics[width=16cm]{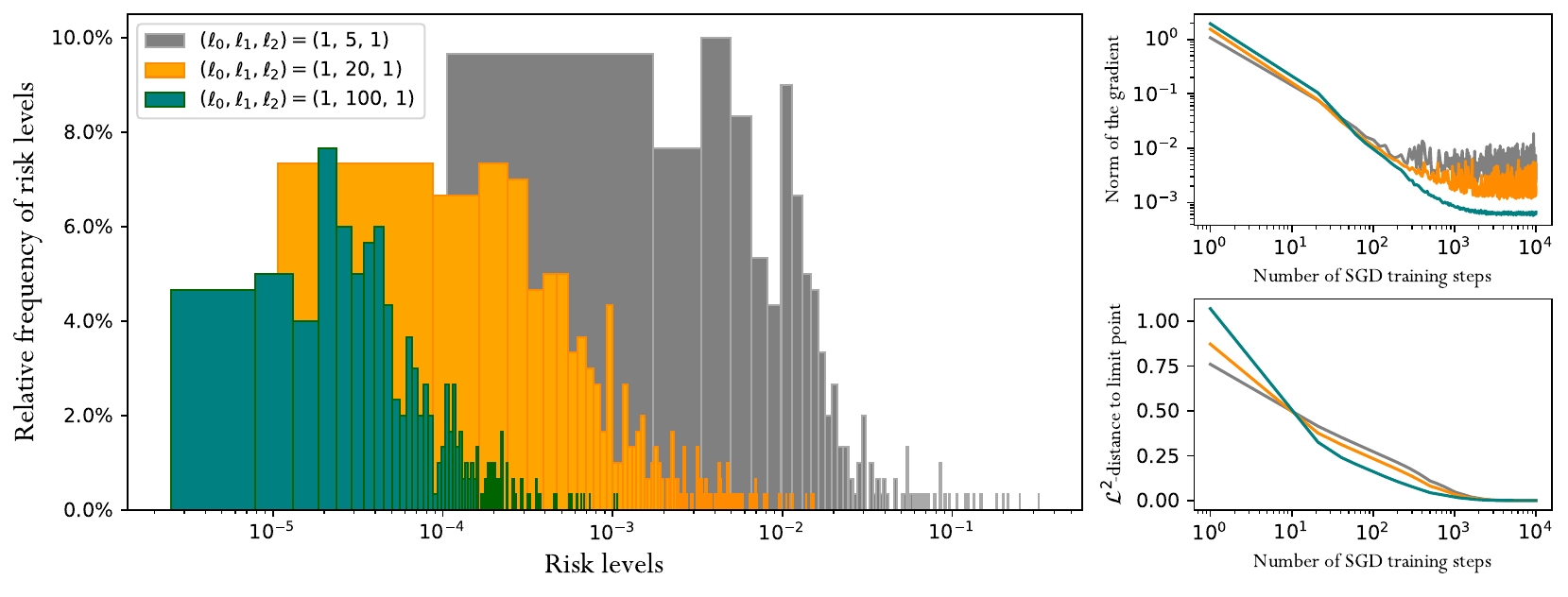}}
\caption{Result of the numerical simulations described in Subsection~\ref{subsec_num_sim_ReLU_inner_bias_1_hidden_layer}.}
\label{fig:inner_biases_ReLU_1_hidden_layer}
\end{figure}

\subsection{SGD training of all parameters of ReLU ANNs (3 layers)}
\label{subsec_num_sim_ReLU_all_params_1_hidden_layer_random_normal}

In this subsection we present numerical simulation results (see \cref{fig:all_parameters_ReLU_1_hidden_layer_random_normal} below) in a supervised learning framework for the training of shallow \ReLU\ \ANNs\ with 3 layers (corresponding to $L = 2$ in \cref{setting:dnn}) with a 1-dimensional input layer (corresponding to $\ell_0 = 1$ in \cref{setting:dnn}), an $\ell_1$-dimensional hidden layer, and a 1-dimensional output layer (corresponding to $\ell_L = 1$ in \cref{setting:dnn}) where $\ell_1 \in \{10, 100, 1000\}$ and where all \ANN\ parameters are modified during the training process. We also refer to \cref{figure_deep_3_layers} above for a graphical illustration of the \ANN\ architectures considered in this subsection.

Assume \cref{setting:dnn}, assume $L = 2$, $\ell_0 = \ell_2 = 1$, $\ell_1 \in \{10, 100, 1000\}$, and $\batchsize = 1024$, assume for all $x \in \R$ that $\act_{\infty} (x) = \allowbreak \max\{x, \allowbreak 0\}$ and $f(x) = \sin(3x)$, assume for all $k \in \{1, \allowbreak 2, \allowbreak \ldots, \allowbreak \batchsize\}$, $n \in \N$, $x \in [0, 1]$ that $\P(X_k^n < x) = x$, assume
\begin{equation}\label{eqn:index_all_params_ReLU_random_normal}
\textstyle \indexset = \N \cap [1, \dimension],
\end{equation}
assume for all $n \in \N$, $g_1, g_2, \ldots, g_n \in \R^{\dimension}$ that $\psi_n(g_1, \allowbreak g_2, \allowbreak \ldots, \allowbreak g_n) = 10^{-3} n^{-\nicefrac{1}{10}} g_n$, and assume for all $j \in \indexset$, $x \in \R$ that 
\begin{equation}
\textstyle \P\big(\Theta_{0}^j \le x\big) = \int_{-\infty}^x [\frac{1}{2 \pi}]^{1/2} \exp(-\frac{y^2}{2}) \, \d y.
\end{equation}
Let us add a few comments regarding the above presented setup. \Nobs that the set $\indexset$ in \cref{eqn:index_all_params_ReLU_random_normal} represents the index set labeling all parameters of the considered shallow \ANNs. Within this subsection all of the parameters of the \ANNs\ are modified during the training procedure. Furthermore, \nobs that the distribution of $\Theta_0$ is nothing else but a standard normal random variable. Moreover, \nobs that the assumption that for all $n \in \N$, $g_1, g_2, \ldots, g_n \in \R^{\dimension}$ it holds that $\psi_n(g_1, g_2, \ldots, g_n) = 10^{-3} n^{-\nicefrac{1}{10}} g_n$ ensures that the stochastic process $\Theta = (\Theta^1, \ldots, \Theta^{\dimension}) \colon \N_0 \times \Omega \to \R^{\dimension}$ in \cref{eqn:setting:dnn:stochastic_pr} satisfies for all $n \in \N$ that
\begin{equation}
\textstyle \Theta_n = \Theta_{n - 1} - 10^{-3} n^{-\nicefrac{1}{10}} \fG_n(\Theta_{n - 1})
\end{equation}
(cf.\ Subsection~\ref{subsubsec_Standard_SGD}). Next \nobs that the assumption that $L = 2$, the assumption that $\ell_0 = \ell_2 = 1$, and the assumption that $\ell_1 \in \{10, \allowbreak 100, \allowbreak 1000\}$ demonstrate that the number $\dimension \in \N$ of the employed real parameters to describe the considered shallow \ANNs\ (see \cref{setting:dnn}) satisfies
\begin{equation}
\textstyle \dimension = \ell_1(\ell_0 + 1) + \ell_2 (\ell_1 + 1) = 3 \ell_1 + 1 = 
\begin{cases}
31 & \colon \ell_1 = 10 \\
301 & \colon \ell_1 = 100 \\
3001 & \colon \ell_1 = 1000.
\end{cases}
\end{equation}
In addition, \nobs that \cref{eqn:setting:dnn:loss}, the assumption that $L = 2$, and the assumption that $\batchsize = 1024$ assure that for all $n \in \N$, $\theta \in \R^{\dimension}$, $\omega \in \Omega$ it holds that $\lossapp_{\infty}^{n} (\theta, \omega) = \frac{1}{1024} \sum_{k = 1}^{1024} \abs{\cN_{\infty}^{2, \theta}(X_k^n (\omega)) - f(X_k^n (\omega))}^2$.

In the left image of Figure~\ref{fig:all_parameters_ReLU_1_hidden_layer_random_normal} we approximately plot 300 samples of the random variable $\lossapp_{\infty}^{10001}(\Theta_{10000})$ in a histogram with 40 distinct subintervals (with 40 bins) on the $x$-axis. In the upper-right image of Figure~\ref{fig:all_parameters_ReLU_1_hidden_layer_random_normal} we approximately plot the $L^2 (\P; \R^{\dimension})$-norm $(\E[\norm{\fG_{n + 1}(\Theta_{n})}]^2)^{1/2}$ of the generalized gradient against the number $n \in (\cup_{k = 0}^{498} \{20 k + 1\})$ of \SGD\ training steps. In the lower-right image of Figure~\ref{fig:all_parameters_ReLU_1_hidden_layer_random_normal} we approximately plot the $L^2 (\P; \R^{\dimension})$-distance $(\E[\norm{\Theta_{n} - \Theta_{9981}}]^2)^{1/2}$ between $\Theta_n$ and $\Theta_{9981}$ against the number $n \in (\cup_{k = 0}^{498} \{20 k + 1\})$ of \SGD\ training steps. In, both, the upper-right image of Figure~\ref{fig:all_parameters_ReLU_1_hidden_layer_random_normal} and the lower-right image of Figure~\ref{fig:all_parameters_ReLU_1_hidden_layer_random_normal} we approximate the expectations in the $L^2 (\P; \R^{\dimension})$-norms by means of 300 independent Monte Carlo samples. The source code used to create Figure~\ref{fig:all_parameters_ReLU_1_hidden_layer_random_normal} can be found at \url{https://github.com/deeplearningmethods/overcome-bad-local-minima} (cf.\ \cref{remark:python_codes} above).

\begin{figure}[H]
\centering
\subfloat{\includegraphics[width=16cm]{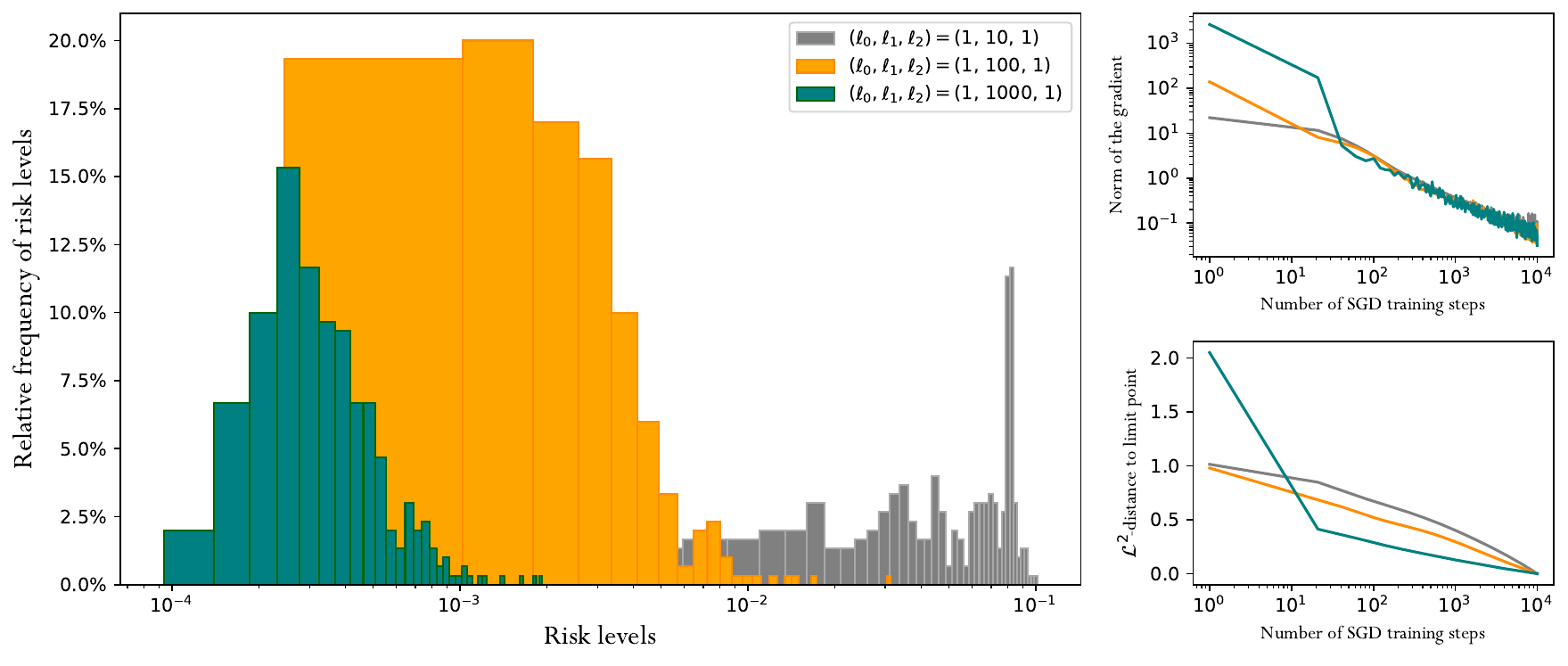}}
\caption{Result of the numerical simulations described in  Subsection~\ref{subsec_num_sim_ReLU_all_params_1_hidden_layer_random_normal}.} 
\label{fig:all_parameters_ReLU_1_hidden_layer_random_normal}
\end{figure}

\subsection{SGD training of all parameters of He initialized ReLU ANNs (3 layers)}
\label{subsec_num_sim_ReLU_all_params_1_hidden_layer_He}

In this subsection we present numerical simulation results (see \cref{fig:all_parameters_ReLU_1_hidden_layer_He} below) in a supervised learning framework for the training of shallow \ReLU\ \ANNs\ with 3 layers (corresponding to $L = 2$ in \cref{setting:dnn}) with a 1-dimensional input layer (corresponding to $\ell_0 = 1$ in \cref{setting:dnn}), a 32-dimensional hidden layer (corresponding to $\ell_1 = 32$ in \cref{setting:dnn}), and a 1-dimensional output layer (corresponding to $\ell_L = 1$ in \cref{setting:dnn}) where all \ANN\ parameters are modified during the training process. We also refer to \cref{figure_deep_3_layers} above for a graphical illustration of the \ANN\ architecture considered in this subsection.

Assume \cref{setting:dnn}, assume $L = 2$, $\ell_0 = \ell_2 = 1$, $\ell_1 = 32$, and $\batchsize = 1024$, assume for all $x \in \R$ that $\act_{\infty} (x) = \allowbreak \max\{x, \allowbreak 0\}$ and $f(x) = x^2 + 2x$, assume for all $k \in \{1, \allowbreak 2, \allowbreak \ldots, \allowbreak \batchsize\}$, $n \in \N$, $x \in [0, 1]$ that $\P(X_k^n < x) = x$, assume
\begin{equation}\label{eqn:index_all_params_ReLU_He}
\textstyle \indexset = \N \cap [1, \dimension],
\end{equation}
assume for all $k \in \{1, 2, \ldots, L\}$, $i \in \{1, 2, \ldots, \ell_k\}$, $j \in \{1, 2, \ldots, \ell_{k - 1}\}$, $x \in \R$ that
\begin{equation}
\textstyle \P\big(\fb_i^{k, \Theta_0} = 0 \big) = 1 \qqandqq \P\big(\fw_{i, j}^{k, \Theta_0} \ge x \big) = \int_x^{\infty} \frac{1}{2} \bigl[\frac{\ell_{k - 1}}{\pi}\bigr]^{1/2} \exp\bigl(- \frac{y^2 \ell_{k - 1}}{4}\bigr) \, \d y,
\end{equation}
and assume for all $n \in \N$, $g_1, g_2, \ldots, g_n \in \R^{\dimension}$ that $\psi_n(g_1, \allowbreak g_2, \allowbreak \ldots, \allowbreak g_n) = 10^{-2} g_n$.

Let us add a few comments regarding the above presented setup. \Nobs that the set $\indexset$ in \cref{eqn:index_all_params_ReLU_He} represents the index set labeling all parameters of the considered shallow \ANNs. Within this subsection all of the parameters of the \ANNs\ are modified during the training procedure. Furthermore, \nobs that the distribution of $\Theta_0$ is nothing else but the He normal initialization (cf., e.g., \cite{tensorflow}) for \ReLU\ \ANNs\ with 3 layers with a 1-dimensional input layer, a 32-dimensional hidden layer, and a 1-dimensional output layer. Moreover, \nobs that the assumption that for all $n \in \N$, $g_1, g_2, \ldots, g_n \in \R^{\dimension}$ it holds that $\psi_n(g_1, g_2, \ldots, g_n) = 10^{-2} g_n$ ensures that the stochastic process $\Theta = (\Theta^1, \ldots, \Theta^{\dimension}) \colon \N_0 \times \Omega \to \R^{\dimension}$ in \cref{eqn:setting:dnn:stochastic_pr} satisfies for all $n \in \N$ that
\begin{equation}
\textstyle \Theta_n = \Theta_{n - 1} - 10^{-2} \fG_n(\Theta_{n - 1})
\end{equation}
(cf.\ Subsection~\ref{subsubsec_Standard_SGD}). Next \nobs that the assumption that $L = 2$, the assumption that $\ell_0 = \ell_2 = 1$, and the assumption that $\ell_1 = 32$ demonstrate that the number $\dimension \in \N$ of the employed real parameters to describe the considered shallow \ANNs\ (see \cref{setting:dnn}) satisfies
\begin{equation}
\textstyle \dimension = \ell_1(\ell_0 + 1) + \ell_2 (\ell_1 + 1) = 3 \ell_1 + 1 = 97.
\end{equation}
In addition, \nobs that \cref{eqn:setting:dnn:loss}, the assumption that $L = 2$, and the assumption that $\batchsize = 1024$ assure that for all $n \in \N$, $\theta \in \R^{\dimension}$, $\omega \in \Omega$ it holds that $\lossapp_{\infty}^{n} (\theta, \omega) = \frac{1}{1024} \sum_{k = 1}^{1024} \abs{\cN_{\infty}^{2, \theta}(X_k^n (\omega)) - f(X_k^n (\omega))}^2$.

In the left image of Figure~\ref{fig:all_parameters_ReLU_1_hidden_layer_He} we approximately plot 300 samples of the random variable $\lossapp_{\infty}^{10001}(\Theta_{10000})$ in a histogram with 80 distinct subintervals (with 80 bins) on the $x$-axis. In the upper-right image of Figure~\ref{fig:all_parameters_ReLU_1_hidden_layer_He} we approximately plot the $L^2 (\P; \R^{\dimension})$-norm $(\E[\norm{\fG_{n + 1}(\Theta_{n})}]^2)^{1/2}$ of the generalized gradient against the number $n \in (\cup_{k = 0}^{998} \{10 k + 1\})$ of \SGD\ training steps. In the lower-right image of Figure~\ref{fig:all_parameters_ReLU_1_hidden_layer_He} we approximately plot the $L^2 (\P; \R^{\dimension})$-distance $(\E[\norm{\Theta_{n} - \Theta_{9991}}]^2)^{1/2}$ between $\Theta_n$ and $\Theta_{9991}$ against the number $n \in (\cup_{k = 0}^{998} \{10 k + 1\})$ of \SGD\ training steps. In, both, the upper-right image of Figure~\ref{fig:all_parameters_ReLU_1_hidden_layer_He} and the lower-right image of Figure~\ref{fig:all_parameters_ReLU_1_hidden_layer_He} we approximate the expectations in the $L^2 (\P; \R^{\dimension})$-norms by means of 300 independent Monte Carlo samples. The source code used to create Figure~\ref{fig:all_parameters_ReLU_1_hidden_layer_He} can be found at \url{https://github.com/deeplearningmethods/overcome-bad-local-minima} (cf.\ \cref{remark:python_codes} above).

\begin{figure}[H]
\centering
\subfloat{\includegraphics[width=16cm]{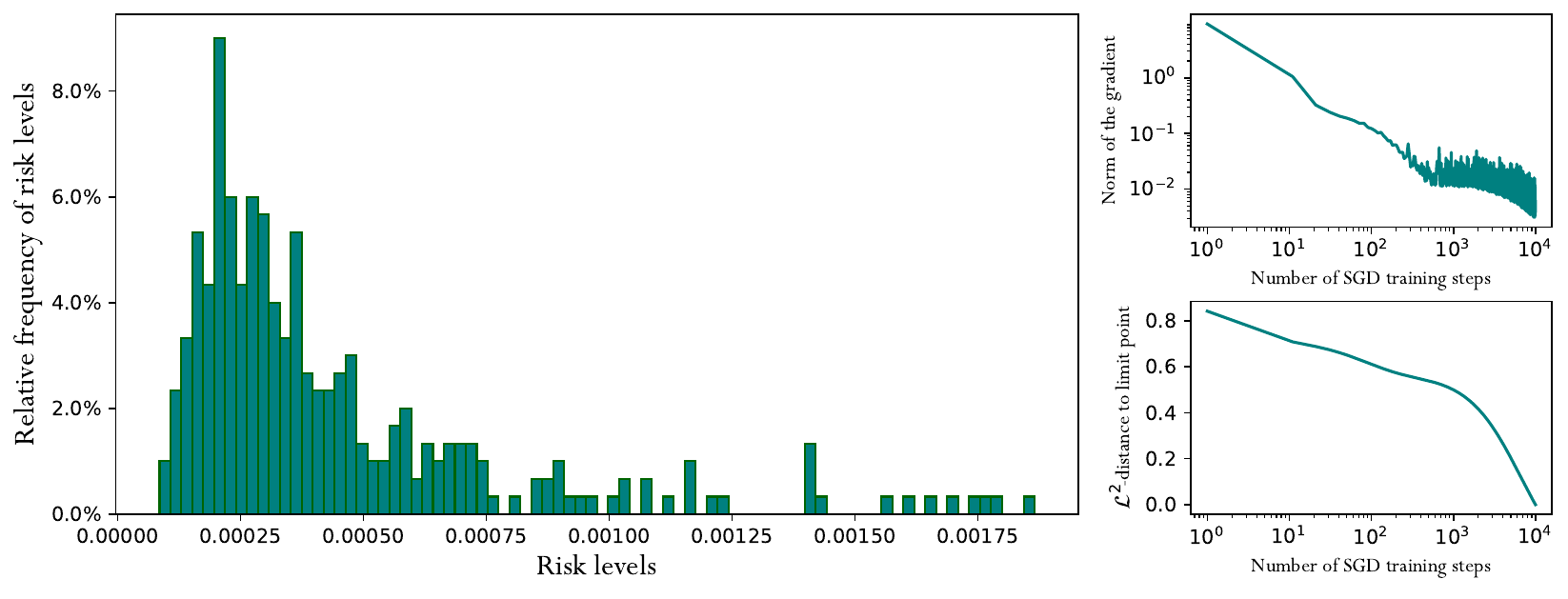}}
\caption{Result of the numerical simulations described in  Subsection~\ref{subsec_num_sim_ReLU_all_params_1_hidden_layer_He}.} 
\label{fig:all_parameters_ReLU_1_hidden_layer_He}
\end{figure}

\subsection{SGD training of all parameters of ReLU ANNs (4 layers)}
\label{subsec_num_sim_ReLU_all_params_2_hidden_layers_random_normal}

In this subsection we present numerical simulation results (see \cref{fig:all_parameters_ReLU_2_hidden_layers_random_normal} below) in a supervised learning framework for the training of deep \ReLU\ \ANNs\ with 4 layers (corresponding to $L = 3$ in \cref{setting:dnn}) with a 1-dimensional input layer (corresponding to $\ell_0 = 1$ in \cref{setting:dnn}), $\ell_1$ and $\ell_2$-dimensional hidden layers, and a 1-dimensional output layer (corresponding to $\ell_L = 1$ in \cref{setting:dnn}) where $\ell_1 = \ell_2 \in \{8, 15, 50\}$ and where all \ANN\ parameters are modified during the training process. We also refer to \cref{figure_deep_4_layers} below for a graphical illustration of the \ANN\ architectures considered in this subsection.

Assume \cref{setting:dnn}, assume $L = 3$, $\ell_0 = \ell_3 = 1$, $\ell_1 = \ell_2 \in \{8, 15, 50\}$, and $\batchsize = 1024$, assume for all $x \in \R$ that $\act_{\infty} (x) = \allowbreak \max\{x, \allowbreak 0\}$ and $f(x) = \sin(3x)$, assume for all $k \in \{1, \allowbreak 2, \allowbreak \ldots, \allowbreak \batchsize\}$, $n \in \N$, $x \in [0, 1]$ that $\P(X_k^n < x) = x$, assume
\begin{equation}\label{eqn:index_all_params_ReLU_random_normal_4_layers}
\textstyle \indexset = \N \cap [1, \dimension],
\end{equation}
assume for all $n \in \N$, $g_1, g_2, \ldots, g_n \in \R^{\dimension}$ that $\psi_n(g_1, \allowbreak g_2, \allowbreak \ldots, \allowbreak g_n) = 0.0005 n^{-\nicefrac{1}{10}} g_n$, and assume for all $j \in \indexset$, $x \in \R$ that 
\begin{equation}
\textstyle \P\big(\Theta_{0}^j \le x\big) = \int_{-\infty}^x [\frac{1}{2 \pi}]^{1/2} \exp(-\frac{y^2}{2}) \, \d y.
\end{equation}

\def\layersep{2.5cm}
\begin{figure}[H]
\centering
\begin{adjustbox}{width=\textwidth}
\begin{tikzpicture}[shorten >=1pt,->,draw=black!50, node distance=\layersep]
\tikzstyle{every pin edge}=[<-,shorten <=1pt]
\tikzstyle{input neuron}=[very thick, circle,draw=red, fill=red!20, minimum size=15pt,inner sep=0pt]
\tikzstyle{output neuron}=[very thick, circle, draw=green,fill=green!20,minimum size=20pt,inner sep=0pt]
\tikzstyle{hidden neuron}=[very thick, circle,draw=blue,fill=blue!20,minimum size=20pt,inner sep=0pt]
\tikzstyle{annot} = [text width=9em, text centered]
\tikzstyle{annot2} = [text width=4em, text centered]

%----------Neuron(s) of input layer----------
\node[input neuron] (I) at (0,-1.5) {$x$};

%----------Neuron(s) of 1st hidden layer----------
%\foreach \name / \y in {1, ..., 8}
\path[yshift = 1.5cm]
%node[hidden neuron] (H0-\name) at (\layersep, -\y cm) {};
node[hidden neuron](H0-1) at (\layersep, -1 cm) {};
\path[yshift = 1.5cm]
node[hidden neuron](H0-2) at (\layersep, -2 cm) {};
\path[yshift = 1.5cm]
node(H0-dots) at (\layersep, -2.9 cm) {\vdots};
\path[yshift = 1.5cm]
node[hidden neuron](H0-49) at (\layersep, -4 cm) {};
\path[yshift = 1.5cm]
node[hidden neuron](H0-50) at (\layersep, -5 cm) {};

%----------Neuron(s) of 2nd hidden layer----------
\path[yshift = 1.5cm]
%node[hidden neuron] (H0-\name) at (\layersep, -\y cm) {};
node[hidden neuron](H1-1) at (2*\layersep, -1 cm) {};
\path[yshift = 1.5cm]
node[hidden neuron](H1-2) at (2*\layersep, -2 cm) {};
\path[yshift = 1.5cm]
node(H1-dots) at (2*\layersep, -2.9 cm) {\vdots};
\path[yshift = 1.5cm]
node[hidden neuron](H1-49) at (2*\layersep, -4 cm) {};
\path[yshift = 1.5cm]
node[hidden neuron](H1-50) at (2*\layersep, -5 cm) {};

%----------Neuron(s) of output layer----------
\path[yshift = 1.5cm]
node[output neuron](O) at (3*\layersep,-3 cm) {$\cN^{\theta}(x)$};

%----------Arrow(s) from 1st to 2nd layer----------
\foreach \dest in {1,2,49,50}
\path[line width = 0.8] (I) edge (H0-\dest);

%----------Arrow(s) from 2nd to 3rd layer----------
\foreach \source in {1,2,49,50}
\foreach \dest in {1,2,49,50}
\path[line width = 0.8] (H0-\source) edge (H1-\dest);

%----------Arrow(s) from 3rd to 4th layer----------
\foreach \source in {1,2,49,50}
\path[line width = 0.8] (H1-\source) edge (O);

% Annotate the layers
\node[annot,above of=H0-1, node distance=1cm, align=center] (hl) {$1^{\text{st}}$ hidden layer\\($2^{\text{nd}}$ layer)};
\node[annot,above of=H1-1, node distance=1cm, align=center] (hl2) {$2^{\text{nd}}$ hidden layer\\($3^{\text{rd}}$ layer)};
\node[annot,left of=hl, align=center] {Input layer\\ ($1^{\text{st}}$ layer)};
\node[annot,right of=hl2, align=center] {Output layer\\($4^{\text{th}}$ layer)};

\node[annot2,below of=H0-50, node distance=1cm, align=center] (sl) {$\ell_1$ \\ neurons};
\node[annot2,below of=H1-50, node distance=1cm, align=center] (sl2) {$\ell_2$ \\ neurons};
\node[annot2,left of=sl, align=center] {$\ell_0=1$ \\ neuron};
\node[annot2,right of=sl2, align=center] {$\ell_3=1$ \\ neuron};
\end{tikzpicture}
\end{adjustbox}
\caption{Graphical illustration of the \ANN\ architectures used in Subsections~\ref{subsec_num_sim_ReLU_all_params_2_hidden_layers_random_normal}, \ref{subsec_num_sim_ReLU_all_params_2_hidden_layers_Xavier}, \ref{subsec_num_sim_ReLU_all_params_2_hidden_layers_He}, \ref{subsec_num_sim_ReLU_all_params_2_hidden_layers_Xavier_Adam}, and \ref{subsec_num_sim_ReLU_all_params_2_hidden_layers_He_Adam}.}
\label{figure_deep_4_layers}
\end{figure}

Let us add a few comments regarding the above presented setup. \Nobs that the set $\indexset$ in \cref{eqn:index_all_params_ReLU_random_normal_4_layers} represents the index set labeling all parameters of the considered deep \ANNs. Within this subsection all of the parameters of the \ANNs\ are modified during the training procedure. Furthermore, \nobs that the distribution of $\Theta_0$ is nothing else but a standard normal random variable. Moreover, \nobs that the assumption that for all $n \in \N$, $g_1, g_2, \ldots, g_n \in \R^{\dimension}$ it holds that $\psi_n(g_1, g_2, \ldots, g_n) = 0.0005 n^{-\nicefrac{1}{10}} g_n$ ensures that the stochastic process $\Theta = (\Theta^1, \ldots, \Theta^{\dimension}) \colon \N_0 \times \Omega \to \R^{\dimension}$ in \cref{eqn:setting:dnn:stochastic_pr} satisfies for all $n \in \N$ that
\begin{equation}
\textstyle \Theta_n = \Theta_{n - 1} - 0.0005 n^{-\nicefrac{1}{10}} \fG_n(\Theta_{n - 1})
\end{equation}
(cf.\ Subsection~\ref{subsubsec_Standard_SGD}). Next \nobs that the assumption that $L = 3$, the assumption that $\ell_0 = \ell_3 = 1$, and the assumption that $\ell_1 = \ell_2 \in \{8, \allowbreak 15, \allowbreak 50\}$ demonstrate that the number $\dimension \in \N$ of the employed real parameters to describe the considered deep \ANNs\ (see \cref{setting:dnn}) satisfies
\begin{equation}
\textstyle \dimension = \ell_1(\ell_0 + 1) + \ell_2 (\ell_1 + 1) + \ell_3 (\ell_2 + 1) =
\begin{cases}
97 & \colon \ell_1 = \ell_2 = 8 \\
286 & \colon \ell_1 = \ell_2 = 15 \\
2701 & \colon \ell_1 =\ell_2 = 50.
\end{cases}
\end{equation}
In addition, \nobs that \cref{eqn:setting:dnn:loss}, the assumption that $L = 3$, and the assumption that $\batchsize = 1024$ assure that for all $n \in \N$, $\theta \in \R^{\dimension}$, $\omega \in \Omega$ it holds that $\lossapp_{\infty}^{n} (\theta, \omega) = \frac{1}{1024} \sum_{k = 1}^{1024} \abs{\cN_{\infty}^{3, \theta}(X_k^n (\omega)) - f(X_k^n (\omega))}^2$.

In the left image of Figure~\ref{fig:all_parameters_ReLU_2_hidden_layers_random_normal} we approximately plot 300 samples of the random variable $\lossapp_{\infty}^{50001}(\Theta_{50000})$ in a histogram with 100 distinct subintervals (with 100 bins) on the $x$-axis. In the upper-right image of Figure~\ref{fig:all_parameters_ReLU_2_hidden_layers_random_normal} we approximately plot the $L^2 (\P; \R^{\dimension})$-norm $(\E[\norm{\fG_{n + 1}(\Theta_{n})}]^2)^{1/2}$ of the generalized gradient against the number $n \in (\cup_{k = 0}^{2498} \{20 k + 1\})$ of \SGD\ training steps. In the lower-right image of Figure~\ref{fig:all_parameters_ReLU_2_hidden_layers_random_normal} we approximately plot the $L^2 (\P; \R^{\dimension})$-distance $(\E[\norm{\Theta_{n} - \Theta_{49981}}]^2)^{1/2}$ between $\Theta_n$ and $\Theta_{49981}$ against the number $n \in (\cup_{k = 0}^{2498} \{20 k + 1\})$ of \SGD\ training steps. In, both, the upper-right image of Figure~\ref{fig:all_parameters_ReLU_2_hidden_layers_random_normal} and the lower-right image of Figure~\ref{fig:all_parameters_ReLU_2_hidden_layers_random_normal} we approximate the expectations in the $L^2 (\P; \R^{\dimension})$-norms by means of 300 independent Monte Carlo samples. The source code used to create Figure~\ref{fig:all_parameters_ReLU_2_hidden_layers_random_normal} can be found at \url{https://github.com/deeplearningmethods/overcome-bad-local-minima} (cf.\ \cref{remark:python_codes} above).

\begin{figure}[H]
\centering
\subfloat{\includegraphics[width=16cm]{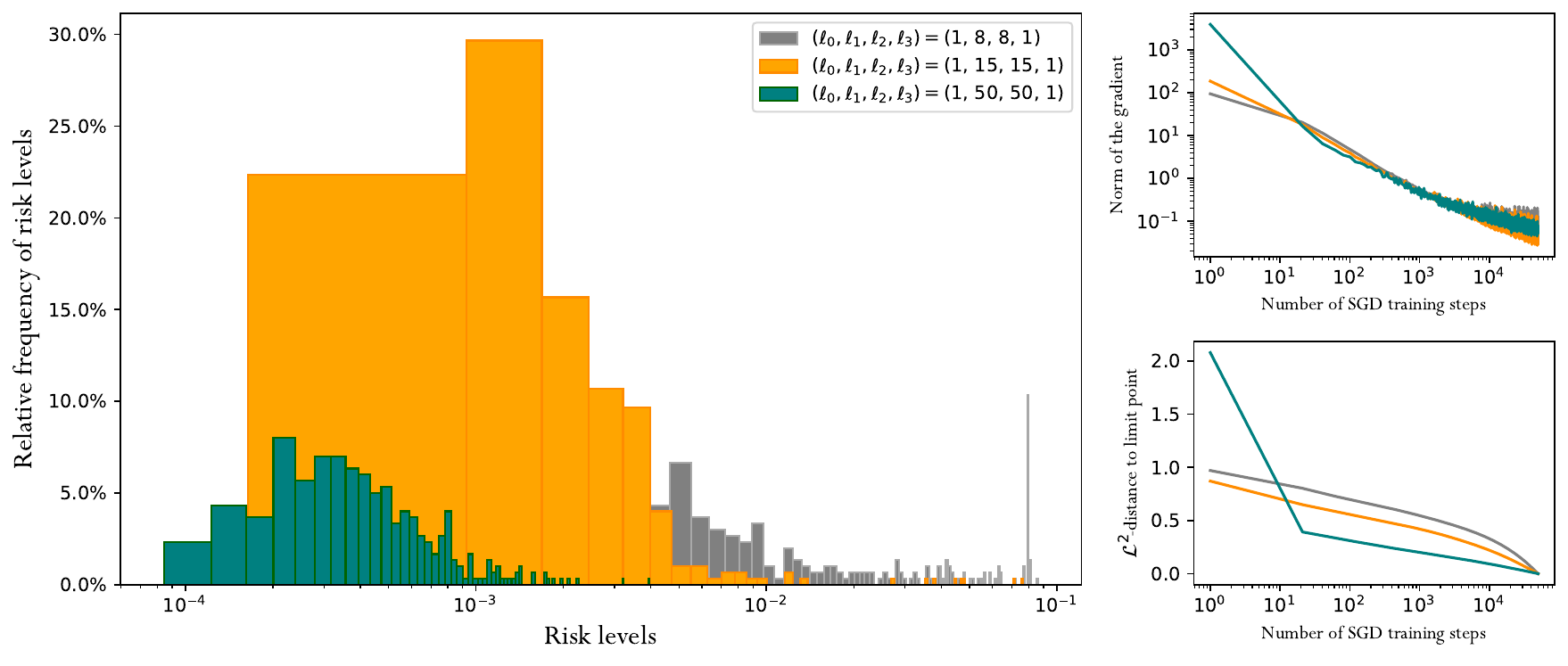}}
\caption{Result of the numerical simulations described in  Subsection~\ref{subsec_num_sim_ReLU_all_params_2_hidden_layers_random_normal}.} 
\label{fig:all_parameters_ReLU_2_hidden_layers_random_normal}
\end{figure}

\subsection{SGD training of all parameters of Xavier initialized ReLU ANNs (4 layers)}
\label{subsec_num_sim_ReLU_all_params_2_hidden_layers_Xavier}

In this subsection we present numerical simulation results (see \cref{fig:all_parameters_ReLU_2_hidden_layers_Xavier} below) in a supervised learning framework for the training of deep \ReLU\ \ANNs\ with 4 layers (corresponding to $L = 3$ in \cref{setting:dnn}) with a 1-dimensional input layer (corresponding to $\ell_0 = 1$ in \cref{setting:dnn}), two 50-dimensional hidden layers (corresponding to $\ell_1 = \ell_2 = 50$ in \cref{setting:dnn}), and a 1-dimensional output layer (corresponding to $\ell_L = 1$ in \cref{setting:dnn}) where all \ANN\ parameters are modified during the training process. We also refer to \cref{figure_deep_4_layers} above for a graphical illustration of the \ANN\ architecture considered in this subsection.

Assume \cref{setting:dnn}, assume $L = 3$, $\ell_0 = \ell_3 = 1$, $\ell_1 = \ell_2 = 50$, and $\batchsize = 1024$, assume for all $x \in \R$ that $\act_{\infty} (x) = \allowbreak \max\{x, \allowbreak 0\}$ and $f(x) = x^2 + 2x$, assume for all $k \in \{1, \allowbreak 2, \allowbreak \ldots, \allowbreak \batchsize\}$, $n \in \N$, $x \in [0, 1]$ that $\P(X_k^n < x) = x$, assume
\begin{equation}\label{eqn:index_all_params_ReLU_Xavier_4_Layers}
\textstyle \indexset = \N \cap [1, \dimension],
\end{equation}
assume for all $k \in \{1, 2, \ldots, L\}$, $i \in \{1, 2, \ldots, \ell_k\}$, $j \in \{1, 2, \ldots, \ell_{k - 1}\}$, $x \in \R$ that
\begin{equation}
\textstyle \P\big(\fb_i^{k, \Theta_0} = 0 \big) = 1 \qqandqq \P\big(\fw_{i, j}^{k, \Theta_0} \ge x \big) = \int_x^{\infty} \frac{1}{2} \bigl[\frac{\ell_{k - 1} + \ell_k}{\pi}\bigr]^{1/2} \exp\bigl(- \frac{y^2 (\ell_{k - 1} + \ell_k)}{4}\bigr) \, \d y,
\end{equation}
and assume for all $n \in \N$, $g_1, g_2, \ldots, g_n \in \R^{\dimension}$ that $\psi_n(g_1, \allowbreak g_2, \allowbreak \ldots, \allowbreak g_n) = 10^{-2} g_n$.

Let us add a few comments regarding the above presented setup. \Nobs that the set $\indexset$ in \cref{eqn:index_all_params_ReLU_Xavier_4_Layers} represents the index set labeling all parameters of the considered deep \ANNs. Within this subsection all of the parameters of the \ANNs\ are modified during the training procedure. Furthermore, \nobs that the distribution of $\Theta_0$ is nothing else but the Xavier normal (Glorot normal) initialization (cf., e.g., \cite{tensorflow}) for \ReLU\ \ANNs\ with 4 layers with a 1-dimensional input layer, two 50-dimensional hidden layers, and a 1-dimensional output layer. Moreover, \nobs that the assumption that for all $n \in \N$, $g_1, g_2, \ldots, g_n \in \R^{\dimension}$ it holds that $\psi_n(g_1, g_2, \ldots, g_n) = 10^{-2} g_n$ ensures that the stochastic process $\Theta = (\Theta^1, \ldots, \Theta^{\dimension}) \colon \N_0 \times \Omega \to \R^{\dimension}$ in \cref{eqn:setting:dnn:stochastic_pr} satisfies for all $n \in \N$ that
\begin{equation}
\textstyle \Theta_n = \Theta_{n - 1} - 10^{-2} \fG_n(\Theta_{n - 1})
\end{equation}
(cf.\ Subsection~\ref{subsubsec_Standard_SGD}). Next \nobs that the assumption that $L = 3$, the assumption that $\ell_0 = \ell_3 = 1$, and the assumption that $\ell_1 = \ell_2 = 50$ demonstrate that the number $\dimension \in \N$ of the employed real parameters to describe the considered deep \ANNs\ (see \cref{setting:dnn}) satisfies
\begin{equation}
\textstyle \dimension = \ell_1(\ell_0 + 1) + \ell_2 (\ell_1 + 1) + \ell_3 (\ell_2 + 1) = 50 (1 + 1) + 50 (50 + 1) + 1 (50 + 1) = 2701.
\end{equation}
In addition, \nobs that \cref{eqn:setting:dnn:loss}, the assumption that $L = 3$, and the assumption that $\batchsize = 1024$ assure that for all $n \in \N$, $\theta \in \R^{\dimension}$, $\omega \in \Omega$ it holds that $\lossapp_{\infty}^{n} (\theta, \omega) = \frac{1}{1024} \sum_{k = 1}^{1024} \abs{\cN_{\infty}^{3, \theta}(X_k^n (\omega)) - f(X_k^n (\omega))}^2$.

In the left image of Figure~\ref{fig:all_parameters_ReLU_2_hidden_layers_Xavier} we approximately plot 300 samples of the random variable $\lossapp_{\infty}^{10001}(\Theta_{10000})$ in a histogram with 80 distinct subintervals (with 80 bins) on the $x$-axis. In the upper-right image of Figure~\ref{fig:all_parameters_ReLU_2_hidden_layers_Xavier} we approximately plot the $L^2 (\P; \R^{\dimension})$-norm $(\E[\norm{\fG_{n + 1}(\Theta_{n})}]^2)^{1/2}$ of the generalized gradient against the number $n \in (\cup_{k = 0}^{998} \{10 k + 1\})$ of \SGD\ training steps. In the lower-right image of Figure~\ref{fig:all_parameters_ReLU_2_hidden_layers_Xavier} we approximately plot the $L^2 (\P; \R^{\dimension})$-distance $(\E[\norm{\Theta_{n} - \Theta_{9991}}]^2)^{1/2}$ between $\Theta_n$ and $\Theta_{9991}$ against the number $n \in (\cup_{k = 0}^{998} \{10 k + 1\})$ of \SGD\ training steps. In, both, the upper-right image of Figure~\ref{fig:all_parameters_ReLU_2_hidden_layers_Xavier} and the lower-right image of Figure~\ref{fig:all_parameters_ReLU_2_hidden_layers_Xavier} we approximate the expectations in the $L^2 (\P; \R^{\dimension})$-norms by means of 300 independent Monte Carlo samples. The source code used to create Figure~\ref{fig:all_parameters_ReLU_2_hidden_layers_Xavier} can be found at \url{https://github.com/deeplearningmethods/overcome-bad-local-minima} (cf.\ \cref{remark:python_codes} above).

\begin{figure}[H]
\centering
\subfloat{\includegraphics[width=16cm]{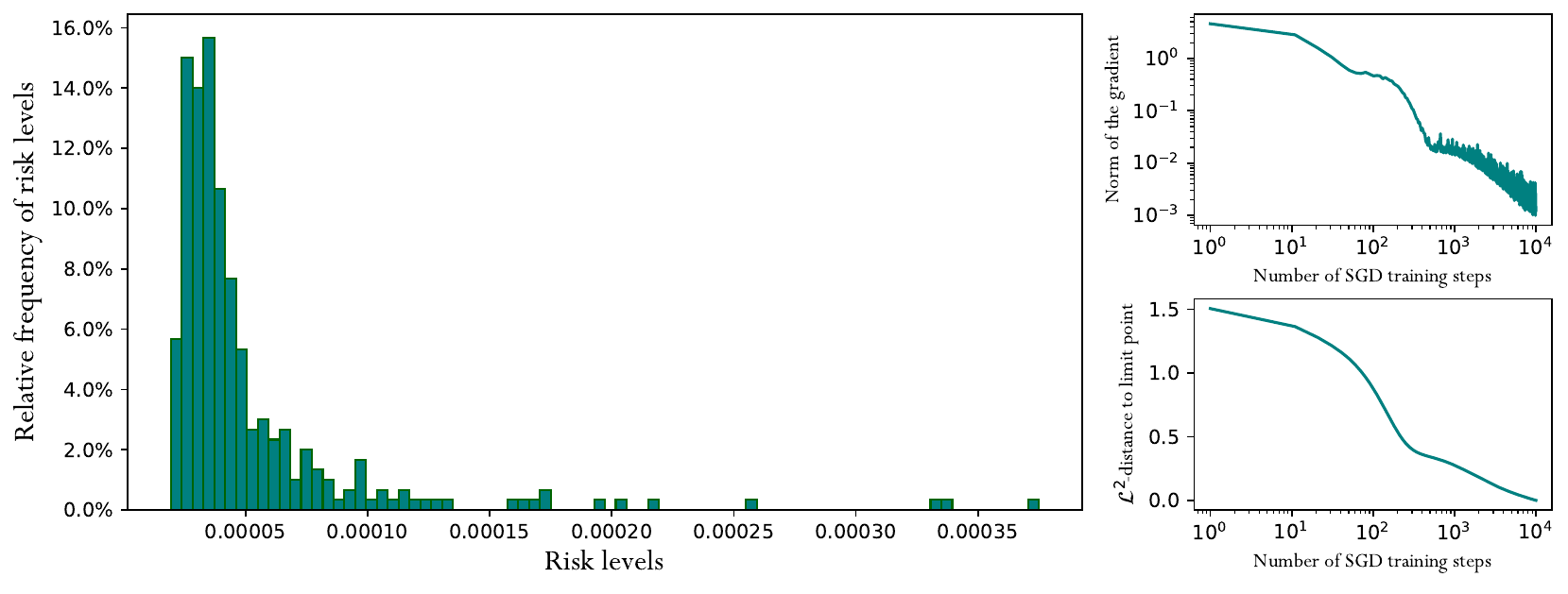}} 
\caption{Result of the numerical simulations described in Subsection~\ref{subsec_num_sim_ReLU_all_params_2_hidden_layers_Xavier}.}
\label{fig:all_parameters_ReLU_2_hidden_layers_Xavier}
\end{figure}

\subsection{SGD training of all parameters of He initialized ReLU ANNs (4 layers)}
\label{subsec_num_sim_ReLU_all_params_2_hidden_layers_He}

In this subsection we present numerical simulation results (see \cref{fig:all_parameters_ReLU_2_hidden_layers_He} below) in a supervised learning framework for the training of deep \ReLU\ \ANNs\ with 4 layers (corresponding to $L = 3$ in \cref{setting:dnn}) with a 1-dimensional input layer (corresponding to $\ell_0 = 1$ in \cref{setting:dnn}), two 50-dimensional hidden layers (corresponding to $\ell_1 = \ell_2 = 50$ in \cref{setting:dnn}), and a 1-dimensional output layer (corresponding to $\ell_L = 1$ in \cref{setting:dnn}) where all \ANN\ parameters are modified during the training process. We also refer to \cref{figure_deep_4_layers} above for a graphical illustration of the \ANN\ architecture considered in this subsection.

Assume \cref{setting:dnn}, assume $L = 3$, $\ell_0 = \ell_3 = 1$, $\ell_1 = \ell_2 = 50$, and $\batchsize = 1024$, assume for all $x \in \R$ that $\act_{\infty} (x) = \allowbreak \max\{x, \allowbreak 0\}$ and $f(x) = x^2 + 2x$, assume for all $k \in \{1, \allowbreak 2, \allowbreak \ldots, \allowbreak \batchsize\}$, $n \in \N$, $x \in [0, 1]$ that $\P(X_k^n < x) = x$, assume
\begin{equation}\label{eqn:index_all_params_ReLU_He_4_Layers}
\textstyle \indexset = \N \cap [1, \dimension],
\end{equation}
assume for all $k \in \{1, 2, \ldots, L\}$, $i \in \{1, 2, \ldots, \ell_k\}$, $j \in \{1, 2, \ldots, \ell_{k - 1}\}$, $x \in \R$ that
\begin{equation}
\textstyle \P\big(\fb_i^{k, \Theta_0} = 0 \big) = 1 \qqandqq \P\big(\fw_{i, j}^{k, \Theta_0} \ge x \big) = \int_x^{\infty} \frac{1}{2} \bigl[\frac{\ell_{k - 1}}{\pi}\bigr]^{1/2} \exp\bigl(- \frac{y^2 \ell_{k - 1}}{4}\bigr) \, \d y,
\end{equation}
and assume for all $n \in \N$, $g_1, g_2, \ldots, g_n \in \R^{\dimension}$ that $\psi_n(g_1, \allowbreak g_2, \allowbreak \ldots, \allowbreak g_n) = 10^{-2} g_n$.

Let us add a few comments regarding the above presented setup. \Nobs that the set $\indexset$ in \cref{eqn:index_all_params_ReLU_He_4_Layers} represents the index set labeling all parameters of the considered deep \ANNs. Within this subsection all of the parameters of the \ANNs\ are modified during the training procedure. Furthermore, \nobs that the distribution of $\Theta_0$ is nothing else but the He normal initialization (cf., e.g., \cite{tensorflow}) for \ReLU\ \ANNs\ with 4 layers with a 1-dimensional input layer, two 50-dimensional hidden layers, and a 1-dimensional output layer. Moreover, \nobs that the assumption that for all $n \in \N$, $g_1, g_2, \ldots, g_n \in \R^{\dimension}$ it holds that $\psi_n(g_1, g_2, \ldots, g_n) = 10^{-2} g_n$ ensures that the stochastic process $\Theta = (\Theta^1, \ldots, \Theta^{\dimension}) \colon \N_0 \times \Omega \to \R^{\dimension}$ in \cref{eqn:setting:dnn:stochastic_pr} satisfies for all $n \in \N$ that
\begin{equation}
\textstyle \Theta_n = \Theta_{n - 1} - 10^{-2} \fG_n(\Theta_{n - 1})
\end{equation}
(cf.\ Subsection~\ref{subsubsec_Standard_SGD}). Next \nobs that the assumption that $L = 3$, the assumption that $\ell_0 = \ell_3 = 1$, and the assumption that $\ell_1 = \ell_2 = 50$ demonstrate that the number $\dimension \in \N$ of the employed real parameters to describe the considered deep \ANNs\ (see \cref{setting:dnn}) satisfies
\begin{equation}
\textstyle \dimension = \ell_1(\ell_0 + 1) + \ell_2 (\ell_1 + 1) + \ell_3 (\ell_2 + 1) = 50 (1 + 1) + 50 (50 + 1) + 1 (50 + 1) = 2701.
\end{equation}
In addition, \nobs that \cref{eqn:setting:dnn:loss}, the assumption that $L = 3$, and the assumption that $\batchsize = 1024$ assure that for all $n \in \N$, $\theta \in \R^{\dimension}$, $\omega \in \Omega$ it holds that $\lossapp_{\infty}^{n} (\theta, \omega) = \frac{1}{1024} \sum_{k = 1}^{1024} \abs{\cN_{\infty}^{3, \theta}(X_k^n (\omega)) - f(X_k^n (\omega))}^2$.

In the left image of Figure~\ref{fig:all_parameters_ReLU_2_hidden_layers_He} we approximately plot 300 samples of the random variable $\lossapp_{\infty}^{10001}(\Theta_{10000})$ in a histogram with 80 distinct subintervals (with 80 bins) on the $x$-axis. In the upper-right image of Figure~\ref{fig:all_parameters_ReLU_2_hidden_layers_He} we approximately plot the $L^2 (\P; \R^{\dimension})$-norm $(\E[\norm{\fG_{n + 1}(\Theta_{n})}]^2)^{1/2}$ of the generalized gradient against the number $n \in (\cup_{k = 0}^{998} \{10 k + 1\})$ of \SGD\ training steps. In the lower-right image of Figure~\ref{fig:all_parameters_ReLU_2_hidden_layers_He} we approximately plot the $L^2 (\P; \R^{\dimension})$-distance $(\E[\norm{\Theta_{n} - \Theta_{9991}}]^2)^{1/2}$ between $\Theta_n$ and $\Theta_{9991}$ against the number $n \in (\cup_{k = 0}^{998} \{10 k + 1\})$ of \SGD\ training steps. In, both, the upper-right image of Figure~\ref{fig:all_parameters_ReLU_2_hidden_layers_He} and the lower-right image of Figure~\ref{fig:all_parameters_ReLU_2_hidden_layers_He} we approximate the expectations in the $L^2 (\P; \R^{\dimension})$-norms by means of 300 independent Monte Carlo samples. The source code used to create Figure~\ref{fig:all_parameters_ReLU_2_hidden_layers_He} can be found at \url{https://github.com/deeplearningmethods/overcome-bad-local-minima} (cf.\ \cref{remark:python_codes} above).

\begin{figure}[H]
\centering
\subfloat{\includegraphics[width=16cm]{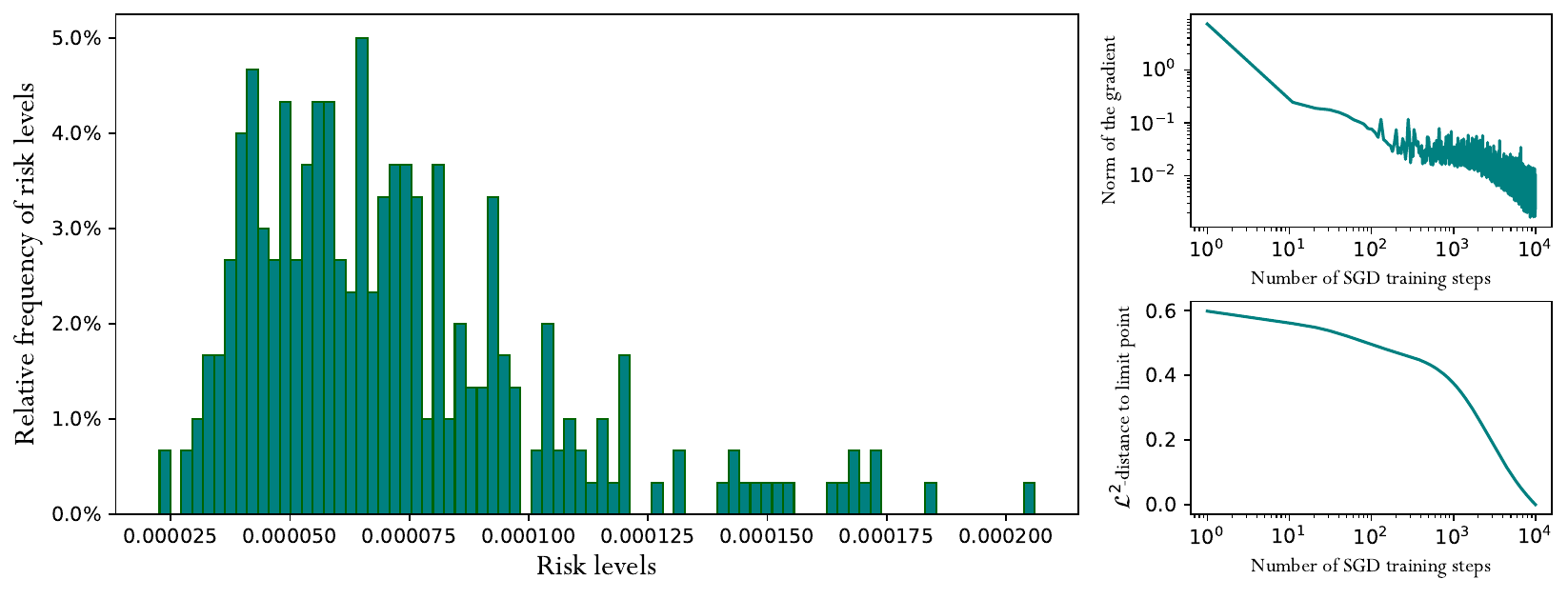}} 
\caption{Result of the numerical simulations described in  Subsection~\ref{subsec_num_sim_ReLU_all_params_2_hidden_layers_He}.}
\label{fig:all_parameters_ReLU_2_hidden_layers_He}
\end{figure}

\subsection{Adam training of all parameters of Xavier initialized ReLU ANNs (4 layers)}
\label{subsec_num_sim_ReLU_all_params_2_hidden_layers_Xavier_Adam}

In this subsection we present numerical simulation results (see \cref{fig:all_parameters_ReLU_2_hidden_layers_Xavier_Adam} below) in a supervised learning framework for the training of deep \ReLU\ \ANNs\ with 4 layers (corresponding to $L = 3$ in \cref{setting:dnn}) with a 1-dimensional input layer (corresponding to $\ell_0 = 1$ in \cref{setting:dnn}), two 50-dimensional hidden layers (corresponding to $\ell_1 = \ell_2 = 50$ in \cref{setting:dnn}), and a 1-dimensional output layer (corresponding to $\ell_L = 1$ in \cref{setting:dnn}) where all \ANN\ parameters are modified during the training process. We also refer to \cref{figure_deep_4_layers} above for a graphical illustration of the \ANN\ architecture considered in this subsection.

Assume \cref{setting:dnn}, assume $L = 3$, $\ell_0 = \ell_3 = 1$, $\ell_1 = \ell_2 = 50$, and $\batchsize = 1024$, assume for all $x \in \R$ that $\act_{\infty} (x) = \allowbreak \max\{x, \allowbreak 0\}$ and $f(x) = x^2 + 2x$, assume
\begin{equation}\label{eqn:index_all_params_ReLU_Xavier_Adam_4_Layers}
\textstyle \indexset = \N \cap [1, \dimension],
\end{equation}
assume for all $n \in \N$, $g_1, g_2, \ldots, g_n \in \R^{\dimension}$, $i \in \{1, 2, \ldots, \dimension\}$ that 
\begin{equation}\label{eqn:Adam_l_rate}
\textstyle \psi_n^{i} (g_1, g_2, \ldots, g_n) \textstyle = 10^{-2} \, 2^{- \floor{\frac{n}{500}}} \left[\frac{\sum_{k = 1}^n 0.9^{n-k} 10^{-1} g_k^{i}}{1 - 0.9^n} \right] \left[10^{-8} + \Bigl[\frac{\sum_{k = 1}^n 0.999^{n-k} 10^{-3} \abs{g_k^{i}}^2}{1 - 0.999^n}\Bigr]^{\nicefrac{1}{2}} \right]^{-1},
\end{equation}
assume for all $k \in \{1, 2, \ldots, L\}$, $i \in \{1, 2, \ldots, \ell_k\}$, $j \in \{1, 2, \ldots, \ell_{k - 1}\}$, $x \in \R$ that
\begin{equation}\label{eqn:Adam_Theta_0}
\textstyle \P\big(\fb_i^{k, \Theta_0} = 0 \big) = 1 \qqandqq \P\big(\fw_{i, j}^{k, \Theta_0} \ge x \big) = \int_x^{\infty} \frac{1}{2} \bigl[\frac{\ell_{k - 1} + \ell_k}{\pi}\bigr]^{1/2} \exp\bigl(- \frac{y^2 (\ell_{k - 1} + \ell_k)}{4}\bigr) \, \d y,
\end{equation}
and assume for all $k \in \{1, \allowbreak 2, \allowbreak \ldots, \allowbreak \batchsize\}$, $n \in \N$, $x \in [0, 1]$ that $\P(X_k^n < x) = x$.

Let us add a few comments regarding the above presented setup. \Nobs that the set $\indexset$ in \cref{eqn:index_all_params_ReLU_Xavier_Adam_4_Layers} represents the index set labeling all parameters of the considered deep \ANNs. Within this subsection all of the parameters of the \ANNs\ are modified during the training procedure. Furthermore, \nobs that the distribution of $\Theta_0$ is nothing else but the Xavier normal (Glorot normal) initialization (cf., e.g., \cite{tensorflow}) for \ReLU\ \ANNs\ with 4 layers with a 1-dimensional input layer, two 50-dimensional hidden layers, and a 1-dimensional output layer. Moreover, \nobs that \cref{eqn:Adam_l_rate} ensures that the stochastic process $\Theta = (\Theta^1, \ldots, \Theta^{\dimension}) \colon \N_0 \times \Omega \to \R^{\dimension}$ in \cref{eqn:setting:dnn:stochastic_pr} satisfies for all $n \in \N$, $i \in \{1, 2, \ldots, \dimension\}$ that
\begin{equation}
\textstyle \Theta_n^i = \Theta_{n - 1}^i - 10^{-2} \, 2^{- \floor{\frac{n}{500}}} \left[\frac{\sum_{k = 1}^n 0.9^{n-k} 10^{-1} \fG_k^i(\Theta_{k - 1})}{1 - 0.9^n} \right] \left[10^{-8} + \Bigl[\frac{\sum_{k = 1}^n 0.999^{n-k} 10^{-3} \abs{\fG_k^i(\Theta_{k - 1})}^2}{1 - 0.999^n}\Bigr]^{\nicefrac{1}{2}} \right]^{-1}
\end{equation}
(cf.\ Subsection~\ref{subsubsec_Adam_SGD}). Next \nobs that the assumption that $L = 3$, the assumption that $\ell_0 = \ell_3 = 1$, and the assumption that $\ell_1 = \ell_2 = 50$ demonstrate that the number $\dimension \in \N$ of the employed real parameters to describe the considered deep \ANNs\ (see \cref{setting:dnn}) satisfies
\begin{equation}
\textstyle \dimension = \ell_1(\ell_0 + 1) + \ell_2 (\ell_1 + 1) + \ell_3 (\ell_2 + 1) = 50 (1 + 1) + 50 (50 + 1) + 1 (50 + 1) = 2701.
\end{equation}
In addition, \nobs that \cref{eqn:setting:dnn:loss}, the assumption that $L = 3$, and the assumption that $\batchsize = 1024$ assure that for all $n \in \N$, $\theta \in \R^{\dimension}$, $\omega \in \Omega$ it holds that $\lossapp_{\infty}^{n} (\theta, \omega) = \frac{1}{1024} \sum_{k = 1}^{1024} \abs{\cN_{\infty}^{3, \theta}(X_k^n (\omega)) - f(X_k^n (\omega))}^2$.

In the left image of Figure~\ref{fig:all_parameters_ReLU_2_hidden_layers_Xavier_Adam} we approximately plot 300 samples of the random variable $\lossapp_{\infty}^{10001}(\Theta_{10000})$ in a histogram with 80 distinct subintervals (with 80 bins) on the $x$-axis. In the upper-right image of Figure~\ref{fig:all_parameters_ReLU_2_hidden_layers_Xavier_Adam} we approximately plot the $L^2 (\P; \R^{\dimension})$-norm $(\E[\norm{\fG_{n + 1}(\Theta_{n})}]^2)^{1/2}$ of the generalized gradient against the number $n \in (\cup_{k = 0}^{998} \{10 k + 1\})$ of \SGD\ training steps. In the lower-right image of Figure~\ref{fig:all_parameters_ReLU_2_hidden_layers_Xavier_Adam} we approximately plot the $L^2 (\P; \R^{\dimension})$-distance $(\E[\norm{\Theta_{n} - \Theta_{9991}}]^2)^{1/2}$ between $\Theta_n$ and $\Theta_{9991}$ against the number $n \in (\cup_{k = 0}^{998} \{10 k + 1\})$ of \SGD\ training steps. In, both, the upper-right image of Figure~\ref{fig:all_parameters_ReLU_2_hidden_layers_Xavier_Adam} and the lower-right image of Figure~\ref{fig:all_parameters_ReLU_2_hidden_layers_Xavier_Adam} we approximate the expectations in the $L^2 (\P; \R^{\dimension})$-norms by means of 300 independent Monte Carlo samples. The source code used to create Figure~\ref{fig:all_parameters_ReLU_2_hidden_layers_Xavier_Adam} can be found at \url{https://github.com/deeplearningmethods/overcome-bad-local-minima} (cf.\ \cref{remark:python_codes} above).

\begin{figure}[H]
\centering
\subfloat{\includegraphics[width=16cm]{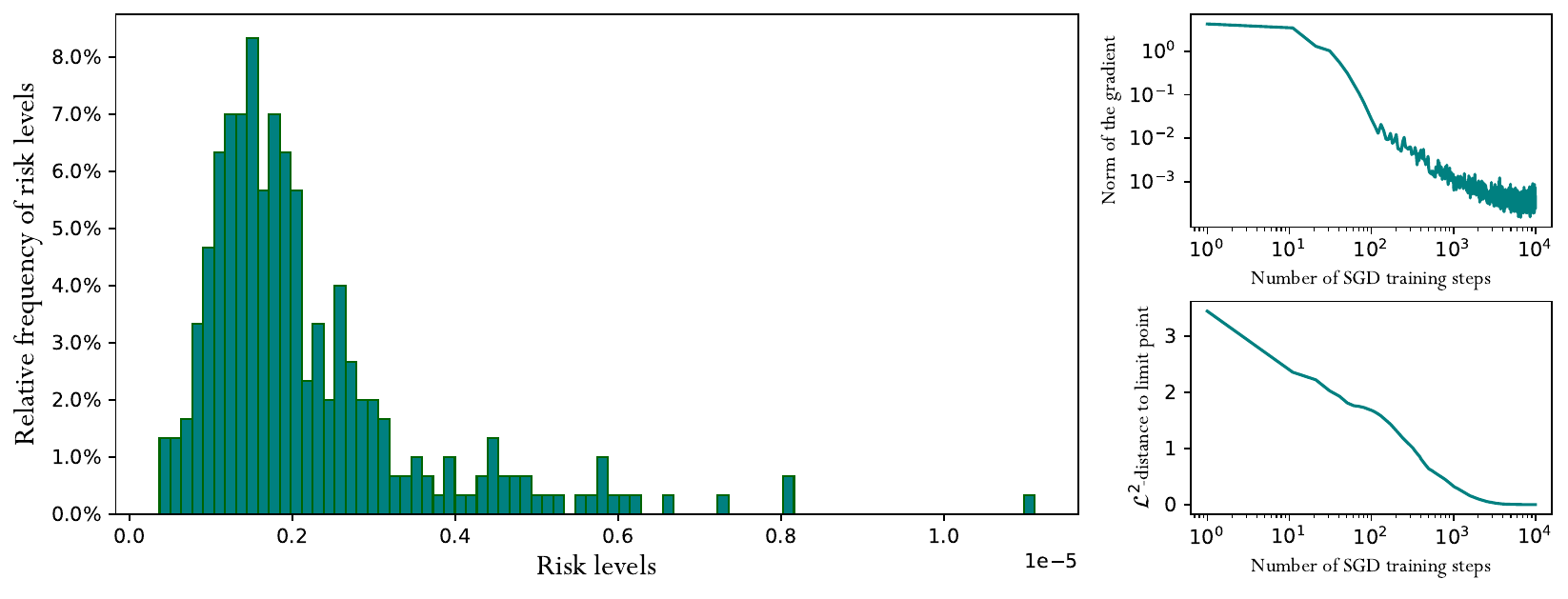}} 
\caption{Result of the numerical simulations described in Subsection~\ref{subsec_num_sim_ReLU_all_params_2_hidden_layers_Xavier_Adam}.}
\label{fig:all_parameters_ReLU_2_hidden_layers_Xavier_Adam}
\end{figure}

\subsection{Adam training of all parameters of He initialized ReLU ANNs (4 layers)}
\label{subsec_num_sim_ReLU_all_params_2_hidden_layers_He_Adam}

In this subsection we present numerical simulation results (see \cref{fig:all_parameters_ReLU_2_hidden_layers_He_Adam} below) in a supervised learning framework for the training of deep \ReLU\ \ANNs\ with 4 layers (corresponding to $L = 3$ in \cref{setting:dnn}) with a 1-dimensional input layer (corresponding to $\ell_0 = 1$ in \cref{setting:dnn}), two 50-dimensional hidden layers (corresponding to $\ell_1 = \ell_2 = 50$ in \cref{setting:dnn}), and a 1-dimensional output layer (corresponding to $\ell_L = 1$ in \cref{setting:dnn}) where all \ANN\ parameters are modified during the training process. We also refer to \cref{figure_deep_4_layers} above for a graphical illustration of the \ANN\ architecture considered in this subsection.

Assume \cref{setting:dnn}, assume $L = 3$, $\ell_0 = \ell_3 = 1$, $\ell_1 = \ell_2 = 50$, and $\batchsize = 1024$, assume for all $x \in \R$ that $\act_{\infty} (x) = \allowbreak \max\{x, \allowbreak 0\}$ and $f(x) = x^2 + 2x$, assume
\begin{equation}\label{eqn:index_all_params_ReLU_He_Adam_4_Layers}
\textstyle \indexset = \N \cap [1, \dimension],
\end{equation}
assume for all $n \in \N$, $g_1, g_2, \ldots, g_n \in \R^{\dimension}$, $i \in \{1, 2, \ldots, \dimension\}$ that 
\begin{equation}\label{eqn:Adam_l_rate_He}
\textstyle \psi_n^{i} (g_1, g_2, \ldots, g_n) \textstyle = 10^{-2} \, 2^{- \floor{\frac{n}{500}}} \left[\frac{\sum_{k = 1}^n 0.9^{n-k} 10^{-1} g_k^{i}}{1 - 0.9^n} \right] \left[10^{-8} + \Bigl[\frac{\sum_{k = 1}^n 0.999^{n-k} 10^{-3} \abs{g_k^{i}}^2}{1 - 0.999^n}\Bigr]^{\nicefrac{1}{2}} \right]^{-1},
\end{equation}
assume for all $k \in \{1, 2, \ldots, L\}$, $i \in \{1, 2, \ldots, \ell_k\}$, $j \in \{1, 2, \ldots, \ell_{k - 1}\}$, $x \in \R$ that
\begin{equation}\label{eqn:Adam_Theta_0_He}
\textstyle \P\big(\fb_i^{k, \Theta_0} = 0 \big) = 1 \qqandqq \P\big(\fw_{i, j}^{k, \Theta_0} \ge x \big) = \int_x^{\infty} \frac{1}{2} \bigl[\frac{\ell_{k - 1}}{\pi}\bigr]^{1/2} \exp\bigl(- \frac{y^2 \ell_{k - 1}}{4}\bigr) \, \d y,
\end{equation}
and assume for all $k \in \{1, \allowbreak 2, \allowbreak \ldots, \allowbreak \batchsize\}$, $n \in \N$, $x \in [0, 1]$ that $\P(X_k^n < x) = x$.

Let us add a few comments regarding the above presented setup. \Nobs that the set $\indexset$ in \cref{eqn:index_all_params_ReLU_He_Adam_4_Layers} represents the index set labeling all parameters of the considered deep \ANNs. Within this subsection all of the parameters of the \ANNs\ are modified during the training procedure. Furthermore, \nobs that the distribution of $\Theta_0$ is nothing else but the He normal initialization (cf., e.g., \cite{tensorflow}) for \ReLU\ \ANNs\ with 4 layers with a 1-dimensional input layer, two 50-dimensional hidden layers, and a 1-dimensional output layer. Moreover, \nobs that \cref{eqn:Adam_l_rate_He} ensures that the stochastic process $\Theta = (\Theta^1, \ldots, \Theta^{\dimension}) \colon \N_0 \times \Omega \to \R^{\dimension}$ in \cref{eqn:setting:dnn:stochastic_pr} satisfies for all $n \in \N$, $i \in \{1, 2, \ldots, \dimension\}$ that
\begin{equation}
\textstyle \Theta_n^i = \Theta_{n - 1}^i - 10^{-2} \, 2^{- \floor{\frac{n}{500}}} \left[\frac{\sum_{k = 1}^n 0.9^{n-k} 10^{-1} \fG_k^i(\Theta_{k - 1})}{1 - 0.9^n} \right] \left[10^{-8} + \Bigl[\frac{\sum_{k = 1}^n 0.999^{n-k} 10^{-3} \abs{\fG_k^i(\Theta_{k - 1})}^2}{1 - 0.999^n}\Bigr]^{\nicefrac{1}{2}} \right]^{-1}
\end{equation}
(cf.\ Subsection~\ref{subsubsec_Adam_SGD}). Next \nobs that the assumption that $L = 3$, the assumption that $\ell_0 = \ell_3 = 1$, and the assumption that $\ell_1 = \ell_2 = 50$ demonstrate that the number $\dimension \in \N$ of the employed real parameters to describe the considered deep \ANNs\ (see \cref{setting:dnn}) satisfies
\begin{equation}
\textstyle \dimension = \ell_1(\ell_0 + 1) + \ell_2 (\ell_1 + 1) + \ell_3 (\ell_2 + 1) = 50 (1 + 1) + 50 (50 + 1) + 1 (50 + 1) = 2701.
\end{equation}
In addition, \nobs that \cref{eqn:setting:dnn:loss}, the assumption that $L = 3$, and the assumption that $\batchsize = 1024$ assure that for all $n \in \N$, $\theta \in \R^{\dimension}$, $\omega \in \Omega$ it holds that $\lossapp_{\infty}^{n} (\theta, \omega) = \frac{1}{1024} \sum_{k = 1}^{1024} \abs{\cN_{\infty}^{3, \theta}(X_k^n (\omega)) - f(X_k^n (\omega))}^2$.

In the left image of Figure~\ref{fig:all_parameters_ReLU_2_hidden_layers_He_Adam} we approximately plot 300 samples of the random variable $\lossapp_{\infty}^{10001}(\Theta_{10000})$ in a histogram with 80 distinct subintervals (with 80 bins) on the $x$-axis. In the upper-right image of Figure~\ref{fig:all_parameters_ReLU_2_hidden_layers_He_Adam} we approximately plot the $L^2 (\P; \R^{\dimension})$-norm $(\E[\norm{\fG_{n + 1}(\Theta_{n})}]^2)^{1/2}$ of the generalized gradient against the number $n \in (\cup_{k = 0}^{998} \{10 k + 1\})$ of \SGD\ training steps. In the lower-right image of Figure~\ref{fig:all_parameters_ReLU_2_hidden_layers_He_Adam} we approximately plot the $L^2 (\P; \R^{\dimension})$-distance $(\E[\norm{\Theta_{n} - \Theta_{9991}}]^2)^{1/2}$ between $\Theta_n$ and $\Theta_{9991}$ against the number $n \in (\cup_{k = 0}^{998} \{10 k + 1\})$ of \SGD\ training steps. In, both, the upper-right image of Figure~\ref{fig:all_parameters_ReLU_2_hidden_layers_He_Adam} and the lower-right image of Figure~\ref{fig:all_parameters_ReLU_2_hidden_layers_He_Adam} we approximate the expectations in the $L^2 (\P; \R^{\dimension})$-norms by means of 300 independent Monte Carlo samples. The source code used to create Figure~\ref{fig:all_parameters_ReLU_2_hidden_layers_He_Adam} can be found at \url{https://github.com/deeplearningmethods/overcome-bad-local-minima} (cf.\ \cref{remark:python_codes} above).

\begin{figure}[H]
\centering
\subfloat{\includegraphics[width=16cm]{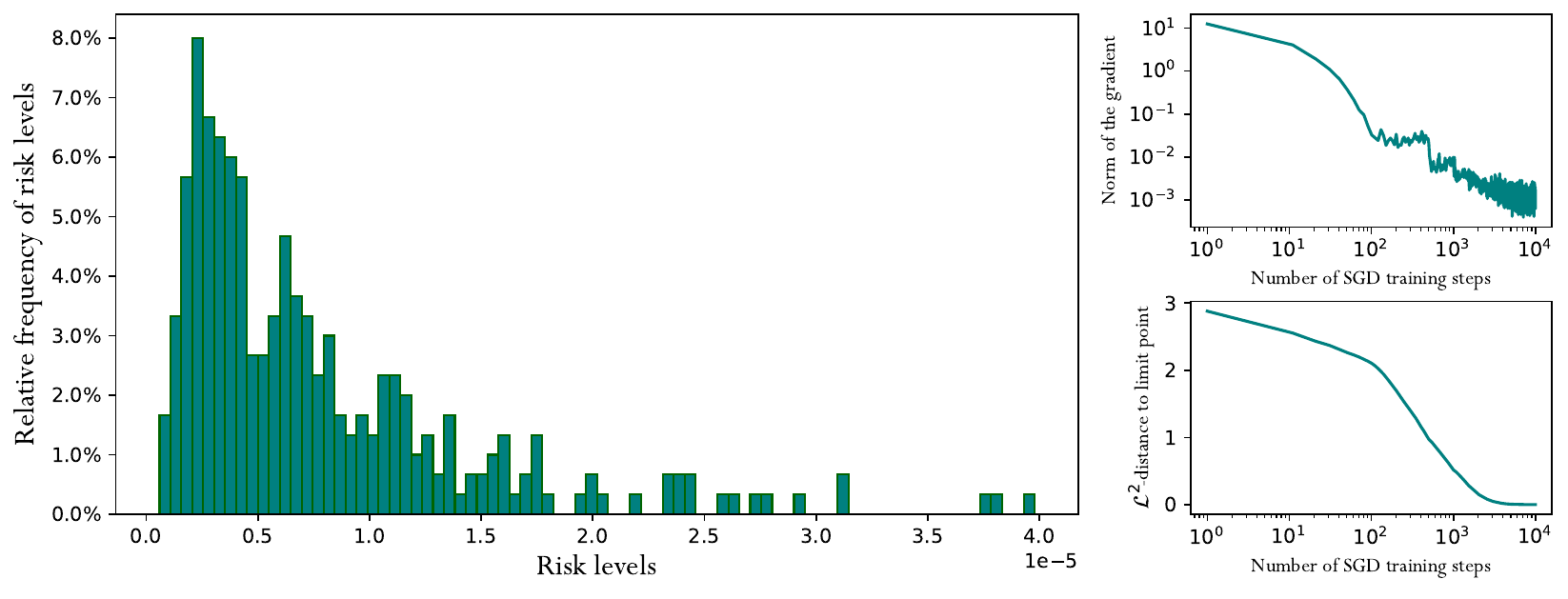}} 
\caption{Result of the numerical simulations described in Subsection~\ref{subsec_num_sim_ReLU_all_params_2_hidden_layers_He_Adam}.}
\label{fig:all_parameters_ReLU_2_hidden_layers_He_Adam}
\end{figure}

\subsection{SGD training of all parameters of Xavier initialized ReLU ANNs (5 layers)}
\label{subsec_num_sim_ReLU_all_params_3_hidden_layers_Xavier}

In this subsection we present numerical simulation results (see \cref{fig:all_parameters_ReLU_3_hidden_layers_Xavier} below) in a supervised learning framework for the training of deep \ReLU\ \ANNs\ with 5 layers (corresponding to $L = 4$ in \cref{setting:dnn}) with a 1-dimensional input layer (corresponding to $\ell_0 = 1$ in \cref{setting:dnn}), three 32-dimensional hidden layers (corresponding to $\ell_1 = \ell_2 = \ell_3 = 32$ in \cref{setting:dnn}), and a 1-dimensional output layer (corresponding to $\ell_L = 1$ in \cref{setting:dnn}) where all \ANN\ parameters are modified during the training process. We also refer to \cref{figure_deep_5_layers} below for a graphical illustration of the \ANN\ architecture considered in this subsection.

Assume \cref{setting:dnn}, assume $L = 4$, $\ell_0 = \ell_4 = 1$, $\ell_1 = \ell_2 = \ell_3 = 32$, and $\batchsize = 1024$, assume for all $x \in \R$ that $\act_{\infty} (x) = \allowbreak \max\{x, \allowbreak 0\}$ and $f(x) = x^2 + 2x$, assume for all $k \in \{1, \allowbreak 2, \allowbreak \ldots, \allowbreak \batchsize\}$, $n \in \N$, $x \in [0, 1]$ that $\P(X_k^n < x) = x$, assume
\begin{equation}\label{eqn:index_all_params_ReLU_Xavier_5_Layers}
\textstyle \indexset = \N \cap [1, \dimension],
\end{equation}
assume for all $k \in \{1, 2, \ldots, L\}$, $i \in \{1, 2, \ldots, \ell_k\}$, $j \in \{1, 2, \ldots, \ell_{k - 1}\}$, $x \in \R$ that
\begin{equation}
\textstyle \P\big(\fb_i^{k, \Theta_0} = 0 \big) = 1 \qqandqq \P\big(\fw_{i, j}^{k, \Theta_0} \ge x \big) = \int_x^{\infty} \frac{1}{2} \bigl[\frac{\ell_{k - 1} + \ell_k}{\pi}\bigr]^{1/2} \exp\bigl(- \frac{y^2 (\ell_{k - 1} + \ell_k)}{4}\bigr) \, \d y,
\end{equation}
and assume for all $n \in \N$, $g_1, g_2, \ldots, g_n \in \R^{\dimension}$ that $\psi_n(g_1, \allowbreak g_2, \allowbreak \ldots, \allowbreak g_n) = 10^{-2} g_n$.

\def\layersep{2.5cm}
\begin{figure}[H]
\centering
\begin{adjustbox}{width=\textwidth}
\begin{tikzpicture}[shorten >=1pt,->,draw=black!50, node distance=\layersep]
\tikzstyle{every pin edge}=[<-,shorten <=1pt]
\tikzstyle{input neuron}=[very thick, circle,draw=red, fill=red!20, minimum size=15pt,inner sep=0pt]
\tikzstyle{output neuron}=[very thick, circle, draw=green,fill=green!20,minimum size=20pt,inner sep=0pt]
\tikzstyle{hidden neuron}=[very thick, circle,draw=blue,fill=blue!20,minimum size=20pt,inner sep=0pt]
\tikzstyle{annot} = [text width=9em, text centered]
\tikzstyle{annot2} = [text width=4em, text centered]

%----------Neuron(s) of input layer----------
\node[input neuron] (I) at (0,-1.5) {$x$};

%----------Neuron(s) of 1st hidden layer----------
%\foreach \name / \y in {1, ..., 8}
\path[yshift = 1.5cm]
%node[hidden neuron] (H0-\name) at (\layersep, -\y cm) {};
node[hidden neuron](H0-1) at (\layersep, -1 cm) {};
\path[yshift = 1.5cm]
node[hidden neuron](H0-2) at (\layersep, -2 cm) {};
\path[yshift = 1.5cm]
node(H0-dots) at (\layersep, -2.9 cm) {\vdots};
\path[yshift = 1.5cm]
node[hidden neuron](H0-31) at (\layersep, -4 cm) {};
\path[yshift = 1.5cm]
node[hidden neuron](H0-32) at (\layersep, -5 cm) {};

%----------Neuron(s) of 2nd hidden layer----------
\path[yshift = 1.5cm]
%node[hidden neuron] (H0-\name) at (\layersep, -\y cm) {};
node[hidden neuron](H1-1) at (2*\layersep, -1 cm) {};
\path[yshift = 1.5cm]
node[hidden neuron](H1-2) at (2*\layersep, -2 cm) {};
\path[yshift = 1.5cm]
node(H1-dots) at (2*\layersep, -2.9 cm) {\vdots};
\path[yshift = 1.5cm]
node[hidden neuron](H1-31) at (2*\layersep, -4 cm) {};
\path[yshift = 1.5cm]
node[hidden neuron](H1-32) at (2*\layersep, -5 cm) {};

%----------Neuron(s) of 3rd hidden layer----------
\path[yshift = 1.5cm]
%node[hidden neuron] (H0-\name) at (\layersep, -\y cm) {};
node[hidden neuron](H2-1) at (3*\layersep, -1 cm) {};
\path[yshift = 1.5cm]
node[hidden neuron](H2-2) at (3*\layersep, -2 cm) {};
\path[yshift = 1.5cm]
node(H1-dots) at (3*\layersep, -2.9 cm) {\vdots};
\path[yshift = 1.5cm]
node[hidden neuron](H2-31) at (3*\layersep, -4 cm) {};
\path[yshift = 1.5cm]
node[hidden neuron](H2-32) at (3*\layersep, -5 cm) {};

%----------Neuron(s) of output layer----------
\path[yshift = 1.5cm]
node[output neuron](O) at (4*\layersep,-3 cm) {$\cN^{\theta}(x)$};

%----------Arrow(s) from 1st to 2nd layer----------
\foreach \dest in {1,2,31,32}
\path[line width = 0.8] (I) edge (H0-\dest);

%----------Arrow(s) from 2nd to 3rd layer----------
\foreach \source in {1,2,31,32}
\foreach \dest in {1,2,31,32}
\path[line width = 0.8] (H0-\source) edge (H1-\dest);

%----------Arrow(s) from 3rd to 4th layer----------
\foreach \source in {1,2,31,32}
\foreach \dest in {1,2,31,32}
\path[line width = 0.8] (H1-\source) edge (H2-\dest);

%----------Arrow(s) from 4th to 5th layer----------
\foreach \source in {1,2,31,32}
\path[line width = 0.8] (H2-\source) edge (O);

% Annotate the layers
\node[annot,above of=H0-1, node distance=1cm, align=center] (hl) {$1^{\text{st}}$ hidden layer\\($2^{\text{nd}}$ layer)};
\node[annot,above of=H1-1, node distance=1cm, align=center] (hl2) {$2^{\text{nd}}$ hidden layer\\($3^{\text{rd}}$ layer)};
\node[annot,above of=H2-1, node distance=1cm, align=center] (hl3) {$3^{\text{rd}}$ hidden layer\\($4^{\text{th}}$ layer)};
\node[annot,left of=hl, align=center] {Input layer\\ ($1^{\text{st}}$ layer)};
\node[annot,right of=hl3, align=center] {Output layer\\($5^{\text{th}}$ layer)};

\node[annot2,below of=H0-32, node distance=1cm, align=center] (sl) {$\ell_1=32$ \\ neurons};
\node[annot2,below of=H1-32, node distance=1cm, align=center] (sl2) {$\ell_2=32$ \\ neurons};
\node[annot2,below of=H2-32, node distance=1cm, align=center] (sl3) {$\ell_3=32$ \\ neurons};
\node[annot2,left of=sl, align=center] {$\ell_0=1$ \\ neuron};
\node[annot2,right of=sl3, align=center] {$\ell_4=1$ \\ neuron};
\end{tikzpicture}
\end{adjustbox}
\caption{Graphical illustration of the \ANN\ architectures used in Subsections~\ref{subsec_num_sim_ReLU_all_params_3_hidden_layers_Xavier} and \ref{subsec_num_sim_ReLU_all_params_3_hidden_layers_He}.}
\label{figure_deep_5_layers}
\end{figure}

Let us add a few comments regarding the above presented setup. \Nobs that the set $\indexset$ in \cref{eqn:index_all_params_ReLU_Xavier_5_Layers} represents the index set labeling all parameters of the considered deep \ANNs. Within this subsection all of the parameters of the \ANNs\ are modified during the training procedure. Furthermore, \nobs that the distribution of $\Theta_0$ is nothing else but the Xavier normal (Glorot normal) initialization (cf., e.g., \cite{tensorflow}) for \ReLU\ \ANNs\ with 5 layers with a 1-dimensional input layer, three 32-dimensional hidden layers, and a 1-dimensional output layer. Moreover, \nobs that the assumption that for all $n \in \N$, $g_1, g_2, \ldots, g_n \in \R^{\dimension}$ it holds that $\psi_n(g_1, g_2, \ldots, g_n) = 10^{-2} g_n$ ensures that the stochastic process $\Theta = (\Theta^1, \ldots, \Theta^{\dimension}) \colon \N_0 \times \Omega \to \R^{\dimension}$ in \cref{eqn:setting:dnn:stochastic_pr} satisfies for all $n \in \N$ that
\begin{equation}
\textstyle \Theta_n = \Theta_{n - 1} - 10^{-2} \fG_n(\Theta_{n - 1})
\end{equation}
(cf.\ Subsection~\ref{subsubsec_Standard_SGD}). Next \nobs that the assumption that $L = 4$, the assumption that $\ell_0 = \allowbreak \ell_4 \allowbreak = 1$, and the assumption that $\ell_1 = \ell_2 = \ell_3 = 32$ demonstrate that the number $\dimension \in \N$ of the employed real parameters to describe the considered deep \ANNs\ (see \cref{setting:dnn}) satisfies
\begin{equation}
\begin{split}
\textstyle \dimension & \textstyle = \ell_1(\ell_0 + 1) + \ell_2 (\ell_1 + 1) + \ell_3 (\ell_2 + 1) + \ell_4 (\ell_3 + 1) \\
& \textstyle  = 32 (1 + 1) + 32 (32 + 1) + 32 (32 + 1) + 1 (32 + 1) = 2209.
\end{split}
\end{equation}
In addition, \nobs that \cref{eqn:setting:dnn:loss}, the assumption that $L = 4$, and the assumption that $\batchsize = 1024$ assure that for all $n \in \N$, $\theta \in \R^{\dimension}$, $\omega \in \Omega$ it holds that $\lossapp_{\infty}^{n} (\theta, \omega) = \frac{1}{1024} \sum_{k = 1}^{1024} \abs{\cN_{\infty}^{4, \theta}(X_k^n (\omega)) - f(X_k^n (\omega))}^2$.

In the left image of Figure~\ref{fig:all_parameters_ReLU_3_hidden_layers_Xavier} we approximately plot 300 samples of the random variable $\lossapp_{\infty}^{10001}(\Theta_{10000})$ in a histogram with 80 distinct subintervals (with 80 bins) on the $x$-axis. In the upper-right image of Figure~\ref{fig:all_parameters_ReLU_3_hidden_layers_Xavier} we approximately plot the $L^2 (\P; \R^{\dimension})$-norm $(\E[\norm{\fG_{n + 1}(\Theta_{n})}]^2)^{1/2}$ of the generalized gradient against the number $n \in (\cup_{k = 0}^{998} \{10 k + 1\})$ of \SGD\ training steps. In the lower-right image of Figure~\ref{fig:all_parameters_ReLU_3_hidden_layers_Xavier} we approximately plot the $L^2 (\P; \R^{\dimension})$-distance $(\E[\norm{\Theta_{n} - \Theta_{9991}}]^2)^{1/2}$ between $\Theta_n$ and $\Theta_{9991}$ against the number $n \in (\cup_{k = 0}^{998} \{10 k + 1\})$ of \SGD\ training steps. In, both, the upper-right image of Figure~\ref{fig:all_parameters_ReLU_3_hidden_layers_Xavier} and the lower-right image of Figure~\ref{fig:all_parameters_ReLU_3_hidden_layers_Xavier} we approximate the expectations in the $L^2 (\P; \R^{\dimension})$-norms by means of 300 independent Monte Carlo samples. The source code used to create Figure~\ref{fig:all_parameters_ReLU_3_hidden_layers_Xavier} can be found at \url{https://github.com/deeplearningmethods/overcome-bad-local-minima} (cf.\ \cref{remark:python_codes} above).

\begin{figure}[H]
\centering
\subfloat{\includegraphics[width=16cm]{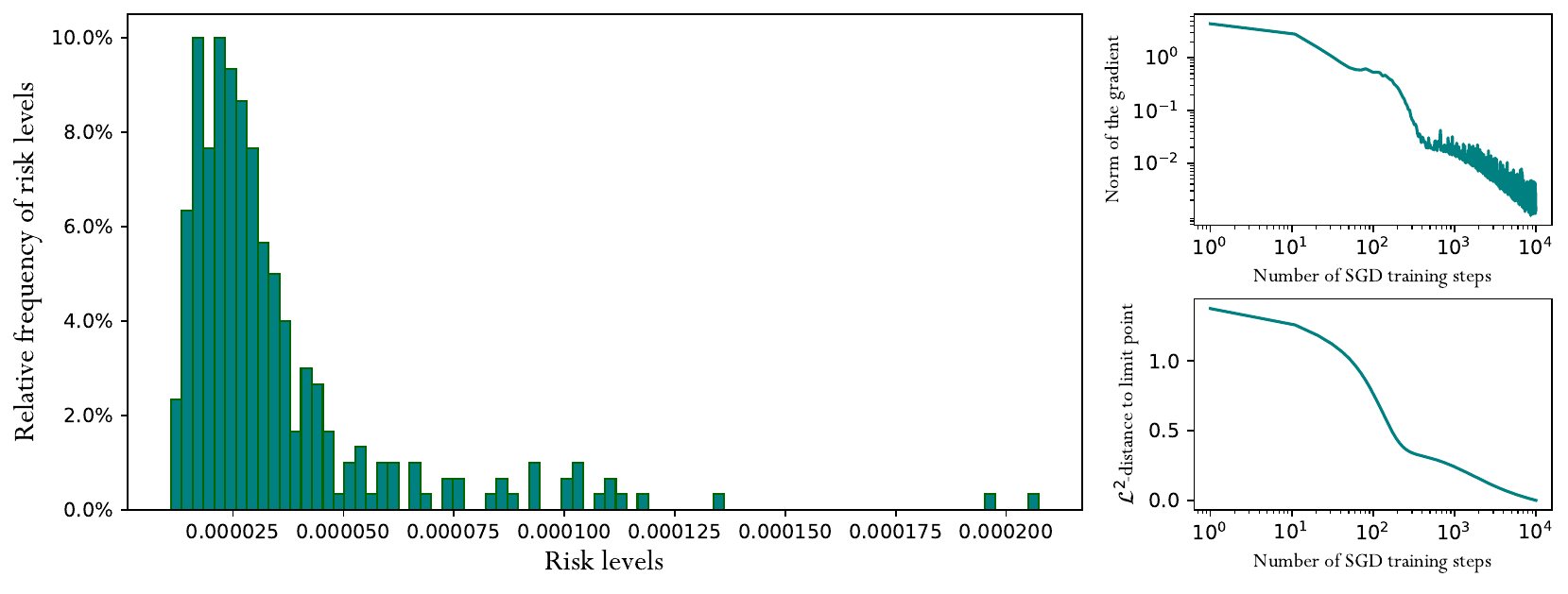}} 
\caption{Result of the numerical simulations described in Subsection~\ref{subsec_num_sim_ReLU_all_params_3_hidden_layers_Xavier}.}
\label{fig:all_parameters_ReLU_3_hidden_layers_Xavier}
\end{figure}

\subsection{SGD training of all parameters of He initialized ReLU ANNs (5 layers)}
\label{subsec_num_sim_ReLU_all_params_3_hidden_layers_He}

In this subsection we present numerical simulation results (see \cref{fig:all_parameters_ReLU_3_hidden_layers_He} below) in a supervised learning framework for the training of deep \ReLU\ \ANNs\ with 5 layers (corresponding to $L = 4$ in \cref{setting:dnn}) with a 1-dimensional input layer (corresponding to $\ell_0 = 1$ in \cref{setting:dnn}), three 32-dimensional hidden layers (corresponding to $\ell_1 = \ell_2 = \ell_3 = 32$ in \cref{setting:dnn}), and a 1-dimensional output layer (corresponding to $\ell_L = 1$ in \cref{setting:dnn}) where all \ANN\ parameters are modified during the training process. We also refer to \cref{figure_deep_5_layers} above for a graphical illustration of the \ANN\ architecture considered in this subsection.

Assume \cref{setting:dnn}, assume $L = 4$, $\ell_0 = \ell_4 = 1$, $\ell_1 = \ell_2 = \ell_3 = 32$, and $\batchsize = 1024$, assume for all $x \in \R$ that $\act_{\infty} (x) = \allowbreak \max\{x, \allowbreak 0\}$ and $f(x) = x^2 + 2x$, assume for all $k \in \{1, \allowbreak 2, \allowbreak \ldots, \allowbreak \batchsize\}$, $n \in \N$, $x \in [0, 1]$ that $\P(X_k^n < x) = x$, assume
\begin{equation}\label{eqn:index_all_params_ReLU_He_5_Layers}
\textstyle \indexset = \N \cap [1, \dimension],
\end{equation}
assume for all $k \in \{1, 2, \ldots, L\}$, $i \in \{1, 2, \ldots, \ell_k\}$, $j \in \{1, 2, \ldots, \ell_{k - 1}\}$, $x \in \R$ that
\begin{equation}
\textstyle \P\big(\fb_i^{k, \Theta_0} = 0 \big) = 1 \qqandqq \P\big(\fw_{i, j}^{k, \Theta_0} \ge x \big) = \int_x^{\infty} \frac{1}{2} \bigl[\frac{\ell_{k - 1}}{\pi}\bigr]^{1/2} \exp\bigl(- \frac{y^2 \ell_{k - 1}}{4}\bigr) \, \d y,
\end{equation}
and assume for all $n \in \N$, $g_1, g_2, \ldots, g_n \in \R^{\dimension}$ that $\psi_n(g_1, \allowbreak g_2, \allowbreak \ldots, \allowbreak g_n) = 10^{-2} g_n$.

Let us add a few comments regarding the above presented setup. \Nobs that the set $\indexset$ in \cref{eqn:index_all_params_ReLU_He_5_Layers} represents the index set labeling all parameters of the considered deep \ANNs. Within this subsection all of the parameters of the \ANNs\ are modified during the training procedure. Furthermore, \nobs that the distribution of $\Theta_0$ is nothing else but the He normal initialization (cf., e.g., \cite{tensorflow}) for \ReLU\ \ANNs\ with 5 layers with a 1-dimensional input layer, three 32-dimensional hidden layers, and a 1-dimensional output layer. Moreover, \nobs that the assumption that for all $n \in \N$, $g_1, g_2, \ldots, g_n \in \R^{\dimension}$ it holds that $\psi_n(g_1, g_2, \ldots, g_n) = 10^{-2} g_n$ ensures that the stochastic process $\Theta = (\Theta^1, \ldots, \Theta^{\dimension}) \colon \N_0 \times \Omega \to \R^{\dimension}$ in \cref{eqn:setting:dnn:stochastic_pr} satisfies for all $n \in \N$ that
\begin{equation}
\textstyle \Theta_n = \Theta_{n - 1} - 10^{-2} \fG_n(\Theta_{n - 1})
\end{equation}
(cf.\ Subsection~\ref{subsubsec_Standard_SGD}). Next \nobs that the assumption that $L = 4$, the assumption that $\ell_0 = \allowbreak \ell_4 \allowbreak = 1$, and the assumption that $\ell_1 = \ell_2 = \ell_3 = 32$ demonstrate that the number $\dimension \in \N$ of the employed real parameters to describe the considered deep \ANNs\ (see \cref{setting:dnn}) satisfies
\begin{equation}
\begin{split}
\textstyle \dimension & \textstyle = \ell_1(\ell_0 + 1) + \ell_2 (\ell_1 + 1) + \ell_3 (\ell_2 + 1) + \ell_4 (\ell_3 + 1) \\
& \textstyle  = 32 (1 + 1) + 32 (32 + 1) + 32 (32 + 1) + 1 (32 + 1) = 2209.
\end{split}
\end{equation}
In addition, \nobs that \cref{eqn:setting:dnn:loss}, the assumption that $L = 4$, and the assumption that $\batchsize = 1024$ assure that for all $n \in \N$, $\theta \in \R^{\dimension}$, $\omega \in \Omega$ it holds that $\lossapp_{\infty}^{n} (\theta, \omega) = \frac{1}{1024} \sum_{k = 1}^{1024} \abs{\cN_{\infty}^{4, \theta}(X_k^n (\omega)) - f(X_k^n (\omega))}^2$.

In the left image of Figure~\ref{fig:all_parameters_ReLU_3_hidden_layers_He} we approximately plot 300 samples of the random variable $\lossapp_{\infty}^{10001}(\Theta_{10000})$ in a histogram with 80 distinct subintervals (with 80 bins) on the $x$-axis. In the upper-right image of Figure~\ref{fig:all_parameters_ReLU_3_hidden_layers_He} we approximately plot the $L^2 (\P; \R^{\dimension})$-norm $(\E[\norm{\fG_{n + 1}(\Theta_{n})}]^2)^{1/2}$ of the generalized gradient against the number $n \in (\cup_{k = 0}^{998} \{10 k + 1\})$ of \SGD\ training steps. In the lower-right image of Figure~\ref{fig:all_parameters_ReLU_3_hidden_layers_He} we approximately plot the $L^2 (\P; \R^{\dimension})$-distance $(\E[\norm{\Theta_{n} - \Theta_{9991}}]^2)^{1/2}$ between $\Theta_n$ and $\Theta_{9991}$ against the number $n \in (\cup_{k = 0}^{998} \{10 k + 1\})$ of \SGD\ training steps. In, both, the upper-right image of Figure~\ref{fig:all_parameters_ReLU_3_hidden_layers_He} and the lower-right image of Figure~\ref{fig:all_parameters_ReLU_3_hidden_layers_He} we approximate the expectations in the $L^2 (\P; \R^{\dimension})$-norms by means of 300 independent Monte Carlo samples. The source code used to create Figure~\ref{fig:all_parameters_ReLU_3_hidden_layers_He} can be found at \url{https://github.com/deeplearningmethods/overcome-bad-local-minima} (cf.\ \cref{remark:python_codes} above).

\begin{figure}[H]
\centering
\subfloat{\includegraphics[width=16cm]{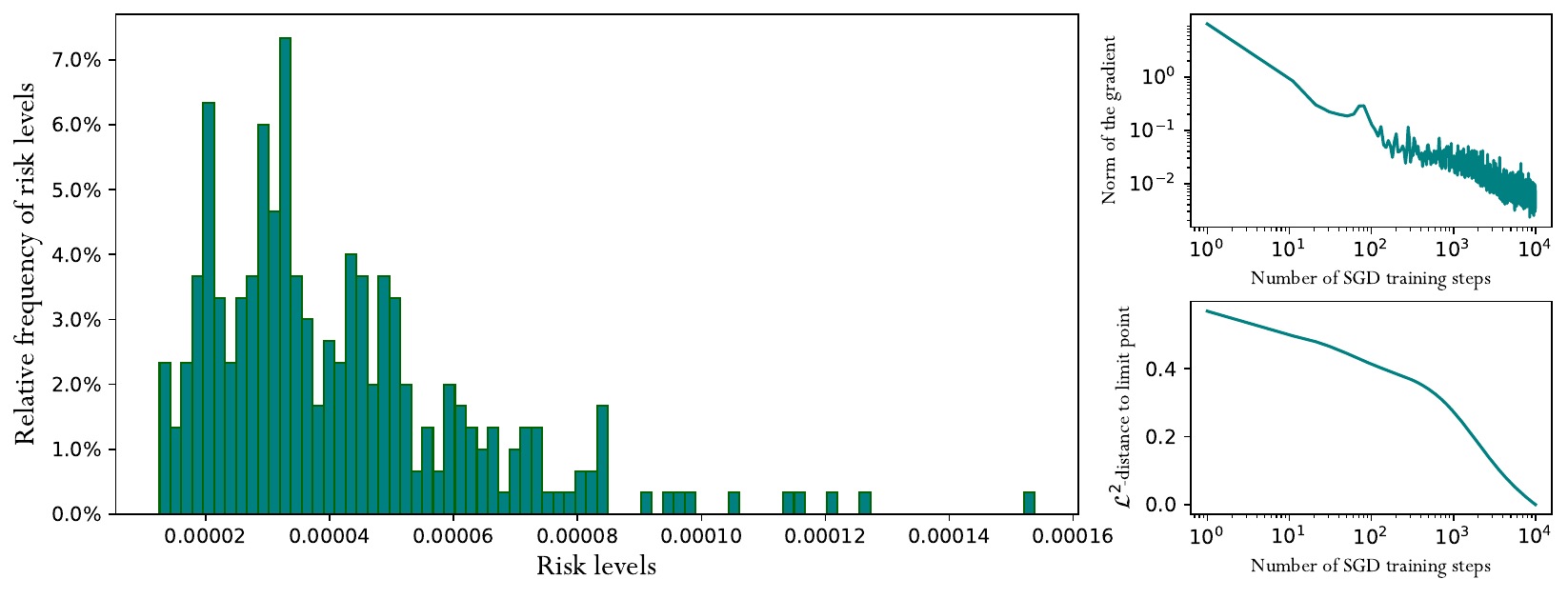}} 
\caption{Result of the numerical simulations described in  Subsection~\ref{subsec_num_sim_ReLU_all_params_3_hidden_layers_He}.}
\label{fig:all_parameters_ReLU_3_hidden_layers_He}
\end{figure}

\subsection*{Acknowledgements}
Special thanks are due to Sebastian Becker for several very useful comments regarding the numerical simulations. This project has been partially funded by the European Union (ERC, MONTECARLO, 101045811). The views and the opinions expressed in this work are however those of the authors only and do not necessarily reflect those of the European Union or the European Research Council (ERC). Neither the European Union nor the granting authority can be held responsible for them. Moreover, we gratefully acknowledge the Cluster of Excellence EXC 2044-390685587, Mathematics M\"{u}nster: Dynamics-Geometry-Structure funded by the Deutsche Forschungsgemeinschaft (DFG, German Research Foundation).

\newpage
%------------------------------------------------------------------------------%
%\printbibliography
\bibliographystyle{acm}
\bibliography{bibfile}

\end{document}